\documentclass[12pt, reqno]{amsart}
\usepackage{mathptmx, enumerate, amsmath, amsthm, amscd, amsfonts, amssymb, graphicx, color}\usepackage[bookmarksnumbered, colorlinks, plainpages]{hyperref}
\hypersetup{colorlinks=true,linkcolor=red, anchorcolor=green, citecolor=cyan, urlcolor=red, filecolor=magenta, pdftoolbar=true}
\input{mathrsfs.sty}
\newtheorem{theorem}{Theorem}[section]
\newtheorem{antitheorem}[theorem]{Anti-theorem}
\newtheorem{lemma}[theorem]{Lemma}
\newtheorem{proposition}[theorem]{Proposition}
\newtheorem{corollary}[theorem]{Corollary}
\theoremstyle{definition}
\newtheorem{definition}[theorem]{Definition}
\newtheorem{example}[theorem]{Example}

\newtheorem{problem}[theorem]{Problem}
\theoremstyle{remark}
\newtheorem{remark}[theorem]{Remark}
\numberwithin{equation}{section}

\def\11b#1{\mathbf{1}_{_{#1}}}

\begin{document}

\title{Metric invariants in Banach and Jordan--Banach algebras}

\author[A.M. Peralta]{Antonio M. Peralta}
\address[A.M. Peralta]{Instituto de Matem{\'a}ticas de la Universidad de Granada (IMAG). Departamento de An{\'a}lisis Matem{\'a}tico, Facultad de Ciencias, Universidad de Granada, 18071 Granada, Spain.}
\email{aperalta@ugr.es}




\begin{abstract} In this note we collect some significant contributions on metric invariants for complex Banach algebras and Jordan--Banach algebras established during the last fifteen years. This note is mainly expository, but it also contains complete proofs and arguments, which in many cases are new or have been simplified. We have also included several new results. The common goal in the results is to seek for ``natural'' subsets, $\mathfrak{S}_{A},$ associated with each complex Banach or Jordan--Banach algebra $A$, sets which when equipped with a certain metric, $d_{A}$, enjoys the property that each surjective isometry from $(\mathfrak{S}_{A},d_A)$ to a similar set, $(\mathfrak{S}_{B},d_B),$ associated with another Banach or Jordan--Banach algebra $B$, extends to a surjective real-linear isometry from $A$ onto $B$. In case of a positive answer to this question, the problem of discussing whether in such a case the algebras $A$ and $B$ are in fact isomorphic or Jordan isomorphic is the subsequent question. The main results presented here will cover the cases in which the sets $(\mathfrak{S}_{A},d_A)$ and $(\mathfrak{S}_{B},d_B)$ are in one of the following situations:\begin{enumerate}[$(\checkmark)$]\item Subsets of the set of invertible elements in a unital complex Banach algebra or in a unital complex Jordan--Banach algebra with the metric induced by the norm. Specially in the cases of unital C$^*$- and JB$^*$-algebras.
\item The sets of positive invertible elements in unital C$^*$- or JB$^*$-algebras with respect to the metric induced by the norm and with respect to the Thompson's metric.
\item Subsets of the set of unitary elements in unital C$^*$- and JB$^*$-algebras.     
	\end{enumerate}	
\end{abstract}

\maketitle

\tableofcontents

\section{Introduction}\label{sec: intro}

This paper, originally thought as an expository essay on the results establishing metric invariants on associative Banach algebras and Jordan--Banach algebras, gathers in a unified perspective the main conclusions obtained during the last fifteen years  in this line. Nevertheless, what do we mean by a metric invariant? The idea behind a metric invariant is to find ``natural'' subsets, $\mathfrak{S}_{A},$ associated with each complex Banach or Jordan--Banach algebra $A$, sets which when equipped with a certain metric, $d_{A}$, enjoys the property that each surjective isometry from $(\mathfrak{S}_{A},d_A)$ to a similar set $(\mathfrak{S}_{B},d_B),$ associated with another Banach or Jordan--Banach algebra $B$, extends to a surjective real-linear isometry from $A$ onto $B$. In case of a positive answer to this question, the problem of discussing whether in such a case the algebras $A$ and $B$ are in fact isomorphic or Jordan isomorphic is the subsequent question. We shall see that examples of metric invariants include many natural subsets, like certain subsets of the set of invertible elements in a unital complex Banach algebra or in a unital complex Jordan--Banach algebra with the metric given by the norm, the set of positive invertible elements in a unital C$^*$- or JB$^*$-algebra with the distance given by the norm or with respect to the Thompson's metric, and certain subsets of the set of unitary elements in unital C$^*$- or JB$^*$-algebras with the distance given by the norm. These are typical problems on preservers.\smallskip

During the elaboration of this paper we struggled for completeness and the final version contains most of the results with complete proofs and references to prior works and basic background. The manuscript actually contains some new results, many new or modified and simplified arguments, and a detailed revision of the similarities and differences between the different conclusions for Banach and Jordan--Banach algebras. All Banach algebras and spaces in this note will be over the complex field.\smallskip

Why metric invariants? A natural candidate as metric invariant is the set, $A^{-1},$ of all invertible elements in a Banach algebra $A$ with the metric given by the norm. However, by a well known result which dates back to the seventies, there are examples of commutative complex Banach algebras $A$ and $B$ whose subgroups of invertible elements are topologically isomorphic as groups but $A$ and $B$ are not isomorphic as Banach algebras nor isometric as Banach spaces (cf. Example~\ref{example no algebra isomorphic but group of invertible topologically isomorphic}). It should be noted that the existence of a surjective linear isometry between the Banach spaces underlying to Banach algebras does not necessarily imply that both algebras are isomorphic (see Example~\ref{example the Banach space structure is not an invariant} and Anti-theorem~\ref{anti-theorem BKU} for an example in the setting of JB$^*$-algebras). \smallskip

Section~\ref{sec: invertible elements as an invariant} is devoted to the results studying the sets of invertible elements as a metric invariant. In \cite{Hat09} Hatori established the first result in this line showing the following:  Let $A$ and $B$ be unital Banach algebras. We additionally assume that $A$ is semisimple and commutative. Let $\mathfrak{A}$ and $\mathfrak{B}$ be open subgroups of $A^{-1}$ and $B^{-1}$, respectively. Let $\Delta : \mathfrak{A}\to \mathfrak{B}$ be a surjective isometry for the distances given by the norms of $A$ and $B$, respectively. Then $B$ is semisimple and commutative, and the mapping $\Delta(\11b{A})^{-1} \Delta$ admits an extension to an isometric real-linear algebra isomorphism from $A$ onto $B$. In particular, $A^{-1}$ and $B^{-1}$ are isometrically isomorphic as metric groups (see Theorem~\ref{thm Hatori semisimple commutative}). The hypotheses concerning commutativity and semisimplicity cannot be relaxed in the final conclusion (see Example~\ref{example Hatori isomerically isomorphic non isomorphic}). If we relax the hypotheses on $A$, another result by Hatori shows that if $\mathfrak{A}$ and $\mathfrak{B}$ are open subgroups of $A^{-1}$ and $B^{-1}$, respectively, and $\Delta : \mathfrak{A}\to \mathfrak{B}$ is a surjective isometry for the distances given by the norms of $A$ and $B$, respectively, then there exist a surjective real-linear isometry $T: A\to B,$ and an element $u_0$ in the Jacobson radical of $B$ which satisfy that $\Delta (a) = T (a)+ u_0,$ for every $a$ in $\mathfrak{A}$ (cf. Theorem~\ref{thm Hatori 2011}). If $\Delta : \mathfrak{A}\to \mathfrak{B}$ is a surjective isometric group isomorphism,  then there exists a real-linear isometric isomorphism $T: A\to B$ satisfying $\Delta (a) = T (a),$ for every $a$ in $\mathfrak{A}$. In particular, $A^{-1}$ and $B^{-1}$ are isometrically isomorphic as metrizable groups (see Corollary~\ref{cor thm group isometric Hatori 2011}).\smallskip

Among the consequences of the above results we conclude in Corollary~\ref{cor 2 thm group isometric Hatori 2011} that the following statements are equivalent for any pair of unital Banach algebras $A$ and $B$: \begin{enumerate}[$(a)$]
	\item $A$ and $B$ are isometrically isomorphic as real Banach algebras;
	\item  There exist two open subgroups $\mathfrak{A}$ and $\mathfrak{B}$ of $A^{-1}$ and $B^{-1}$, respectively, which are isometrically isomorphic as groups;
	\item $A^{-1}$ and $B^{-1}$ are isometrically isomorphic as metrizable groups.
\end{enumerate}\smallskip

We shall obtain the previous results from the corresponding versions in the wider setting of Jordan--Banach algebras. This will be a characteristic of this note because we try to emphasize the links and divergences between both settings. To facilitate the accessibility for novice readers, section~\ref{sec: taste of radicals in Jordan} is devoted to present the basic notions, definitions and facts in the theory of Jordan--Banach algebras with a complete list of references. We shall pay especial attention to the McCrimmon radical of Jordan--Banach algebras and its different characterizations (see Theorem~\ref{t Aupetit c 1} and Proposition~\ref{p Hatori McCrimmon radical}), the connected components of the invertible elements in a unital Jordan--Banach algebra (see Theorem~\ref{t Loos thm characterization of the principal component}), and the characterization of the principal component of the set of all invertible elements in a Jordan--Banach algebra $M$ as the least quadratic subset of $M^{-1}$ containing $\exp(M)$ (see Lemma~\ref{l principal component as the least quadratic subset of invertible elements}).\smallskip

Continuing with the results about subsets of invertible elements in a unital complex Jordan--Banach algebra $M$, the main conclusion known until now is treated in Theorem~\ref{t surjective isometries between invertible clopen} and asserts the following:  Let $M$ and $N$ be unital Jordan-Banach algebras. Suppose that $\mathfrak{M}\subseteq M^{-1}$ and $\mathfrak{N}\subseteq N^{-1}$ are clopen quadratic subsets which are closed for powers and inverses. Let $\Delta : \mathfrak{M}\to \mathfrak{N}$ be a surjective isometry. Then there exist a surjective real-linear isometry $T_0: M\to N$ and an element $u_0$ in the McCrimmon radical of $N$ such that $\Delta (a) = T_0(a) +u_0,$ for all $a\in \mathfrak{M}$. The Jordan product of two invertible elements is not, in general, an invertible element. This is the reason for considering quadratic subsets. We recall that a subset $\mathfrak{M}$ of $M^{-1}$ is called a \emph{quadratic subset} if it is non-empty and $U_{a} (b) \in \mathfrak{M}$ for all $a,b\in \mathfrak{M}$ (see Definition~\ref{def cuadratic subset}). The main results concerning Jordan--Banach algebras in this section are inspired by those in \cite{Pe2021}.\smallskip

Among the new results in section~\ref{sec: invertible elements as an invariant} we establish that if $A$ is a semisimple and commutative unital Banach algebra, $N$ is a unital Jordan--Banach algebra, $\mathfrak{A}$ and $\mathfrak{N}$ are clopen quadratic subsets closed for powers and inverses of $A^{-1}$ and $N^{-1},$ respectively, and  $\Delta : \mathfrak{A}\to \mathfrak{N}$ is a surjective isometry for the distances given by the norms of $A$ and $N$, then $N$ is semisimple and associative, and the mapping $\Delta(\11b{A})^{-1} \Delta$ admits an extension to an isometric real-linear algebra isomorphism from $A$ onto $N$. In particular, $A^{-1}$ and $N^{-1}$ are isometrically isomorphic as metric groups (see Theorem~\ref{t real-linear isometries from a semisimple commutative Banach algebra}). In the proof we employ a version of the Gleason-Kahane-\.{Z}elazko theorem for Jordan--Banach algebras from \cite{EscPeVill2023} and a new subtle variant for linear maps valued in the exponential spectrum (see Lemma~\ref{l GKZ for exponential spectrum valued real-linear}). \smallskip

Another new result derived from Theorem~\ref{t real-linear isometries from a semisimple commutative Banach algebra} reads as follows: Let $A$ be a unital and commutative Banach algebra and let $N$ be a semisimple unital Jordan-Banach algebra. Suppose that $\mathfrak{A}\subseteq A^{-1}$ and $\mathfrak{N}\subseteq N^{-1}$ are clopen quadratic subsets which are closed for powers and inverses. Let $\Delta : \mathfrak{A}\to \mathfrak{N}$ be a surjective isometry. Then $A$ is semisimple and $N$ is associative, and the mapping $\Delta(\11b{A})^{-1} \Delta$ admits an extension to an isometric real-linear algebra isomorphism from $A$ onto $N$. In particular, $A^{-1}$ and $N^{-1}$ are isometrically isomorphic as metric groups (see  Corollary~\ref{t real-linear isometries from a semisimple to a commutative Banach algebra}).\smallskip

In the setting of Jordan--Banach algebras preservers of Jordan products do not make any sense, the natural hypothesis is the preservation of expressions of the form $U_a(b)$. In Proposition~\ref{p surjective isometries preserving quadratic expressions} it is established that if $M$ and $N$ are unital Jordan-Banach algebras,  $\mathfrak{M}\subseteq M^{-1}$ and $\mathfrak{N}\subseteq N^{-1}$ are clopen quadratic subsets which are closed for powers and inverses, and $\Delta : \mathfrak{M}\to \mathfrak{N}$ is a surjective isometry satisfying $$\Delta (U_a (b)) = U_{\Delta(a)} (\Delta(b)), \hbox{ for all } a,b\in \mathfrak{M},$$ then setting $u = \Delta(\11b{M})$,  there exists a real-linear isometric Jordan isomorphism $J$ from $M$ onto the $u$-isotope $N_{(u)}$ such that $\Delta (a) = J(a) $ for all $a\in \mathfrak{M}$.\smallskip

The final conclusion in the setting of Jordan--Banach algebras is surveyed in Corollary~\ref{cor 2 thm group isometric Hatori 2011 JB} and shows that for any pair $M,N$ of unital complex Jordan--Banach algebras the following statements are equivalent:\begin{enumerate}[$(a)$]
	\item There exists an isometric real-linear Jordan isomorphism from $M$ onto $N$.
	\item  There exist two clopen quadratic subsets $\mathfrak{M}\subseteq M^{-1}$ and $\mathfrak{N}\subseteq N^{-1}$ which are closed for powers and inverses and a surjective isometry $\Delta : \mathfrak{M}\to \mathfrak{N}$ satisfying $\Delta(\11b{M})=\11b{N}$ and $\Delta (U_a (b)) = U_{\Delta(a)} (\Delta(b)),$  for all $a,b\in \mathfrak{M}.$
	\item There exists a surjective isometry $\Delta: M^{-1}\to N^{-1}$ satisfying $\Delta(\11b{M}) = \11b{N}$ and $\Delta (U_a (b)) = U_{\Delta(a)} (\Delta(b))$, for all $a,b\in {M}^{-1}$.
\end{enumerate}

We have already commented the existence of non-isomorphic Banach algebras which are isometrically isomorphic as Banach spaces (see Example~\ref{example Hatori isomerically isomorphic non isomorphic}). A celebrated result of Kadison proves that this pathology cannot happen if both algebras are C$^*$-algebras (see Theorem~\ref{t Kadison 55}). \smallskip

We devote the first part of section~\ref{sec: geometric goodness of C* and JB-algebras} to survey the results on surjective isometries between subsets of invertible elements in two unital C$^*$-algebras obtained by Hatori and Watanabe in \cite{HatWan2012}. It is worth recalling that every C$^*$-algebra is semisimple.\smallskip

Let $A$ and $B$ be unital C$^*$-algebras, and let $\mathfrak{A}$ and $\mathfrak{B}$ be open subgroups of $A^{-1}$ and $B^{-1},$ respectively. Under these hypotheses, a bijection $\Delta: \mathfrak{A}\to \mathfrak{B}$ is an isometry if and only if $\Delta(\11b{A})$ is unitary in $B$ and there exist a projection $p$ in the centre of $B$, and a complex-linear Jordan $^*$-isomorphism $\tilde{J}$ from $A$ onto $B$ such that $$ \Delta (a) = \Delta(\11b{A}) p \tilde{J}(a) + \Delta(\11b{A}) (\11b{N}-p) \tilde{J} (a)^*, \hbox{ for all } a \in  \mathfrak{A}.$$ Therefore $\Delta$ admits an extension to a surjective real-linear isometry from $A$ onto $B$ (cf. Theorem~\ref{t HatWat}).\smallskip

JB$^*$-triple theory with Kaup's Banach--Stone theorem (Theorem~\ref{t Kaup-Banach-Stone}) and the theory of isotopes for unital JB$^*$-algebras are subsequently employed to establish the following result: Let $M$ and $N$ be unital JB$^*$-algebras. Suppose that $\mathfrak{M}\subseteq M^{-1}$ and  $\mathfrak{N}\subseteq N^{-1}$ are clopen quadratic subsets which are closed for powers and inverses. Let $\Delta : \mathfrak{M}\to \mathfrak{N}$ be a surjective isometry. Then $\Delta(\11b{M}) =u$ is a unitary element in $N$ and there exists a projection $p$ in the centre of $M$ and a complex-linear Jordan $^*$-isomorphism $J$ from $M$ onto the $u^*$-homotope $N_{(u^*)} = (N,\circ_{u},*_u)$ such that $$\Delta (a) = J(p\circ a) + J ((\11b{M}-p) \circ a^*), \hbox{ for all } a\in \mathfrak{M}.$$ Assuming that there exists a unitary $\omega_0$ in $N$ such that the identity $U_{\omega_0} (\Delta(\mathbf{1})) = \mathbf{1}$ holds, then there exist a central projection $p\in M$ and a complex-linear Jordan $^*$-isomorphism $\Phi$ from $M$ onto $N$ such that $$\Delta (a) = U_{w_0^{*}} \left(\Phi (p\circ a) + \Phi ((\11b{M}-p) \circ a^*)\right),$$ for all $a\in \mathfrak{M}$ (see Theorem~\ref{t surjective isometries between invertible clopen JBstar}).\smallskip

The final affirmation can be better understood after a couple of comments. It is well known that the group of surjective linear isometries in a unital C$^*$-algebra $A$ acts transitively on the set of unitaries of $A$. However, this fact is, in general, false in the setting of JB$^*$-algebras, actually there exists a unital JC$^*$-algebra $M$ containing two unitaries $u_1$ and $u_2$ satisfying $T (u_1)\neq u_2$ for all $T\in Iso (M)$. Furthermore, each unitary produces a unital JB$^*$-algebra $M_{(u_1^*)}= (M, \circ_{u_1},*_{u_1})$ and $M_{(u_2^*)}= (M, \circ_{u_2},*_{u_2})$ which are not Jordan $^*$-isomorphic but have the same underlying Banach space (cf. Anti-theorem~\ref{anti-theorem BKU}). In a very informal language, the same complex Banach space $M$ can be equipped with many non-Jordan $^*$-isomorphic structures of JB$^*$-algebra, although by Kaup's Banach--Stone theorem (Theorem~\ref{t Kaup-Banach-Stone}) there is only a triple product making $M$ a JB$^*$-triple.\smallskip

The idea of considering the set of all positive invertible elements in a unital C$^*$-algebra as a metric invariant was considered by Hatori and Moln{\'a}r in \cite[page 167]{HatMol2014}. If we equip this set with the usual metric given by the norm, the conclusion can be achieved via a key theorem due to Mankiewicz (see Theorem~\ref{t Mankiewicz}). Here we extend the result to the setting of unital JB$^*$-algebras by showing that each surjective isometry between sets of positive invertible elements in two unital JB$^*$-algebras $M$ and $N$ admits an extension to a Jordan $^*$-isomorphism $J: M\to N$ (see Proposition~\ref{p cone of positive elements as metric invariant usual norm}).\smallskip

The first part of section~\ref{sec: geometric tool of Hatori Hirasawa Miura Molnar} is devoted to present an outstanding geometric tool due to Hatori, Hirasawa, Miura and Moln{\'a}r in \cite{HatHirMiuMol2012}, which provides a fix point theorem type for maps between abstract spaces preserving a certain mild relation $d$ between elements, which is a far-reaching generalization of the usual notion of distance (cf. Theorem~\ref{t main geometric tool}). The arguments are related to certain ideas in V\"{a}is{\"a}l{\"a}'s simplified proof of the celebrated Mazur--Ulam theorem, and the statement of this theorem is so general that permits many different applications of this exceptional geometric tool. Here we employ it in the study of surjective isometries between sets of positive invertible elements in two unital JB$^*$-algebras when equipped with the Thompson's metric (subsection~\ref{subsec: positive invertible with Thompson metric as invariant}), and in the study of surjective isometries between sets of unitary elements in two unital JB$^*$-algebras with respect to the metric given by the norm (section~\ref{sec: unitaries as invariant}). \smallskip

Thompson's metric on the interior of the positive cone $A^{+}$ (i.e. the set of all Archimedean order units) in an order unit space $A$ equipped with an Archimedean order unit is defined as follows:  For $a, b \in \hbox{int} (A^{+}),$ we set $$ M(a/b) := \inf\{ \beta > 0 \colon a \leq  \beta b\} (<\infty),$$ and Thompson's metric is defined by $$d_T (a, b) = \log\left(\max\{M(a/b), M(b/a)\}\right).$$ One of the most fruitful sources of examples in this paper is coming from the cones of positive elements in a C$^*$-algebra $A$, or more generally, in a JB$^*$-algebra $M$. The interior of $M^+$  is non-empty if and only if $M$ is unital and in such a case $\hbox{int}\left(M^+\right)= M^{++}$ (see Proposition~\ref{p when a JB-algebra admits an Archimedean order unit}). Andruchow, Corach, and Stojanoff gave in \cite{Andru2000} an explicit formula to compute Thompson's metric on the set positive invertible elements in a unital C$^*$-algebra $A$ via differential geometric methods. An alternative explanation for $B(H)^{++}$ was given by Moln{\'a}r in \cite{Mol2009}. The concrete expression affirms that for every $a, b \in A^{++}$ the Thompson's metric satisfies $$ d_{T} (a,b) = \left\| \log\left( a^{-1/2} b a^{-1/2} \right)\right\| \ \ \hbox{ (cf. Proposition~\ref{p Andrichoff Corach Stojanoff geomesic and Thompson's metric}).}$$

By applying the celebrated Shirshov--Cohn theorem (see Theorem~\ref{t Shirshov-Cohn}) we shall see that an appropriate Jordan version of the formula proved by Andruchow, Corach, and Stojanoff remains valid for unital JB$^*$-algebras. The Jordan version of this formula, which was originally established by Lemmens, Roelands, and Wortel in \cite{LemmRoeWor2018}, asserts that for every unital JB$^*$-algebra $M,$ the Thompson's metric on $M^{++}$ satisfies $$ d_{T} (a,b) = \left\| \log\left( U_{a^{-1/2}} (b) \right)\right\|, \hbox{ for every } a, b \in M^{++}.$$ Furthermore, $$ \begin{aligned}
	d_T \left(U_b (a^{-1}) , a \right) &= 2 d_T \left(b, a \right), \hbox{ and }\\
	d_T(U_c (a^{-1}), U_c (b^{-1})) &= d_T( a^{-1}, b^{-1}) = d_T (a,b),
\end{aligned}$$ for every $a, b, c \in M^{++}$ (see Proposition~\ref{p Lemmens on the Thompson distance for positive invertible in JB*-algebras}).\smallskip

In subsection \ref{subsec: positive invertible with Thompson metric as invariant} we revisit the results exploring the set of positive invertible elements with the Thompson's metric as a metric invariant for unital C$^*$-algebras and JB$^*$-algebras. A change in the metric produces a important difference in the arguments with respect to those employed in subsection~\ref{subsection positive invertible with the usual metric}, where this set was equipped with the metric given by the norm.\smallskip

In the case of C$^*$-algebras, a result due to Hatori and Moln{\'a}r \cite{HatMol2014} proves that if $A$ and $B$ are unital C$^*$-algebras, then a mapping $\Delta : (A^{++},d_T) \to (B^{++},d_T)$ is a surjective isometry if and only if there is a central projection $p$ in $B$ and a Jordan $^*$-isomorphism $J : A\to B$ such that $\Delta$ is of the form
$$\Delta (a) = \Delta(\11b{A})^{1/2} \left( p J(a) + (1 - p) J(a^{-1}) \right) \Delta(\11b{A})^{1/2}, \ \ a\in A^{++}$$ (see Theorem~\ref{t Hatori Molnar Thompson isometries}). A forerunner to this result in the case that $A = B= B(H)$ was previously established by Moln{\'a}r (see \cite[Theorem 1]{Mol2009}).\smallskip

The reader has probably realized that, since the question makes perfectly sense in the case of unital JB$^*$-algebras, it also seems interesting to determine if the set of positive invertible elements in a unital JB$^*$-algebra is a metric invariant when equipped with the Thompson's metric. This problem was studied, and solved, by Lemmens, Roelands, and Wortel via techniques based on differential geometry and geodesics. Here, we present an adaptation of the ideas by Hatori and Moln{\'a}r in \cite{HatMol2014} to conclude in Theorem~\ref{t Lemmens Roelands Wortel Thompson isometries} that if $M$ and $N$ are unital JB$^*$-algebras, then a mapping $\Delta : (M^{++},d_T) \to (N^{++},d_T)$ is a surjective isometry if and only if there is a central projection $p$ in $N$ and a Jordan $^*$-isomorphism $J : M\to N$ such that $\Delta$ is of the form
$$\Delta (a) = U_{\Delta(\11b{M})^{1/2}} \left( p\circ  J(a) + (\11b{N} - p)\circ J(a^{-1}) \right), \ \ a\in M^{++}.$$

We are slowly approaching to the contains in the final section~\ref{sec: unitaries as invariant}, where we study the set $\mathcal{U} (A)$, of all unitary elements in a unital C$^*$-algebra $A$, or more generally in a unital JB$^*$-algebra, as a metric invariant when equipped with the metric given by the norm. In strict chronological order, Hatori and Moln{\'a}r made the first contribution in the case that $A = B(H)$ in 2012 (see \cite{HatMol2012}). Concretely, if $H$ is a complex Hilbert space and $\Delta : \mathcal{U}(B(H)) \to  \mathcal{U}(B(H))$ is a surjective isometry with respect to the metric given by the operator norm, then there exists a surjective real-linear isometry (equivalently, a triple automorphism) $T$ on $B(H)$ such that $\Delta (u) = T (u)$ for all $u\in \mathcal{U}(B(H))$ (cf. Theorem~\ref{t Hatori Molnar unitaries BH}). In this special case we can further conclude that there exist unitaries $v,w\in \mathcal{U}(B(H))$ such that $\Delta$ is of one of the following forms:
\begin{enumerate}[$(a)$]
	\item $\Delta (u) = v u w^*,\ $ for all $u\in \mathcal{U}(B(H))$,
	\item $\Delta (u) = v u^* w^*,\ $ for all $u\in \mathcal{U}(B(H))$,
	\item $\Delta(u ) = v u^t w^*,\ $ for all $u\in \mathcal{U}(B(H)),$
	\item $\Delta(u) = v \overline{u} w^*, \ $ for all $u\in \mathcal{U}(B(H))$,
\end{enumerate}  where the maps $a\mapsto \overline{a}, a^t$ denote an abstract conjugation and transposition on $B(H)$, respectively.\smallskip

In a general unital C$^*$-algebra the set $\mathcal{U} (A)$ is a multiplicative subgroup which need not be connected, and the extension of a surjective isometry between the unitary subgroups of two unital C$^*$-algebras is not always possible (cf. \cite{Hatori14}). However, in a subsequent paper,  Hatori and Moln{\'a}r established that for each  surjective isometry $\Delta  : \mathcal{U} (A) \to \mathcal{U} (B)$, where $A$ and $B$ are unital C$^*$-algebras and we consider the metrics given  the C$^*$-norms, we have $\Delta (e^{i A_{sa}}) = \Delta(\11b{A})\  e^{ i B_{sa}}$ and there is a central projection $p \in B$ and a Jordan $^*$-isomorphism $J : A\to B$ such that
\begin{equation}\label{eq 2608} \Delta (e^{ix}) = \Delta(\11b{A}) (p J(e^{ix}) + (\11b{B} - p) J(e^{ix})^*),\  x\in A_{sa}.
\end{equation} Under the additional assumptions that $A$ or $B$ is a von Neumann algebra, then both algebras are von Neumann algebras and $\Delta$ admits an extension to a real-linear isometry from $A$ onto $B$. More precisely, there exists a central projection $p \in B$ and a Jordan $^*$-isomorphism $J : A\to B$ such that the conclusion in \eqref{eq 2608} holds when $e^{ix}$ is replaced with an arbitrary unitary $u$ in $A$ (cf. Theorem~\ref{t Hatori Molnar surjecitve isometries between unitaries vN}).\smallskip

We can therefore arrive to the striking conclusion that two unital C$^*$-algebras are isomorphic as Jordan $^*$-algebras if and only if their unitary groups are isometric as metric spaces (see Corollary~\ref{c Hatori Molnar isometrically isomorphic}).\smallskip

We shall also study in this note the set of unitary elements in a unital JB$^*$-algebra $M$ as a metric invariant with respect to the metric given by the norm. Subsections~\ref{subsec: connected components of the unitay set}, \ref{subsect: isometries between connected components} and \ref{subsec: when is a surjective isometry extendable} contain the results in this direction obtained by Cueto, Enami, Hirota, Miura and the author of this note in \cite{CuEnHirMiuPe2022}.\smallskip

As in the case of invertible elements, the Jordan product of two unitaries in a unital JB$^*$-algebra $M$ need not be, in general, a unitary. The key property in the setting of unital JB$^*$-algebras assures that given $u,v\in \mathcal{U} (M)$, the element $U_{u} (v)$ lies in $\mathcal{U} (M)$ (cf. Lemma~\ref{l first properties invertible Jordan}).  So, the key hypotheses in this setting are given in terms of expressions of the form $U_{u} (v)$ and quadratic subsets. As in the case of C$^*$-algebras, the set $\mathcal{U} (M)$ is not, in general, connected. The connected component of $ \mathcal{U} (M)$ containing the unit element is called the \emph{principal connected component} or simply the \emph{principal component} and will be denoted by $\mathcal{U}^{0}(M)$.\smallskip

The principal component of the unitary elements in a unital C$^*$-algebra $A$ is a well studied object, whose description can be found in books (cf. \cite[Exercises 4.6.6 and 4.6.7]{KR1}, \cite[Exercise 2, page 56]{Tak}, and \cite[Lemma 3.2]{Hatori14}. It is known that $$ \begin{aligned}
	\mathcal{U}^0(A) &= \{ e^{i h_1} \cdot \ldots \cdot e^{i h_n} : n\in \mathbb{N}, \ h_1, \ldots, h_n\in A_{sa} \}\\ &= \{ e^{i h_1} \cdot \ldots \cdot e^{i h_n} \11b{A} e^{i h_n} \cdot \ldots \cdot e^{i h_1} : n\in \mathbb{N}, \ h_1, \ldots, h_n\in A_{sa} \}
\end{aligned}$$ is open, closed, and path connected in the (relative) norm topology on $\mathcal{U}(A)$. Furthermore, by the continuity of the module mapping $x\mapsto |x|:=(x^* x)^\frac12$, we have $$A^{-1}_{\mathbf{1}} \cap \mathcal{U}(A) =\mathcal{U}^0(A),$$ where $A^{-1}_{\mathbf{1}}$ stands for the principal component of the invertible elements in $A$ (see Theorem~\ref{t principal connected component in unital C$^*$-algebras}).\smallskip

An algebraic characterization of the principal component of the set of unitaries in a unital JB$^*$-algebra has not been available until the recent reference \cite{CuEnHirMiuPe2022}. Here the result is presented in Theorem~\ref{t principal component as product of exp} and Proposition~\ref{p self adjoint quadratic subset princ comp}, where it is shown that for each unital JB$^*$-algebra $M$ with principal component of $\mathcal{U}(M)$ denoted by $\mathcal{U}^0(M)$, the following statements are equivalent for every unitary $u$ in $M$:
\begin{enumerate}[$(a)$]\item $u\in M^{-1}_{\mathbf{1}}\cap\mathcal{U} (M)$;
	\item There exists a continuous path $\Gamma : [0,1]\to \mathcal{U} (M)$ with $\Gamma (0) = \11b{M}$ and $\Gamma(1) = u$;
	\item $u\in \mathcal{U}^0(M)$;
	\item $u = U_{e^{i h_n}} \cdots U_{e^{i h_1}}(\11b{M}),$ for some $n\in \mathbb{N}$, $h_1,\ldots,h_n\in M_{sa}$;
	\item There exists $w\in \mathcal{U}^0(M)$ such that $\|u-w\|<2$.
\end{enumerate}
Consequently,
\begin{equation}\label{eq algebraic charact principal component unitaries JBstar}\begin{aligned} \mathcal{U}^0(M) &=  M^{-1}_{\mathbf{1}}\cap\mathcal{U} (M)\\
		&=\left\lbrace  U_{e^{i h_n}}\cdots U_{e^{i h_1}}(\11b{M}) \colon n\in \mathbb{N}, \ h_j\in M_{sa} \ \forall\ 1\leq j \leq n  \right\rbrace \\
		&= \left\lbrace  u\in \mathcal{U} (M) : \hbox{ there exists } w\in \mathcal{U}^0(M) \hbox{ with } \|u-w\|<2 \right\rbrace
\end{aligned}\end{equation} is analytically arc-wise connected. Moreover, $\mathcal{U}^{0} (M)$ is a self-adjoint quadratic subset of $\mathcal{U} (M)$. $U_{\mathcal{U}^{0} (M)}$ is precisely the smallest quadratic subset of $\mathcal{U} (M)$ containing the set $e^{i M_{sa}}$. \smallskip

After a series of technical results, we present a characterization of all surjective isometries between the principal components of the unitary sets in two unital JB$^*$-algebras in Theorem~\ref{t surj isom principal components}. Suppose $\Delta: \mathcal{U}^0 (M)\to \mathcal{U}^0 (N)$ is a surjective isometry between the principal components of two unital JB$^*$-algebras. Then there exist self-adjoint elements $k_1,\ldots,k_n$ in $N$, a central projection $p$ in $N,$ and a Jordan $^*$-isomorphism $\Phi:M\to N$ such that $$\begin{aligned}\Delta(u) &= p\circ U_{e^{i k_n}} \ldots  U_{e^{i k_1}}  \Phi (u) + (\11b{N}-p)\circ \left( U_{e^{-i k_n}} \ldots  U_{e^{-i k_1}} \Phi(u)\right)^*,
\end{aligned}$$ for all $u\in \mathcal{U}^0(M)$. Moreover, the elements $k_1,\ldots,k_n\in N_{sa}$ can be replaced by any finite collection of self-adjoint elements $\tilde{k}_1,\ldots,\tilde{k}_m$ such that $ \Delta(\11b{M}) = U_{e^{i \tilde{k}_m}} \ldots  U_{e^{i \tilde{k}_1}} (\11b{N})$. Under these hypotheses $M$ and $N$ are Jordan $^*$-isomorphic, and there exists a surjective real-linear isometry {\rm(}i.e., a real-linear triple isomorphism{\rm)} from $M$ onto $N$ whose restriction to $\mathcal{U}^0 (M)$ is $\Delta$.\smallskip

If under the previous hypotheses, we also assume that there exists $v\in \mathcal{U} (N)$ such that $U_v (\Delta(\11b{M})) = \11b{N}$. Then there exist a central projection $p$ in $N$ and a Jordan $^*$-isomorphism $\Phi:M\to N$ such that $$\begin{aligned}\Delta(w) &= p\circ U_{v^*}  \Phi (w) + (\11b{N}-p)\circ \left( U_{v} \Phi(w)\right)^*,
\end{aligned}$$ for all $w\in \mathcal{U}^0(M)$ (cf. Corollary~\ref{c surj isom principal components simplified}).\smallskip

In Theorem~\ref{t HM for unitary sets of unitary JBW*-algebra to a JB*-algebra} we obtain an strengthened  version of a result from the recent reference \cite{CuEnHirMiuPe2022}, showing that if $M$ and $N$ are two unital JB$^*$-algebras such that $M$ or $N$ is a JBW$^*$-algebra, and $\Delta: \mathcal{U} (M)\to \mathcal{U} (N)$ is a surjective isometry, then $M$ and $N$ are JBW$^*$-algebras and there exist a self-adjoint element $k_1$ in $N,$ a central projection $q\in N$, and a Jordan $^*$-isomorphism $\Phi:M\to N$ such that $$\begin{aligned}\Delta( u ) &= q \circ  U_{e^{i k_1}}\left(  \Phi(u) \right) + (\11b{N}-q)  \circ U_{e^{i k_1}}\left( \Phi( u)^*\right),\\
\end{aligned}$$ for all $u\in \mathcal{U}(M)$. If $\Delta$ is unital we can take $k_1=0$.  Consequently, $\Delta$ admits a (unique) extension to a surjective real-linear isometry from $M$ onto $N$. A version of this result assuming that $M$ and $N$ are von Neumann algebras was previously obtained by Hatori and Moln{\'a}r in \cite{HatMol2014}. \smallskip   

The consequences of the results in subsection~\ref{subsec: connected components of the unitay set} are summarized in Corollary~\ref{c two unital JB*algebras are isomorphic iff their unitaries are isometric iff their principal components are isometric}. To present the conclusion, let $M$ and $N$ be unital JB$^*$-algebras. Consider the following statements 
\begin{enumerate}[$(1)$]
	\item $M$ and $N$ are isometrically isomorphic as (complex) Banach spaces. 
	\item $M$ and $N$ are isometrically isomorphic as real Banach spaces. 
	\item $\mathcal{U}(M)$ and $\mathcal{U}(N)$ are isometrically isomorphic as metric spaces.
\end{enumerate}
\begin{enumerate}[$(a)$]
	\item $M$ and $N$ are Jordan $^*$-isomorphic;
	\item There exists a surjective isometry $\Delta: \mathcal{U}(M)\to \mathcal{U}(N)$ satisfying $\Delta (\11b{M})\in \mathcal{U}^0(N)$;
	\item There exists a surjective isometry $\Delta: \mathcal{U}^0(M)\to \mathcal{U}^0(N).$
\end{enumerate} Then $(1)\Leftrightarrow (2) \Leftrightarrow (3) \Leftarrow (a) \Leftrightarrow (b) \Leftrightarrow (c)$. The implication $(3)\Rightarrow (a)$ holds when $M$ and $N$ are unital C$^*$-algebras (cf. Corollary~\ref{c two unital C*algebras are isomorphic iff their unitaries are isometric iff their principal components are isometric}). We shall see in Corollary~\ref{r existence of unitaries with non isomorphic nor isometric connected components in the Jordan setting} that the implication $(3)\Rightarrow (a)$ is, in general, false for unital JB$^*$-algebras. \smallskip

We shall see in Example~\ref{example existence of non-extendible surjective isometries between unitaries} that a surjective isometry between the sets of unitary elements in two unital JB$^*$-algebras is not always extendable to a surjective real-linear isometry between the algebras. A complete characterization of those surjective isometries between sets of unitary elements in two unital JB$^*$-algebras admitting an extension to a surjective real-linear isometry is obtained in Corollary~\ref{c extendibility of surjective isometries}.\smallskip

The concluding result in this note gives a sufficient condition to guarantee that a surjective isometry between the unitary elements in two unital JB$^*$-algebras $M$ and $N$ is extendable to a surjective real-linear isometry. Indeed, let $\Delta :  \mathcal{U} (M)\to \mathcal{U} (N)$ be a surjective isometry and let $\{\mathcal{U}^0 (M)\}\cup \{ \mathcal{U}^{j} (M) : j\in \Lambda\}$ be the collection of all different connected components of $\mathcal{U} (M)$. Assume that for each $j\in \Lambda$ there exists $u_j \in \mathcal{U} (M)$ with $u_j^2 \in \mathcal{U}^{j} (M)$. Suppose additionally that \begin{equation}\label{eq preservation of quadratic products extra hypo} \Delta\{u,v,u\}= \Delta U_u (v^*) = U_{\Delta(u)} \left(\Delta(v)^* \right) =\{\Delta(u), \Delta(v), \Delta(u)\},
\end{equation} for all $u,v\in \mathcal{U} (M)$. Then $\Delta$ admits an extension to a surjective {\rm(}real-linear{\rm)} isometry from $M$ onto $N$ (see Proposition~\ref{p extra algebraic hypotheses}). 

\section{A taste of radicals in Banach and Jordan--Banach algebras}\label{sec: taste of radicals in Jordan}

We devote this section to recall the basic notions and results for Jordan--Banach algebras and the McCrimmon radical. Our aim is to present the main results needed in later sections following a parallelism between associative Banach and Jordan--Banach algebras.  Perhaps the notion of Jordan-Banach algebra has already reached to a wide enough audience, after an almost a century of developments,  for not needing a detailed presentation. Anyway, in order to make this paper accessible to newcomers (or reader who are less familiarized with the theory of Jordan algebras) we shall equip the reader with full information and bibliographic sources. For a full basic-background on Jordan--Banach algebras and any result or notation not explicitly presented here, the reader is referred to \cite{McCrimm66,McCrimm69,McCrimm71,HogbMcCrimm81, Jacobson81, HOS, Aupetit93, AlfsenShultz2003} and the most updated monograph \cite{Cabrera-Rodriguez-vol1}.\smallskip

By quoting a phrase from Jacobson's book \cite{Jacobson81}, ``The Jordan theory had its birth in an attempt by Jordan and subsequently by Jordan, von Neumann and Wigner to formulate the foundation of quantum mechanics''. These renowned mathematicians gave a complete determination of the finite-dimensional Jordan algebras over the reals which are formally real \cite{JorvonNeuWigner34}.\smallskip 

A \emph{Jordan algebra} over $\mathbb{K} = \mathbb{R}$ or $\mathbb{C}$, is a linear space $M$ over $\mathbb{K}$ equipped with a bilinear product $M\times M\to M$, $(a,b)\mapsto a\circ b$ ($a,b\in M$)  satisfying the following two axioms \begin{enumerate}[$(1)$]
\item (commutativity) $a\circ b = b \circ a$ ($a,b\in M$).
\item (\emph{Jordan identity}) $(a \circ b)\circ a^2 = a\circ (b \circ a^2)$ ($a,b\in M$). 	
\end{enumerate} By replacing the usual  associativity of the product by the weaker assumption in the Jordan identity we get a wider class of algebras. A Jordan algebra $M$ is called \emph{unital} if there exists an element $\mathbf{1} = \11b{M}\in M$ (called the unit of $M$) satisfying $\11b{M}\circ a= a$ for all $a\in M$. A subalgebra $N$ of $M$ is said to be a \emph{strongly associative subalgebra} if $$ (b\circ x)\circ a - b\circ (x\circ a) = 0,\hbox{ for all } a,b \in N \hbox{ and } x  \in M.$$
A \emph{normed Jordan algebra} is a Jordan algebra $M$ equipped with a norm, $\|.\|$, satisfying $\|a\circ b\| \leq \|a\| \ \|b\|$ ($a,b\in M$). A \emph{Jordan-Banach algebra} is a normed Jordan algebra $M$ whose norm is complete. Unless otherwise stated, all Banach spaces and algebras in this paper are over the complex field. If $M$ is a unital Jordan--Banach algebra with unit $\11b{M}$, we shall also require that $\|\11b{M}\|=1$. When $M$ lacks of a unit, it can be always embedded in a unital Jordan-Banach algebra.\smallskip

Along this paper, we shall employ the following notation for the powers of an element $a$ in a Jordan algebra $M$:
$$ a^1 =a;\ \  a^{n+1} =a\circ a^n, \quad n\geq 1.$$ If $M$ is unital, we set $a^0 =\11b{M}$. The Jordan identity assures that each Jordan algebra $M$ is \emph{power associative}, that is, each subalgebra generated by a single element is associative,  equivalently, $ a^m\circ a^n=a^{m+n},$ for all $a\in M$, $m,n\in \mathbb{N}$ (see \cite[Lemma 2.4.5]{HOS}).
\smallskip

From a purely algebraic point of view, elements $a, b$ in a Jordan algebra $M$ are said to \emph{operator commute} if $$(a\circ c)\circ b= a\circ (c\circ b),$$ for all $c\in M$ (cf. \cite[4.2.4]{HOS}). The \emph{centre} of $M$ is, by definition, the set of all elements of $M$ which operator commute with any other element in $M$, and will be denoted by $Z(M)$. Elements in the centre of $M$ are called \emph{central}.\smallskip

By a \emph{Banach algebra} we mean an associative Banach algebra. There is a natural way to regard every real or complex Banach algebra as a real or complex Jordan-Banach algebra by just considering the natural Jordan product defined by $a\circ b: = \frac12 (a b +ba)$. Norm closed Jordan subalgebras of associative Banach algebras are called \emph{special}. Jordan algebras which are not special will be called exceptional. Examples of the latter can be found in \cite[Chapter IX]{Jacobson81}. So, the class of Jordan algebras is strictly wider than the mere collection of all special Jordan algebras (more concrete examples will appear later). \smallskip   

Along this paper, we shall assume that our Jordan-Banach algebras are unital. A (real-linear) Jordan homomorphism between Jordan Banach algebras is a (real-linear) linear map preserving Jordan products (equivalently, squares of elements). \smallskip

Invertibility in Jordan--Banach algebras requires some extra attention.  We recall that an element $a$ in a unital Jordan--Banach algebra $M$ is said to be \emph{invertible} if there exists $b\in M$ satisfying $a \circ b = \11b{M}$ and $a^2 \circ b = a.$ In such a case, the element $b$ is unique, it is called the inverse of $a$, and will be denoted by $a^{-1}$ (cf. \cite[3.2.9]{HOS} and \cite[Definition 4.1.2]{Cabrera-Rodriguez-vol1}). The symbol $M^{-1}$ will stand for the set of all invertible elements in $M$.\smallskip

In an associative Banach algebra $A$, each element $a$ gives rise to the left and right multiplication operators defined by $L_a (x) =  a x$ and $R_a (x) = x a$ ($x\in A$). When $A$ is unital, the invertibility of $a$ is equivalent to the invertibility of $L_a$ or $R_a$. In a Jordan Banach algebra $M$ left and right multiplication operators are essentially indistinguishable. The (Jordan) multiplication operator by an element $a\in M$ (denoted by $M_a$) is defined by $M_a (x)= a \circ x$ ($x\in M$), while the $U_a$ operator is given by $$U_{a} (x) = (2 L_a^2 -L_{a^2}) (x) = 2 (a\circ x) \circ a  - a^2\circ x \ \  (x\in M).$$ Furthermore, given $a,b\in M$, we shall write $U_{a,b}:M\to M$ for the linear mapping on $M$ defined by $$U_{a,b} (x)=(a\circ x) \circ b + (b\circ x)\circ a - (a\circ b)\circ x,\ \ (x\in M).$$ The "fundamental formula" in Jordan algebras asserts that \begin{equation}\label{eq fundamental identity UaUbUa} \ \  U_a U_b U_a = U_{U_a(b)}, \end{equation} for all $a,b$ in $M$ (cf. \cite[2.4.18]{HOS} or \cite{McCrimm66} or \cite[Proposition 3.4.15]{Cabrera-Rodriguez-vol1}). We shall also employ the following version of the fundamental formula, which can be easily obtained via induction, \begin{equation}\label{eq fundamental identity n elements} \ \   U_{U_{a_1} U_{a_2} \ldots U_{a_n} (b)} =U_{a_1} U_{a_2}\ldots  U_{a_n} U_b U_{a_n} \ldots U_{a_2} U_{a_1}, \end{equation} for all $n\in \mathbb{N}$, $a_1,\ldots, a_n,b$ in $M$.\smallskip

We gather in the next lemma some starting properties of invertible elements in unital Jordan--Banach algebras. From now on we shall write $L(X)$ for the space of all linear maps on a linear space $X$.

\begin{lemma}\label{l first properties invertible Jordan}\cite[Theorem 4.1.3, Proposition 4.1.6, and Theorem 4.1.7]{Cabrera-Rodriguez-vol1} Let $M$ be a Jordan--Banach algebra. Then the following statements hold:
\begin{enumerate}[$(a)$]
\item An element $a\in M$ is invertible if and only if $U_a$ is a bijective mapping (i.e. invertible in $L(M)$), and in such a case $U_a^{-1} = U_{a^{-1}}$.
\item The set $M^{-1}$ of all invertible elements in $M$ is open, and consequently $M^{-1}$ is locally connected and all its connected components are open.
\item An element of the form $U_x (y)$ is invertible if and only if both $x$ and $y$ belong to $M^{-1}$, and such a case, $(U_x (y))^{-1} = U_x^{-1} (y^{-1})= U_{x^{-1}} (y^{-1})$. 
\item If $a\in M^{-1}$, then $a^n$ and $(a^{-1})^m$ operator commute for all $n,m\in \mathbb{N}$. Furthermore, by defining $a^{-n} := (a^{-1})^n$ for all natural $n$ and $a^0 := \11b{M}$, then $a^k\circ a^l = a^{k+l}$ for all integers $k,l$.
\item If $A$ is a unital complex Banach algebra regarded as a Jordan--Banach algebra with respect to the natural Jordan product, and element $a\in A$ is invertible in the Jordan sense if and only if it is invertible in the usual sense, and in such a case the inverse is the same in both senses.
\item The mapping $M^{-1}\to M^{-1}$, $a\mapsto a^{-1}$ is a homeomorphism.
\end{enumerate}
\end{lemma}

The reader should be warned that the set invertible elements in a Jordan--Banach algebra is not, in general, stable under Jordan products. \smallskip

As in the setting of Banach algebras, the spectrum of an element $a$ of a complex Jordan--Banach algebra $M$, $Sp(a),$ $Sp(M,a)$ or J-$\sigma (a)$, is defined as the set of all $\lambda\in \mathbb{C}$ for which $a-\lambda \11b{M}$ is not invertible. An appropriate Gelfand theorem for Jordan--Banach algebras assures that $Sp(a)$ is a non-empty compact subset of the complex plane (see, for example,  \cite[Theorem 4.1.17]{Cabrera-Rodriguez-vol1} where even a Jordan version of the celebrated Gelfand-Beurling formula is stated). The multi-valued mapping $a\mapsto Sp(a)$ is upper semi-continuous on $M$. We shall write  $\rho(a)$ for the \emph{spectral radius} of $a$, that is, $\rho(a)= \max\{|\lambda|: \lambda \in Sp(a)\}.$ It should be noted that if $A$ is a complex Banach algebra, the spectrum of each element in $A$ coincides with its spectrum when $A$ is regarded as a unital complex Jordan--Banach algebra.  \smallskip

From now on, given a Banach space $X$, the open (respectively, closed) ball of radius $\rho$ and center $a\in X$ will be denoted by $B_{_X}(a,\rho)$ or by $B(a,\rho)$ (respectively, $\overline{B}_{_X}(a,\rho)$). \smallskip

Continuous and holomorphic functional calculus for elements in a complex Jordan--Banach algebra $M$ are relatively simple and natural to understand. Indeed, the closed subalgebra, $M_a$, of a unital Jordan-Banach algebra $M$ generated by $\11b{M}$ and an element $a\in M$ is an associative and commutative Banach algebra (cf. \cite[Theorem 2.3]{Aupetit95} and \cite[Theorems 4.1.88 and 4.1.93]{Cabrera-Rodriguez-vol1}). Therefore, the standard holomorphic functional calculus in complex Banach algebras can be perfectly understood in the setting of complex Jordan-Banach algebras. We gather in the next result some of the basic properties of this holomorphic functional calculus.

\begin{theorem}[Holomorphic functional calculus]\label{t holomorphic functional calculus} {\rm(\cite[Theorem 2.3]{Aupetit95} and \cite[Theorems 4.1.88 and 4.1.93]{Cabrera-Rodriguez-vol1})} Let $a$ be an element in a unital complex Jordan-Banach algebra $M$. Suppose $\Omega$ is an open neighbourhood of Sp$(a)$, and let $h_0$ denote the function defined by the inclusion mapping $\Omega\hookrightarrow \mathbb{C}$. Then there is a unique continuous unit-preserving algebra homomorphism $f \mapsto f (a)$ from the complex algebra $\mathcal{H} (\Omega),$ of all complex-valued holomorphic functions on $\Omega,$ into $M$ taking $h_0$ to $a$ and satisfying the following properties:\begin{enumerate}[$(i)$]\item $\displaystyle f(a) = \frac{1}{2 \pi i} \int_{\Gamma} f(\lambda) (\lambda \11b{M}-a)^{-1} d\lambda$, for any positively orientated curve included in $\Omega$ and surrounding Sp$(a)$;
		\item $f(a)$ is contained in the smallest closed strongly associative subalgebra of $M$ containing $\11b{M}$ and $a$;
		\item For each $f\in \mathcal{H} (\Omega)$ we have
		$\hbox{Sp}(M, f (a))$ $=$ $f (\hbox{Sp}(M,a))$ \ \  (Spectral mapping theorem);
		\item The set $$M_{\Omega} : =\{ x\in M : \hbox{Sp} (M,x) \subseteq \Omega\}$$ is a non-empty open subset of $M$, and the mapping  $\widetilde{f} : x \mapsto f (x)$ from $M_{\Omega}$ to $M$ is holomorphic {\rm(}i.e. there exists a bounded linear operator $T: M\to M$ such that $\displaystyle \lim_{h\to 0} \frac{\|f(x + h) - f(x)- T(h)\|}{\|h\|} = 0$ for any $x\in M_{\Omega}${\rm)}.
		\item Let $\Omega_1$ be another open subset in $\mathbb{C}$, $f\in\mathcal{H} (\Omega)$, $g\in \mathcal{H} (\Omega_1)$ such that $f(\Omega)\subseteq \Omega_1$. Then $(g\circ f ) (a) = g(f(a))$.
	\end{enumerate}
\end{theorem}

\subsection{Jacobson and McCrimmon Radicals}\label{subsec Jacobson and McCrimmon radicals}\ \smallskip

We begin this subsection by recalling the notion of Jacobson radical for Banach algebras (see \cite[\S 12]{Zelaz73}). Given a Banach algebra $A$ and a complex vector space $X$, a \emph{representation} of $A$ in the space $X$ is a homomorphism $\pi: A\to L(X)$, $x\mapsto  \pi(x)$. A subspace $X_0$ of $X$ is called \emph{invariant} for the representation $\pi$ if $\pi(x) (X_0)\subseteq X_0$ for all $x \in A$.  A representation is called irreducible if the only invariant subspaces are $\{0\}$ and $X$. A two-sided ideal $I\subseteq A$ is said to be a \emph{left-} (respectively, \emph{right-} or \emph{modular}) \emph{modular ideal} when there exists an element $e\in A$ such that $x-xe\in I$ (respectively, $x-e x \in I$ or $x- x e, x- e x\in I$) for any $x\in A$. Finally, $I$ is called \emph{primitive} if it is of the form $I = \{x\in A : x A \subseteq J\}$, where $J$ is a maximal left-modular ideal in $A$. A celebrated result by Rickart shows that a two-sided ideal is primitive if and only if it is the kernel of an irreducible representation. Moreover, every maximal two-sided modular ideal is primitive. If $A$ is commutative, primitive ideals are precisely the maximal modular ideals (see \cite[Chapter II]{RickartBook}).\smallskip 

Assuming that $A$ is non-unital, an element $x$ in $A$ is called quasi-regular if there exists an element $y \in A$ (called the quasi-inverse of $x$) such that $x+ y -x y  =x + y - y x = 0$, i.e., $\11b{A}-x$ is invertible in the unitization of $A$ with inverse $\11b{A}-y$.

\begin{definition}\label{def Jacobson Radical}(\cite[\S 12]{Zelaz73}, \cite[Chapter II]{RickartBook} or \cite[\S 4.3]{PalmerBook}) The Jacobson radical of any Banach algebra $A$, rad$(A)$, is defined as the intersection of all primitive ideals of $A$. If $A\neq  \hbox{rad} (A)$, the Jacobson radical coincides with the intersection of all maximal left-modular ideals, or, equivalently, with the intersection of all right-modular ideals. Furthermore 
	$$\begin{aligned}
	\hbox{rad}(A) &= \{x\in A : \forall \lambda\in \mathbb{C} \hbox{ the element } \lambda x + y x \hbox{ is quasi-regular } \forall y\in A\} \\
	&= \{x \in A : \forall \lambda \in \mathbb{C} \hbox{ the element } \lambda x +x y \hbox{ is quasi-regular } \forall y \in A\}.
	\end{aligned}$$ We say that $A$ is semisimple if rad$(A)=\{0\}$. 
\end{definition}

An element $a$ in a Banach algebra $A$ is called \emph{quasi-nilpotent} if its spectral radius, $\rho (a),$ is zero. It is known that the Jacobson radical of $A$ is contained in the set of all quasi-nilpotents elements of $A$, however the Jacobson radical is, in general, smaller. More analytic characterizations of the Jacobson radical in terms of the spectral radius were obtained by  S{\l}odkowski, Wojty{\'n}ski, and Zemz{\'a}nek, and subsequently extended by Aupetit. 

\begin{theorem}\label{t characterization of the Jacobson radical} {\rm(\cite[Theorem 2.3.3]{PalmerBook} and \cite[Theorem 1]{Aupetit93})} For each complex Banach algebra $A$ we have  $$\begin{aligned} \hbox{{\rm rad}}(A) &= \{x\in A : \rho(a + x) =\rho(a),\  \forall a\in A\} = \{x \in A : \rho(a x) = 0,\ \forall a \in A\} \\
&=  \{x \in A : \sigma(a + x) = \sigma(a) ,\ \forall a \in A\}\\
& = \left\{ x\in A : \hbox{ for each } a\in A \hbox{ the set  } \{\rho(a+\lambda x) : \lambda\in \mathbb{C}\} \hbox{ is bounded}\right\}.
	\end{aligned}$$
\end{theorem}

We know, by an observation due to Hatori \cite{Hat2011}, that the second  characterization of the radical presented in Theorem~\ref{t characterization of the Jacobson radical} admits a weaker version in case that $A$ is unital. We shall present here an alternative proof to that given by Hatori. Let us write  $A_{\mathbf{1}}^{-1}$ for the  principal connected component of $A^{-1}$. 

\begin{lemma}\label{l Hatori associative radical}{\rm\cite[Lemma 3.1]{Hat2011}} Let $A$ be a unital (associative) complex Banach algebra. Then $$\begin{aligned} \hbox{{\rm rad}}(A) &= \left\{x \in A : \rho(a x) = 0,\ \forall a \in A^{-1}\right\} = \left\{x\in A : \rho(a x) =0,\  \forall a\in A_{\mathbf{1}}^{-1}\right\} .
	\end{aligned}$$
\end{lemma} 

\begin{proof} The containing $\subseteq$ is clearly true in both cases. So, we only need to show that the set on the right hand side is inside the Jacobson radical. Suppose that $\rho(a x) =0,$  for all $a\in A_{\mathbf{1}}^{-1}$. Fix an arbitrary $b\in A$ and $\mu \in \mathbb{C}$ with $|\mu| > \rho (b)$. Clearly, $\11b{A}-\mu^{-1} b$ is invertible, and since $\lim_{|\mu|\to \infty} (\11b{A}-\mu^{-1} b )^{-1}=\11b{A}$, we deduce that $(\11b{A}-\mu^{-1} b )^{-1} \in A_{\mathbf{1}}^{-1}$, which implies that $\beta (\11b{A}-\mu^{-1} b )^{-1} \in A_{\mathbf{1}}^{-1}$ for all $\beta \in \mathbb{C}\backslash\{0\}$.  It follows from the assumptions that $\rho\left( \beta (\11b{A}-\mu^{-1} b )^{-1} x \right) =0$, for all $\beta \in \mathbb{C}\backslash\{0\}$. \smallskip

By applying the previous conclusion we get $$\begin{aligned}
\sigma \left((\11b{A}-\mu^{-1} b )^{-1} (\11b{A} -\mu^{-1} b- \beta x)  \right) &= \sigma (\11b{A} - \beta  (\11b{A} -\mu^{-1} b )^{-1} x) \\ 
&= 1- \sigma (\beta (\11b{A}-\mu^{-1} b )^{-1} x )= \{1\},
	\end{aligned}$$ and thus $(\11b{A} -\mu^{-1} b- \beta x)$ must be invertible for all $\beta\in \mathbb{C}\backslash\{0\}$. However, it is clear that the also have an invertible element for $\beta =0$. This is equivalent to say that $-\mu \11b{A} + b +\lambda x) = -\mu ((\11b{A} -\mu^{-1} b- \lambda \mu^{-1} x))\in A^{-1},$ equivalently, $\mu \notin \sigma (b +\lambda x),$ for all $\lambda \in \mathbb{C}$. We have shown that $$\mathbb{C}\backslash B_{\mathbb{C}} (0,\rho(b) +\varepsilon) \subseteq \mathbb{C}\backslash \sigma(b +\lambda x),$$ for all $\lambda\in \mathbb{C}$ and $\varepsilon>0$, and consequently, $\rho (b +\lambda x) \leq \rho (b),$ for all $\lambda\in \mathbb{C}$. The arbitrariness of $b$ and the last characterization in Theorem~\ref{t characterization of the Jacobson radical} prove that $x\in \hbox{rad} (A)$.   
\end{proof}

We turn now to Jordan--Banach algebras. An element $a$ in a unital Jordan-Banach algebra $M$ is \emph{quasi-invertible} with \emph{quasi-inverse} $b \in M$ if $\11b{M} -a$ is invertible with inverse $\11b{M} - b$. A subset of $M$ is called quasi-invertible if all its elements are quasi-invertible in $M$.\smallskip

An \emph{outer ideal} of a Jordan-Banach algebra $M$ is a subspace $J$ such that $U_{M} (J)\subseteq J$; an \emph{inner ideal} or a \emph{strict inner ideal} is a subspace $I$ satisfying $U_{I} (M) \subseteq I$ and $I^2 \subseteq I$ (if $M$ is unital the second condition clearly follows from the first one). An ideal $J$ is a subspace that is both an outer and a strict inner ideal (i.e. $J^2 \subseteq J$). It is known that if $J$ is an ideal or an inner ideal, then the quasi-inverse of any quasi-invertible element in $J$ also belongs to $J$ (see \cite[page 672]{McCrimm69}).\smallskip

The existence of an appropriate notion of radical in the setting of Jordan--Banach algebras was proved by McCrimmon \cite{McCrimm69}.

\begin{theorem}\label{t existence of McCrimmon radical}{\rm(\cite[Theorem 1]{McCrimm69}, \cite[Theorem 1.1]{HogbMcCrimm81}, \cite[Proposition 4.4.11 and Definition 4.4.12]{Cabrera-Rodriguez-vol1})}  Let $M$ be a Jordan-Banach algebra $M$. Then there exists a unique maximum quasi-invertible ideal $\mathcal{R}ad(M)$ of $M,$ called the \emph{McCrimmon radical} of $M$, which contains all quasi-invertible ideals. Actually the McCrimmon radical of a unital Jordan--Banach algebra $M$ coincides with the intersection of all maximal inner ideals in $M$.
\end{theorem}

A Jordan-Banach algebra $M$ is called \emph{Jordan-semisimple}, \emph{J-semisimple}, or simply \emph{semisimple} if $\mathcal{R}ad(M)= \{0\}$.

\begin{remark}\label{r Jordan radical and McCrimmon radical coincide when it makes sense}\cite[Theorem 3.6.21]{Cabrera-Rodriguez-vol1}  Let $A$ be a Banach algebra, and let $M$ denote the Jordan-Banach algebra obtained when $A$ is equipped with its natural Jordan product. Then the McCrimmon radical of $M$ coincides with the usual Jacobson radical of $A$.
\end{remark}

Our next result, originally established by Aupetit in \cite{Aupetit93}, is a Jordan version of the characterizations of the radical in terms of the spectral radius presented in Theorem~\ref{t characterization of the Jacobson radical}. 

\begin{theorem}\label{t Aupetit c 1}{\rm \cite[Theorem 2, Corollaries 1 and 2]{Aupetit93}} Let $x$ be an element in a Jordan-Banach algebra $M$. Then the following statements are equivalent:\begin{enumerate}[$(a)$]\item $x$ is in the McCrimmon radical of $M$;
\item For every $a\in M$ the set  $\{\rho(a+\lambda x) : \lambda\in \mathbb{C}\}$ is bounded;
\item $\rho(U_a(x)) = 0,$ for every $a$ in $M$;
\item There exists $C\geq 0$ such that $\rho (a) \leq C \|a-x\|,$ for all $a$ in a neighbourhood of $x$.
\end{enumerate}
\end{theorem}

Our next goal is a Jordan version of  Lemma~\ref{l Hatori associative radical}, which improves the characterization established by Aupetit in the previous theorem. We also include a slightly different proof from the respect to the original source. According to the usual notation, we shall write $M_{\mathbf{1}}^{-1}$ for the principal connected component of the invertible elements in a unital Jordan--Banach algebra $M$.  

\begin{proposition}\label{p Hatori McCrimmon radical}{\rm\cite[Proposition 2.3]{Pe2021}} Let $M$ be a unital complex Jordan--Banach algebra. Then $$\begin{aligned} \mathcal{R}ad(M) &= \left\{x \in M : \rho(U_a (x)) = 0,\ \forall a \in M^{-1}\right\} \\
		&= \left\{x\in M : \rho(U_a (x))  =0,\  \forall a\in M_{\mathbf{1}}^{-1}\right\} .
	\end{aligned}$$
\end{proposition} 

\begin{proof} As in the associative case, it suffices to prove that the last set is contained in the McCrimmon radical of $M$. Assume that $x\in M$ satisfies $\rho(U_a (x))  =0$,  for all $a\in M_{\mathbf{1}}^{-1}$.\smallskip
	
Fix an arbitrary $b\in M$ and $\mu \in \mathbb{C}$ with $|\mu|>\rho (b).$ In such a case $b - \mu \11b{M}\in M^{-1},$ equivalently, $-\mu^{-1} b +\11b{M} \in M^{-1}.$  \smallskip

Let us make some observations on the spectrum of the elements of the form $-\mu^{-1} b +  \11b{M}\in M_{\mathbf{1}}^{-1}.$ Clearly $$\hbox{Sp}\left(- \mu^{-1} b +  \11b{M}\right) =1- \mu^{-1} \hbox{Sp}(y) \subseteq 1 - {B}_{\mathbb{C}} (0,1) \subseteq \{z\in \mathbb{C} : \Re\hbox{e}(z) >0\}= \Omega.$$ We are in a position to apply the holomorphic functional calculus (cf. Theorem~\ref{t holomorphic functional calculus}), via the principal branch of the square root $f: \Omega \to \mathbb{C},$ $f(\lambda) = \frac{1}{\sqrt{\lambda}}$, to define the element $\left(\11b{M} - \mu^{-1} b \right)^{-\frac12} : = f \left(\11b{M} -  \mu^{-1} b \right)$ ($|\mu|>\rho (b)$).  By the continuity of the holomorphic functional calculus, $\lim_{|\mu|\to \infty} \left(\11b{M} - \mu^{-1} b \right)^{-\frac12} = \11b{M}$, and consequently $\beta \left(\11b{M} - \mu^{-1} b \right)^{-\frac12} \in M_{\11b{M}}^{-1},$ for all $\beta \in \mathbb{C}\backslash\{0\}$.\smallskip

By the assumptions $\rho \left(U_{\beta \left(\11b{M} - \mu^{-1} b \right)^{-\frac12}} (x)\right)=0,$ for all non-zero $\beta\in \mathbb{C}$, $|\mu|>\rho (b)$). Since $\hbox{Sp} \left(U_{\beta \left(\11b{M} - \mu^{-1} b \right)^{-\frac12}} (x)\right)=\{0\}$, the identity $$U_{\left(\11b{M} - \mu^{-1} b \right)^{-\frac12}} \Big(\11b{M} -  \mu^{-1} b- \beta^2 x\Big) = \11b{M} - U_{\beta \left(\11b{M} - \mu^{-1} b \right)^{-\frac12}} (x)$$ implies that $\displaystyle U_{\left(\11b{M} - \mu^{-1} b \right)^{-\frac12}} \Big(\11b{M} -  \mu^{-1} b- \beta^2 x\Big)$ is an invertible element for all $\beta\in \mathbb{C}\backslash\{0\}.$ But the same conclusion trivially holds for $\beta =0$. We therefore deduce that the element $\displaystyle U_{\left(\11b{M} - \mu^{-1} b \right)^{-\frac12}} \Big(\11b{M} -  \mu^{-1} b- \alpha x\Big)$ is invertible for all $\alpha\in \mathbb{C}$. The fundamental identity \eqref{eq fundamental identity UaUbUa} shows that $$U_{\left(\11b{M} - \mu^{-1} b \right)^{-\frac12}} \Big(\11b{M} -  \mu^{-1} b- \alpha x\Big) = U_{\left(\11b{M} - \mu^{-1} b \right)^{-\frac12}} U_{\11b{M} -  \mu^{-1} b- \alpha x} U_{\left(\11b{M} - \mu^{-1} b \right)^{-\frac12}},$$ and since $\left(\11b{M} - \mu^{-1} b \right)^{-\frac12}$ is clearly invertible, the element $\11b{M} -  \mu^{-1} b- \alpha x,$ equivalently $- \mu \11b{M} +  b+ \alpha x,$ must be invertible for all $\alpha\in \mathbb{C}$ (cf. Lemma~\ref{l first properties invertible Jordan}$(c)$). This shows that $\mu\notin \hbox{Sp} (b+ \alpha x)$ ($\alpha\in \mathbb{C}$, $|\mu|>\rho(b)$), and hence $\mathbb{C}\backslash \overline{B}_{\mathbb{C}} (0,\rho(b))\subseteq \mathbb{C}\backslash \hbox{Sp} (b+ \alpha x)$, for all $\alpha \in \mathbb{C}$. We have therefore shown that $\rho (b+\alpha x)\leq \rho (b)$ for all $\alpha\in \mathbb{C}$. The arbitrariness of $b$ and Theorem~\ref{t Aupetit c 1} imply that $x$ lies in the McCrimmon radical of $M$.	 
\end{proof}

\subsection{On the connected components of the set of invertible elements}\label{subsect: connected components of invertible elements}\ \smallskip

The aim of this subsection is to revisit the topological properties of the set of invertible elements in a unital Jordan--Banach algebra and its connected components. We have already made use of the principal connected component in the latest characterizations of the Jacobson and McCrimmon radicals in unital Banach and Jordan--Banach algebras.  \smallskip 

Well by employing holomorphic functional calculus, or by analytic series, for each element $a$ in a unital complex Jordan--Banach algebra $M$, can define the exponential of $a,$ $\displaystyle \exp(a) =\sum_{n=0}^\infty \frac{a^n}{n!}$. Clearly, the mapping $a\mapsto \exp (a)$ is analytic. 
If an associative Banach algebra $A$ is regarded as a Jordan-Banach algebra with its natural Jordan product, the exponential of each element $a\in A$ coincides with the usual notion in the setting of associative Banach algebras.\smallskip  

Since the Jordan-Banach subalgebra $M_a$ of $M$ generated by $a$ and $\mathbf{1}$ is a commutative unital Banach algebra with respect to the inherited Jordan product \cite[1.1]{Jacobson81} and contains $\exp(a)$, we deduce that $\exp(s a )\  \exp(t a) = \exp((s + t) a)$
for all $s, t \in \mathbb{C}$. We also have $U_a (b) = a^2 \circ b$ for $a, b \in  C$, which yields $$U_{\exp( t x)} (\exp(-tx)) = \exp(tx),\hbox{ and } U_{\exp(tx)} (\exp( -tx)^2) = \mathbf{1}.$$ Invertibility is actually a local notion, i.e., $a$ is invertible in $M$ if and only if it is invertible in $M_a$, $\exp(a)\in M^{-1}$ with inverse $\exp(-a)$ for all $a\in M$ (cf. \cite[1.6.1]{Jacobson81}).\smallskip

We gather in the next result a list of topological and algebraic properties of the connected components of the invertible elements in a unital complex Banach algebra. 

\begin{proposition}\label{p properties of components in Ainverse}{\rm(\cite[Propositions 8.6 and 8.7]{BonDunCNA73} or \cite[\S 7]{Zelaz73})}  Let $A$ be a complex unital Banach algebra, and let $A_{\mathbf{1}}^{-1}$ denote the principal connected  component of $A^{-1}$. Then the following statements hold:
\begin{enumerate}[$(a)$]
	\item $A_{\mathbf{1}}^{-1}$ is the least subgroup of $A^{-1}$ containing $\exp(A)$.  
	\item $A_{\mathbf{1}}^{-1}$ is a normal subgroup {\rm(}i.e. $a b a^{-1}\in A_{\mathbf{1}}^{-1}$, for all $a\in A^{-1}$, $b\in A_{\mathbf{1}}^{-1}${\rm)} and is a maximal connected subgroup of $A^{-1}$. 
	\item  If $A$ is commutative, $A_{\mathbf{1}}^{-1}$ is $\exp(A)$.
\end{enumerate}
\end{proposition}

\begin{remark}\label{r open subgroups of Ainverse are clopen} Suppose $A$ is a unital complex Banach algebra. Let $\mathfrak{A}$ be an open subgroup of $A^{-1}$. It is well known that $\mathfrak{A}$ must be also closed, and consequently $\mathfrak{A}$ contains the principal component of $A^{-1}$. To see this it suffices to show that $A^{-1}\backslash \mathfrak{A}$ is open. Indeed, for each $b\in A^{-1}\backslash \mathfrak{A}$, the set $\mathfrak{A} b$ is an open neighbourhood of $b$, whose intersection with $\mathfrak{A}$ must be empty, since if $a\in \mathfrak{A}\cap \mathfrak{A} b$, then $a = c b$ with $c\in \mathfrak{A}$ and, by applying that $\mathfrak{A}$ is a subgroup, we have $b = c^{-1}(c b)\in \mathfrak{A}$, which is impossible.\smallskip  
	
Furthermore, since $A^{-1}$ is open, and hence locally connected and its connected components are open, and $\mathfrak{A}$ is clopen a classic result in topology {\rm(}see, for example, \cite[Theorem 5]{Niem77}{\rm)} asserts that $\mathfrak{A}$ must coincide with the union of some connected components of $A^{-1}$.\smallskip

If $\mathcal{C}$ is another connected component of $A^{-1}$ and we pick $b\in \mathcal{C}$, it is easy to see that $\mathcal{C} = b A^{-1}_{\mathbf{1}} = A_{\mathbf{1}}^{-1} b$ {\rm(}just apply that the left and right multiplication operators by $b$ are homeomorphisms mapping $1$ to $b${\rm)}. The component $\mathcal{C}$ need not be a subgroup, however, if we replace the product of $A$ with the one defined by $x \cdot_{b^{-1}} y = x b^{-1} y$, we get another another associative complete normed algebra $A_{b^{-1}}$, called the \emph{$b^{-1}$-homotope} of $A$, with unit $b$. The main connected component of $A_{b^{-1}}$ is precisely $\mathcal{C}$, which is therefore a subgroup of $A_{b^{-1}}$.
\end{remark}

The topological and algebraic properties of the principal connected component $M_{\mathbf{1}}^{-1}$ of the invertible elements in a unital complex Jordan--Banach algebra $M$ have not been completely clarified until quite recently. For example, Aupetit asked in \cite[page 484]{Aupetit95} whether $M_{\mathbf{1}}^{-1}$ is connected by
analytic arcs, as it is in Banach algebras. An algebraic characterization was finally found by Loos in 1996.\smallskip 

\begin{theorem}\label{t Loos thm characterization of the principal component}\cite{Loos96} For each unital complex Jordan--Banach algebra $M$ we have $$M_{\mathbf{1}}^{-1} = \left\{ U_{\exp(a_1)}\cdots U_{\exp(a_n)} (\11b{M}) : a_i\in M, \ n\geq 1 \right\}.$$ Consequently, $M_{\mathbf{1}}^{-1}$ is open and closed in $M^{-1}$. Moreover each connected component of $M^{-1}$ is analytically arc-wise connected.  
\end{theorem}

It should be noted that in Loos's proof of the previous result the norm is not required to satisfy $\|x\circ y\|\leq \|x\| \ \|y\|$ and $\|\mathbf{1}\|=1$.\smallskip

The arguments in Loos' paper rely on the fact that passing to an appropriate isotope we can avoid the use of left or right multiplication by an invertible element like in an associative Banach algebra, operation which has no direct Jordan analogue. Suppose $c$ is an invertible element in a unital Jordan-Banach algebra $M$ then the vector space $M$ becomes a complete normed Jordan algebra $M_{(c)}$ with unit element $c^{-1}$ and Jordan product $a\circ_c b := U_{a,b} (c)$, for all $a,b\in M$. This algebra is called the \emph{$c$-isotope}\label{def of u-isotope} of $M$. It is known that the quadratic operator in the $u$-isotope satisfies \begin{equation}\label{eq quadratic operator in the isotope} U^{(c)}_{a} = U_a U_c, \hbox{ for all } a,b\in M.
\end{equation}

\begin{remark}\label{r invertible in the c-isotope} Let $c$ be an invertible element in a unital complex Jordan--Banach algebra $M$, and let $M_{(c)}$ stand for the $c$-isotope of $M$. Then $M^{-1} = M_{(c)}^{-1}$. If $a\in M^{-1},$ 
having in mind that the corresponding $U_{a}^{(c)}$ operator in the $c$-isotope is given by $U_{a}^{(c)} = U_{a} U_{c},$ where the two operators on the right-hand-side are invertible, we deduce that $a\in M_{(c)}^{-1}$. The  reciprocal inclusion follows from the identity $U_a^{(c)} U_{c}^{-1}= U_a$. Moreover, $\left(U_{a}^{(c)} \right)^{-1} = U_{c}^{-1} U_{a^{-1}},$ and by the fundamental identity, $U^{(c)}_{U_{c^{-1}}(a^{-1})} = U_{U_{c^{-1}}(a^{-1})} U_{c} = U_{c^{-1}} U_{a^{-1}} U_{c^{-1}} U_{c} = U_{c^{-1}} U_{a^{-1}},$ which shows that $ U_{c^{-1}}(a^{-1})$ is the inverse of $a$ in $M_{(c)}$.   
\end{remark}

The transfer and repercussion between the results for Banach algebras and Jordan--Banach algebras is not always smooth, especially from the Jordan setting to the associative case. A similar conclusion to that in Loos' theorem was independently discovered by Hatori for the principal component of the unitary group in a unital C$^*$-algebra \cite[Lemma 3.2]{HatMol2014}.\smallskip

The Jordan product of two invertible elements in a unital Jordan-Banach algebra $M$ is not, in general, invertible. The best we can conclude is that $U_a (b)\in M^{-1}$ for all $a,b\in M^{-1}$ (cf. Lemma~\ref{l first properties invertible Jordan}). This motivates us to introduce the following concept. 

\begin{definition}\label{def cuadratic subset} A subset $\mathfrak{M}$ of $M^{-1}$ will be called a \emph{quadratic subset} if it is non-empty and $U_{a} (b) \in \mathfrak{M}$ for all $a,b\in \mathfrak{M}$.
\end{definition}

We can now present a recent algebraic characterization of the principal component of the invertible elements in a unital Jordan--Banach algebra.

\begin{lemma}\label{l principal component as the least quadratic subset of invertible elements}\cite[Lemma 2.5]{Pe2021} Let $M$ be a unital complex Jordan-Banach algebra. Then the principal component of $M^{-1}$ is precisely the least quadratic subset of $M^{-1}$ containing $\exp(M)$. Consequently, for each natural number $n$ and each $a\in M_{\mathbf{1}}^{-1}$ the element $a^{n}$ lies in $M_{\mathbf{1}}^{-1}$. 
\end{lemma}

\begin{proof} Obviously, $M_{\mathbf{1}}^{-1}$ contains $\exp(M)$. To prove that $M_{\mathbf{1}}^{-1}$ is a quadratic subset of $M^{-1}$ Loos' Theorem~\ref{t Loos thm characterization of the principal component} plays a central role, since we can restrict to elements of the form $a = U_{\exp(a_1)}\cdots U_{\exp(a_n)} (\mathbf{1})$ and $b = U_{\exp(b_1)}\cdots U_{\exp(b_m)} (\mathbf{1})$ in $M_{\mathbf{1}}^{-1}$. By the fundamental formula \eqref{eq fundamental identity n elements} we have $$\begin{aligned}
U_{a} (b) &= U_{U_{\exp(a_1)}\cdots U_{\exp(a_n)} (\mathbf{1})} (U_{\exp(b_1)}\cdots U_{\exp(b_m)} (\mathbf{1})) \\
 &= U_{\exp(a_1)} \cdots U_{\exp(a_n)} U_{\mathbf{1}} U_{\exp(a_n)}\cdots U_{\exp(a_1)} U_{\exp(b_1)}\cdots U_{\exp(b_m)} (\mathbf{1}) \in M_{\mathbf{1}}^{-1}.
	\end{aligned}$$ Given $a\in M_{\mathbf{1}}^{-1}$, it follows from the above that $a^2 = U_{a} (\mathbf{1})\in M_{\mathbf{1}}^{-1}$, $a^3 = U_{a} (a)\in M_{\mathbf{1}}^{-1}$, 
	and by an induction argument $a^{n} = U_{a} (a^{n-2})\in M_{\mathbf{1}}^{-1}$ for all $n\geq 3$. \smallskip
	
We finally assume that $\mathfrak{M}$ is a \emph{quadratic subset} of $M^{-1}$ containing $\exp(M)$. Since, clearly $\mathbf{1}\in \mathfrak{M}$, we have $U_{\exp(a)} (\mathbf{1})\in \mathfrak{M}$. Suppose that $U_{\exp(a_n)}\cdots U_{\exp(a_1)} (\mathbf{1})\in \mathfrak{M}$ for $n\geq 2$, $a_i\in M$.  Then, by the assumption on $\mathfrak{M}$ and the induction hypothesis, we have $$U_{\exp(a_{n+1})} U_{\exp(a_n)}\cdots U_{\exp(a_1)} (\mathbf{1}) \in \mathfrak{M}.$$ A new application of Loos' Theorem~\ref{t Loos thm characterization of the principal component} implies that $M_{\mathbf{1}}^{-1} \subseteq \mathfrak{M}.$
\end{proof}

\begin{remark}\label{r the principal component is closed for inverses} Let $M$ be a unital complex Jordan--Banach algebra. Since the inverse mapping $a\mapsto a^{-1}$ is a homeomorphism on $M^{-1}$ (cf. Lemma~\ref{l first properties invertible Jordan}$(f)$), and maps $\11b{M}$ to itself, it fixes the principal component of $M^{-1}$, that is, $a^{-1}\in M_{\mathbf{1}}^{-1}$ for all $a\in M_{\mathbf{1}}^{-1}$.   Therefore, $a^{n}\in M_{\mathbf{1}}^{-1}$ for all $a\in M_{\mathbf{1}}^{-1}$ and every integer number $n$. 
\end{remark}

\begin{problem}\label{r open subgroup} For a Jordan-Banach algebra $M$, it would be interesting to know if every open quadratic subset $\mathfrak{M}$ of $M^{-1}$ containing the unit must be also closed, and hence $M_{\mathbf{1}}^{-1}\subseteq \mathfrak{M}$.
\end{problem}

The set $M^{-1}$ is a clopen quadratic subset of itself closed for powers and inverses.\smallskip

We shall conclude this section by showing that clopen quadratic subsets of the invertible elements in a unital Jordan--Banach algebra which are closed for powers and inverses are also stable under translations by elements in the McCrimmon radical. 

\begin{lemma}\label{l N-u0 in N} Let $M$ be a unital Jordan-Banach algebra, and let $\mathfrak{M}\subseteq M^{-1}$ be a clopen quadratic subset of $M^{-1}$ which is closed for powers and inverses. Then for each element $u_0$ in the McCrimmon radical of $M$ the identity $\mathfrak{M}-u_0 =\mathfrak{M}$ holds.
\end{lemma}

\begin{proof} Clearly, $\mathfrak{M}$ contains the unit of $M$ because $U_{a^{-1}} (a^2) = \11b{M}\in \mathfrak{M}$ for all $a\in \mathfrak{M}$. Since $\mathfrak{M}$ is clopen and the connected components of $M^{-1}$ are open, we deduce that $\mathfrak{M}$ coincides with the union of some of the connected components of $M^{-1}$ and $M^{-1}_{\mathbf{1}}\subseteq \mathfrak{M}$. \smallskip
	
Fix an arbitrary $b\in \mathfrak{M}$. A result by McCrimmon assures that, the McCrimmon radical of the $b$-homotope algebra $M_{(b)}$ satisfies $$\mathcal{R}ad(M_{(b)}) = \{z\in M : U_{b} (z)\in \mathcal{R}ad(M)\} \hbox{ (see \cite[Proposition 3 in page 110]{McCrimm71})}.$$ Consequently, $U_{b^{-1}} (u_0) = U_{b}^{-1} (u_0) \in \mathcal{R}ad(M_{(b)})$. It follows that $$\rho_{_{M_{(b)}}}(U^{(b)}_{x} (U_{b^{-1}} (u_0)))=0, \hbox{ for all $x\in M_{(b)}$ (cf. Theorem~\ref{t Aupetit c 1}).}$$ Therefore $\rho_{_{M_{(b)}}}(t U^{(b)}_{b^{-1}} (U_{b^{-1}} (u_0)))=0$, which implies that $$t U^{(b)}_{b^{-1}} (U_{b^{-1}} (u_0)) +b^{-1}\in M_{(b)}^{-1},$$ for all $t\in \mathbb{R}_0^+$ (because $b^{-1}$ is the unit of the homotope $M_{(b)}$). Moreover, since $U^{(b)}_{b^{-1}}$ is the identity mapping because $b^{-1}$ is the unit of the homotope $M_{(b)}$, the previous conclusion implies that $t U_{b^{-1}} (u_0) +b^{-1}\in (M_{(b)})^{-1} = M^{-1},$ for all $t\in \mathbb{R}_0^+$. This proves that $U_{b^{-1}} (u_0) +b^{-1}$ is in the connected component of $M^{-1}$ containing $b^{-1}$, and thus the element $$u_0 +b = U_b(  U_{b^{-1}} (u_0) +b^{-1}) $$ is in the connected component of $M^{-1}$ containing $b,$ and hence $u_0 +b \in \mathfrak{M}$. The arbitrariness of $b$ and the fact that $-b = U_{i b} (b^{-1})\in \mathfrak{M}$ give the desired conclusion.
\end{proof}

\section{Sets of invertible elements as a metric invariant}\label{sec: invertible elements as an invariant}

In the seeking of a certain subset of unital Banach algebras which can serve as an invariant to identify or distinguish them, the first natural candidate is the set of invertible elements. If $A$ is a (associative) unital Banach algebra, the set $A^{-1}$ of all invertible elements in $A$ is a multiplicative subgroup which contains the unit element $\11b{A}$. In the next example we illustrate well known fact in spectral theory assuring that there exist non-isomorphic unital commutative Banach algebras whose subgroups of invertible elements are topologically homeomorphic as topological groups.

\begin{example}\label{example no algebra isomorphic but group of invertible topologically isomorphic}\cite[Remark I.7.8]{Zelaz73} The commutative Banach algebras $A	= C([-1, -\frac12 ]\cup [\frac12, 1])$ and $B = C([0, 1]\cup\{2\})$ are not isomorphic since $I = \{b\in B : b|_{[0,1]} \equiv 0\}$ is a minimal ideal in $B,$ and $A$ lacks of minimal ideals. Actually, the Banach spaces $A$ and $B$ are not isometric. The mapping $a\mapsto b$ with $$b(t) :=\left\{\begin{array}{ll}
		a(t-1), & \hbox{if $t\in [0,\frac12]$;} \\
		{a\left(-\frac12\right)}/{a\left(\frac12\right)} \, a(t), & \hbox{if $t\in [\frac12,1]$;} \\
		{a\left(\frac12\right)}/{a\left(-\frac12\right)}, & \hbox{if $t=\frac12,$} \\
	\end{array} \right.$$ is a group homeomorphism from $A^{-1}$ onto $B^{-1}$.
\end{example} 

Roughly speaking, the topological aspect of the subgroup of invertible elements in a unital Banach is not an invariant to distinguish the whole algebra. \smallskip

Hatori introduced in \cite{Hat09} the problem of determining whether the set of invertible elements regarded as a metric space with the induced distance is an invariant. The first success appears in the setting of unital semisimple commutative Banach algebras.   

\begin{theorem}\label{thm Hatori semisimple commutative}\cite{Hat09} Let $A$ and $B$ be unital Banach algebras. We additionally assume that $A$ is semisimple and commutative. Let $\mathfrak{A}$ and $\mathfrak{B}$ be open subgroups of $A^{-1}$ and $B^{-1}$, respectively. Let $\Delta : \mathfrak{A}\to \mathfrak{B}$ be a surjective isometry for the distances given by the norms of $A$ and $B$, respectively. Then $B$ is semisimple and commutative, and the mapping $\Delta(\11b{A})^{-1} \Delta$ admits an extension to an isometric real-linear algebra isomorphism from $A$ onto $B$. In particular, $A^{-1}$ and $B^{-1}$ are isometrically isomorphic as metric groups.
\end{theorem}

There are certain aspects of the previous result that should be commented more in detail.  We begin with a handicap, the Banach space structure of a complex Banach algebra is not, in general,  an invariant by itself, that is, there exist (unital commutative and semisimple) Banach algebras which are isometrically isomorphic as Banach spaces but not isomorphic as Banach algebras.

\begin{example}\label{example the Banach space structure is not an invariant}\cite[Example 3.2]{Hat09} The existence of a surjective linear isometry between two Banach algebras does not imply that these algebras are isomorphic as Banach algebras, even if we assume that they are unital, semisimple, and commutative. For example, the Wiener algebra, ${A}(\mathbb{T})$, is usually defined as the Banach space of absolutely convergent Fourier series, that is, $${A}(\mathbb{T}): =\left\{f \in C(\mathbb{T}) : \| f\| :=\sum_{k=-\infty}^{\infty} |\widehat{f}(n)|<\infty \right\},$$ equipped with the pointwise product, where $\mathbb{T}$ stands for the unit sphere of $\mathbb{C}$ and $\widehat{f}(n)$ denotes the $n$th Fourier coeficient of $f$. It is known from Gelfand theory that ${A}(\mathbb{T})$ is isometrically isomorphic, via the isomorphism given by the Fourier transform, to the the discrete group Banach algebra $\ell_{1}(\mathbb{Z})$ with convolution product (cf. \cite[page 30]{ArvesonBook}). Furthermore, the spectrum of ${A}(\mathbb{T})$ identifies with $\mathbb{T}$.\smallskip
	
On the other hand, we can consider the subalgebra $$ \ell_{1}(\mathbb{Z}^+):=\left\{f \in {A}(\mathbb{T}) :   \widehat{f}(n) =0, \hbox{ for all } n <0\right\},$$ whose spectrum is the closed unit disc $\overline{\mathbb{D}}$ (cf. \cite{MoriniRenteln89}).  Since $\overline{\mathbb{D}}$ is not homeomorphic to $\mathbb{T}$,  ${A}(\mathbb{T})$ and $\ell_{1}(\mathbb{Z}^+)$ are not isomorphic as complex Banach algebras, neither are as Jordan--Banach algebras.\smallskip
	
The mapping $T_W : {A}(\mathbb{T})\to 	\ell_{1}(\mathbb{Z}^+)$ defined by 
$$(T_W (f))(e^{it}) =\sum_{n=0}^{\infty} \widehat{f} (n) e^{2nit} +
\sum_{n=1}^{\infty} \widehat{f}(-n) e^{(2n-1)i t}$$ is a surjective linear isometry.  
\end{example}

Another example, also due to Hatori, shows that a surjective unital isometry between two unital Banach algebras need not be, in general, multiplicative. 

\begin{example}\label{example Hatori isomerically isomorphic non isomorphic}(\cite[Example 3.4]{Hat09} and \cite[Example 5.3]{Hat2011}) Let $B$ denote the unitization of the algebra of all $3\times 3$ upper triangular matrices with the usual matrix product, that is, $$B = \left\{
\begin{pmatrix}
	\alpha & a & b\\
	0 & \alpha & c \\
	0 & 0 & \alpha
\end{pmatrix}  \colon 	\alpha,a,b,c\in \mathbb{C} \right\}.$$ Let $A$ denote the algebra obtained from the same linear space of $B$ equipped with the trivial product given by $(\alpha \mathbf{1} + y) (\beta \mathbf{1} + x) = \alpha \beta \mathbf{1} +\alpha x + \beta y$ for all $\alpha, \beta\in \mathbb{C}$ and every pair of upper triangular matrices $M,N$. Both algebras are equipped with the usual operator norm. Clearly $A$ and $B$ are unital but not semisimple with $A$ being commutative.  Thus, the identity mapping $Id: x\mapsto x$ is a unital surjective linear isometry from $A$ to $B$. It is easy to check that $Id$ maps invertible elements to invertible elements, where $A^{-1} = B^{-1} =  \left\{
\begin{pmatrix}
\alpha & a & b\\
0 & \alpha & c \\
0 & 0 & \alpha
\end{pmatrix}  \colon 	\alpha,a,b,c\in \mathbb{C}, \alpha\neq 0\right\} .$ Actually, $A^{-1}$ and $B^{-1}$ cannot be isomorphic as groups because one is commutative and the other one is not.  
\end{example}

We shall see in the subsequent Anti-theorem~\ref{anti-theorem BKU} the existence of two unital JB$^*$-algebras which are not Jordan $^*$-isomorphic but share the same underlying Banach space. \smallskip

The weakest hypotheses that one could consider for Theorem~\ref{thm Hatori semisimple commutative} in the associative setting are those relaxing the assumptions that $A$ is commutative and semisimple.  

\begin{theorem}\label{thm Hatori 2011}\cite{Hat2011} Let $A$ and $B$ be unital Banach algebras. Assume that $\mathfrak{A}$ and $\mathfrak{B}$ are open subgroups of $A^{-1}$ and $B^{-1}$, respectively. Let $\Delta : \mathfrak{A}\to \mathfrak{B}$ be a surjective isometry for the distances given by the norms of $A$ and $B$, respectively. Then there exist a surjective real-linear isometry $T: A\to B,$ and an element $u_0$ in the Jacobson radical of $B$ which satisfy that $\Delta (a) = T (a)+ u_0,$ for every $a$ in $\mathfrak{A}$. Furthermore, the element $u_0$ can be obtained as $\displaystyle \lim _{\mathfrak{A}\ni a\to 0} \Delta (a).$
\end{theorem} 

Instead of revisiting the original arguments leading to the proofs of Theorems~\ref{thm Hatori semisimple commutative} and \ref{thm Hatori 2011}, we shall obtain them as a consequence of a more general conclusion in the wider setting of Jordan--Banach algebras. First, we shall obtain some corollaries. \smallskip

It should be observed now that in Theorem~\ref{thm Hatori 2011}, the hypotheses on the mapping $\Delta$ do not require that it preserves the group structures of $\mathfrak{A}$ and $\mathfrak{B}$. If we assume this extra hypothesis the conclusion is stronger.   

\begin{corollary}\label{cor thm group isometric Hatori 2011}\cite{Hat2011} Let $A$ and $B$ be unital Banach algebras. Assume that $\mathfrak{A}$ and $\mathfrak{B}$ are open subgroups of $A^{-1}$ and $B^{-1}$, respectively. Let $\Delta : \mathfrak{A}\to \mathfrak{B}$ be a surjective isometric group isomorphism. Then there exist a real-linear isometric isomorphism $T: A\to B$ satisfying $\Delta (a) = T (a),$ for every $a$ in $\mathfrak{A}$. In particular, $A^{-1}$ and $B^{-1}$ are isometrically isomorphic as metrizable groups.
 \end{corollary}  

\begin{proof} By Theorem~\ref{thm Hatori 2011} we can find a surjective real-linear isometry $T: A\to B$ and $u_0 =\displaystyle \lim_{\mathfrak{A}\ni a\to 0} \Delta (a)\in \hbox{rad} (B)$ satisfying $\Delta (a) = T(a) + u_0$ for all $a\in \mathfrak{A}$. By the surjectivity of $\Delta$ we can find $a_0\in \mathfrak{A}$ such that $\Delta (a_0) = 2\ \11b{B}$. Since $\mathfrak{A}$ is a subgroup and $\Delta$ is multiplicative, we have $$u_0 =  \lim_{\mathfrak{A}\ni a\to 0} \Delta (a) = \lim_{\mathfrak{A}\ni a\to 0} \Delta (a_0 a)  = \lim_{\mathfrak{A}\ni a\to 0} \Delta (a_0) \Delta(a) =  \lim_{\mathfrak{A}\ni a\to 0} 2\ \11b{B}\ \Delta(a)  = 2 u_0,$$ which shows that $u_0=0$, and hence $\Delta$ agrees with $T$ on $\mathfrak{A}$. In particular, $T(\11b{A}) =\11b{B}$. Finally, fix $a,b\in A$, for $\delta >0$ small enough, the elements $\11b{A} + \delta a$ and  $\11b{A} + \delta b$ lie in $\mathfrak{A}$, and hence $$\begin{aligned}
(\11b{B} + \delta T(a))  (\11b{B} + \delta T(b)) &= T(\11b{A} + \delta a)  T(\11b{A} + \delta b) = \Delta(\11b{A} + \delta a)  \Delta(\11b{A} + \delta b) \\
&= \Delta \left((\11b{A} + \delta a)  (\11b{A} + \delta b)\right)= T \left((\11b{A} + \delta a)  (\11b{A} + \delta b)\right) \\
&= \11b{B}  + \delta T(a)  + \delta T(b) + \delta^2 T(a b), 
	\end{aligned} $$ implying that $T(a) T(b) = T(a b)$ for all $a,b\in A$, by the arbitrariness of $a$ and $b$. Therefore, $T$ is a real-linear algebra isometric algebra isomorphism, whose restriction to $A^{-1}$ gives an isometric group isomorphism onto $B^{-1}$.
\end{proof}

The next corollary is clear.

\begin{corollary}\label{cor 2 thm group isometric Hatori 2011}\cite{Hat2011} Let $A$ and $B$ be unital Banach algebras. The following statements are equivalent:\begin{enumerate}[$(a)$]
 		\item $A$ and $B$ are isometrically isomorphic as real Banach algebras.
 		\item  There exist two open subgroups $\mathfrak{A}$ and $\mathfrak{B}$ of $A^{-1}$ and $B^{-1}$, respectively, which are isometrically isomorphic as groups.
 		 \item $A^{-1}$ and $B^{-1}$ are isometrically isomorphic as metrizable groups.
 	\end{enumerate}
\end{corollary}  

The last main result in this section is an extension of Theorem~\ref{thm Hatori 2011} to the setting of Jordan--Banach algebras. The reader who is familiar with the theory of Jordan--Banach algebras, or has previously taken a look to section~\ref{sec: taste of radicals in Jordan}, will easily appreciate the similarities between the corresponding statements and conclusions up to the necessary adaptations required in the Jordan setting. We result reads as follows (see section~\ref{sec: taste of radicals in Jordan} for any non-clear notion).

\begin{theorem}\label{t surjective isometries between invertible clopen}{\rm\cite{Pe2021}} Let $M$ and $N$ be unital Jordan-Banach algebras. Suppose that $\mathfrak{M}\subseteq M^{-1}$ and $\mathfrak{N}\subseteq N^{-1}$ are clopen quadratic subsets which are closed for powers and inverses. Let $\Delta : \mathfrak{M}\to \mathfrak{N}$ be a surjective isometry. Then the following statements hold: \begin{enumerate}[$(a)$]
\item For each $w_0$ in the McCrimmon radical of $N$ the mapping $a\mapsto \Delta (a) - w_0$ is a surjective isometry from $\mathfrak{M}$ to $\mathfrak{N}.$
\item The norm limit $\displaystyle \lim_{\mathfrak{M}\ni a\to 0} \Delta (a)$ exists and defines an element $u_0$ in the McCrimmon radical of $N$.
\item There exists a surjective real-linear isometry $T_0: M\to N$ such that $\Delta (a) = T_0(a) +u_0,$ for all $a\in \mathfrak{M}$.
\end{enumerate}
\end{theorem} 

The proof will be obtained in a series of results in the next subsection~\ref{subsection proof of Hatori thm for Jordan Banach}. Before presenting it, we shall see how Theorems~\ref{thm Hatori semisimple commutative} and \ref{thm Hatori 2011} can be obtained from the previous Theorem~\ref{t surjective isometries between invertible clopen}. Theorem~\ref{thm Hatori 2011} is easier to obtain since each open subgroup of the invertible elements of a unital Banach algebra $A$ is a clopen, quadratic subset which is closed for powers and inverses (cf. Remark~\ref{r open subgroups of Ainverse are clopen}) and the McCrimmon radical of $A$ regarded as a Jordan--Banach algebra is precisely the Jacobson radical (see Remark~\ref{r Jordan radical and McCrimmon radical coincide when it makes sense}).\label{label page mathfrakA satisfies mathfrakM} \smallskip

The proof of Theorem~\ref{thm Hatori semisimple commutative} is a bit more laborious. We shall actually prove a bit more general statement. Let us remark that the hypotheses are fully justified in view of Example~\ref{example the Banach space structure is not an invariant}. 

\begin{theorem}\label{t real-linear isometries from a semisimple commutative Banach algebra} Let $A$ be a semisimple and commutative unital Banach algebra, and let $N$ be a unital Jordan--Banach algebra. Let $\mathfrak{A}$ and $\mathfrak{N}$ be clopen quadratic subsets closed for powers and inverses of $A^{-1}$ and $N^{-1},$ respectively. Let $\Delta : \mathfrak{A}\to \mathfrak{N}$ be a surjective isometry for the distances given by the norms of $A$ and $N$, respectively. Then $N$ is semisimple and associative, and the mapping $\Delta(\11b{A})^{-1} \Delta$ admits an extension to an isometric real-linear algebra isomorphism from $A$ onto $N$. In particular, $A^{-1}$ and $N^{-1}$ are isometrically isomorphic as metric groups.
\end{theorem}

For the proof we shall need some technical tools. Let us begin with the celebrated Gleason--Kahane--\.{Z}elazko theorem. 

\begin{theorem}\label{t GKZ}{\rm \cite[Theorem 2.4.13]{PalmerBook}} Let ${A}$ be a complex unital Banach algebra, and let $\varphi : {A}\to \mathbb{C}$ be a non-zero linear functional. Then the following statements are equivalent:\begin{enumerate}[$(a)$]\item $\varphi$ is spectral-valued, that is, $\varphi (a)\in \sigma(a)$ for all $a\in {A}$.
\item $\varphi$ is unital and maps invertible elements in ${A}$ to non-zero complex numbers.
\item $\varphi$ is multiplicative.
\end{enumerate}
\end{theorem}

An appropriate version of the Gleason--Kahane--\.{Z}elazko theorem for complex Jordan--Banach algebras has been recently established in \cite{EscPeVill2023} via local theory.

\begin{theorem}\label{t GKZ JB}{\rm\cite[Theorem 3.2]{EscPeVill2023}} Let $M$ be a complex Jordan--Banach algebra, and let $\varphi : M\to \mathbb{C}$ be a non-zero linear functional. Then the following statements are equivalent:\begin{enumerate}[$(a)$]\item $\varphi (a)\in \hbox{Sp}(a)$ for all $a\in M$.
\item $\varphi$ is unital or admits a unital extension to the unitization of $M$ and maps invertible elements in $M$ to non-zero complex numbers.
\item $\varphi$ is {\rm(}Jordan-{\rm)}multiplicative.
	\end{enumerate}
\end{theorem}

As noted by Hatori in \cite{Hat09}, the arguments in most of available proofs of the Gleason--Kahane--\.{Z}elazko theorem are also valid for linear maps valued in the exponential spectrum \cite{Harte1976}. We consider this exponential spectrum in the setting of Jordan--Banach algebras. 

\begin{definition}\label{def eponential spectrum} Let $a$ be an element in a unital complex Jordan--Banach algebra $M$. The exponential spectrum of $a$ is defined as $$\hbox{Sp}_{exp} (a) := \{\lambda \in \mathbb{C} : a-\lambda\11b{M}\notin M^{-1}_{\11b{M}}\}.$$
\end{definition}

Clearly, $\hbox{Sp} (a)\subseteq \hbox{Sp}_{exp} (a)$ is a non-empty compact subset of $\mathbb{C}$ bounded by $\|a\|$.

\begin{lemma}\label{l GKZ for exponential spectrum valued} Let $\varphi : M\to \mathbb{C}$ be a non-zero linear functional, where $M$ is a complex unital Jordan--Banach algebra. Suppose that  $\varphi (a)\in \hbox{Sp}_{exp}(a)$ for all $a\in M$. Then $\varphi$ is (Jordan-)multiplicative.
\end{lemma}

\begin{proof} We observe that for each $a\in M$ and $\lambda\in \mathbb{C}$ the element $e^{\lambda a}$ lies in $M_{\mathbf{1}}^{-1}$, and hence $\varphi(a)\in \hbox{Sp}_{exp}(a)\in \mathbb{C}\backslash\{0\}$. We also know that $\varphi (\11b{M})\subseteq \hbox{Sp}_{exp}(\11b{M}) =\{1\}$. Therefore, the mapping $\varphi : \mathbb{C} \to \mathbb{C}$, $f (\lambda) = \varphi(e^{\lambda a})$ never vanishes and $f(0) =1$. We can now follow the standard reasoning. Having in mind that $\varphi$ is linear and continuous we deduce that  $\displaystyle f (\lambda) = \sum_{n=0}^{\infty} \frac{\varphi(a^n)}{n!} \lambda^n$ for all $\lambda\in \mathbb{C}.$ Clearly, $f$ is an entire function that never vanishes admits, and hence it admits an entire ``logarithm'', that is, there exists an entire function $g :\mathbb{C}\to \mathbb{C}$ satisfying $f (\lambda) = e^{g(\lambda)}$ ($\lambda\in \mathbb{C}$). The inequality $$ |f (\lambda)| \leq \sum_{n=0}^{\infty} \frac{|\varphi(a^n)|}{n!} |\lambda|^n \leq \|\varphi\| \ e^{|\lambda| \ \|a\|}, \ (\lambda\in \mathbb{C})$$ implies that  $f$ has exponential type bounded by one. It follows from Hadamard's Factorization theorem (see \cite[Theorem 2.7.1]{Boas}) that $g (\lambda ) = \alpha \lambda +\beta$ for suitable complex numbers $\alpha$ and $\beta$. By applying that $f(0)=1$, we get $\beta =0$. Consequently, $$ \sum_{n=0}^{\infty} \frac{\varphi(a^n)}{n!} \lambda^n = f (\lambda) = e^{\alpha \lambda }= \sum_{n=0}^{\infty} \frac{\alpha^n}{n!} \lambda^n, \hbox{ for all } \lambda\in \mathbb{C},$$ which proves that $\varphi (a^n) = \alpha^n $ for all natural $n$, and thus $ \varphi (a) = \alpha$ and $\varphi(a^2) = \alpha^2 = \varphi(a)^2$.
\end{proof}

The next result was established by Hatori in the case of Banach algebras. The idea comes from a result by Kowalski and S{\l}odkowski (cf. \cite[Lemma 2.1]{KoSlod}).

\begin{lemma}\label{l GKZ for exponential spectrum valued real-linear} Let $\varphi : M\to \mathbb{C}$ be a non-zero real-linear functional, where $M$ is a complex unital Jordan--Banach algebra. Suppose that  $\varphi (a)\in \hbox{Sp}_{exp}(a)$ for all $a\in M$. Then $\varphi$ is complex linear and {\rm(}Jordan-{\rm)}multiplicative.
\end{lemma}

\begin{proof} The desired conclusion can be obtained by fixing an element $a$ in $M$ and considering the Jordan--Banach subalgebra $M_{a}$ generated by $a$ and $\11b{M}$, which is a commutative unital complex Banach algebra. Since $\left( M_{a}\right)^{-1}_{\mathbf{1}}\subseteq M_{\mathbf{1}}^{-1}$, and hence $$\hbox{Sp$_{exp}(M,b)\subseteq$Sp$_{exp} (M_a,b),$}$$ for every $b\in M_a$, we can apply \cite[Lemma 2.2]{Hat09} to the restricted mapping $\varphi|_{M_a}:M_a\to \mathbb{C}$. For the sake of self-containment, we revisit all details.\smallskip
		
As in the proof of \cite[Lemma 2.1]{KoSlod} and \cite[Lemma 2.2]{Hat09} (see also \cite[Lemma 3.3]{LiPeWanWan2019}), the functionals $$\varphi_1 (x) := \Re\hbox{e} \varphi(x) + i \Im\hbox{m} (-i \varphi(ix)) = \Re\hbox{e} \varphi(x) - i \Re\hbox{e} \varphi(ix)$$
$$\varphi_2 (x) := \Re\hbox{e} (-i \varphi(i x)) + i \Im\hbox{m} \varphi(x) = \Im\hbox{m} \varphi(ix) + i \Im\hbox{m} \varphi(x)$$ are continuous, complex linear,  $\Re\hbox{e} (\varphi(x)) = \Re\hbox{e} (\varphi_1(x))$ and $\Im\hbox{m} (\varphi(x)) = \Im\hbox{m} (\varphi_2(x))$ for all $x\in M,$ and $\varphi_j (x) \in \hbox{Sp}_{exp}(x),$ for every $x\in M$, $j=1,2$. By Lemma~\ref{l GKZ for exponential spectrum valued}, $\varphi_1$ and $\varphi_2$ are Jordan multiplicative functionals on $M$. The desired conclusion will follow as soon as we prove that $\varphi_1 = \varphi_2$. If that is not the case, there exists $a\in M$ with $\varphi_1 (a) =1$ and $\varphi_2(a) = 0$. The element $b=\exp(\frac{ \pi i a}{2})-\11b{M}\in M$ satisfies that $b+\11b{M} = \exp(\frac{ \pi i a}{2})\in M_{\mathbf{1}}^{-1},$ and thus $-1\notin \hbox{Sp}_{exp} (b)$. On the other hand by applying that $\varphi_1$ and $\varphi_2$ are continuous homomorphisms we deduce that $\varphi_1 (b) = i-1,$ $\varphi_2 (b) = 0$, and $\varphi (b) = -1$. The hypotheses on $\varphi$ imply that $-1\in \hbox{Sp}_{exp} (b)$, which  is impossible.
\end{proof}

We can now present the proof of Theorem~\ref{t real-linear isometries from a semisimple commutative Banach algebra} which follows the lines in \cite[Theorem 3.3]{Hat09}. The proof is a bit technical, but perhaps the reader will appreciate a complete access to all details.

\begin{proof}[Proof of Theorem~\ref{t real-linear isometries from a semisimple commutative Banach algebra}] Let $\Delta : \mathfrak{A}\to \mathfrak{N}$ be a surjective isometry satisfying the hypotheses of the theorem. We observe that under these assumptions $\mathfrak{A}\supseteq A_{\mathbf{1}}^{-1}$ and $\mathfrak{N}\supseteq N_{\mathbf{1}}^{-1}$. 
By  Theorem~\ref{t surjective isometries between invertible clopen} applied to $\Delta^{-1}$, there exists a real-linear surjective isometry $T: N\to A$ such that $\Delta^{-1} =T|_{\mathfrak{N}}$. We need to prove that $T$ is a Jordan homomorphism. Let $\Omega_{A}$ denote the character space of $A$, and let $G_A : A\to C(\Omega_A)$ be the Gelfand transform, which is an isomorphism of Banach algebras. 
	
In a first step we shall prove that 
\begin{equation}\label{eq images of the unitaries we know} |G_{A} (T(\lambda \11b{N})) (t) | =1, \hbox{ for all } \lambda \in \mathbb{T}, t\in \Omega_{A}. 
\end{equation}	

Clearly, $|G_{A} (T(\lambda \11b{N})) (t) | \leq \| G_{A} (T(\lambda \11b{N}))\| \leq \|T(\lambda \11b{N}) \| = \|\lambda \11b{N} \|=1.$ If there exists $t\in \Omega_A$ with $|G_{A} (T(\lambda \11b{N})) (t) | <1$, for $\alpha =  G_{A} (T(\lambda \11b{N})) (t)$ we have $$1>|\alpha| = \|T(\lambda \11b{N}) - (T(\lambda \11b{N})- \alpha \11b{A})\| = \| \lambda \11b{N} - T^{-1} (T(\lambda \11b{N})- \alpha \11b{A})\|,$$ and consequently $T^{-1} (T(\lambda \11b{N})- \alpha \11b{A})\in N_{\mathbf{1}}^{-1}\subseteq \mathfrak{N}.$ Since $T(\mathfrak{N})= \Delta^{-1} (\mathfrak{N}) = \mathfrak{A}\subseteq A^{-1}$, we deduce that $T(\lambda \11b{N})- \alpha \11b{A}\in \mathfrak{A}\subseteq A^{-1}$, which is impossible since $$G_{A}\left( T(\lambda \11b{N})- \alpha \11b{A} \right) (t) = G_{A}\left( T(\lambda \11b{N}) \right) (t) - G_{A}\left( T(\lambda \11b{N})\right) (t) =0.$$ This finishes the proof of \eqref{eq images of the unitaries we know}.\smallskip

Define now a new real-linear bijection $S: N\to A$ defined by  $S(b) := T(\11b{N})^{-1} T(b) = \left(\Delta^{-1}(\11b{N})\right)^{-1} T(b)$. Clearly $S(\11b{N}) = \11b{A}$ and $\|b\|= \|T(b)\| \leq \| T(\11b{N}) \|\  \| T(\11b{N})^{-1} T(b)\| = \|S(b)\|.$ Furthermore \begin{equation}\label{eq S maps frak N to frak A} S(\mathfrak{N}) = T(\11b{N})^{-1} T(\mathfrak{N}) = T(\11b{N})^{-1} \mathfrak{A}\subseteq {A}^{-1}.
\end{equation}
	
Next we show that \begin{equation}\label{eq S(i) 0 pm i} G_{A} (S(i\ \11b{N})) (t) \in \{\pm i\}, \hbox{ for all } t\in \Omega_{A}. 
\end{equation}	This will be simpler from \eqref{eq images of the unitaries we know}, since by this result we get $$\begin{aligned}
|1 \pm G_{A} (S(i \11b{N})) (t)| &= |G_{A} (T(\11b{N})) (t) \pm G_{A} (T(i \11b{N})) (t)| \\
&\leq \| T(\11b{N}\pm i \11b{N})\| = \|(1+i) \11b{N}\| = \sqrt{2}.
\end{aligned}$$ Having in mind that $G_{A} (S(i \ \11b{N})) (t)\in \mathbb{T}$, we arrive to the desired conclusion in \eqref{eq S(i) 0 pm i}. \smallskip

Therefore, the set $\Omega_{A}$ decomposes as the disjoint union of the (possibly empty) clopen sets $$\Omega_{A}^{+} = \{t\in \Omega_{A} : G_{A} (S(i \ \11b{N})) (t) =i \} \hbox{ and }\Omega_{A}^{-} = \{t\in \Omega_{A} : G_{A} (S(i\ \11b{N})) (t) =- i \}.$$ 	  

We shall next prove the following statement:
\begin{equation}\label{eq deltat S is a complex linear Jordan homo for +}\begin{aligned}&\hbox{For each $t\in \Omega_{A}^+$ the mapping } \varphi_t := \delta_{t}\circ G_{A}\circ S : N\to \mathbb{C} \hbox{ is a real-linear} \\
& \hbox{functional with $\varphi_t (b) \in \hbox{Sp}_{exp} (b)$ for each  element $b\in N$, and hence } \\
& \hbox{it is a complex linear Jordan homomorphism by Lemma~\ref{l GKZ for exponential spectrum valued real-linear}.}		
\end{aligned}
\end{equation} To get the desired conclusion, let $\mu = \varphi_t (b)$ for some $b\in N$. Clearly, $G_A (S(b) -\mu \11b{A}) (t) =0$. If $\mu\notin \hbox{Sp}_{exp} (N,b)$ the element $b- \mu \11b{N}\in N_{\mathbf{1}}^{-1}\subseteq \mathfrak{N}$. Having in mind \eqref{eq S maps frak N to frak A}, the previous properties imply that $G_{A}\left(S(b-\mu \11b{N})\right) \in G_{A} \left(S(\mathfrak{N}\right) = G_{A}\left( T(\11b{N})^{-1} \mathfrak{A}\right)\subseteq G_{A}\left( A^{-1}\right) \subseteq C(\Omega)^{-1}.$ However, since $t\in \Omega_{A}^+$ we get $$ \begin{aligned}
G_{A}\left(S(b-\mu \11b{N})\right) (t) &= G_{A}\left(S(b-\Re\hbox{e}(\mu) \11b{N}- \Im\hbox{m}(\mu) i \11b{N})\right)  (t) \\
&= \left(G_{A}\left(S(b)\right) - \Re\hbox{e}(\mu) G_{A}\left(S( \11b{N} )\right) - \Im\hbox{m}(\mu) i G_{A}\left(S( \11b{N} )\right) \right) (t) \\
& = \left(G_{A}\left(S(b)\right) - \mu G_{A}\left(S( \11b{N} )\right) \right) (t) = G_A (S(b) -\mu \11b{A} ) (t),
\end{aligned}   $$ contradicting that $G_A (S(b) -\mu \11b{A}) (t) =0$.\smallskip

Similar arguments to those given in the proof of \eqref{eq deltat S is a complex linear Jordan homo for +} lead to \begin{equation}\label{eq deltat S is a complex linear Jordan homo for -}\begin{aligned}&\hbox{For each $t\in \Omega_{A}^-$ the mapping } \varphi_t := \overline{\delta_{t}\circ G_{A}\circ S} : N\to \mathbb{C} \hbox{ is a real-linear} \\
& \hbox{functional with $\varphi_t (b) \in \hbox{Sp}_{exp} (b)$ for each  element $b\in N$, and hence } \\
& \hbox{it is a complex linear Jordan homomorphism by Lemma~\ref{l GKZ for exponential spectrum valued real-linear}.}		
	\end{aligned}
\end{equation}

Now, by combining \eqref{eq deltat S is a complex linear Jordan homo for +} and \eqref{eq deltat S is a complex linear Jordan homo for -} and the fact that $G_{A}$ is a Jordan isomorphism, we deduce that $S$ is a Jordan isomorphism. This implies that $N$ is a semisimple associative unital complex Jordan--Banach algebra (i.e. a semisimple commutative unital Banach algebra), and $S$ is an expansive isomorphism of Banach algebras.\smallskip

Finally, having in mind that $N$ is a unital commutative semsimple complex Banach algebra, the arguments in the first part also assure that the mapping $R = \Delta(\11b{A})^{-1} T^{-1} = \left(T^{-1}(\11b{A})\right)^{-1} T^{-1} : A\to N$ is an isomorphism of Banach algebras and expansive (i.e. $\|R(a)\| \geq \|a\|$ for all $a\in A$). It is not hard to check that $R =S^{-1}$. Namely, 
$$\begin{aligned}
R \left(S(b)\right) &= R \left( \left(\Delta^{-1}(\11b{N})\right)^{-1} T(b) \right)  = R \left( \left(\Delta^{-1}(\11b{N})\right)^{-1}\right) R \left( T(b) \right)  \\ &= R \left( \Delta^{-1}(\11b{N})\right)^{-1} R \left( T(b) \right) = \left[\Delta(\11b{A})^{-1} T^{-1} T (\11b{N}) \right]^{-1} \Delta(\11b{A})^{-1} T^{-1} T(b) = b,
\end{aligned} $$ for all $b\in N$, and similarly $S\left(R(a)\right) = a,$ for all $a\in A.$ All together gives $\|b\|\leq \|S(b)\| \leq \| R\left(S(b)\right)\| =\|b\|,$ for all $b$ in $N$, which proves that $S$ is an isometry.  
\end{proof}

The next result is essentially due to Hatori {\rm(see \cite[Corollary 5.2]{Hat2011})}, here we state it in a wider setting.  

\begin{corollary}\label{t real-linear isometries from a semisimple to a commutative Banach algebra} Let $A$ be a unital and commutative Banach algebra and let $N$ be a semisimple unital Jordan-Banach algebra. Suppose that $\mathfrak{A}\subseteq A^{-1}$ and $\mathfrak{N}\subseteq N^{-1}$ are clopen quadratic subsets which are closed for powers and inverses. Let $\Delta : \mathfrak{A}\to \mathfrak{N}$ be a surjective isometry. Then $A$ is semisimple and $N$ is associative, and the mapping $\Delta(\11b{A})^{-1} \Delta$ admits an extension to an isometric real-linear algebra isomorphism from $A$ onto $N$. In particular, $A^{-1}$ and $N^{-1}$ are isometrically isomorphic as metric groups.
\end{corollary}

\begin{proof} The desired conclusion will follow from Theorem~\ref{t real-linear isometries from a semisimple commutative Banach algebra} as soon as we prove that $A$ is semisimple. Suppose, on the contrary that there exists $u_0\in \mathcal{R}ad(A)\backslash\{0\}$.  \smallskip
	
Since $N$ is semisimple its McCrimmon radical reduces to $\{0\}$.	Then an application of Theorem~\ref{t surjective isometries between invertible clopen} gives the existence of a surjective real-linear isometry $T: A\to N$ such that $\Delta =T|_{\mathfrak{A}}$. By Lemma~\ref{l N-u0 in N}, and the hypotheses, the mapping $a\mapsto \Delta (a-u_0)$ is a surjective isometry from $\mathfrak{A}$ onto $\mathfrak{N}$. A new application of Theorem~\ref{t surjective isometries between invertible clopen} guarantees the existence of a surjective real-linear isometry $T_{u_0}: A\to N$ such that $\Delta (a-u_0) =T_{u_0} (a)$ for all $a\in \mathfrak{A}$. Having in mind that $\11b{A}\in A_{\mathbf{1}}^{-1} \subseteq \mathfrak{A},$ we have $$T_{u_0} \left(\frac{1}{n} \11b{A} \right) = \Delta \left(\frac{1}{n} \11b{A} -u_0\right) = T \left(\frac{1}{n} \11b{A} -u_0\right), \hbox{ for all natural } n.$$ The continuity of the surjective real-linear isometries $T$ and $T_{u_0}$ assures that $T(-u_0) =0,$ and thus $u_0=0$, which is impossible.
\end{proof}

The set of invertible elements in a unital complex Jordan--Banach algebra is not stable under Jordan products. So, there is no hope to expect a result in the line of Corollary~\ref{cor thm group isometric Hatori 2011}, instead of that we can consider preservers of expressions of the form $U_a(b)$.

\begin{proposition}\label{p surjective isometries preserving quadratic expressions}\cite[Proposition 3.9]{Pe2021} Let $M$ and $N$ be unital Jordan-Banach algebras. Suppose that $\mathfrak{M}\subseteq M^{-1}$ and $\mathfrak{N}\subseteq N^{-1}$ are clopen quadratic subsets which are closed for powers and inverses. Let $\Delta : \mathfrak{M}\to \mathfrak{N}$ be a surjective isometry and set $u = \Delta(\11b{M})$. We shall also assume that $\Delta$ satisfies the following property: \begin{equation}\label{property preservation of quadratic expressions} \Delta (U_a (b)) = U_{\Delta(a)} (\Delta(b)), \hbox{ for all } a,b\in \mathfrak{M}.
	\end{equation} Then there exists a real-linear isometric Jordan isomorphism $J$ from $M$ onto the $u$-isotope $N_{(u)}$ such that $\Delta (a) = J(a) $ for all $a\in \mathfrak{M}$.
\end{proposition}

\begin{proof} Recall that $N_{(u)}^{-1} = N^{-1}$ (cf. Remark~\ref{r invertible in the c-isotope}). Then $\mathfrak{N}\subseteq N_{(u)}^{-1}$. Let us prove that $\mathfrak{N}$ satisfies the same hypotheses as a subset of $N_{(u)}^{-1}$ in the $u$-homotope. Indeed, let us take $a,b\in \mathfrak{N}$. It is known that $U_{u^{-1}} (a^{-1})$ is the inverse of $a$ in the $u$-homotope $N_u$ (see Remark~\ref{r invertible in the c-isotope}) and belongs to $\mathfrak{N}$ by the hypotheses. Furthermore, by definition
	$U_{a}^{(u)} (b) =  U_{a} U_{u} (b) \in \mathfrak{N},$ because, by the hypotheses, $U_{u} (b) \in \mathfrak{N}$. Furthermore $a\circ_{u} a = U_a (u)\in \mathfrak{N}$ by assumptions, and a simple induction argument shows that $\mathfrak{N}$ is closed for powers in $N_{(u)}$.  \smallskip

We shall next show that $\Delta$ also satisfies property \eqref{property preservation of quadratic expressions} for the Jordan product in the $u$-homotope. Namely, by the hypotheses,
	\begin{equation}\label{eq 9 new}\begin{aligned}U^{(u)}_{\Delta(a)} (\Delta(b)) &= U_{\Delta(a)} U_{\Delta(1)} (\Delta(b)) = U_{\Delta(a)} \Delta\left( U_{1} (b)\right)\\
			& = U_{\Delta(a)} \Delta\left(b\right) = \Delta\left(U_{a}(b)\right),
		\end{aligned}
	\end{equation} for all $a,b\in \mathfrak{M}$, which proves the claim. Therefore $\Delta: \mathfrak{M}\subseteq M^{-1}\to \mathfrak{N}\subseteq (N_{(u)})^{-1}$ is a unital mapping satisfying the same hypotheses in our statement. By applying Theorem~\ref{t surjective isometries between invertible clopen}$(c)$ to the latter mapping, we deduce the existence of a surjective real-linear isometry $T_0: M\to N$ and an element $u_0$ in the McCrimmon radical of $N_{(u)}$ such that $\Delta (a) = T_0(a) + u_0$ for all $a\in \mathfrak{M}$.\smallskip
	
Pick $a\in \mathfrak{M}$ with $\Delta(a) = 2 u$. Theorem~\ref{t surjective isometries between invertible clopen}$(b)$ and \eqref{property preservation of quadratic expressions} imply that for each norm-null sequence $(a_n)_n$ in $\mathfrak{M}$ we have $$u_0 = \lim_{n} \Delta\left( U_a (a_n)\right) =\lim_{n} U^{(u)}_{\Delta(a)} (\Delta(a_n)) = \lim_{n} U^{(u)}_{2 u} (\Delta(a_n))  = \lim_{n} 2 \Delta(a_n) = 2 u_0,$$ which assures that $u_0=0$. Therefore $\Delta (a) = T_0(a)$ for all $a\in \mathfrak{M}$.\smallskip
	
It is now time to explore the properties of the real-linear isometry $T_0: M\to N_u$. Fix an arbitrary $a\in M$. Since for $r\in\mathbb{R}$ with $|r|$ large enough we have $a+ r \ \11b{M}\in M_{\mathbf{1}}^{-1}\subseteq \mathfrak{M},$ we deduce from \eqref{property preservation of quadratic expressions} and the conclusions in the previous paragraphs that \begin{equation}\label{equation 1 quadratic morphisms}\begin{aligned}U_{T_0(a+ r \11b{M})}^{(u)} \left( T_0(a+ r \11b{M}) \right) &= U_{\Delta(a+ r \11b{M})}^{(u)} \left( \Delta(a+ r \11b{M}) \right) = \Delta\left( U_{a+ r \11b{M}} (a+ r \11b{M}) \right)\\
			&= T_0\left( U_{a+ r \11b{M}} (a+ r \11b{M}) \right).
		\end{aligned}
	\end{equation} By expanding the extreme terms in the previous identity we get
	$$\begin{aligned}U_{T_0(a+ r \11b{M})}^{(u)} \left( T_0(a+ r \mathbf{1}) \right)&= U_{T_0(a) + r u }^{(u)} \left( T_0(a)+ r u \right)
		= U_{T_0(a) }^{(u)} \left( T_0(a) \right) + r U_{T_0(a) }^{(u)} (u) \\
		&+ r^2  T_0(a) + r^3 u + 2 r U_{T_0(a), u }^{(u)} (T_0(a)) + 2 r^2 U_{T_0(a),u }^{(u)} (u) \\
		&= U_{T_0(a) }^{(u)} \left( T_0(a) \right) + 3 r T_0(a) \circ_{u} T_0(a) + 3 r^2  T_0(a) + r^3 u
	\end{aligned}  $$ and $$U_{a+ r \11b{M}} (a+ r \11b{M}) = U_a(a) + 3 r a^2 + 3 r^2 a + r^3 \11b{M},$$ which combined with \eqref{equation 1 quadratic morphisms} lead to $$\begin{aligned}
		& U_{T_0(a) }^{(u)} \left( T_0(a) \right) + 3 r T_0(a) \circ_{u} T_0(a) + 3 r^2  T_0(a) + r^3 u \\
		&= T_0(U_a(a)) + 3 r T_0(a\circ a) + 3 r^2 T_0(a) + r^3 u,
	\end{aligned}$$ for all $r\in \mathbb{R}$ with $|r|$ big enough. Replacing $r$ with $-r$ and subtracting we get $T_0(a) \circ_{u} T_0(a) = T_0(a\circ a)$ for all $a\in M$. It is standard to check that $T_0$ is a Jordan isomorphism from $M$ onto $N_{(u)}$.
\end{proof}

The last result in this section is a clear consequence of the previous proposition. 

\begin{corollary}\label{cor 2 thm group isometric Hatori 2011 JB} Let $M$ and $N$ be unital complex Jordan--Banach algebras. The following statements are equivalent:\begin{enumerate}[$(a)$]
\item There exists an isometric real-linear Jordan isomorphism from $M$ onto $N$.
\item  There exist two clopen quadratic subsets $\mathfrak{M}\subseteq M^{-1}$ and $\mathfrak{N}\subseteq N^{-1}$ which are closed for powers and inverses and a surjective isometry $\Delta : \mathfrak{M}\to \mathfrak{N}$ satisfying $\Delta(\11b{M})=\11b{N}$ and $\Delta (U_a (b)) = U_{\Delta(a)} (\Delta(b)),$  for all $a,b\in \mathfrak{M}.$
\item There exists a surjective isometry $\Delta: M^{-1}\to N^{-1}$ satisfying $\Delta(\11b{M}) = \11b{N}$ and $\Delta (U_a (b)) = U_{\Delta(a)} (\Delta(b))$, for all $a,b\in {M}^{-1}$.
	\end{enumerate}
\end{corollary}

\subsection{Proof of Theorem~\ref{t surjective isometries between invertible clopen}}\label{subsection proof of Hatori thm for Jordan Banach}\ \smallskip

The proof presented in this section relies, among other tools, on a powerful result established by Mankiewicz in \cite{Mank1972} on the affine behaviour of certain surjective isometries between convex sets. All the arguments have been borrowed from \cite{Pe2021}. We recall that a mapping $F$ from a convex subset $\mathcal{C}$ of a Banach space $X$ to another Banach space $Y$ is said to be \emph{real affine}, or simply \emph{affine}, if for each $t\in [0,1]$, and $x,y\in \mathcal{C}$ we have $F(t x + (1-t) y ) = t F(x) + (1-t) F(y)$. A \emph{convex body} in a normed space is a closed convex set with non-empty interior.  \smallskip

The mentioned result by Mankiewicz, which is  generalization of the celebrated Mazur--Ulam theorem, is the following:

\begin{theorem}\label{t Mankiewicz}{\rm\cite[Theorem 5]{Mank1972}}
Every surjective isometry between convex bodies in two arbitrary normed spaces can be uniquely extended to an affine function between the spaces.	
\end{theorem} 

It is worth note that the above theorem by Mankiewicz motivated the introduction of the so-called strong Mankiewicz property in \cite{MoriOza2018}, a property which plays a fundamental role in the positive solution to Tingley's problem on the extension of surjective isometries from the unit spehere of a unital C$^*$-algebra (respectively,  a JBW$^*$-triple, a compact C$^*$-algebra and unital JB$^*$-algebra) onto the unit sphere of another Banach space \cite{MoriOza2018} (respectively, \cite{BeCuFerPe2021,KaPe2021,Pe2023,PeSvard2022}). We shall see further applications in a subsequent section.\smallskip

In the rest of this section $M$, $N$, $\mathfrak{M},$ $\mathfrak{N},$ and $\Delta$ will have the meaning in the hypotheses of Theorem~\ref{t surjective isometries between invertible clopen}. Our arguments start with an observation on the local affine behaviour of $\Delta$. 

\begin{lemma}\label{l Delta is locally affine}\cite{Pe2021} Let $\Delta : \mathfrak{M}\to \mathfrak{N}$ be a surjective isometry. Then $\Delta$ is a local affine mapping, concretely, for each $a\in \mathfrak{M}$ there exists a positive $\delta_a$, depending on $a$, such that $B(a,\delta_a)\subseteq \mathfrak{M}$, $B(\Delta(a), \delta_a)\subseteq \mathfrak{N}$, $\Delta(B(a,\delta_a)) = B(\Delta(a), \delta_a)$ and there exists a surjective affine isometry $F_{a, \delta_a} : M\to N$ such that $\Delta|_{B(a,\delta_a)} = F_{a, \delta_a}|_{B(a,\delta_a)}$ {\rm(}and hence $\Delta|_{B(a,\delta_a)}$ is an affine mapping{\rm)}. Furthermore, suppose that $\gamma : [0,1]\to \mathfrak{M}$ is a continuous path. Then there exists a surjective affine isometry $F: M\to N$ and an open neighborhood $U$ of $\gamma([0,1])$ such that $U\subseteq \mathfrak{M}$ and $F|_{U} = \Delta$.
\end{lemma}

\begin{proof} Since $\Delta$ is a surjective isometry, for each $\delta >0$, $a\in \mathfrak{M}$, we have $$\Delta\left( B(a,\delta) \cap \mathfrak{M} \right) =B(\Delta(a),\delta) \cap \mathfrak{N}.$$ Since $\mathfrak{M}$ and $\mathfrak{N}$ are open we can find $\delta_a>0$ satisfying $B(a,\delta_a)\subset \mathfrak{M}$, $B(\Delta(a), \delta_a)\subset \mathfrak{N}$ and $\Delta(B(a,\delta_a)) \subseteq B(\Delta(a), \delta_a)$. Clearly, the restricted mapping $$\Delta|_{B(a,\delta_a)}: B(a,\delta_a) \to B(\Delta(a), \delta_a)$$ is a surjective isometry. The existence of a surjective affine isometry $F_{a, \delta_a} : M\to N$ extending the mapping $\Delta|_{B(a,\delta_a)}$ (and thus proving that the latter is an affine mapping) follows from Mankiewicz theorem (Theorem~\ref{t Mankiewicz}).\smallskip 
	
In order to prove the final statement we employ an argument based on compactness. Suppose $\gamma : [0,1]\to \mathfrak{M}$ is continuous path. By the conclusion in the first part, for each $t\in [0,1],$  there exist $\delta_t>0$ and a surjective affine isometry $F_{t} : M\to N$ satisfying $B(\gamma(t),\delta_t)\subseteq \mathfrak{M}$, $B(\Delta(\gamma(t)), \delta_t)\subseteq \mathfrak{N}$, $\Delta(B(\gamma(t),\delta_t)) = B(\Delta(\gamma(t)), \delta_t),$ and $\Delta|_{B(\gamma(t),\delta_t)} = F_{t}|_{B(\gamma(t),\delta_t)}$. \smallskip

By applying that $\gamma([0,1])$ is compact we deduce the existence of $t_0=0<t_1<\ldots<t_{n-1}<t_n =1,$ such that $\displaystyle \gamma ([0,1]) \subseteq \bigcup_{k=0}^n B(\gamma(t_{k}),\delta_{t_k})$. The existence of the surjective affine isometry $F$ will be clear if we show that the corresponding maps $F_{t_{j}}$ satisfy $F_{t_{k}} = F_{t_{l}}$ for all $k,l\in \{0,\ldots,n\}$. A very basic identity principle assures that two surjective affine isometries between normed spaces are the same if and only if they coincide on an open subset. So, if the open set $B(\gamma(t_{k}),\delta_{t_k})\cap B(\gamma(t_{l}),\delta_{t_l})$ is non-empty we clearly have $F_{t_{k}} = F_{t_{l}}$, otherwise, by the connectedness of $\gamma([0,1])$, we can find a finite collection $t_k= t_{j_1},\ldots,t_{j_m}=t_{l}\in \{t_0=0,t_1,\ldots,t_{n-1},t_n\}$ such that $$B(\gamma(t_{j_k}),\delta_{t_{j_k}})\cap B(\gamma(t_{j_{k+1}}),\delta_{t_{j_{k+1}}})\neq \emptyset,$$ for all $k\in \{1,\ldots, m-1\}$. It then follows that $F_{t_{k}} = F_{t_{j_1}} = F_{t_{j_{2}}} = \ldots= F_{t_{j_m}} = F_{t_{l}}$, which gives the desired statement. To conclude the proof it suffices to take $F = F_{t_0}$ and the open neighbourhood of $\gamma ([0,1])$ given by $\displaystyle U = \bigcup_{k=0}^n B(\gamma(t_{k}),\delta_{t_k})$.
\end{proof}

The next tool, probably known in the folklore of affine maps, was designed ad hoc to play a specific role in the proof of Theorem~\ref{t surjective isometries between invertible clopen}. Here we present a more general version with a simplified proof. 

\begin{lemma}\label{l affine on segements} Let $X$ and $Y$ be real normed spaces with dim$(Y)\geq 2$. Suppose $G:Y \to X$ is a continuous mapping satisfying:
\begin{enumerate}[$(a)$]\item $G(0) =0;$
\item The restriction of $G$ to each segment not containing zero is an affine map.
\end{enumerate} Then $G$ is real-linear.
\end{lemma}

\begin{proof} It suffices to show that $G$ is affine on each segment containing $0$. If $0\in [x,y]$, the vectors $x,y$ are linearly dependent. We may assume that $x\neq 0$. Take $z\in Y$ such that $\{x,z\}$ is a linearly independent set. Clearly $0\notin [x,y]+\varepsilon z = [x+\varepsilon z,y+\varepsilon z]$ for all $\varepsilon>0$. By hypothesis $(b)$, the restriction of $G$ to $[x+\varepsilon z,y+\varepsilon z]$ is an affine map, meaning that $$G (t (x+\varepsilon z) + (1-t) (y+\varepsilon z) ) = t G ( (x+\varepsilon z)) + (1-t) G(y+\varepsilon z)$$ for all $t\in[0,1].$ Letting $\varepsilon\to 0$ we deduce from the continuity of $G$ that $G (t (x) + (1-t) (y) ) = t G ( x) + (1-t) G(y),$ for all $t\in[0,1].$
\end{proof}

Let us observe that the conclusion of the above lemma is not true when dim$(Y) =1$, the absolute value on $\mathbb{R}$ is a counterexample.\smallskip

Before entering into the details of the proof of Theorem~\ref{t surjective isometries between invertible clopen} we recall the definition and some basic facts on numerical ranges. By a \emph{numerical range space} we mean a (real or complex) Banach space $X$ with a distinguished element $u$ in its unit sphere. The \emph{state space} of the numerical range space $(X,u)$ is the defined as the set $$D(X) = D(X,u) = \{\phi \in X^* : \|\phi \| = \phi (u) =1 \},$$ and the \emph{numerical range} of an element $x\in X$ is the non-empty compact and convex set given by $$V(X,x) = V(x) =\{\phi (x) : \phi \in D(X,u)\}.$$ Finally, the \emph{numerical radius} of $x\in X$ is the number given by $$v(x) = \max\{ |\lambda| : \lambda\in V(x)\},$$ while the \emph{numerical index} of $X$ is defined by $$n(X,u)= n(X) = \inf\{ v(x) : x\in X, \ \|x\|=1 \} $$ $$= \max\{ \alpha \geq 0 : \alpha \|x\| \leq v(x) \hbox{ for all } x\in X\}.$$ It may occur that $n(X,u)=0$. Those elements $u$ in the unit sphere of $X$ satisfying $n(X, u) > 0$ are called \emph{geometrically unitaries}. The Bohnenblust--Karlin theorem is a milestone result in the theory of complex Banach algebras. 

\begin{theorem}[Bohnenblust--Karlin theorem]\label{t Bohnenblust--Karlin}{\rm(\cite[Theorem 2.6.4]{PalmerBook} and \cite[Proposition 2.1.11]{Cabrera-Rodriguez-vol1})} Let $A$ be a unital Banach algebra with unit $\11b{A}$ {\rm(}and $\|\11b{A}\|=1${\rm)}. Then the numerical radius of the numerical range space $(A,\11b{A})$ is a norm on $A$ which is equivalent to the original norm of $A$. Moreover, $n(A,\11b{A}) \geq \frac1e$ and thus $v(a) \leq \|a\|\leq e \ v(a)$ for all $a\in A$. Furthermore, the same conclusion holds when $A$ is replaced with a $M$ is a unital (non-necessarily associative) normed complex algebra $M$ with $\|\11b{A}\|=1$. 
\end{theorem}  

We proceed next to present the promised proof of Theorem~\ref{t surjective isometries between invertible clopen}. 

\begin{proof}[Proof of Theorem~\ref{t surjective isometries between invertible clopen}] $(a)$ Let $w_0$ be an element in the McCrimmon radical of $N$.  Lemma~\ref{l N-u0 in N} proves that $\mathfrak{N}-w_0 = \mathfrak{N}$, and hence $a\mapsto \Delta(a)-w_0$ is a surjective isometry from $\mathfrak{M}$ onto $\mathfrak{N}$.\smallskip
	
$(b)$ Observe first that $0$ always lies in the norm adherence of $\mathfrak{M}$. Take any 
sequence $(a_n)_n\subseteq \mathfrak{M}$ converging to $0$ in norm. Since $\Delta$ is an isometry, the sequence $(\Delta(a_n))_n$ is Cauchy in $N$, and hence there exists $u_0^a\in N$, depending on $(a_n)_n$, such that $(\Delta(a_n))_n\to u_0^{a}.$ Actually the limit $u_0^{a}$ does not depend on the chosen sequence.  Namely, if $(b_n)_n\subseteq \mathfrak{M}$ with $(\Delta(b_n))_n\to u_0^{b},$ by the assumptions on $\Delta$, we have $$\begin{aligned}\|u_0^{a}-u_0^{b}\| &\leq \|u_0^{a}- \Delta(a_n) \| + \|\Delta (b_n)- \Delta(a_n) \| + \| \Delta(b_n) - u_0^{b} \| \\
&= \|u_0^{a}- \Delta(a_n) \| + \| b_n- a_n \| + \| \Delta(b_n) - u_0^{b} \|, 
\end{aligned}$$ which implies that $u_0^{a}=u_0^{b}$ because the sequences on the right hand side of the above inequality tend to zero. We have therefore proven the existence of the limit $\displaystyle \lim_{\mathfrak{M}\ni a\to 0} \Delta (a) = u_0\in N.$ \smallskip

We shall finally show that $u_0$ lies in the McCrimmon radical of $N$. The key result is the characterization given in Proposition~\ref{p Hatori McCrimmon radical}. By the just quoted result, it suffices to show that $Sp(U_b(u_0))\neq \{0\}$ for every $b \in N_{\mathbf{1}}^{-1}\subseteq \mathfrak{N}$ (observe that under our assumptions $\mathfrak{M}$ must coincide with the union of several connected components of $M^{-1}$ including $M^{-1}_{\mathbf{1}}$ among them). \smallskip

Arguing by contradiction, we assume the existence of $b\in N_{\mathbf{1}}^{-1}\subseteq \mathfrak{N}$ and a non-zero $\lambda \in Sp(U_b(u_0))$. By applying that $\mathfrak{N}$ is a quadratic subset of $N^{-1}$ with $\11b{N}\in \mathfrak{N}$, we have $\alpha^2 c = U_{\alpha \11b{N}} (c) \in \mathfrak{N}$ for all $\alpha\in \mathbb{C}\backslash\{0\}$, $c\in \mathfrak{N}$, that is, $\left(\mathbb{C}\backslash\{0\}\right) \mathfrak{N} = \mathfrak{N}$. Now, having in mind that $\mathfrak{N}$ is closed for powers and inverses, and the previous comment, we get  $-\lambda b^{-2}\in \mathfrak{N} = \Delta\left(\mathfrak{M} \right)$. Therefore, there exists (a unique) $c_{\lambda}\in \mathfrak{M}$ such that $\Delta(c_{\lambda}) = - \lambda b^{-2}$. Actually, by Lemma~\ref{l principal component as the least quadratic subset of invertible elements} and Remark~\ref{r the principal component is closed for inverses} we know that $-\lambda b^{-2}\in {N}^{-1}_{\mathbf{1}}$.\smallskip 

By the arguments above, for each $0< s<\frac12,$ the segment $[s c_{\lambda},(1-s) c_{\lambda}]$ lies in $\mathfrak{M}$. By Lemma~\ref{l Delta is locally affine} we can find an open neighbourhood $U$ of $[s c_{\lambda},(1-s) c_{\lambda}]$ with $U \subseteq \mathfrak{M}$ and a surjective affine isometry $F: M\to N$ such that $\Delta|_{U} = F|_{U}$. The properties of these maps assure that  \begin{equation}\label{eq identities before limit in s} \begin{aligned}\Delta \left(\frac{c_{\lambda}}{2}\right) &= F \left(\frac12 s c_{\lambda} + \frac12 (1-s) c_{\lambda}\right) = \frac{F((1-s) c_{\lambda}) +F (s c_{\lambda})}{2} \\
		&= \frac{\Delta((1-s) c_{\lambda}) +\Delta (s c_{\lambda})}{2}.
	\end{aligned}
\end{equation} 
Having in mind that $\displaystyle \lim_{\mathfrak{M}\ni a\to 0} \Delta (a) = u_0$ and the continuity of $\Delta$, we allow $s\to 0$ in \eqref{eq identities before limit in s} to deduce that $$ \Delta \left(\frac{c_{\lambda}}{2}\right) = \frac{\Delta(c_{\lambda}) + u_0}{2} = \frac{- \lambda b^{-2} + u_0}{2},$$ and consequently,  \begin{equation}\label{eq lambda not in the spectrum} - \lambda \11b{N} + U_{b} (u_0) = -\lambda U_b (b^{-2}) + U_{b} (u_0) = U_b \left( - \lambda b^{-2} + u_0 \right)= 2 U_b \left( \Delta \left(\frac{c_{\lambda}}{2}\right)\right). 
\end{equation}
Since $b$ and $ \Delta \left(\frac{c_{\lambda}}{2}\right) $ belong to $\mathfrak{N},$ and the latter is a quadratic subset of $N^{-1}$, the element $2 U_b \left( \Delta \left(\frac{c_{\lambda}}{2}\right)\right)$ lies in $U_{\mathfrak{N}} (\mathfrak{N})\subseteq N^{-1}.$ The identity in \eqref{eq lambda not in the spectrum} implies that $\lambda\notin Sp (U_b(u_0)),$ which contradicts our initial assumption. \smallskip

$(c)$ Let $u_0=\displaystyle \lim_{\mathfrak{M}\ni a\to 0} \Delta (a) \in \mathcal{R}ad (N)$ be the element given by $(b)$. The conclusion in $(a)$ implies that the mapping $\Delta_0 : \mathfrak{M}\to \mathfrak{N},$ $\Delta_0 (a) = \Delta (a) - u_0$ is a surjective isometry. Furthermore, $\displaystyle \exists\lim_{\mathfrak{M}\ni a\to 0} \Delta_0(a) = 0.$\smallskip

The set $\mathfrak{M}$ is not a linear subspace. However, as we have seen in the proof of $(b)$, $\left(\mathbb{C}\backslash\{0\}\right) \mathfrak{M} = \mathfrak{M}$. In particular, for each $a\in\mathfrak{M}$, $\left(\mathbb{C}\backslash\{0\}\right) a\subseteq  \mathfrak{M}$ and hence the mapping $\lambda \mapsto G_a(\lambda) = \Delta_0 (\lambda a)$ is well defined on $
\mathbb{C}\backslash\{0\}$. We extend the definition to the whole complex plane by setting $G_a(0)=0$. 
\begin{equation}\label{eq Delta0 is real-linear} \hbox{We claim that for each $a\in \mathfrak{M}$ the mapping } G_a :\mathbb{C}\to N \hbox{ is real-linear.}
\end{equation} $G_a$ is continuous since $\displaystyle \lim_{\mathfrak{M}\ni a\to 0} \Delta_0(a) = 0$ and $\Delta_0$ is an isometry. Moreover, if $[\lambda,\mu]$ is a segment in $\mathbb{C}$ not containing zero, Lemma~\ref{l Delta is locally affine} applied to $\Delta_0$ on $[\lambda a, \mu a ]$, implies that $\Delta_0 |_{[\lambda a ,\mu a]}$ is an affine mapping (note that $[\lambda a,\mu a]$ is even convex). Therefore, $G_a$ is affine on each segment not containing $0$. Lemma~\ref{l affine on segements} now gives the desired statement.\smallskip

We have actually proved a better property for $\Delta_0$. Suppose $a,b\in \mathfrak{M}$ with $[a,b]\subset \mathfrak{M}$. A new application of Lemma~\ref{l Delta is locally affine} implies that $\Delta_0|_{[a,b]}$ is affine. By combining this property with \eqref{eq Delta0 is real-linear} we obtain \begin{equation}\label{eq Delta0 is additive on elements whose segment is in mathfrakM} \frac12 \Delta_0 (a+b) = \Delta_0\left(\frac{a+b}{2}\right) = \frac{\Delta_0(a)+\Delta_0(b)}{2},\ \forall a,b\in M \hbox{ with } [a,b]\subset \mathfrak{M}.
\end{equation}

For the final step we consider the set $$\Omega = \bigcup_{\alpha>0} B(\alpha \11b{M}, \alpha) = \{x\in M : \|x - \alpha \11b{M}\|< \alpha\hbox{ for some } \alpha >0\} \subseteq M_{\mathbf{1}}^{-1}\subset \mathfrak{M}.$$ Clearly $\Omega$ is open since it is the union of a family of open balls. Given $x,y\in \Omega$ there exist positive $\alpha,\beta $ such that $\|x - \alpha \mathbf{1}\|< \alpha$ and $\|y - \beta \mathbf{1}\|< \beta$, and consequently, for each $0<t<1$ we have $$\left\|
t x + (1-t) y - (\alpha t-\beta (1-t))\11b{M} \right\|\leq t \left\|
x - \alpha \11b{M} \right\| + (1-t) \left\|y - \beta \11b{M} \right\| < t \alpha t + (1-t) \beta, $$ which implies that $\Omega$ is convex. The set $\Omega$ also enjoys the following property: $\mathbb{R}^+ \Omega = \Omega$ (just observe that for $t>0$, we have $\|x-\alpha \11b{M}\|<\alpha $ $\Rightarrow$ $\|t x- t \alpha \11b{M}\|<t \alpha$).\smallskip

The next step consists in showing that $\Delta_0 (\Omega)$ also is convex. Given $a,b\in \Omega$ and $0<t<1$, the elements $t a $ and $(1-t) b$ belong to $\Omega$, and hence $[t a , (1-t) b]\subset \Omega$ by convexity. Applying \eqref{eq Delta0 is real-linear} we arrive to $$\begin{aligned}
 t \Delta_0 (a) + (1-t) \Delta_0 (b) &= t G_a (1) + (1-t) G_b (1)= G_a (t) + G_b (1-t) \\
 & = \Delta_0 ( t a)  +\Delta_0( (1-t) b).
\end{aligned} $$ Having in mind that $[t a , (1-t) b]\subset \Omega\subseteq \mathfrak{M},$ the conclusion in \eqref{eq Delta0 is additive on elements whose segment is in mathfrakM}  leads to 
$\Delta_0 ( t a)  +\Delta_0( (1-t) b) =  \Delta_0 ( t a + (1-t) b),$ and thus $ t \Delta_0 (a) + (1-t) \Delta_0 (b) =  \Delta_0 ( t a + (1-t) b).$  This guarantees that $\Delta_0 (\Omega)$ is convex. \smallskip 

The mapping $\Delta_0|_{\Omega} : \Omega \to \Delta_0(\Omega)$ is a surjective isometry between two open convex sets. Similar arguments to those in the proof of $(b)$, show that $\Delta_0|_{\Omega}$ admits an extension to a surjective isometry from the closure of $\Omega$ onto the closure of $\Delta_0 (\omega)$.  Mankiewicz theorem (see Theorem~\ref{t Mankiewicz}) applied to the extension of $\Delta_0$ to the corresponding norm closures, implies the existence of a surjective affine isometry $T_0 : M\to N$ such that $\Delta_0|_{\Omega} = T_0|_{\Omega}$.  As in $(b)$, $\displaystyle 0= \lim_{\Omega\ni a\to 0} \Delta_0(a) = \lim_{\Omega\ni a\to 0} T_0(a) = T_0 (0)$. Therefore,  $T_0$ is a surjective real-linear isometry.\smallskip

The desired statement in $(c)$ will be proved if we show that $T_0|_{\mathfrak{M}} = \Delta_0$. Pick an arbitrary $a\in \mathfrak{M}$. Since $\|\pm a +2 \|a \| \11b{M} - 2 \|a\| \11b{M} \| = \|a\|< 2 \|a\|,$ we deduce that  $\pm a +2 \|a \| \11b{M}\in \Omega,$ and thus $T_{0} (\pm a +2 \|a \| \11b{M}) = \Delta_0 (\pm a +2 \|a \| \11b{M})$. \smallskip

Now, having in mind that $T_0$ is a surjective real-linear isometry, $\Delta_0$ is an isometry and $\Delta_0|_{\Omega} = T_0|_{\Omega}$, we arrive to \begin{equation}\label{eq final identity 1 0711} \begin{aligned} 2 \|a\| &= \| a+ 2\|a\| \11b{M} - a  \| = \|\Delta_0( a+ 2\|a\| \11b{M} ) -\Delta_0 (a) \| \\
&= \|  a+ 2\|a\| \11b{M} -T_0^{-1} \Delta_0 (a) \|,
\end{aligned}		
\end{equation} and if we also apply \eqref{eq Delta0 is real-linear} we get
\begin{equation}\label{eq final identity 2 0711} \begin{aligned}
		2 \|a\| &= \| - a+2 \|a\| \11b{M} + a  \| = \|\Delta_0( - a+2 \|a\| \11b{M}) -\Delta_0 (-a) \| \\ &= \| - a+2 \|a\| \11b{M} -T_0^{-1} \Delta_0 (-a) \|= \| - a+2 \|a\| \11b{M} +T_0^{-1} \Delta_0 (a) \|.
	\end{aligned}
\end{equation} Let us take a state $\phi\in D(M,\11b{M}) = \{\phi \in M^* : \|\phi \| = \phi (\11b{M}) =1 \}$. Now, combining \eqref{eq final identity 1 0711} and \eqref{eq final identity 2 0711} we deduce that $$\begin{aligned}|2 \|a\| \pm \phi (a -T_0^{-1} \Delta_0 (a) )| &= |\phi (2 \|a\| \11b{M} \pm ( a -T_0^{-1} \Delta_0 (a)))| \\ &\leq \| 2 \|a\| \11b{M} \pm (a - T_0^{-1} \Delta_0 (a)) \|= 2 \ \|a\|,
\end{aligned}$$ which proves that $\phi (a -T_0^{-1} \Delta_0 (a) ) =0$. The arbitrariness of $\phi$ proves that the numerical radius of $a -T_0^{-1} \Delta_0 (a)$ is zero, that is,  $v(a -T_0^{-1} \Delta_0 (a))=0$. The Bohnenblust--Karlin theorem (see Theorem~\ref{t Bohnenblust--Karlin}) implies that $a =T_0^{-1} \Delta_0 (a),$ and hence $\Delta (a) = T_0 (a) + u_0$ for all $a\in \mathfrak{M}$, as desired.
\end{proof}     

\section{The geometric goodness of C$^*$- and JB$^*$-algebras}\label{sec: geometric goodness of C* and JB-algebras}

C$^*$-algebras form a very well known and deeply studied model in mathematics, particularly in functional analysis and physic, whose origins go back to their use in quantum mechanics too. A C$^*$-algebra is a complex Banach algebra $A$ equipped with an algebra involution $a\mapsto a^*$ satisfying the Gelfand--Naimark axiom $\|aa^* \| = \|a\|^2$ for all $a\in A$. The Gelfand--Naimark theorem (1943) affirms that each C$^*$-algebra admits a concrete representation as a norm closed self-adjoint subalgebra of the space $B(H)$ of all bounded linear operators on a complex Hilbert space $H$. These objects have been intensively studied since the forties of the last century, a wide collection of books and structure results have been completely devoted to study them from many different points of view. For the sake of brevity we shall not employ too much time on introducing examples nor basic results, which seem to be widely known. We refer to the monographs \cite{BonDunCNA73, KR1,PalmerBook,Ped, RickartBook,Tak} as reference sources on C$^*$-algebras.

It was soon realised, since the early stages of the theory of C$^*$-algebras, that these objects enjoy very special geometric properties. We can trace back the good geometric properties in a result earlier on time. The Banach--Stone theorem asserts that given two compact Hausdorff spaces $K_1$ and $K_2$ and a surjective linear isometry $T : C(K_1) \to C(K_2)$, then there exists a homeomorphism $\varphi : K_2 \to K_1$ and a unimodular function $u \in C(K_2)$ satisfying $$(Tf) (s) =u (s) f(\varphi (s)) \hbox{ for all } s\in K_2,\ f\in C(K_1).$$ The composition operator $f\mapsto f\circ \varphi$ is a $^*$-isomorphism from $C(K_1)$ onto $C(K_2)$. An appropriate version of the Banach--Stone theorem works for non-unital commutative C$^*$-algebras (i.e. $C_0(L)$ spaces where $L$ is a locally compact Hausdorff space, see \cite[Corollary 2.3.12]{FleJa2003v1}). In particular, the counterexample in Example~\ref{example the Banach space structure is not an invariant} can never occur when both algebras are commutative C$^*$-algebras. This is just the starting point of a series of deep geometric results. The first step was the following result proved by Kadison in the unital case, and by Paterson and Sinclair in the non-unital setting. We first recall that the \emph{multiplier algebra}, $M(A),$ of a C$^*$-algebra $A$ is the idealizer of $A$ in its bidual $A^{**},$ that is, the largest
C$^*$-subalgebra of the von Neumann algebra $A^{**}$ containing $A$ as an ideal). In other words $M(A) =\left\{ x\in A^{**} : x A, A x\subseteq A\right\}$ \cite[\S 3.12]{Ped}.

\begin{theorem}\label{t Kadison 55}{\rm(\cite[Theorem 7]{Kad51} and \cite[Theorem 1]{PatSin72})} Let $A$ and $B$ be C$^*$-algebras. Then for each surjective linear isometry $T: A\to B$ there exists a unitary element $u$ in the multiplier algebra of $B$ inside $B^{**}$, and a Jordan $*$-isomorphism $\Phi :A\to B$  such that $$T (x) = u \ \Phi (x), \ (x\in A).$$
\end{theorem}

The above result shows that two C$^*$-algebras are Jordan isomorphic if and only if their underlying Banach spaces are isometrically isomorphic, so there is no place for counterexamples like Example~\ref{example the Banach space structure is not an invariant} in the setting of C$^*$-algebras.\smallskip

Theorem~\ref{t Kadison 55} also implies that  every surjective linear isometry between C$^*$-algebras preserves triple products of the form \begin{equation}\label{eq triple produt Cstar algebras} \{ x,y,z\} = 2^{-1} ( x y^* z + z y^* x).
\end{equation} Another consequence of the theorem proves that every unital surjective isometry between unital C$^*$-algebras is a Jordan $^*$-isomorphism. \smallskip

Hatori and Watanabe explored in \cite{HatWan2012} an appropriate version of Theorem~\ref{thm Hatori 2011} in the case that $A$ and $B$ are unital C$^*$-algebras and they found that a stronger conclusion is valid under these extra hypotheses. It should be recalled that every C$^*$-algebra is semisimple (cf. \cite[Propositions 2.3.17 and 2.3.18]{PalmerBook}). So, the mapping $\Delta$ in Theorem~\ref{thm Hatori 2011} extends to a surjective real-linear isometry $T$ from $A$ into $B$. We shall see that $T$ admits a particular form.

\begin{theorem}\label{t HatWat}{\rm(\cite[Theorem 2.2]{HatWan2012})} Let $A$ and $B$ be unital C$^*$-algebras, and let $\mathfrak{A}$ and $\mathfrak{B}$ be open subgroups of $A^{-1}$ and $B^{-1},$ respectively. Suppose $\Delta$ is a bijection from $\mathfrak{A}$ onto $\mathfrak{B}$.
	Then $\Delta$ is an isometry if and only if $\Delta(\11b{A})$ is unitary in $B$ and there exist a projection $p$ in the centre of $B$, and a complex-linear Jordan $^*$-isomorphism $\tilde{J}$ from $A$ onto $B$ such that $$ \Delta (a) = \Delta(\11b{A}) p \tilde{J}(a) + \Delta(\11b{A}) (\11b{N}-p) \tilde{J} (a)^*, \hbox{ for all } a \in  \mathfrak{A}.$$ 
	
\noindent Furthermore the operator $ \Delta(\11b{A}) p \tilde{J}(\cdot) + \Delta(\11b{M}) (\11b{N}-p) \tilde{J}(\cdot)^*$ defines a surjective real-linear isometry from $A$ onto $B$.
\end{theorem} 

As in the previous section, we shall obtain Theorem~\ref{t HatWat} above from its version for JB$^*$-algebras that will be presented in Theorem~\ref{t surjective isometries between invertible clopen JBstar}.\smallskip  

JB$^*$-algebras are the Jordan alter ego of C$^*$-algebras. The notion was by proposed by Kaplansky during his celebrated final lecture to the 1976 St. Andrews Colloquium of the Edinburgh Mathematical Society. Concretely, a \emph{JB$^*$-algebra} is a complex Jordan-Banach algebra $M$ equipped with an algebra involution $^*$ satisfying the following version of the Gelfand--Naimark axiom: \begin{equation}\label{eq GN Jordan} \|\{ a,{a},a\} \|= \|a\|^3, \hbox{ for all  }a\in M,
\end{equation} (where $\{a,{a},a \}= U_a (a^*) = 2 (a\circ a^*) \circ a - a^2 \circ a^*$). A basic result in the theory of JB$^*$-algebras asserts that the involution of every JB$^*$-algebra is an isometry (cf. \cite[Lemma 4]{youngson1978vidav} or \cite[Proposition 3.3.13]{Cabrera-Rodriguez-vol1}).\smallskip

Among the starring examples of JB$^*$-algebras we find all C$^*$-algebra when they are equipped with the same norm and involution and the natural Jordan product. In such a case $U_a ({a^*}) = aa^* a,$ and hence $\| aa^* a\| =\|U_a ({a^*}) \|= \|a\|^3,$ which is an equivalent reformulation of the Gelfand--Naimark axiom. JB$^*$-subalgebras of C$^*$-algebras are called \emph{JC$^*$-algebras}. We should note the existence of exceptional JB$^*$-algebras which cannot be embedded as Jordan $^*$-subalgebras of a C$^*$-algebra (see, for example,  \cite[Corollary 2.8.5]{HOS}).\smallskip

JB$^*$-algebras are indissolubly linked with JB-algebras. A \emph{JB-algebra} is a real Jordan algebra $\mathfrak{J}$ equipped with a complete norm satisfying \begin{equation}\label{eq axioms of JB-algebras} \|a^{2}\|=\|a\|^{2}, \hbox{ and } \|a^{2}\|\leq \|a^{2}+b^{2}\|\ \hbox{ for all } a,b\in \mathfrak{J}.
\end{equation} Wright proved in \cite{Wri77} that every JB-algebra $\mathfrak{J}$ corresponds uniquely to the self-adjoint part $M_{sa}=\{x\in M : x^* =x\}$ of a JB$^*$-algebra $M$, answering in this way a conjecture by Kaplansky.\smallskip 

A JBW$^*$-algebra (respectively, a JBW-algebra) is a JB$^*$-algebra (respectively, a JB-algebra) which is also a dual Banach space. Every JBW$^*$-algebra (respectively, each JBW-algebra) contains a unit element (see \cite[\S 4]{HOS} or \cite{AlfsenShultz2003}). It is worth to note that JBW-algebras are precisely the self-adjoint parts of JBW$^*$-algebras (see \cite[Theorems 3.2 and 3.4]{Ed80JBW} or \cite[Corollary 2.12]{MarPe00}). In particular every von Neumann algebra is a JBW$^*$-algebra. \smallskip

Concerning surjective linear isometries, Wright and Youngson established that every unital and surjective linear isometry between two unital JB$^*$-algebras is a Jordan
$*$-isomorphism  \cite[Theorem 6]{WriYo}. A result by Isidro and Rodr\'{\i}guez Palacios proves that every surjective linear
isometry $T$ between two JB-algebras $\mathfrak{J}_1$ and $\mathfrak{J}_2$ is of the form $T (x) = s \circ \Phi(x),$ where $s$ is a central symmetry in the multiplier algebra of $\mathfrak{J}_2$ and $\Phi: \mathfrak{J}_1\to \mathfrak{J}_2$ is a surjective Jordan isomorphism (cf. \cite[Theorem 1.9]{IsRo}). We recall that an element $s$ in a unital JB-algebra is called a \emph{symmetry} if we can write it in the form $s = p-q$ with $p q =0$ and $p+q =\mathbf{1}$ (equivalently, $s^2 = \mathbf{1}$). We also recall that the \emph{(Jordan) multiplier algebra} of a JB$^*$-algebra (respectively, a JB-algebra) $M$ is defined as $$M(M):=\{x\in M^{**}: x\circ  M \subseteq  M\}.$$ According to the result by Isidro and Rodr\'{\i}guez Palacios, a bijective linear map $\Phi$ between two JB-algebras is an isometry if and only if
$\Phi$ is a triple-isomorphism with respect to the triple product \begin{equation}\label{eq triple product JBalgebras}
\{ x,y,z\} = (x \circ y) \circ z + (z \circ y) \circ x -
(x \circ z) \circ y.\end{equation}

There is no better way to understand surjective isometries between JB$^*$-algebras than when they are regarded inside the strictly wider class of JB$^*$-triples. The complex Banach spaces in this class are characterized by the holomorphic property that their open unit balls are bounded symmetric domains \cite{Ka83}, a property already observed by Harris for the open unit ball of each C$^*$-algebra \cite{Harris74}. The contribution of Kaup gives an axiomatic description of these complex Banach spaces.  A JB$^*$-triple is a complex Banach space $E$ equipped with a continuous triple product $\{.,.,.\} : E\times E\times E \to E,$ $(a,b,c)\mapsto \{a,b,c\},$ which is bilinear and symmetric in $(a,c)$ and conjugate linear in $b$,
and satisfies the following axioms for all $a,b,x,y\in E$:
\begin{enumerate}[{\rm (a)}] \item (Jordan identity) $$L(a,b) L(x,y) = L(x,y) L(a,b) + L(L(a,b)x,y)
	- L(x,L(b,a)y),$$ where $L(a,b):E\to E$ is the operator defined by $L(a,b) x = \{ a,b,x\};$
	\item $L(a,a)$ is a hermitian operator with non-negative spectrum;
	\item $\|\{a,a,a\}\| = \|a\|^3$.\end{enumerate}

The triple product of every JB$^*$-triple is a continuous mapping with norm $1$, that is,
\begin{equation}\label{eq triple product non-expansive} \|\{a,b,c\}\|\leq \|a\| \|b\| \|c\|,\ \hbox{ for all } a,b,c \hbox{ (cf. \cite[Corollary 3]{FriRu86}).}
\end{equation}

The class of JB$^*$-triples is rich of interesting examples, complex Hilbert spaces, spaces of operators of the form $B(H,K)$, C$^*$-algebras and JB$^*$-algebras are among them. The triple products for C$^*$-algebras and JB$^*$-algebras are given by  \begin{equation}\label{eq product operators} \{ x,y,z\} =\frac12 (x y^* z +z y^*
	x),\end{equation}  and \begin{equation}\label{eq product jordan}\{x,y,z\} = (x\circ y^*) \circ z + (z\circ y^*)\circ x -
	(x\circ z)\circ y^*, \end{equation} respectively (cf. \cite{Harris74} and \cite[Theorem 3.3]{BraKaUp78}). Note that we have already employed the triple product given in \eqref{eq product operators} after the Kadison--Paterson-Sinclair theorem in \eqref{eq triple produt Cstar algebras}. Additional examples can be found by recalling that each element in the widely studied class of Hilbert C$^*$-modules is a JB$^*$-triple (cf. \cite{Isidro2003,Kusu2014}).\smallskip

For each element $a$ in a JB$^*$-triple $E$, we shall write $Q(a)$ for the conjugate linear operator on $E$ defined by $Q(a) (x)= \{a,x,a\}$ ($x\in E$). Another important connection with JB$^*$-algebras and $U$ operators affirms that if $M$ is a JB$^*$-algebra regarded as a JB$^*$-triple with the triple product given in \eqref{eq product jordan}, for each $a\in M$, we have $Q(a) (x) = U_a (x^*)$ ($x\in M$).\smallskip

Let $E$ and $F$ denote two JB$^*$-triples. A linear mapping $T:E\to F$ is a \emph{triple homomorphism} if it preserves triple products. A triple homomorphism which is also a linear bijection is called a \emph{triple isomorphism}. If the corresponding maps only are real-linear we call them \emph{real-linear triple homomorphisms} and \emph{real-linear triple isomorphisms}, respectively.\smallskip

There is an undoubtedly attractiveness in the geometric properties of JB$^*$-triples expressed in the following Kaup's Banach--Stone theorem. 

\begin{theorem}\label{t Kaup-Banach-Stone}{\rm(\cite[Proposition 5.5]{Ka83} or \cite[Theorem 5.6.57]{Cabrera-Rodriguez-vol2})} 
Let $E$ and $F$ be JB$^*$-triples. Then a linear bijection $T: E\to F$ is an isometry if and only if it is a triple isomorphism.
\end{theorem} 

Alternative proofs and additional information can be found in \cite{DanFriRu1990}, \cite[Theorem 2.2]{FerMarPe2004}, \cite[Corollary 3.4]{BecPe2004}, and \cite{FerMarPe2004iso}. Kaup's Banach--Stone theorem reaffirms the deep connection between the triple product and the geometric properties of the underlying Banach space.    

In order to understand the implications of Kaup's Banach--Stone Theorem, we need to recall the notion of tripotent and some variants like unitary and complete tripotents. An element $e$ in a JB$^*$-triple $E$ is said to be a \emph{tripotent} if $\{e,e,e\}=e$. In such a case, the tripotent $e$ produces a \emph{Peirce decomposition}\label{eq Peirce decomposition} of $E,$ $$E= E_{2} (e) \oplus E_{1} (e) \oplus E_0 (e),$$ in terms of the eigenspaces of the mapping $L(e,e)$, that is,  $$E_j (e)=\{ x\in E : \{e,e,x\} = \frac{j}{2}x \}\ \ (j\in\{0,1,2\}.$$ Triple products among elements in different Peirce subspaces follow the following patterns known under the name of \emph{Peirce arithmetic} or \emph{Peirce rules}. $$\begin{aligned} \{ {E_{i}(e)},{E_{j} (e)},{E_{k} (e)}\} &\subseteq E_{i-j+k} (e)\  \hbox{ if $i-j+k \in \{ 0,1,2\},$ }\\
	\{ {E_{i}(e)},{E_{j} (e)},{E_{k} (e)} \}& =\{0\} \hbox{ if $i-j+k \notin \{ 0,1,2\},$ }
\end{aligned}$$ and $\{ {E_{2} (e)},{E_{0}(e)},{E} \} = \{ {E_{0} (e)},{E_{2}(e)},{E} \}=0.$ Consequently, each Peirce subspace $E_j(e)$ is a JB$^*$-subtriple of $E$.\smallskip

The natural projecton of $E$ onto $E_{k} (e)$ is called the \emph{Peirce $k$-projection}, and is denoted by $P_{k_{}}(e)$. Among their properties, we highlight that Peirce projections are contractive (cf. \cite[Corollary 1.2]{FriRu85}), and can be algebraically recovered by the following identities: $P_{2}(e) = Q(e)^2,$ $P_{1}(e) =2(L(e,e)-Q(e)^2),$ and $P_{0}(e) =\hbox{Id}_E - 2 L(e,e) + Q(e)^2,$ where $Q(e):E\to E$ is the conjugate linear map defined by $Q(e) (x) =\{e,x,e\}$.\smallskip

In case that a C$^*$-algebra $A$ is regarded  as a JB$^*$-triple respect to the triple product in \eqref{eq product operators}, the identity $e=\{e,e,e\} = e e^* e$ is an equivalent reformulation of the fact that $e$ is a partial isometry. Therefore partial isometries and tripotents in $A$ are the same elements. In this particular case, Peirce projections are given by the left and right projections of $e$, that is, $P_2 (e) (a) = ee^* a e^*e$, $P_1 (e) (a) = (\mathbf{1}-ee^*) a e^*e + ee^* a (\mathbf{1}-e^*e),$ and $P_0 (e) (a) = (\mathbf{1}-ee^*) a (\mathbf{1}-e^*e)$ ($a\in A$).\smallskip

For our purposes here, we note that for each tripotent $e$ in a JB$^*$-triple $E$, the Peirce 2-subspace $E_2 (e)$ is a unital JB$^*$-algebra with unit $e$,
Jordan product $x\circ_e y := \{ x,e,y\}$ and involution $x^{*_e} := \{e,x,e\}$, respectively (see \cite[Theorem 4.1.55]{Cabrera-Rodriguez-vol1}).\label{eq Peirce-2 is a JB-algebra} 
\smallskip

Following standard notation, a tripotent $e$ in in a JB$^*$-triple $E$ is called \emph{unitary} (respectively, \emph{complete} or \emph{maximal})
if $E_2 (e) = E$. (respectively, $E_0 (e) =\{0\}$). There exists complete tripotents which are not unitaries, for example, a non-surjective isometry in $B(H)$. 
For each unitary $u$ in $E$, the space $E$ is a unital JB$^*$-algebra with unit $u,$ Jordan product $\circ_{u}$, and involution $*_{u}$.  
\smallskip

We could have produce a certain polysemy for the word unitary since the term is already in use in unital C$^*$-algebras and unital JB$^*$-algebras. We recall that an element $u$ in a unital C$^*$-algebra $A$ is called \emph{unitary} if $u u^* = u^* u =\11b{A}$, equivalently, $u$ is invertible with inverse $u^*$. The latter reformulation literally makes sense in unital JB$^*$-algebras.  Since the class of JB$^*$-triples englobes all unital JB$^*$-algebras, we might have a possible conflict of meanings. Fortunately, the concerns are completely clarified by a result of Braun, Kaup, and Upmeier showing that unitary elements in a unital JB$^*$-algebra $M$ are precisely the unitary tripotents in $M$ when the latter is regarded as a JB$^*$-triple {\rm(}cf. \cite[Proposition 4.3]{BraKaUp78} or \cite[Theorem 4.2.24, Definition 4.2.25 and Fact 4.2.26]{Cabrera-Rodriguez-vol1}{\rm)}. The set of all unitary elements in a unital C$^*$-algebra $A$ or in a unital JB$^*$-algebra $M$ will be denoted by $\mathcal{U} (A)$ and  $\mathcal{U} (M),$ respectively. Note that given a unitary $u$ in a unital JB$^*$-algebra $M$, the unital JB$^*$-algebra $(M,\circ_{u},*_{u})$ is precisely the $u^*$-isotope $M_{(u^*)}$ of $M$ defined in page \pageref{def of u-isotope}. The next lemma also is a straight consequence of \cite[Proposition 4.3]{BraKaUp78}.\label{comments unitaries as JB*-algebra an triple coincide}

\begin{lemma}\label{l unitaries in the homotope coincide}\cite[Proposition 4.3]{BraKaUp78} Let $M$ be a unital JB$^*$-algebra and let $u$ be a unitary element in $M$. The unitary elements of the JB$^*$-algebras $M$ and $M_{(u^*)}=(M,\circ_u,*_u)$ are the same, and they are precisely the unitary tripotents of $M$ when the latter is regarded as a JB$^*$-triple. 
\end{lemma}

A \emph{JBW$^*$-triple} is a {JB$^*$-triple} which is also a dual Banach space. Examples of JBW$^*$-triples include von Neumann algebras and JBW$^*$-algebras. More examples can be obtained by having in mind that the second dual of a JB$^*$-triple $E$ is a JBW$^*$-triple under a triple product extending the product of $E$ \cite{Di86}. Among the additional properties of JBW$^*$-triples we note that every JBW$^*$-triple admits a unique isometric predual, and its triple product is separately weak$^*$ continuous \cite{BarTi86}. The celebrated Dixmier-Ng theorem\label{comments on Dixmier-Ng theorem} characterizes when a normed space is in fact a dual Banach space. Concretely, a normed space $X$ is a dual Banach space if and only if there exists a Hausdorff locally convex topology $\tau$ on $X$ so that the closed unit ball of $X$ is $\tau$-compact. It follows that if two JB$^*$-triples are isometrically isomorphic as Banach spaces (equivalently, JB$^*$-triple isomorphic by Theorem~\ref{t Kaup-Banach-Stone}), then one is a JBW$^*$-triple if and only if the other one is.\smallskip

Let us go back to surjective linear isometries between JB$^*$-algebras. Suppose $T:M\to N$ is a surjective linear isometry between two JB$^*$-algebras. Clearly, $T^{**}: M^{**}\to N^{**}$ is a surjective linear isometry between two JBW$^*$-algebras (which are always unital, cf. \cite[Proposition 2.7]{AlfsenShultz2003}). Kaup's Banach--Stone theorem (Theorem~\ref{t Kaup-Banach-Stone}) assures that $T^{**}$ is a triple isomorphism and hence $T^{**} (\11b{M^{**}})=u$ must be a unitary in $N^{**}$ and $T^{**}: M\to (N^{**}, \circ_{u},*_{u})$ is a unital Jordan $^*$-isomorphism (cf. \cite{WriYo}). The restriction $T^{**}|_{M} :M\to (N,\circ_u,*_{u})$ is a Jordan $^*$-isomorphism. We have obtained a Jordan $^*$-isomorphism at the cost of replacing the replacing the Jordan product and involution in the JB$^*$-algebra $N$. \smallskip

Since Theorem~\ref{t surjective isometries between invertible clopen} only assures the existence of a real-linear surjective isometry extending the mapping $\Delta$, we shall also require the following result. Henceforth, we shall say that a JB$^*$-algebra $M$ is the (orthogonal) direct sum of two JB$^*$-subalgebras $M_1$ and $M_2$ (denoted by $M= M_1\oplus^{\ell_{\infty}} M_2$) if $M = M_1\oplus M_2$ and $a b = 0$ for all $a\in M_1$, $b\in M_2$.\smallskip

A more precise notion of orthogonality between general elements in a JB$^*$-triple has been widely study in the setting of preservers (cf. \cite{BurGarPe2011C,BurGarPe2011JB,OikPeRa2011,OikPePu2013,OikPe2013}). Elements $a,b$ in a JB$^*$-triple $E$ are called \emph{orthogonal}\label{def orthogonality} (written $a\perp b$) if $L(a,b) =0$. It is known that $a\perp b$ $\Leftrightarrow$ $\{ a,a,b\} =0$ $\Leftrightarrow$ $\{b,b,a\}=0$ $\Leftrightarrow$ $b\perp a$ (see \cite[Lemma 1]{BurFerGarMarPe08}).

\begin{theorem}\label{t Dang real-linear surjective isometries JBstar algebras}{\rm(\cite[Corollary 3.2]{Dang92} or by \cite[Corollary 3.4]{FerMarPe2004iso})} Every surjective real-linear isometry $T$ between JB$^*$-algebras $M$ and $N$ is a real-linear triple homomorphism. Moreover, there exist two JB$^*$-subalgebras $M_1$ and $M_2$ of $M$ such that $M$ is the (orthogonal) direct sum $M= M_1\oplus^{\ell_{\infty}} M_2$, $T|_{M_1}  : M_1 \to N$ is a complex-linear triple homomorphism and $T|_{M_2}  : M_2 \to N$ is a conjugate-linear triple homomorphism. 
\end{theorem}

For each Banach space $X$ we shall write $Iso(X)$ for the group of all surjective linear isometries on $X$. It is well known that for each unital C$^*$-algebra $A$ the group $Iso(A)$ acts transitively on $\mathcal{U} (A)$. Namely, if $u,v\in \mathcal{U} (A)$, the mapping $x\mapsto vu^*x$ is a surjective linear isometry on $A$ mapping $u$ to $v$. At this point we must remark one of the main divergences in C$^*$-algebra theory and JB$^*$-algebra theory. When $A$ is replaced with a unital JB$^*$-algebra $M$, this conclusion is no longer true. 

\begin{antitheorem}\label{anti-theorem BKU}{\rm(\cite[Example 5.7]{BraKaUp78}) and \cite[Antitheorem 3.4.34]{Cabrera-Rodriguez-vol1}}  There exists a unital JC$^*$-algebra $M$ containing two unitaries $u_1$ and $u_2$ satisfying $T (u_1)\neq u_2$ for all $T\in Iso (M)$. Furthermore, the unital JC$^*$-algebras $M_{(u_1^*)}= (M, \circ_{u_1},*_{u_1})$ and $M_{(u_2^*)}= (M, \circ_{u_2},*_{u_2})$ are linearly isometric but not Jordan $^*$-isomorphic.
\end{antitheorem}

We have seen after Kadison--Paterson--Sinclair theorem that there is no place for counterexamples like Example~\ref{example the Banach space structure is not an invariant} in the setting of C$^*$-algebras. The previous anti-theorem shows that these pathologies are possible for JB$^*$-algebras.\smallskip

Anti-theorem~\ref{anti-theorem BKU} also proves that the same complex Banach space $M$ may admit two different structures of JB$^*$-algebras which are not Jordan $^*$-isomorphic. However, by Kaup's Banach--Stone theorem (Theorem~\ref{t Kaup-Banach-Stone}) there is only one structure of JB$^*$-triple. 

\begin{corollary}\label{c uniqueness of triple product} Let $u$ be a unitary in a unital JB$^*$-algebra $M$, then there is a unique triple product $\{.,.,\}$ on $M$ making it a JB$^*$-triple, and it satisfies $$\begin{aligned} \{x,y,z\} &= (x\circ y^*) \circ z + (z\circ y^*)\circ x -
		(x\circ z)\circ y^*\\
		&= (x\circ_{u} y^{*_{u}}) \circ_{u} z + (z\circ_{u} y^{*_{u}})\circ_{u} x -
		(x\circ_{u} z)\circ_{u} y^{*_{u}},
	\end{aligned}$$ for all $x,y,z\in M$.
\end{corollary} 

We are now in a position to determine the concrete form of a surjective isometry between two clopen quadratic subsets of the invertible elements in two unital JB$^*$-algebras which are closed for powers and inverses.

\begin{theorem}\label{t surjective isometries between invertible clopen JBstar}\cite[Theorem 4.1]{Pe2021} Let $M$ and $N$ be unital JB$^*$-algebras. Suppose that $\mathfrak{M}\subseteq M^{-1}$ and  $\mathfrak{N}\subseteq N^{-1}$ are clopen quadratic subsets which are closed for powers and inverses. Let $\Delta : \mathfrak{M}\to \mathfrak{N}$ be a surjective isometry. Then $\Delta(\11b{M}) =u$ is a unitary element in $N$ and there exists a projection $p$ in the centre of $M$ and a complex-linear Jordan $^*$-isomorphism $J$ from $M$ onto the $u^*$-homotope $N_{(u^*)} = (N,\circ_{u},*_u)$ such that $$\Delta (a) = J(p\circ a) + J ((\11b{M}-p) \circ a^*), \hbox{ for all } a\in \mathfrak{M}.$$
	
\noindent If we additionally assume that there exists a unitary $\omega_0$ in $N$ such that the identity $U_{\omega_0} (\Delta(\mathbf{1})) = \mathbf{1}$ holds, then there exist a central projection $p\in M$ and a complex-linear Jordan $^*$-isomorphism $\Phi$ from $M$ onto $N$ such that $$\Delta (a) = U_{w_0^{*}} \left(\Phi (p\circ a) + \Phi ((\11b{M}-p) \circ a^*)\right),$$ for all $a\in \mathfrak{M}$.
\end{theorem}

\begin{proof} We recall that every JB$^*$-algebra $M$ is (Jordan-)semisimple, that is $\mathcal{R}ad (M) =\{0\}$ (cf. \cite[Lemma 4.4.28$(iii)$]{Cabrera-Rodriguez-vol1}). Thus, by Theorem~\ref{t surjective isometries between invertible clopen} there exists a surjective real-linear isometry $T_0: M\to N$ such that $\Delta (a) = T_0(a) $ for all $a\in \mathfrak{M}$.  Theorem~\ref{t Dang real-linear surjective isometries JBstar algebras} implies that $T_0$ is a real-linear triple isomorphism, that is, $T_0$ preserves triple products of the form $\{a,b,c\} = (a\circ b^*)\circ c + (c\circ b^*) \circ a - (a\circ c) \circ b^*$ ($a,b,c\in M$). Actually, there exist two JB$^*$-subalgebras $M_1$ and $M_2$ of $M$ such that $M$ is the (orthogonal) direct sum $M= M_1\oplus^{\ell_{\infty}} M_2$, $T_0|_{M_1}  : M_1 \to N$ is a complex-linear triple homomorphism and $T_0|_{M_2}  : M_2 \to N$ is a conjugate-linear triple homomorphism. \smallskip
	
It is now easy to check that the element $u= T_0(\11b{M})$ is a unitary in $N$ (observe that $\{u,u,u\} = T_0 (\{\11b{M},\11b{M},\11b{M}\}) = T(\11b{M}) =u$, and, by the surjectivity of $T_0$, for each $z$ in $N$ we can find $x\in M$ with $T_0(x) =z,$ and hence$\{u,u,z\} = T_0(\{\11b{M},\11b{M},x\}) = T_0(x) = z$).\smallskip
	
As we have seen above the $u^*$-isotope $N_{(u^*)} = (N,\circ_u, *_u)$ is a unital JB$^*$-algebra with unit $u$. Since $$T_0(x\circ y ) = T_0(\{x,\11b{M},y\}) = \{T_0(x),T_0(\11b{M}),T_0(y)\} = T_0(x)\circ_{u^*} T_0(y)$$ and $$T_0(x^*) = T_0(\{\11b{M},x,\11b{M}\}) = \{T_0(\11b{M}),T_0(x),T_0(\11b{M})\} = T_0(x)^{*_u},\ (x,y\in M),$$ we deduce that $T_0: M\to N_{u^*}$ is an isometric real-linear unital Jordan $^*$-isomorphism.\smallskip
	
Let $p$ denote the unit of $M_1$.  Clearly, $p$ is a central  projection in $M$ and $\11b{M}-p$ is precisely the unit of $M_2$. The arguments in the previous paragraphs show that $T_0|_{M_1}  : M_1 \to N_u$ is a complex-linear Jordan $^*$-monomorphism and $T_0|_{M_2}  : M_2 \to N_u$ is a conjugate-linear Jordan $^*$-monomorphism. Furthermore, bearing in mind that $M_1$ and $M_2$ are orthogonal in $M$ and that $T_0: M\to N_u$ is a real-linear Jordan $^*$-isomorphism, we deduce that $N_1 = T_0(M_1)$ and $N_2 = T_0(M_2)$ are orthogonal JB$^*$-subalgebras of $N_{(u^*)}$ with $N_{(u^*)} = N_1\oplus^{\ell_{\infty}} N_2$. Moreover, the mapping $J: M\to N_{(u^*)}$, $$J(x) = T_0(p\circ x) + T_0((\11b{M}-p) \circ x)^{*_u} = T_0(p\circ x) + T_0((\11b{M}-p) \circ x^*)\ \ (x\in M)$$ is a complex-linear Jordan $^*$-isomorphism (just apply that $x\mapsto p\circ x$ and $x\mapsto (\11b{M}-p) \circ x$ are the natural projections of $M$ onto $M_1$ and $M_2$, respectively), and the identity $$\Delta (a) = J(p\circ a) + J ((\11b{M}-p) \circ a^*),$$ holds for all $a\in \mathfrak{M}$.\smallskip
	
To prove the final statement, we assume the existence of a unitary $w_0\in N$ such that $U_{w_0} (\Delta(\11b{M}))= \11b{N}$. Under this assumption, the mapping $U_{w_0} : N\to N$ is a surjective complex-linear isometry and a triple isomorphism mapping $u$ to $\11b{N}$ (cf. \cite[Theorem 4.2.28]{Cabrera-Rodriguez-vol1}), therefore $U_{w_0} : N_{(u^*)}\to N$ is a complex-linear Jordan $^*$-isomorphism. Let $p\in M$ and $J: M\to N_{(u^*)}$ be the central projection and the complex-linear Jordan $^*$-isomorphism given by our first conclusions. Clearly, $\Phi= U_{w_0} \circ J: M\to N$ is a Jordan $^*$-isomorphism and the equality $$\Delta (a) = U_{w_0^{-1}} \left(\Phi (p\circ a) + \Phi ((\11b{M}-p) \circ a^*)\right)= U_{w_0^{*}} \left(\Phi (p\circ a) + \Phi ((\11b{M}-p) \circ a^*)\right),$$ holds for all $a\in \mathfrak{M}$.
\end{proof}

We shall next obtain Theorem~\ref{t HatWat} as a consequence of the just proved theorem.

\begin{proof}[Proof of Theorem~\ref{t HatWat}] Let $\Delta : \mathfrak{A}\to \mathfrak{B}$ be a surjective isometry satisfying the hypothesis of the theorem. We observe that $\mathfrak{A}$ and  $\mathfrak{B}$ clearly satisfy the hypothesis of Theorem~\ref{t surjective isometries between invertible clopen JBstar} (cf. page~\pageref{label page mathfrakA satisfies mathfrakM}). By Theorem~\ref{t surjective isometries between invertible clopen JBstar},  $\Delta(\11b{A}) =u$ is a unitary element in $N$ and there exists a projection $p$ in the centre of $A$ and a complex-linear Jordan $^*$-isomorphism $J$ from $A$ onto the $u^*$-homotope $B_{(u^*)} = (B,\circ_{u},*_u)$ such that $$\Delta (a) = J(p\circ a) + J ((\11b{A}-p) \circ a^*), \hbox{ for all } a\in \mathfrak{A}.$$ The mapping $L_{u^*}: B\to B,$ $L_{u^*} (x) = u^* x$ is a surjective isometry mapping $u$ to $\11b{A}$, and hence it is a Jordan $^*$-isomorphism from $B_{(u^*)} = (B,\circ_{u},*_u)$ onto $B$. Therefore the mapping $\tilde{J} =   L_{u^*} \circ J$ is a Jordan $^*$-isomorphism from $A$ onto $B$ with $J = L_{u} \circ \tilde{J} = L_{\Delta(\11b{A})} \circ \tilde{J}$, and since $p$ is central,  $$\Delta (a) = \Delta(\11b{A}) \left( \tilde{J} (p  a) + \tilde{J} ((\11b{A}-p)  a)^* \right), \hbox{ for all } a\in \mathfrak{A}.$$ 
\end{proof}

\subsection{The cone of positive elements as a metric invariant}\label{subsection positive invertible with the usual metric} \ \smallskip

An element $a$ in a JB$^*$-algebra $M$ is called positive if $a^*= a$ and the spectrum of $a$ is contained in $\mathbb{R}_0^{+}$. The set $M^{+}$ of all positive elements in $M$ is a closed convex cone \cite[Lemma 3.3.7]{HOS}. If $M$ is a C$^*$-algebra, $M^{+}$ is the usual cone of positive elements. If $M^{+}$ is regarded inside $M_{sa}$ and $M$ is unital, we can also conclude that $M^{+}$ has non-empty interior, and hence it is a convex body. In this case, we shall write $M^{++}$ for the set of all positive and invertible elements in $M$. \smallskip

The idea of considering the cone of positive elements or the set of positive invertible elements in a unital C$^*$-algebra as metric invariants with respect to the original norm was considered by Hatori and Moln{\'a}r in \cite[page 167]{HatMol2014}, the same ideas work for unital JB$^*$-algebras.

\begin{proposition}\label{p cone of positive elements as metric invariant usual norm} Let $M$ and $N$ be two unital JB$^*$-algebras, and let $\Delta: M^{++}\to N^{++}$ be a surjective isometry with respect to the distances induced by the JB$^*$-norms. Then there exits a Jordan $^*$-isomorphism $J: M\to N$ whose restriction to $M^{++}$ is $\Delta$.  
\end{proposition}

\begin{proof} We can extend $\Delta$ to a surjective isometry $\tilde\Delta$ from $M^{+}$ onto $N^{+}$. Since $M^{+}$ and $N^+$ are closed convex bodies of $M_{sa}$ and $N_{sa}$ respectively, we can apply Mankiewicz theorem (see Theorem~\ref{t Mankiewicz}) to deduce that $\tilde\Delta$ extends to a surjective affine isometry $F: M_{sa}\to N_{sa}$. Since $0$ is the unique extreme point of $M^+$ and $N^+$, we have $\tilde{\Delta} (0) = F(0)=0$, and thus $F$ is a linear isometry. Now applying \cite[Theorem 1.9]{IsRo}, and having in mind that $N_{sa}$ is unital, it follows that there exists a central symmetry $s\in N$ and a surjective Jordan isomorphism $J: M_{sa}\to N_{sa}$ such that $T (x) = s \circ J(x),$ for all $x\in M_{sa}$. Since $s= F(\11b{M})=\Delta (\11b{M})  \in N^{+}$, we obtain $s=\11b{M}$ and $F (x) = J(x)$.
\end{proof}	

As we already commented, Mankiewicz theorem (Theorem~\ref{t Mankiewicz}) motivated the introduction of the so-called Mankiewicz property by Mori and Ozawa, who established the following result.

\begin{theorem}\label{t Mori-Ozawa strong Mankiewicz property}{\rm\cite[Theorem 2]{MoriOza2018}}. Let $X$ be a Banach space such that the closed convex hull of the set $\partial_e \left({B}_X\right)$, of all extreme points of the closed unit ball, ${B}_X$, of $X$ has non-empty interior in $X$. Then every convex body $C \subseteq  X$ has the strong Mankiewicz property, that is, every surjective isometry from $C$ onto an arbitrary convex subset $L$ in a normed space $Y$ is affine.
\end{theorem}

Mori and Ozawa applied this result to show that every convex body in a unital C$^*$-algebra, or in a (unital) real von Neumann algebra satisfies the strong Mankiewicz property \cite[Corollary 2]{MoriOza2018}. The same occurs for convex bodies in JBW$^*$-triples \cite[Corollary 2.2]{BeCuFerPe2021}, convex bodies in $C(\Omega,\mathcal{H})$, where $\Omega$ is a compact Hausdorff space, and $\mathcal{H}$ is a real Hilbert space with dim$(\mathcal{H}) \geq 2$ \cite{CuePer19}, and convex bodies in unital JB$^*$-algebras and real JBW$^*$-algebras \cite[pages 15 and 20]{PeSvard2022}. \smallskip

It is worth recalling that a \emph{real JB$^*$-algebra} is a closed $^*$-invariant real subalgebra of a (complex) JB$^*$-algebra. Unital real JB$^*$-algebras were introduced by K. Alvermann in \cite{Alver86}. Besides the just quoted reference an axiomatic intrinsic definition of unital real JB$^*$-algebras can be found in \cite{PeAxiom2003}. Additional examples of real JB$^*$-algebras include all JB-algebras, which are precisely those real JB$^*$-algebras whose involution is the identity mapping, and all real C$^*$-algebras (i.e., closed $^*$-invariant real subalgebras of C$^*$-algebras) equipped with their natural Jordan product \cite{Li2003,GoodearlBook82}. A real JBW$^*$-algebra is a real JB$^*$-algebra which is also a dual Banach space. It follows from the results commented above that every JBW-algebra $\mathfrak{J}$ is a real JBW$^*$-algebra, every convex body in $\mathfrak{J}$ has the strong Mankiewicz property, and in particular $\mathfrak{J}^+$ enjoys this property. We have just discovered the following corollary.

\begin{corollary}\label{c mankiewicz from J+ to a convex subset} Let $M$ be a JBW$^*$-algebra and let $L$ be a convex subset in a normed space $Y$. Suppose $\Delta: M^{+}\to L$ is a surjective isometry.  Then $\Delta$ extends to a surjective isometry from $M_{sa}$ onto a complete subspace of $Y$. 
\end{corollary}

\begin{proof} 
Since $M^+$ satisfies the strong Mankiewicz property as a convex body inside $M_{sa}$, the mapping $\Delta$ is affine. It is part of the folklore of affine maps that in such a case $\Delta$ can be uniquely extended to an affine isometry from $M_{sa}$ onto a complete subspace of $Y$.
\end{proof}

\section{A geometric tool: preservers of inverted Jordan triple products}\label{sec: geometric tool of Hatori Hirasawa Miura Molnar}

In \cite{Vaisala2003} V\"{a}is{\"a}l{\"a} gave an astonishing-simplified proof of the celebrated Mazur--Ulam theorem based on the idea that a surjective isometry between normed real-linear spaces preserve algebraic midpoints, i.e. $T\left(\frac{x+y}{2}\right) = \frac{T(x)+T(y)}{2}$ (cf. the forerunner due to Vogt in \cite{Vogt1973}). By replacing preservation of midpoints with preservation of inverted Jordan triple products of the form $y x^{-1} y$, Hatori, Hirasawa, Miura and Moln{\'a}r established in \cite{HatHirMiuMol2012} a very abstract and general result concerning a local property for maps between abstract spaces preserving a certain mild relation $d$ between elements, which is a far-reaching generalization of the usual notion of distance.\smallskip

The simplicity of the statement and the hypotheses are striking. Let us begin by recalling some definitions. A  \emph{metric} on a set $X$ is a function $d\,\colon X\times X\to \mathbb{R}$ satisfying the following axioms: 
\begin{enumerate}[$(1)$]
\item (Semi-positivity) $d(x,y)\geq 0$;
\item (Identity of indiscernibles) $x= y \Leftrightarrow d(x,y)=0;$
\item (Symmetry)  $d(x,y)=d(y,x);$
\item (Triangle inequality) $d(x,z)\leq d(x,y)+d(y,z)$,
\end{enumerate} for all $x,y,z\in X$. There is a wide range of generalizations of the notion of metric (see, for instance, \cite[\S 1.1]{DezaDeza2016}). For example, a \emph{distance} or a \emph{dissimilarity} on $X$ is a function $d : X \times X\to \mathbb{R}$ which is semi-positive, symmetric and that instead of the identity of indiscernibles, it satisfies the following weaker version 
\begin{enumerate}[$(2')$]
\item (Reflexivity)  $d(x,x)=0$ for all $x\in X$.
\end{enumerate} A \emph{semi-metric} is a mapping $d$ which is semi-positive, symmetric, reflexive and satisfies the triangular inequality. The mildest set of axioms leads us to the notion of \emph{quasi-distance} or \emph{parametric} which is a mapping $d : X \times X \to \mathbb{R}$ which is merely semi-positive and reflexive. Here, we shall work with an even weaker notion, a \emph{non-reflexive quasi-distance} or \emph{non-reflexive parametric} (\emph{nr-quasi-distance} and \emph{nr-parametric} in short) on $X$ is a semi-positive mapping $d : X \times X \to \mathbb{R}$. Here, the term non-reflexive actually means non-necessarily reflexive.\smallskip  

Along this section, $X$ will be a non-empty set equipped with a nr-quasi-distance $d_{X}$. The pair $(X,d_X)$ will be called a nr-quasi-distance space.  If $(Y,d_Y)$ is another nr-quasi-distance space, a mapping $f : X \to Y$ is said to be \emph{$d$-preserving} if $$d_X(x, y) = d_Y (f (x), f (y)),
\hbox{ for every } (x, y) \in  X \times X.$$ Clearly every metric space is a nr-quasi-distance space. Furthermore, a mapping between metric spaces with respect to the $d$ maps given by the distances is $d$-preserving precisely when it is an isometry. If no confusion arises, we shall omit the subindex in $d_{X}$. \smallskip

Our first theorem goes back to a result by Vogt \cite[Theorem 1.2]{Vogt1973} under stronger hypotheses which were later simplified by Hatori, Hirasawa, Miura and Moln{\'a}r in \cite[Lemma 2.3]{HatHirMiuMol2012}, from where we have borrowed the proof.

\begin{theorem}\label{t Vogt}{\rm\cite{Vogt1973,HatHirMiuMol2012}} Let $(X, d)$ be a nr-quasi-distance space. Suppose that there
exists a point $c$ in $X$, a bijective $d$-preserving mapping $\varphi : X\to X$, and a constant $K >1$ satisfying: 
\begin{enumerate}[$(a)$]
\item $d(\varphi(x), x) \geq  K d(x, c),$ for every $x$ in $X$;
\item the function $d(\cdot,c): X\to \mathbb{R}_0^+,$ $x\mapsto d(x,c)$ is bounded. 
\end{enumerate}
Then we have $d(f (c), c) = 0$ for every bijective $d$-preserving map $f: X\to X$.
\end{theorem}

\begin{proof} It follows from hypothesis $(b)$ that the set $$\left\{d(f(c),c) : f:X\to X \hbox{ bijective and $d$-preserving}\right\}$$ is bounded. Let $\lambda$ denote the supreme of this set. For each bijective $d$-preserving mapping $f:X\to X$, the map $\tilde{f} = f^{-1}\circ  \varphi \circ f$ is a $d$-preserving bijection. Therefore, by $(a),$ we have $$\lambda \geq d(\tilde{f} (x), c ) = d(\varphi f(x), f(c) ) \geq  K d(f(x), c).$$   The arbitrariness of $f$ shows that  $\lambda \geq K \lambda,$ which leads to $\lambda =0$. 
\end{proof}

If we assume that in Theorem~\ref{t Vogt} above the mapping $d$ satisfies the additional hypothesis $d(x,y) =0 \Rightarrow x=y$, like, for example, when $d$ is a metric, the previous theorem proves that $c$ is a fixed point for every $d$-preserving bijection $f:X\to X$. In \cite[Theorem 1.2]{Vogt1973}, $X$ is a bounded metric space for which there exists a surjective isometry $\varphi : X\to X$ and $K>1$ satisfying hypothesis $(a)$ in Theorem~\ref{t Vogt} for certain point $c\in X$, then the conclusion is that every surjective isometry $f$ on $M$ fixes $c$ (note that hypothesis $(b)$ clearly follows from the boundedness of $X$).\smallskip

It should be also added that a weaker form of Vogt's result \cite[Theorem 1.2]{Vogt1973}, asserting that if a bounded metric space $(X,d)$ is symmetric with respect to some $c\in X$-- i.e., there exists a surjective isometry $g$ on $X$ such that $d(g(x),x) = 2 d(x,c)$ for all $x\in X$-- then $c$ must be a fixed point of every surjective isometry on $X$, was employed by Mankiewicz in his celebrated extension of the Mazur--Ulam theorem on the extension of surjective isometries between arbitrary convex subsets with non-empty interior in real normed spaces to affine surjective isometries between the spaces (see \cite{Mank1972}).\smallskip 

We are now in a position to revisit the main result of this section which is due to Hatori, Hirasawa, Miura and Moln{\'a}r \cite{HatHirMiuMol2012}.

\begin{theorem}\label{t main geometric tool}{\rm\cite[Theorem 2.4]{HatHirMiuMol2012}} Let $(X, d_X)$ and $(Y, d_Y )$ be two nr-quasi-distance spaces. Let $a$ and $c$ be two arbitrary points in $X$. Suppose we can find a period-2 {\rm(}i.e. $\varphi^2 = \varphi \circ \varphi = Id_{X}${\rm)}, $d$-preserving mapping $\varphi : X \to X$ admitting $c$ as a fixed point. Let $L$ denote the set of all points in $X$ such that the nr-quasi-distance from $a$ and from $\varphi(a)$ to them is equal to $d_X(a, c)$, that is, $$L = \{ x \in X : d_X (a, x) = d_X (\varphi(a), x) = d_X(a, c)\}.$$ We shall additionally assume that the set $\{d_{X} (x, c) : x \in L\}$ is bounded, and there exists a constant $K > 1$ such that $d_X (\varphi(x), x) \geq  K d_X (x, c)$ for every $x \in L$. Then, if $T$ is a bijective $d$-preserving map from $X$ onto $Y$, and $\psi$ is a bijective $d$-preserving map from $Y$ onto itself such that one of the next statements holds \begin{equation}\label{eq hypotheses for T varphi and psi} \psi (T (a)) = T (\varphi(a)), \hbox{ and } \psi(T (\varphi(a))) = T (a),
\end{equation} 
\begin{equation}\label{eq hypotheses for T varphi and psi 2} \psi (T (a)) = T (a), \hbox{ and } \psi(T (\varphi(a))) = T (\varphi(a)),
\end{equation}
we have $d_Y (\psi(T (c)), T (c)) = 0$.
\end{theorem} 

\begin{proof} The hypotheses $\varphi (c) =c$ and $\varphi$ being $d$-preserving imply that $c\in L\neq \emptyset$. By applying that $\varphi$ is a period-2, $d$-preserving mapping, it is easy to check that $\varphi (L)\subseteq L$, and thus $\varphi (L) =L$ because $\varphi^2 = Id_{X}$.\smallskip 	
	
The alter ego of $L$ in $Y$ is the set 
$$L_Y = \{ y \in Y : d_Y (T(a), y) = d_Y (T(\varphi(a)), y) = d_X(a, c)\}.$$ 

We claim that $T(L) = L_Y$. Indeed, for each $x\in L$ the assumptions on $T$ give 
$$d_Y (T(a),T(x)) = d_X (a,x) = d_X (a,c), \hbox{ and }$$ 
$$d_Y (T(\varphi(a)),T(x)) = d_X (\varphi(a),x) = d_X (a,c),$$ proving that $T(x)\in L_Y$ and hence $T(L)\subseteq L_Y$. For the reciprocal, we observe that $\tilde{\varphi} := T \circ \varphi \circ T^{-1} : Y\to Y$ is a period-2, $d$-preserving mapping admitting $T(c)$ as a fixed point, $T^{-1} : Y\to X$ is a $d$-preserving bijection for which $L_Y$ plays the same role of $L$ for $T$. So, it follows from the previous arguments that $T^{-1} (L_Y) \subseteq L$, and hence $T(L)= L_Y$.\smallskip  

The stability of $L_{Y}$ under the mapping $\psi$ follows from the hypotheses in \eqref{eq hypotheses for T varphi and psi} or \eqref{eq hypotheses for T varphi and psi 2}. We shall only prove the case in which \eqref{eq hypotheses for T varphi and psi} holds, the other case is similar. In such a case, for each $y\in L_Y$ we have 
$$d_Y (T(a),\psi(y)) = d_Y (\psi (T(\varphi(a))), \psi(y) ) = d_Y ( T(\varphi(a)), y ) =d_X (a, c ), \hbox{ and }$$
$$d_Y (T(\varphi(a)), \psi(y) ) = d_Y (\psi (T(a)), \psi(y) ) = d_Y (T(a), y ) =d_X (a,c).$$ The last two identities show that $\psi (L_{Y}) \subseteq L_Y$. On the other hand, given $y\in L_Y$, by applying that $\psi^{-1}$ also is $d$-preserving we get 
$$d_{X} (a,c)= d_{Y} (T(a), y) = d_{Y} (\psi^{-1}T(a), \psi^{-1} (y)) = d_{Y} (T(\varphi(a)), \psi^{-1} (y)), \hbox{ and }$$ 
$$d_{X} (a,c) = d_{Y} (T(\varphi(a)), y) = d_{Y} (\psi^{-1} T(\varphi(a)), \psi^{-1} (y)) =d_{Y} (T(a), \psi^{-1} (y)),$$ proving that $\psi^{-1} (y)\in L_Y,$ and thus $\psi (L_Y) = L_Y$.\smallskip

Finally, the mapping $\tilde{\psi} = T^{-1} \circ \psi\circ T$ is a $d$-preserving bijection on $X$ with $\tilde{\psi} (L) =L,$ and $\varphi|_{L}: L\to L$ is a $d$-preserving bijection. Since $d_X (\varphi(x), x) \geq  K d_X (x, c)$ for every $x \in L$ (with $K > 1$) and the set $\{d_{X} (x, c) : x \in L\}$ is bounded, Theorem~\ref{t Vogt}, applied to $\tilde{\psi}|_{L}$ and $\varphi|_{L},$ assures that $0=d_{X} (\tilde{\psi}(c),c) =d_{Y} (\psi(T(c)),T(c))$, as desired.  
\end{proof}

\begin{remark}\label{remark bounded is free for metric} Note that if $d_X$ is a metric, the triangle inequality guarantees that the set $\{d_{X} (x, c) : x \in L\}$ in Theorem~\ref{t main geometric tool} is bounded.
\end{remark}

The refinement in the previous theorem is to consider the set $L$ to guarantee the hypothesis concerning  boundedness in Theorem~\ref{t Vogt}. \smallskip

As the reader might guess, Vogh and Väisälä arguments in \cite{Vogt1973,Vaisala2003} leading to an optimal proof of the Mazur--Ulam theorem, can be now rediscovered from the previous results. Namely, suppose $F:X\to Y$ is a surjective isometry between real normed spaces. We consider the metrics given by the norms on $X$ and $Y$, respectively. For each $c\in X$, let $\varphi_c: X\to X$ denote the mapping given by the reflection of $X$ in $c$, that is, $\varphi_c (x) := 2 c -x$ ($x\in X$). Fix $a,b\in X$ and denote by $\varphi_c$ and $\tilde{\varphi}_{\frac{F(a)+F(b)}{2}}$ the reflections of $X$ and $Y$ in $c= \frac{a+b}{2}$ and $\frac{F(a)+F(b)}{2}$, respectively. We consider the mapping $f = \varphi_c \circ F^{-1} \circ \tilde{\varphi}_{\frac{F(a)+F(b)}{2}}\circ F: X\to X$ which is a $d$-preserving bijection. Clearly, $\varphi_c$ is a period-$2$, $d$-preserving bijection with $\varphi_c(c) =c$, $\varphi_c (a) =b,$ and $\varphi_c(b) =a$. Furthermore, $$d_{X} (\varphi_c(x),x) = \|2 c -x -x \| = 2 \|c-x\|, \hbox{ for all } x\in X.$$ In principle, the function $d_X (\cdot,c)$ is not bounded on $X$. To avoid this obstacle we restrict our maps to a suitable subset $L=\{x\in X : d_{X} (a,x) =d_X (\varphi_c(a),x) = d_{X} (a,c)\}.$ As we saw in the previous theorem, or with simpler arguments, $\varphi_c (L) = L$. By the triangular inequality $\sup\{d_X(x,c) : x\in L\}<\infty $.\smallskip

The mapping $\psi = F\circ f \circ F^{-1} = F\circ \varphi_c \circ F^{-1} \circ \tilde{\varphi}_{\frac{F(a)+F(b)}{2}}$ is a $d$-preserving bijection on $Y$ with $\psi(F(a)) = F\circ \varphi_c \circ F^{-1} \circ \tilde{\varphi}_{\frac{F(a)+F(b)}{2}} (F(a)) = F(a),$  and $\psi(F(\varphi_{c}(a))) = F\circ \varphi_c \circ F^{-1} \circ \tilde{\varphi}_{\frac{F(a)+F(b)}{2}} (F(\varphi_{c}(a))) =  F(b) = F(\varphi_c (a)).$  Theorem~\ref{t main geometric tool} implies that $ F(c) = \psi (F(c)) =F(f(c)) = F\circ \varphi_c \circ F^{-1} \circ \tilde{\varphi}_{\frac{F(a)+F(b)}{2}}\circ F (c),$ or equivalently, $\tilde{\varphi}_{\frac{F(a)+F(b)}{2}}\circ F (c) = F(c)$. Therefore $F(c)$ is the unique fixed point for $\tilde{\varphi}_{\frac{F(a)+F(b)}{2}}$, that is, $\frac{F(a)+F(b)}{2} = F(\frac{a+b}{2})$. A standard argument shows that $F$ is affine.\smallskip 
 
After a first application to rediscover the classical Mazur--Ulam theorem, we should not forget that Theorem~\ref{t main geometric tool} is a device thought to prove the preservation of inverted Jordan triple products of the form $y x^{-1} y$. We shall see different applications in various contexts. We shall adopt the notation from \cite{HatHirMiuMol2012}. \smallskip

Along this section, we shall assume that $\mathcal{G}$ is a group. Unless another assumptions are made, we shall work with a nr-quasi-distance space $(X,d)$ satisfying the mere assumptions that $X$ is a subset of $\mathcal{G}$ and \begin{equation}\label{eq condition p HatoriMiuraMolnar}\hbox{ $y x^{-1} y\in X$ for all $x,y\in X$.}
\end{equation} The hypotheses do not necessarily imply that $X$ is a subgroup of $\mathcal{G}$. We continue by recalling certain extra algebraic-geometric properties.

\begin{definition}\label{def condition B} Let us fix $a,b$ in $X$. We shall say that $(X,d)$ satisfies  property $B(a,b)$  if the following three statements hold:\begin{enumerate}[$(B.1)$]\item  $d(b x^{-1} b, b y^{-1} b) = d(x,y),$ for all $x,y\in X$.
\item $\sup\{d(x,b): x\in L_{a,b}\}<\infty,$ where we set $$L_{a,b} =\{x\in X : d(a,x) = d(b a^{-1} b, x) = d(a,b)\}.$$
\item There exists a constant $K>1$ satisfying $$d(b x^{-1} b,x)\geq K d(x,b),$$ for all $x\in L_{a,b}.$
	\end{enumerate}
\end{definition}

In case that $d$ satisfies the triangular inequality (for example, when $d$ is a semi-metric or a metric), condition $(B.2)$ above is trivially true by the triangular inequality. Following the notation of the pioneering paper \cite{HatHirMiuMol2012}, by an \emph{algebraic mid-point} of a couple of elements $x,y\in \mathcal{G}$ we mean any element $c\in \mathcal{G}$ satisfying $y = c x^{-1} c$. Let us note that $b$ is an algebraic point of $a$ and $b a^{-1} b$, and in case that $X$ satisfies property $B(a,b)$ the element $b$ lies in the set $L_{a,b}$ by definition.    

\begin{definition}\label{def condition C1} Let us fix $a,b\in X$. We shall say that $(X,d)$ satisfies property $C_1(a,b)$ if the following statements hold:\begin{enumerate}[$(C.1)$]\item $a x^{-1} b, b x^{-1} a\in X,$ for every $x\in X$.
\item $d(a x^{-1} b, a y^{-1} b) = d(x,y)$, for all $x,y\in X$.
	\end{enumerate}
\end{definition}

\begin{definition}\label{def condition C2} Let us fix $a,b\in X$. We shall say that $(X,d)$ has property $C_2(a,b)$ if there exists $c\in X$ such that $c a^{-1} c =b$ and $d(c x^{-1} c, c y^{-1} c) = d(x,y)$ for all $x,y\in X$.
\end{definition}

Let us remark that properties $B(a,b)$, $C_1(a,b)$ and $C_2(a,b)$ above depend of the points $a,b$.\smallskip

An element $x\in X$ is called \emph{{2-divisible}} if there exists $y\in X$ such that $y^2 =x$. $X$ is called {\emph{2-divisible}} if every element in $X$ is {2-divisible}. We observe that $X$ contains the unit of $\mathcal{G}$ whenever it is 2-divisible. Namely, if we take $y\in X,$ by assumptions there is $x \in X$ with $x^2 = y$ and hence $e = x y^{-1} x \in X$. Moreover, let us take $a, b \in X,$ by the property on $X$, there exist $d,f \in X$ with $d^2 = a$ and $f^2 = d b^{-1} d$, and hence setting $c = d f^{-1} d$ we have $c \in X$ and $c a^{-1} c = d f^{-1} d a^{-1} d f^{-1} d = d (f^{-1})^2 d = b$.\smallskip

Finally, $X$ is a called {\emph{2-torsion free}} if it contains the unit of $\mathcal{G}$ and the condition $x^2 =1$ with $x\in X$ implies $x =1$.

\begin{definition}\label{def inverse and translation invariant} Let us additionally assume that $X$ is a subgroup of $\mathcal{G}$. We say that $d$ is inverse and translation invariant if $d(x^{-1}, y^{-1}) =  d (a x, a y) = d(x, y)$	for every triplet $x, y, a$ in $X$.
\end{definition}

If $X$ is actually a subgroup of $\mathcal{G}$, then $d$ is inverse and translation invariant if
and only if $(X,d)$ satisfies $C_1(a, b)$ for every $a, b \in X$. If $X$ is a subgroup of $\mathcal{G}$, 2-divisible, and $d$ is inverse and translation invariant, then $(X,d)$ satisfies $C_2(a, b)$ for every $a, b \in X$.\smallskip

The assumptions until now are so flexible that may spaces can serve as examples. The first one is a basic example,  $X$ being an additive subgroup of a normed space $(Z,\|.\|)$ with the metric given by the norm, then $X$ satisfies properties $B(a,b)$ and $C_1(a, b)$ for every $a,b$ in $X$, moreover $X$ satisfies $C_2(a,b)$ in the case that $X$ is 2-divisible.\smallskip

The next three examples have been borrowed from \cite{HatHirMiuMol2012}.

\begin{example}\label{example Tn} Let $X$ denote the n-dimensional torus $\mathbb{T}^n$ with the usual pointwise multiplication and metric defined by $d(x, y) = \max\{ |x_j - y_j| : 1 \leq j \leq n\}$ for $x = (x_1, \ldots,x_n),$ $y =(y_1, \ldots,y_n)$ in $\mathbb{T}^n$. It is shown 
	in \cite[Example 3.7]{HatHirMiuMol2012} that for any $a, b \in \mathbb{T}^n$ with $d(a, b) < \sqrt{2}$, $(X,d)$ satisfies $B(a, b)$ with $K =
\sqrt{2}$. $(X,d)$ does not satisfy property $B(a, -a)$ (the third condition fails).  It is also shown that $X=\mathbb{T}^n$ is 2-divisible and $d$ is inverse and translation invariant, and hence $(X,d)$
satisfies $C_1(a, b)$ and $C_2(a, b)$ for all $a, b\in X$.
\end{example}

\begin{example}\label{example CK spaces with a new metric}\cite[Example 3.8]{HatHirMiuMol2012} Let $K$ be a compact Hausdorff space and $E$ any real-linear subspace of the space $C(K)$ of all complex-valued continuous functions on $K$. Define 
$$ d(f, g) = \max\left\{\left\|{f}/{g}-1\right\|_{\infty}, \left\|{g}/{f}- 1\right\|_{\infty} \right\}$$ for all $f,g \in \exp(E)$, where $\|\cdot\|_{\infty}$ stands for the supremum norm $C(K)$. Then $\exp(E)$ is a (multiplicative) group and $(\exp(E), d)$ is a nr-quasi-distance space. However $d$ is not, in general, a metric, it is semi-positive, symmetric and satisfies the identity of indiscernibles but it fails, in general, the triangle inequality. It is known that $\exp(E)$ is 2-divisible, and $d$ is inverse and translation invariant. The nr-quasi-distance space $(\exp(E), d)$ satisfies properties $C_1(a, b)$ and $C_2(a, b)$ for every $a, b$ in $\exp(E)$.  
\end{example}

The next examples are based on metrics and distances which are not directly related with the norm of a Banach space. In our approach to the Hilbert projective metric and to the Thompson's metric we shall consider them in the setting of {order unit spaces} as in \cite{LemmNussBook2012, HatHirMiuMol2012, RoeWort2019, LemmRoeWor2018, LemmRoeWor2019}. We refer to \cite{Al,AsEll,HOS,LemmRoeWor2018} and the survey \cite[\S 1]{BaRo} for the basic theory any undefined terms and details used on ordered Banach spaces. By a \emph{a partially ordered space} we mean a real vector space $A$ with a \emph{proper cone} $A^+$ (i.e., $A^+$ is convex, $\mathbb{R}^+ A^+ \subseteq A^+$, and  $A^+ \cap \left(- A^+\right) = \{0\}$). The \emph{partial ordering} ``$\leq$'' on $A$ is defined by $a \leq b$ if $b-a \in  A^+$. Given $a,b\in A$ the \emph{order interval} $[a,b]$ is the set of all $x\in A$ such that $a\leq x\leq b$. An element $u\in A^+$ is called an \emph{order unit}, if it satisfies $A = \cup \left\{ [- \lambda \ u, \lambda \ u] : \lambda \geq 0 \right\}.$ For each order unit $u$ in $A^{+}$, the expression $$\|x\|_{u} = \inf \{ \lambda >0 : - \lambda u \leq x \leq \lambda u \} \ \ \  (x\in A),$$ defines a seminorm on $A$. The order unit $u$ is called \emph{Archimedean} if $n a \leq  u$ for all $n \in \mathbb{N}$, then $-a \in A^+$. If $u$ is an Archimedean order unit, the seminorm $\|\cdot\|_u$ is actually a norm, $A^{+}$ is a closed convex cone, and the closed unit ball with respect to this norm is the order interval $[-u,u]$ (cf. \cite[Theorem 2.2.5]{AsEll}). It is further known that $u$ is a norm interior point of $A^{+}$ \cite[Proposition 2.2.1]{AsEll}. Actually, by the just quoted results, the set $\hbox{int} (A^{+})$ of norm interior points of $A^{+}$ coincides with the set of \emph{order units} in $A$ (this was also rediscovered in \cite[\S 1.4]{BaRo} and \cite[Lemma 2.1]{LemmRoeWor2018}).\smallskip 

In the case that we have an order space $A$ which is also a Banach space with respect to a given norm $\|\cdot\|$ (usually called an \emph{ordered Banach space}), for each Archimedean order unit $u$ in $A$ we find a new norm $\|\cdot\|_u$. It is natural to ask what is the relation of these two norms. \smallskip

For our goals in this note a typical example of an ordered Banach space can be given by the self-adjoint part, $A_{sa}$, of a C$^*$-algebra $A$ with the natural C$^*$-norm, or the self-adjoint part, $\mathfrak{A}_{sa},$ of a JB$^*$-algebra $\mathfrak{A}$ with the natural JB$^*$-norm, with the cone of positive elements. We recall that the cone of positive elements in a JB$^*$-algebra is a proper closed convex cone (see \cite[Theorem I.6.1]{Tak} and \cite[\S 3.3]{HOS}). \smallskip

It is shown in \cite[Example 2.15]{GluckWeb2020} that the cone of positive elements in a C$^*$-algebra $A$ admits an interior point with respect to the topology generated by the norm if and only if $A$ is unital, and in such a case, the interior of $A^+$ with respect to the norm topology is the set of all positive invertible elements in $A$.\smallskip
 
A first natural question now is to determine when the self-adjoint part of a C$^*$-algebra (respectively, a JB$^*$-algebra) admits an Archimedean order unit $u$. Let $u$ be an Archimedean order unit in a C$^*$-algebra $A$, and let $\|\cdot\|_u$ denote the corresponding order norm on $A_{sa}$ generated by $u$. Let us take $a\in A_{sa}$. It is well known that $a$ writes in the form $a= a^* - a^-$ with $a^+,a^-\geq 0$, $a^+ a^- =0$ and $\|a \|  =\|a^+\|\vee \|a^-\|.$ There is no loos of generality in assuming that $\|a\| = \|a^*\|.$ If $\delta>0$ satisfies $-\delta u  \leq a \leq \delta u$, by applying that the operator $x\mapsto a^+ x a^+$ is positive, we get $$0\leq \left(a^+\right)^3 = U_{a^+} (a^+) \leq \delta U_{a^+} (u),$$ and hence $$ \|a\|^3 = \|a^+\|^3 = \| \left(a^+\right)^3 \| \leq \delta \|a^+\|^2 \|u\|, \hbox{ and } \|a\| = \|a^+\| \leq \|u\| \delta. $$   Taking infimum in $\delta$ we arrive to \begin{equation}\label{eq C*-norm bounded by order norm}  \|a\| \leq \|u\| \|a\|_{u}, \hbox{ for all } a\in A_{sa}.
\end{equation} Similar arguments to those given above, but applying that for each hermitian element $h$ in a JB$^*$-algebra the operator $u_h$ is positive \cite[Proposition 3.3.6]{HOS}, prove that \eqref{eq C*-norm bounded by order norm} holds for every self-adjoint element in a JB$^*$-algebra $\mathfrak{A}$. However, the inequality in \eqref{eq C*-norm bounded by order norm} nor its version for JB$^*$-algebras allow us to apply the above commented \cite[Example 2.15]{GluckWeb2020} because it does not imply that an interior point of $A^*$ for the norm topology is an interior point for the norm $\|\cdot\|_u$. 

\begin{proposition}\label{p when a JB-algebra admits an Archimedean order unit} Let $M$ be a JB$^*$-algebra (in particular a C$^*$-algebra). Then $M_{sa}$ admits an Archimedean order unit $u$ if and only if it is unital. In such a case $u$ is a positive invertible element in $M$ and $\|\cdot\|_u$ is equivalent to the JB$^*$-norm on $M_{sa}$.  Furthermore, the set of all Archimedean order units in $M$ is precisely the set $M^{++}$ of all positive invertible elements in $M,$ which coincides with the interior of $M^+$ with respect to the JB$^*$-norm or with respect to $\|.\|_u$.  
\end{proposition}

\begin{proof} In the first statement the ``if'' implication is clear. Suppose $u$ is an Archimedean order unit in $M_{sa}$. Let $N$ denote the JB$^*$-subalgebra of $M$ generated by $u$. It is known that $M_{sa}$ is isometrically isomorphic as JB-algebras to $C_0(Sp(u)\backslash\{0\},\mathbb{R})$ and under this identification $u$ corresponds to the inclusion of $ Sp(u)$ into $\mathbb{C}$, where $ Sp(u)$ denotes the spectrum of $u$ (cf. \cite[Theorem 3.2.4]{HOS}). Since $u$ must be also an order unit in $N_{sa}$, it follows that there exits a positive $\delta$ such that $u^{\frac12} \leq \delta u$, which implies that $0$ is isolated in $Sp(u)$, and hence $N_{sa}\cong C(Sp(u),\mathbb{R})$ is an associative unital JB-algebra and $u$ is positive and invertible in $N$. Let $p$ denote the unit of $N$. Let us show that $p$ is the unit in $M$. Given $a\in M_{sa}$, the element $a- U_p(a)$ lies in $M_{sa},$ and by \cite[Lemma 2.4.21]{HOS} $$U_p (a- U_p (a)) = U_p (a) - U_p^2 (a)  =U_p (a) - U_{p^2} (a)=0.$$ Having in mind that $u$ is an order unit in $M_{sa}$ and $u\leq \|u\| p$, we deduce that $p$ also is an order unit in  $M_{sa}$, and so  $\beta p \leq  a- U_p(a) \leq \beta p$ for some $\beta>0$. An standard argument in spectral theory shows that $U_p (a- U_p(a)) = a- U_p(a),$ and thus $a= U_p(a)$ for all $a\in M_{sa}$. Therefore $p$ is the unit in $M.$ This concludes the proof of the first statement.\smallskip 
	
For the second statement, suppose that $u$ is an Archimedean order unit in $M$ and the latter is unital with unit $\11b{M}$. In such a case $\11b{M}\leq \delta u$ for some $\delta >0$, which clearly implies that $u$ is invertible, positive, and lies in the interior of $M^+$ with respect to the JB$^*$-norm.  Let us finally observe that for each $a\in M_{sa}$ we have  $$\frac{1}{\|u^{-1}\|} \ \{\beta > 0 : -\beta \11b{M} \leq a\leq \beta \11b{M}\} \subseteq \{\delta > 0 : -\delta u \leq a\leq \delta u \},$$ and consequently  $\|a\|_{u} \leq \frac{1}{\|u^{-1}\|} \| a\|,$ which combined with \eqref{eq C*-norm bounded by order norm} concludes the proof. 
\end{proof}

The set of invertible elements in a unital C$^*$-algebra $A$ is a multiplicative subgroup. However, the cone of positive elements in $A$ is not, in general, stable under products. The same difficulty occurs to the set of invertible elements and the set of positive elements in a unital JB$^*$-algebra $M$. However, as we have commented before, the sets $M^{-1}$ of all invertible elements, $M^{+}$ of positive elements, and the set $M^{++}$ of all positive and invertible elements in $M$ enjoy another  stability property \begin{equation}\label{eq stability of Uab for positive and invertible} U_a (b)\in M^{\sigma}, \hbox{ for all } a,b\in M^{\sigma} \hbox{ and } \sigma \in \{+,-1, ++\}. 
\end{equation} If we only consider invertible elements, the JB$^*$-algebra $M$ can be replaced with any unital complex Jordan--Banach algebra.  Furthermore, the conclusion is also true when $M^{\sigma}$ is replaced by the set $\mathcal{U}(M)$ of all unitary elements in $M$. The stability exhibited in \eqref{eq stability of Uab for positive and invertible} reaffirms the usefulness of the inverted Jordan triple product. \smallskip

Let us fix a C$^*$-algebra $A$. We equip the sets $A^{-1}$, $A^{++}$, and $\mathcal{U} (A),$ of invertible, positive and invertible and unitary elements in $A$, respectively, with the metric given by the norm, and we consider the multiplicative structure of $A$. Clearly, by \eqref{eq stability of Uab for positive and invertible}, these three sets satisfy the condition in \eqref{eq condition p HatoriMiuraMolnar}. It is easy to check that $A^{-1}$ and $A^{++}$ fail the properties $B(a,b)$, $C_1(a,b)$ and $C_2(a,b)$ in definitions \ref{def condition B}, \ref{def condition C1} and \ref{def condition C2} for most pairs $a,b$ in these sets. As we shall see later, unitaries exhibit a better behaviour. To avoid the difficulties mathematicians have considered alternative metrics on $A^{++}$.\smallskip 
 
It seems an appropriate moment to introduce the Hilbert projective metric and the Thompson's metric on the interior of the positive cone $A^{+}$ (i.e. the set of all Archimedean order units) in an order unit space equipped with an Archimedean order unit. For $a, b \in \hbox{int} (A^{+}),$ by applying the properties of order units, we set $$ M(a/b) := \inf\{ \beta > 0 \colon a \leq  \beta b\} (<\infty).$$ On $\hbox{int} (A^{+})$ Hilbert's projective metric and Thompson's metric are given by \begin{equation}\label{def Tompson and Hilbert distances}
 d_H (a, b) = \log M(a/b) M(b/a), \hbox{ and }  d_T (a, b) = \log\left(\max\{M(a/b), M(b/a)\}\right),	
\end{equation} respectively. The mapping $d_T$ satisfies all properties required for a metric on $\hbox{int} (A^{+})$, while Hilbert's projective metric is a semi-metric (i.e., it is semi-positive, reflexive, symmetric, and satisfies the triangular inequality), however it has a certain handicap in what respect the identity of indiscernibles, since $d_H (\lambda a, \mu b)= d_H (a, b)$ for all $\lambda,\mu > 0$ and $a, b \in\hbox{int} (A^+)$ and $d_H (a, b) = 0$ with $a, b \in \hbox{int} (A^+)$ if and only if $a = \lambda b$ for some $\lambda > 0$ (see \cite[Proposition 2.1.1 and comments in page 30]{LemmNussBook2012}).\smallskip

One of the most fruitful source of examples in this paper is coming from the cones of positive elements in a C$^*$-algebra $A$, or more generally, in a JB$^*$-algebra $M$. The interior of $M^+$  is non-empty if and only if $M$ is unital and in such a case $\hbox{int}\left(M^+\right)= M^{++}$ (see Proposition~\ref{p when a JB-algebra admits an Archimedean order unit}). Do we know a spectral formula for the Thompson's metric?\smallskip

We borrowed some paragraphs from a paper by Andruchow, Corach, and Stojanoff \cite{Andru2000} on the differential geometry of the cone $A^{++}$ of positive invertible elements in a unital C$^*$-algebra $A$ (see \cite{Andru2000, Cor1994, CorachPortaRecht1992,CorachPorRech1993, CorachPorRech1993b,CorachPorRech1994}). $A^{++}$ is an open convex subset of $A_{sa},$ and hence it can be considered as an open submanifold of $A_{sa},$ and the tangent space $(TA^+)_a$ at a point $a\in A^{++}$ will be identified with $A_{sa}$. There is a natural action of $A^{-1}$ on $A^{++}$ given by $(x, a) \mapsto  x a x^*$ ($x \in A^{-1}$, $a \in A^{++}$). This is a transitive action, that is, given $a, c \in  A^{++}$, $x = c^{1/2} a^{-1/2}$ verifies $x a x^* = c$. For every $a \in A^{++}$, the map $\tau_a : A^{-1} \to A^{++}$, $\tau_a (x) = x a x^*$ is a principal fibre bundle with a natural connection. The unique geodesic $\gamma$ such that $\gamma (0) = a$ and $\dot{\gamma} (0) = x \in A_{sa}$ is $$ \gamma (t) = e^{(t/2) x a^{-1} } = a^{1/2} e^{t a^{-1/2} x a^{-1/2}} a^{1/2},$$  and for every $a, b \in A^{++}$ there is a unique geodesic $\gamma_{a,b}$ such that $\gamma_{a,b} (0) = a$ and $\gamma_{a,b}(1) = b$, which is explicitly given by the expression $$\gamma_{a,b} (t) = a^{1/2} (a^{-1/2} b a^{-1/2})^t a^{1/2},$$ corresponding to the vector $x = a^{1/2}  \log\left( a^{-1/2} b a^{-1/2}\right) a^{1/2}.$ 
	
Although $A^{++}$ is not a Riemannian manifold, there is a natural Finsler structure on it. The \emph{geodesic or rectifiable distance} on $A^{++}$ defined by $$d_{g} (a, b) = \inf\{ \hbox{ length } \gamma  : \gamma \hbox{ is a smooth curve in } A^{++} \hbox{ joining $a$ and $b$} \},$$ satisfies $d_g(a, b) = \hbox{length} \gamma_{a,b} = \left\| \log\left( a^{-1/2} b a^{-1/2} \right)\right\|.$ If $a$ and $b$ commute we have $d_g(a, b) = \left\| \log\left( b a^{-1} \right)\right\|.$\smallskip

Andruchow, Corach, and Stojanoff proved in \cite{Andru2000} that the geodesic distance obtained by differential geometric methods coincides with the Thompson's metric. An alternative explanation for $B(H)^{++}$ was given by Moln{\'a}r in \cite{Mol2009}.

\begin{proposition}\label{p Andrichoff Corach Stojanoff geomesic and Thompson's metric}\cite[Proposition in page 1035]{Andru2000} Let $A$ be a unital C$^*$-algebra. For every $a, b \in A^{++}$ the Thompson's metric satisfies \begin{equation}\label{eq fla Thompsosn distance Andrichoff} d_{T} (a,b) = \left\| \log\left( a^{-1/2} b a^{-1/2} \right)\right\|.
	\end{equation}
\end{proposition}	\vspace*{0.25mm}

\begin{example}\label{example positive elements in a C*-algebra and thomson metric}\cite[Example 3.6]{HatHirMiuMol2012} Let $d_T$ denote the Thompson's metric on the set $A^{++}$ of all positive invertible elements in a unital C$^*$-algebra $A$. The condition $B(a, b)$ holds for every pair $a, b \in A^{++}$. Given $a, b, c \in A^{++}$, since the mapping $z\mapsto b z b$ ($z\in A$) is positive, we can easily see that $a\leq t c$ for some $t>0$ if and only if $b a b \leq t b c b,$ if and only if, $t^{-1} c^{-1}\leq a^{-1} $ \cite[Proposition 1.3.6]{Ped}. Thus the equality $d_T(b a^{-1} b, b c^{-1} b) = d_T( a^{-1}, c^{-1}) = d_T (a,c)$ follows directly from the definition of the Thompson's metric in \eqref{def Tompson and Hilbert distances}.  For the next property we take an argument from \cite[Proof of Theorem 1]{Mol2009}. Since $(b a^{-1} b)^{-1} a = (b^{-1} a)^2$, we are dealing with positive elements, by the basic properties of the spectrum we have  $$\begin{aligned} \sigma\left( (b a^{-1} b)^{-1/2} a (b a^{-1} b)^{-1/2} \right) &= \sigma\left( (b a^{-1} b)^{-1} a \right) = \sigma \left( (b^{-1} a)^2 \right) \\
&= \sigma \left( b^{-1} a \right)^2 = \sigma \left(b^{-1/2} a b^{-1/2} \right)^2,	
\end{aligned}$$	and hence $$ \begin{aligned}
&\sigma\left(\log\left( (b a^{-1} b)^{-1/2} a (b a^{-1} b)^{-1/2} \right)\right) =\log \sigma\left( (b a^{-1} b)^{-1/2} a (b a^{-1} b)^{-1/2} \right)\\
&= 2 \log \sigma \left(b^{-1/2} a b^{-1/2} \right) = 2 \sigma \left( \log \left(b^{-1/2} a b^{-1/2} \right)\right).
\end{aligned} $$ 
Taking suprema in the above sides of the above identity, by \eqref{eq fla Thompsosn distance Andrichoff}, we get 
\begin{equation}\label{eq B3 Molnar C*}\begin{aligned}
		d_T \left(b a^{-1} b , a \right) &= \left\|\log\left( (b a^{-1} b)^{-1/2} a (b a^{-1} b)^{-1/2}\right) \right\| \\
		&= 2 \left\| \log\left( b^{-1/2} a b^{-1/2} \right)\right\| = 2 d_T \left(b, a \right).
	\end{aligned}
\end{equation} Since $d_T$ is a metric, property $B(a, c)$ holds for every $a, c \in A^{++}$. Moreover, since $A^{++}$ is $2$-divisible, the conclusions in the first part of this example assure that property $C_2(a, c)$ also holds for every $a, c \in A^{++}$.		
\end{example}

\begin{remark}\label{r Tomhpson metric does not change in subalgebras} The Thompson's metric does not change when computed among positive invertible elements in a unital JB$^*$-subalgebra. That is, if $M$ is a unital JB$^*$-algebra and $N$ is a JB$^*$-subalgebra containing the unit of $M$, then for all $a,b\in N^{++}$ we have $d_{T}^N (a,b) = d_{T}^M (a,b)$.	This is clear from the definition in \eqref{def Tompson and Hilbert distances}. 
\end{remark}

We deal next with the case of JB$^*$-algebras. There is a central result in the theory of JB$^*$-algebras which has not been considered in this note yet. The so-called Shirshov--Cohn theorem gives a detailed information about the structure of the JB$^*$-subalgebra generated by two self-adjoint elements.

\begin{theorem}[Shirshov--Cohn theorem]\label{t Shirshov-Cohn}{\rm(\cite[Theorem 7.2.5]{HOS}, \cite[Proposition 2.1]{Wri77})} Let $M$ be a unital JB$^*$-algebra. Then the unital JB$^*$-subalgebra $M_{a,b,\mathbf{1}}$ generated by two self-adjoint elements $a,b$ and the unit in $M$ is a JC$^*$-algebra, that is, there exists a unital C$^*$-algebra $A$ such that $M_{a,b,\mathbf{1}}$ is a JB$^*$-subalgebra of $A$.
\end{theorem}

The first consequence of this result is a Jordan version of the formula proved by Andruchow, Corach, and Stojanoff in Proposition~\ref{p Andrichoff Corach Stojanoff geomesic and Thompson's metric}. The result was established by Lemmens, Roelands, and Wortel in  \cite{LemmRoeWor2018}, here we present a slightly different proof.

\begin{proposition}\label{p Lemmens on the Thompson distance for positive invertible in JB*-algebras}\cite[Proposition 2.4]{LemmRoeWor2018} Let $M$ be a unital JB$^*$-algebra. Then the Thompson's metric on $M^{++}$ satisfies \begin{equation}\label{eq fla Thompsosn distance JB*algebras} d_{T} (a,b) = \left\| \log\left( U_{a^{-1/2}} (b) \right)\right\|, \hbox{ for every } a, b \in M^{++}.  
\end{equation} Furthermore, \begin{equation}\label{eq conditions Jordan for B and C2} \begin{aligned}
d_T \left(U_b (a^{-1}) , a \right) &= 2 d_T \left(b, a \right), \hbox{ and }\\
 d_T(U_c (a^{-1}), U_c (b^{-1})) &= d_T( a^{-1}, b^{-1}) = d_T (a,b),
\end{aligned}
\end{equation} for every $a, b, c \in M^{++}$. Therefore $(M^{**},d_T)$ satisfies the Jordan analogues of properties $B(a,b)$ and $C_2(a,b)$ for all $a,b\in M^{++}$.  
\end{proposition}

Let us observe, before dealing with the proof, that a direct proof of the fact that the formula on the right-hand-side of \eqref{eq fla Thompsosn distance JB*algebras} is a distance on $M^{++}$ does not seem an easy task. 

\begin{proof} Let $N_{a,b,\mathbf{1}}$ denote the JB$^*$-subalgebra of $M$ generated by $a,b$ and $\mathbf{1}$. Clearly, $a^{-1},b^{-1}\in N_{a,b,\mathbf{1}}$. By the Shirshov-Cohn theorem (Theorem~\ref{t Shirshov-Cohn}) there exists a unital C$^*$-algebra $A$ such that $M_{a,b,\mathbf{1}}$ is a JB$^*$-subalgebra of $A$. By Proposition~\ref{p Andrichoff Corach Stojanoff geomesic and Thompson's metric} we have 
$$d^{M}_{T} (a,b) = d^{N_{a,b,\mathbf{1}}}_{T} (a,b) = d^{A}_{T} (a,b) = \left\| \log\left( U_{a^{-1/2}} (b) \right)\right\|_{A} = \left\| \log\left( U_{a^{-1/2}} (b) \right)\right\|_{M}.$$	Finally, by \eqref{eq B3 Molnar C*}, $$d_T^{M} \left(U_b (a^{-1}) , a \right) = d_T^{N_{a,b,\mathbf{1}}} \left(U_b (a^{-1}) , a \right)= d_T^{A} \left(b a^{-1} b , a \right)= 2 d_T^{A} \left(b, a \right) = 2 d_T^{M} \left(b, a \right).$$ The rest is clear by the properties of the cone of positive elements in $M$ via the same arguments given in Example~\ref{example positive elements in a C*-algebra and thomson metric} (observe that $M^{++}$ is 2-divisible).	
\end{proof}

Besides Theorem~\ref{t main geometric tool}, the reference \cite{HatHirMiuMol2012} contains other appropriate tools to determine preservers of inverted Jordan triple products. 

\begin{theorem}\label{t second geometric main tool for preservers of inverted triple products A}\cite[Corollary 3.9]{HatHirMiuMol2012}  Let $(X, d_X)$ and $(Y, d_Y )$ be two non-empty nr-quasi-distance spaces with the additional property that $d_{Y} (y_1,y_2)=0\Rightarrow y_1=y_2$. We also assume that $X$ and $Y$ are subsets of two groups $\mathcal{G}_1$ and 	$\mathcal{G}_2$, respectively, satisfying \eqref{eq condition p HatoriMiuraMolnar}, that is, $y x^{-1} y\in X$ for all $x,y\in X$ and  similarly for $Y$. Let $\Delta : X\to Y$ be a bijective $d$-preserving map, and fix $a,b\in X$. Suppose that the following hypotheses hold: \begin{enumerate}[$(a)$]
		\item $X$ satisfies property $B(a, b)$. 
		\item $Y$ satisfies property $C_1(\Delta(a), \Delta(b a^{-1} b))$. 
\end{enumerate} Then we have $ \Delta (b a^{-1} b) = \Delta(b) \Delta(a)^{-1} \Delta(b).$
\end{theorem}

\begin{theorem}\label{t second geometric main tool for preservers of inverted triple products}\cite[Corollary 3.10]{HatHirMiuMol2012}  Let $(X, d_X)$ and $(Y, d_Y )$ be two non-empty nr-quasi-distance spaces, and suppose additionally that $d_{Y} (y_1,y_2)=0\Rightarrow y_1=y_2$. We also assume that $X$ and $Y$ are subsets of two groups $\mathcal{G}_1$ and 	$\mathcal{G}_2$, respectively, satisfying \eqref{eq condition p HatoriMiuraMolnar}, that is, $y x^{-1} y\in X$ for all $x,y\in X$ and similarly for $Y$. Let $\Delta : X\to Y$ be a bijective $d$-preserving map, and fix $a,b\in X$. Suppose that the following hypotheses hold: \begin{enumerate}[$(a)$]
\item $X$ satisfies property $B(a, b)$. 
\item $Y$ satisfies $C_2(\Delta(a), \Delta(b a^{-1} b))$. 
\item $Y$ is 2-divisible and 2-torsion free. 
\end{enumerate} Then we have $ \Delta (b a^{-1} b) = \Delta(b) \Delta(a)^{-1} \Delta(b).$
\end{theorem}

The previous theorem was employed by Hatori and  Moln{\'a}r in \cite[Theorem 9]{HatMol2012} to prove that every surjective isometry $\Delta: (A^{++},d_{T}^{A})\to (B^{++},d_{T}^{B})$, where $A$ and $B$ are unital C$^*$-algebras, preserves inverted Jordan triple products. In this case $A^{++}$ can be regarded inside $A^{-1}$ which is multiplicative subgroup of $A$. However, for a unital JB$^*$-algebra, the set $M^{++}$ cannot be easily embedded inside a group, and Theorem~\ref{t second geometric main tool for preservers of inverted triple products} is not useful (there are subtle difficulties for properties $B(a,b)$, $C_1(a,b)$ and $C_2(a,b)$ in the Jordan setting). Despite of this handicap, we establish next a Jordan version of \cite[Theorem 9]{HatMol2012}.

\begin{theorem}\label{t Thompson isometries preserve inverted Jordan triple products Jordan} Let $M$ and $N$ be two unital JB$^*$-algebras. Let $\Delta: (M^{++},d_{T}^{M})\to (N^{++},d_{T}^{N})$ be a surjective isometry. Then $\Delta$ preserves inverted Jordan triple products, that is, $ \Delta (U_b(a^{-1})) = U_{\Delta (b)} (\Delta (a)^{-1}),$ for all $a,b\in M^{++}$.
\end{theorem}

\begin{proof} Consider the mapping $\varphi: M^{++}\to M^{++}$ defined by $\varphi (x) := U_{b} (x^{-1})$. We claim that $\varphi$ satisfies the hypotheses in Theorem~\ref{t main geometric tool}. Namely,  $$\varphi^2 (x) = \varphi \left( U_{b} (x^{-1}) \right) = U_{b} \left((U_{b} (x^{-1}))^{-1}\right) = U_{b} (U_{b^{-1}} (x))=x, \ (x\in X),$$
 and hence $\varphi$ is a period-2 bijection. We deduce that $\varphi$ is $d_{T}^{M}$-preserving by \eqref{eq conditions Jordan for B and C2} in Proposition~\ref{p Lemmens on the Thompson distance for positive invertible in JB*-algebras}. It is also clear that $\varphi (b) = b$. \smallskip
 
Pick $c\in N^{++}$ with $U_c (\Delta (a)^{-1}) = \Delta\left(U_b (a^{-1}) \right)$. Clearly, $d_{T}^{N} (U_c (x^{-1}), U_c (y^{-1})) = d_{T}^{N} (x,y)$ for all $x,y\in N^{++}$ (cf. Proposition~\ref{p Lemmens on the Thompson distance for positive invertible in JB*-algebras}). Define a mapping $\psi: N^{++} \to N^{++}$ given by $\psi (y) := U_{c} (y^{-1})$. As before, $\psi$ is period-2 and $d$-preserving. \smallskip

It is not difficult to check that $$\psi (\Delta (a)) = U_{c} (\Delta(a)^{-1}) = \Delta\left(U_b (a^{-1}) \right)= \Delta (\varphi(a))$$ and $$\psi (\Delta (\varphi(a))) = U_c \left( \Delta\left( U_{b} (a^{-1}) \right)^{-1} \right) = U_c \left(\left( U_c (\Delta (a)^{-1}) \right)^{-1} \right) = \Delta (a).$$ 

The boundedness of the set $\{d_{T}^{M} (x, c) : x \in L\}$ in Theorem~\ref{t main geometric tool} follows from the fact that $d_{T}^{M}$ is a metric (cf. Remark~\ref{remark bounded is free for metric}). Having in mind \eqref{eq conditions Jordan for B and C2} we have  $d_{T}^{M} (\varphi (x), x) = d_T^{M} \left(U_b (x^{-1}) , x \right) = 2 d_T^{M} \left(b, x \right)$ for all $x\in M^{++}$. \smallskip

We are in a position to apply Theorem~\ref{t main geometric tool} to deduce that $$d_{T}^{N} (\psi(\Delta (b)), \Delta (b)) = 0,$$ and hence $U_c (\Delta (b)^{-1}) = \Delta(b)$, because $d_{T}^{N}$ is a metric. In such a case, by the fundamental identity we arrive to $$\begin{aligned}
\11b{N} &= U_{\Delta(b)^{-1/2}} (\Delta(b)) =U_{\Delta(b)^{-1/2}} U_c (\Delta (b)^{-1}) = U_{\Delta(b)^{-1/2}} U_c U_{\Delta(b)^{-1/2}} (\11b{N}) \\
&= U_{U_{\Delta(b)^{-1/2}} (c)} (\11b{N}) = \left(U_{\Delta(b)^{-1/2}} (c) \right)^2,
\end{aligned}$$ which implies that $U_{\Delta(b)^{-1/2}} (c) =\11b{N},$ since $U_{\Delta(b)^{-1/2}} (c) \in N^{++}.$ However, $U_{\Delta(b)^{-1/2}} (c) =\11b{N} = U_{\Delta(b)^{-1/2}} (\Delta(b)),$ combined with the invertibility of $U_{\Delta(b)^{-1/2}}$ implies that $c = \Delta(b)$. Finally, it follows from the definition of $c$ that
$$ \Delta (U_b(a^{-1})) = U_{\Delta (b)} (\Delta (a)^{-1}).$$
\end{proof}

\subsection{Positive invertible elements with the Thompson's metric}\label{subsec: positive invertible with Thompson metric as invariant} \ \smallskip

We have already seen in subsection~\ref{subsection positive invertible with the usual metric} that, thanks to Mankiewicz theorem (Theorem~\ref{t Mankiewicz}), for each unital JB$^*$-algebra $M$, the metric space $M^{++}$ equipped with the distance given by the norm is an invariant to identify unital JB$^*$-algebras (cf. Theorem~\ref{p cone of positive elements as metric invariant usual norm}). In this subsection we shall consider the metric space $(M^{++},d_{T})$.\smallskip

The first result in this line was established by Moln{\'a}r in \cite[Theorem 1]{Mol2009}, who proved that if $H$ is a complex Hilbert space with dim$(H) \geq 3$, and $\Delta : (B(H)^{++},d_T) \to (B(H)^{++},d_T)$ is a bijective isometry, then $\Delta$ is either of the form
$\Delta (a) = b a b^*,$ or of the form  $\Delta(a) = b a b^*,$  ($a\in B(H)^{++}$), where $b$ is an invertible bounded linear or conjugate-linear operator on $H$. The result was latter extended to general unital C$^*$-algebras by Hatori and Moln{\'a}r \cite{HatMol2014}.

\begin{theorem}\label{t Hatori Molnar Thompson isometries}\cite[Theorem 9]{HatMol2014} Let $A$ and $B$ be unital C$^*$-algebra. Then a mapping $\Delta : (A^{++},d_T) \to (B^{++},d_T)$ is a surjective isometry if and only if there is a central projection $p$ in $B$ and a Jordan $^*$-isomorphism $J : A\to B$ such that $\Delta$ is of the form
	$$\Delta (a) = \Delta(\11b{A})^{1/2} \left( p J(a) + (1 - p) J(a^{-1}) \right) \Delta(\11b{A})^{1/2}, \ \ a\in A^{++}.$$
\end{theorem}

A similar problem was studied by Lemmens, Roelands, and Wortel in the wider setting of unital JB$^*$-algebras. 

\begin{theorem}\label{t Lemmens Roelands Wortel Thompson isometries}\cite[Theorem 3.2]{LemmRoeWor2018} Let $M$ and $N$ be unital JB$^*$-algebras. Then a mapping $\Delta : (M^{++},d_T) \to (N^{++},d_T)$ is a surjective isometry if and only if there is a central projection $p$ in $N$ and a Jordan $^*$-isomorphism $J : M\to N$ such that $\Delta$ is of the form
	\begin{equation}\label{eq forma a Thompson isometry in the Jordan setting} \Delta (a) = U_{\Delta(\11b{M})^{1/2}} \left( p\circ  J(a) + (\11b{N} - p)\circ J(a^{-1}) \right), \ \ a\in M^{++}.
		\end{equation}
\end{theorem}

The proof of the above result given by Lemmens, Roelands, and Wortel in \cite{LemmRoeWor2018} is based on differential geometry and geodesics. Here we present an adaptation of the ideas by Hatori and Moln{\'a}r in \cite{HatMol2014}.  

\begin{proof} Let us begin by showing that expressions of the form in \eqref{eq forma a Thompson isometry in the Jordan setting} preserve Thompson's metric. It is clear from the definition in \eqref{def Tompson and Hilbert distances} that every Jordan $^*$-isomorphism  $J: M\to N$ preserves Thompson's metric. The second identity in \eqref{eq conditions Jordan for B and C2} in Proposition~\ref{p Lemmens on the Thompson distance for positive invertible in JB*-algebras} assures that maps of the form $a\mapsto a^{-1}$ and $a\mapsto U_{c} (a)$ (with $c$ fixed in $M^{++}$) preserve Thompson's metric. All these statements together prove the ``if'' implication.\smallskip
	
Let us now deal with the ``only if''  implication. By replacing $\Delta$ with $\Delta_0 (x) = U_{\Delta(\11b{M})^{-1/2}} \left(\Delta(x)\right),$ ($x\in M^{++}$), we can assume that $\Delta$ is unital. Theorem~\ref{t Thompson isometries preserve inverted Jordan triple products Jordan} assures that $\Delta$ preserves inverted Jordan triple products, that is, \begin{equation}\label{eq Delta preserves inverted Jordan triple products in demo} \Delta(U_a (b^{-1})) = U_{\Delta(a)} (\Delta(b)^{-1}), \hbox{ for all } a,b\in M^{++}.
\end{equation} Combining this property with the fact that $\Delta$ is unital we arrive to $\Delta (a^2) = \Delta (a)^2$ for all $a\in M^{++}$, and by induction and a standard argument \begin{equation}\label{eq Delta preserves powers and roots in Lemmens thm} \Delta (a^n) = \Delta (a)^n, \hbox{ and } \Delta (a^{\frac{1}{n}}) = \Delta (a)^{\frac{1}{n}}, \hbox{ for all } a\in M^{++},\ n\in \mathbb{N}.
\end{equation}
	  
Now, as in the proof of \cite[Theorem 9]{HatMol2014} and \cite[Theorem 3.2]{LemmRoeWor2018} we consider the next mapping $S: M_{sa}\to N_{sa}$, $S(x) = \log\Delta (e^{x})$. Clearly, $S(0 ) = \log\Delta (\11b{M}) = 0$.  We claim that $S$ is a surjective isometry. \smallskip

To see the claim, let us fix $x,y$ in $M_{sa}$. By the Lie--Trotter formula established in \cite[Corollary 2.2]{EscPeVill2023} (which is valid in wider setting of unital complex Jordan--Banach algebras), we have $$\lim_{n\to \infty} \Big( {U_{e^{-\frac{x}{2n}} } \left(e^{\frac{y}{n}}\right)}\Big)^{n} = e^{y-x}, \hbox{ and hence } \lim_{n\to \infty} \frac{ \log\Big(U_{e^{-\frac{x}{2n}} } \left(e^{\frac{y}{n}}\right) \Big)}{n^{-1}} ={y-x}.$$ It then follows that \begin{equation}\label{eq Thompson distance is asymptotically the usual distance} \lim_{n\to \infty} \frac{d_T^{M} \left( e^{\frac{x}{n}}, e^{\frac{y}{n}}\right)}{n^{-1}} =  \lim_{n\to \infty} \frac{\left\| \log\Big(U_{e^{-\frac{x}{2n}} } \left(e^{\frac{y}{n}}\right) \Big)\right\|}{n^{-1}} = \|x-y\|, 
\end{equation} which asserts that we can asymptotically approximate the usual distance between self-adjoint elements by Thompson's distances of some exponentials.\smallskip

By applying \eqref{eq Delta preserves powers and roots in Lemmens thm} we deduce that  $e^{\frac{S(x)}{n}} = e^{\frac{\log\Delta (e^{x})}{n}} = \Delta (e^{\frac{x}{n}}),$ for all natural $n$.  Having in mind that $\Delta$ is an isometry for the Thompson metrics we also have $$d_{T}^{N}\left(e^{\frac{S(x)}{n}} , e^{\frac{S(y)}{n}} \right) = d_{T}^{N}\left( \Delta (e^{\frac{x}{n}}) , \Delta (e^{\frac{x}{n}}) \right)  =  d_{T}^{M}\left( e^{\frac{x}{n}} , e^{\frac{x}{n}} \right),$$ identity which combined with the conclusion in \eqref{eq Thompson distance is asymptotically the usual distance} leads to $\|S(x)-S(y)\| = \|x-y\|$, witnessing that $S$ is an isometry. The basic properties of the holomorphic functional calculus show that the exponential is a bijection from $M_{sa}$ onto $M^{++}$, and hence $S: M_{sa}\to N_{sa}$ is a surjective isometry with $S(0) =0$. The Mazur--Ulam theorem now affirms that $S$ is a surjective real-linear isometry between two unital JB-algebras. The result by Isidro and Rodr\'{\i}guez Palacios in \cite[Theorem 1.9]{IsRo} proves the existence of a central symmetry $s$ in $N_{sa}$ and surjective Jordan isomorphism $J: M_{sa}\to N_{sa}$ such that $S (x) = s \circ J(x),$ for all $x\in M_{sa}$. It is known that $ s = p - (\11b{N}-p),$ where $p$ is a central projection in $N$. Therefore, for each $x\in M_{sa}$, since elements of the form $p\circ c$ and $p\circ d$ with $c,d\in N_{sa}$ operator commute and have zero Jordan product, we have $$\begin{aligned}
\Delta (e^{x}) &= e^{S(x)} = e^{ (p - (\mathbf{1}-p)) \circ J (x)} =   e^{ p \circ J (x)} \circ e^{ - (\mathbf{1}-p)\circ J (x)} \\
&= (\11b{N}-p + p \circ   e^{ J (x)}) \circ (p+ (\11b{N}-p)  \circ e^{ - J (x)})
\\
&= (\11b{N}-p + p \circ J \left(  e^{ x}\right) )\circ  (p+ (\11b{N}-p)  \circ J \left( e^{- x}\right)) \\
&= p \circ J \left(  e^{ x}\right)  + (\11b{N}-p)  \circ J \left( e^{- x}\right),
\end{aligned}$$ for all $x\in M_{sa}$, where we have used that $J$ is a Jordan isomorphism. This proves that $\Delta (a) = p \circ J(a) + (\11b{N}-p)  \circ J \left( a^{- 1}\right)$ for all $a\in M^{++}.$ 
\end{proof} 

It is clear that Theorem~\ref{t Hatori Molnar Thompson isometries} can be rediscovered as a consequence of the result for JB$^*$-algebras in Theorem~\ref{t Lemmens Roelands Wortel Thompson isometries}.

\section{Unitary elements as a metric invariant}\label{sec: unitaries as invariant}

After having explored the sets of invertible elements and positive invertible elements as metric invariants (with respect to the usual norm and with respect to the Thompson metric) in unital C$^*$- and JB$^*$-algebras in sections~\ref{sec: geometric goodness of C* and JB-algebras} and \ref{sec: geometric tool of Hatori Hirasawa Miura Molnar}, our goal in this section is to explore whether the subset of unitary elements can be also employed to identify these algebras. \smallskip

The first precedent in the literature is a result by Hatori and Moln{\'a}r for $B(H)$. We first recall the notion of abstract transposition on $B(H)$. Let $H$ be a complex Hilbert space. Suppose we fix a conjugation (i.e. conjugate linear isometry of period-$2$ on $H$). For each $a\in B(H)$ we define $\overline{a},a^t\in B(H)$ by $\overline{a} (\xi):= \overline{a(\overline{\xi})},$ and $a^t := \overline{a^*}$ ($\xi\in H$). We can define an abstract transposition and a conjugation on $B(H)$ defined by $a\mapsto a^t$ and $a\mapsto \overline{a}$, respectively. The transposition is complex linear, satisfies $a^{tt} = a$ and $(a b)^t = b^t a^t$, while the conjugation is conjugate-linear and satisfies $\overline{\overline{a}} =a$ and $\overline{a b} = \overline{a} \overline{b}$, for all $a,b\in B(H)$. 

\begin{theorem}\label{t Hatori Molnar unitaries BH}\cite[Theorem 8]{HatMol2012} Let $H$ be a complex Hilbert space. Suppose that $\Delta : \mathcal{U}(B(H)) \to  \mathcal{U}(B(H))$ is a surjective isometry with respect to the metric given by the operator norm. Then there exists a surjective real-linear isometry (equivalently, a triple automorphism) $T$ on $B(H)$ such that $\Delta (u) = T (u)$ for all $u\in \mathcal{U}(B(H))$. Furthermore, there exist unitaries $v,w\in \mathcal{U}(B(H))$ such that $\Delta$ is of one of the following forms:
\begin{enumerate}[$(a)$]
	\item $\Delta (u) = v u w^*,\ $ for all $u\in \mathcal{U}(B(H))$,
	\item $\Delta (u) = v u^* w^*,\ $ for all $u\in \mathcal{U}(B(H))$,
	\item $\Delta(u ) = v u^t w^*,\ $ for all $u\in \mathcal{U}(B(H)),$
	\item $\Delta(u) = v \overline{u} w^*, \ $ for all $u\in \mathcal{U}(B(H))$. 
\end{enumerate}
\end{theorem}

It is worth to note that the concrete description of the mapping $\Delta$ in $(a)$ to $(d)$ above is a consequence of the first conclusion and a classic result affirming that every surjective real-linear isometry (equivalently, triple isomorphism cf. Theorem~\ref{t Dang real-linear surjective isometries JBstar algebras}) $T$ on $B(H)$ is of the one of the forms given in $(a)$ to $(d)$ for all $u\in B(H)$ in this case (see \cite[page 162]{FleJa2003v2}).\smallskip 

The same authors of the previous theorem also determined the precise form of all surjective isometries between the unitary groups of two von Neumann algebras. 

\begin{theorem}\label{t Hatori Molnar surjecitve isometries between unitaries vN}\cite[Theorem 1]{HatMol2014} Let $A$ and $B$ be unital C$^*$-algebras. Assume that $\Delta  : \mathcal{U} (A) \to \mathcal{U} (B)$ is a surjective isometry between the corresponding groups of unitary elements (with respect to the C$^*$-norms). Then we have $\Delta (e^{i A_{sa}}) = \Delta(\11b{A})\  e^{ i B_{sa}}$ and there is a central projection $p \in B$ and a Jordan $^*$-isomorphism $J : A\to B$ such that
\begin{equation}\label{eq fundamental conclusion on the form of a si on unitaries} \Delta (e^{ix}) = \Delta(\11b{A}) (p J(e^{ix}) + (\11b{B} - p) J(e^{ix})^*),\  x\in A_{sa}.
\end{equation} Furthermore, if $A$ or $B$ is a von Neumann algebra, then both algebras are von Neumann algebras and $\Delta$ admits an extension to a real-linear isometry from $A$ onto $B$, more precisely, there exists a central projection $p \in B$ and a Jordan $^*$-isomorphism $J : A\to B$ such that the conclusion in \eqref{eq fundamental conclusion on the form of a si on unitaries} holds when $e^{ix}$ is replaced with an arbitrary unitary $u$ in $A$. 
\end{theorem}   

\begin{corollary}\label{c Hatori Molnar isometrically isomorphic}\cite[Corollary 2]{HatMol2014}. Two unital C$^*$-algebras are isomorphic as Jordan $^*$-algebras if and only if their unitary groups are isometric as metric spaces. 
\end{corollary}

As in previous sections, we shall include in this note versions of these results in the wider setting of unital JB$^*$-algebras with a complete proof (see subsection~\ref{subsect: isometries between connected components}). Before working with JB$^*$-algebras, we shall turn back our view to the geometric tools in section~\ref{sec: geometric tool of Hatori Hirasawa Miura Molnar} with the aim of establishing a certain local properties of surjective isometries over subgroups of surjective linear isometries which will be useful in subsequent results.  \smallskip

The powerful geometric tool developed in  Theorem~\ref{t main geometric tool} and the subsequent combinations with properties $B(a,b)$, $C_1(a,b)$ and $C_2 (a,b)$ in Theorems~\ref{t second geometric main tool for preservers of inverted triple products A} and \ref{t second geometric main tool for preservers of inverted triple products} played a fundamental role in the study of surjective isometries between sets positive invertible elements in unital C$^*$-algebras \cite{Mol2009,HatMol2014} and JB$^*$-algebras \cite{LemmRoeWor2018}. These tools will be also very important in the results presented in this subsection.\smallskip

Keeping the notation in previous sections, for each Banach space $X$, the symbol Iso$(X)$ will stand for the surjective linear isometries on $X$, which is clearly a multiplicative subgroup of $B(X)$. The group Iso$(X)$ serves as a test bench for properties $B(a,b)$ and $C_1(a,b)$. It will be checked in the proof of the next result, whose proof has been slightly simplified. 

\begin{theorem}\label{t Theorem 6 Haotri Molnar 2012}\cite[Theorem 6]{HatMol2012} Let $X$ and $Y$ be complex Banach spaces and let $\mathcal{G}$ and $\tilde{\mathcal{G}}$ be a subgroup of Iso$(X)$ and Iso$(Y)$, respectively, equipped with the metrics coming from the operator norm. Then $\mathcal{G}$ (respectively, $\tilde{\mathcal{G}}$) always satisfies property $C_1(U,V)$ for all $U,V\in \mathcal{G}$ (respectively, for all $U,V\in \tilde{\mathcal{G}}$). Suppose we take $U, V \in \mathcal{G}$ with $d(U, V ) < 1/2$. Then $\mathcal{G}$ satisfies property $B(U,V)$. Furthermore, every surjective isometry  $\Delta : \mathcal{G}\to  \tilde{\mathcal{G}}$ satisfies 
$$ \Delta (V U^{-1} V ) = \Delta (V ) \Delta(U)^{-1} \Delta(V )$$ for all $U, V \in \mathcal{G}$ with $d(U, V ) < 1/2$.
\end{theorem}

\begin{proof} The three statements are intrinsically related. Since $\mathcal{G}$ is a subgroup of Iso$(X)$ and the operator norm satisfies $$\| U R^{-1} V- U S^{-1} V\| = \| R^{-1}- S^{-1}\| = \|R-S\|$$ for all $R,S,U,V\in \hbox{Iso}(X)$, it follows that $\mathcal{G}$ satisfies property $C_1 (U,V)$ for all $U,V \in \mathcal{G}$. The final conclusion will follow from Theorem~\ref{t second geometric main tool for preservers of inverted triple products A} as soon as we prove that $\mathcal{G}$ satisfies property $B(U,V)$ for all $U,V\in \mathcal{G}$ with $\|U-V\|<\frac12$. \smallskip
	
For $(B.1)$ we observe that  
$$\| V W^{-1} V - V R^{-1} V \| =\|W^{-1}-R^{-1}\| = \| I- W R^{-1}\| = \|R- W\|,$$ for all $U,V,R,W\in \hbox{Iso} (X)$. Property $(B.2)$ trivially holds because the distance given by the operator norm is a metric. The non-trivial property is $(B.3)$. For this purpose take $W\in L_{U,V}$ (that is, $\|U-W\| = \|V U^{-1} V -W\| = \|U-V\|$). The constant $K = 2-2\|U-V\| >1$ by the hypothesis. \smallskip

We shall first prove that $V+W$ is invertible. Namely \begin{equation}\label{eq distance from V to W} \|W-V\|\leq \|W-U\| + \|U-V\| = 2 \|U-V\| <\frac12, 
\end{equation} and consequently, $\frac12 > \left\|\frac{V+W}{2} -V \right\| = \left\|\frac{V+W}{2} V^{-1} -I  \right\|,$ which gives the desired statement. \smallskip

We claim next that \begin{equation}\label{eq norm of V+W inverse} \|(W+V)^{-1}\| \leq \frac1K.
\end{equation} Indeed, 
$$\begin{aligned}
	1 &= \| (W+V) (W+V)^{-1}\|  = \| (2 W - (W-V)) (W+V)^{-1}\| \\
	&\geq  \| 2 W (W+V)^{-1}\| - \| (W-V) (W+V)^{-1}\| \\
	&\geq  2 \| (W+V)^{-1}\| - \| W-V \| \| (W+V)^{-1}\| \\
	&\geq  \hbox{(by \eqref{eq distance from V to W})} \geq 2 \| (W+V)^{-1}\| - 2 \| U-V \| \| (W+V)^{-1}\| \\
	& = (2 - 2 \| U-V \|) \| (W+V)^{-1}\|  = K \| (W-V)^{-1}\|.
\end{aligned}$$
 
Now, we employ \eqref{eq norm of V+W inverse} and the identity $\left( W (W+V)^{-1}\right)^{-1} = I + V W^{-1}$ to deduce that \begin{equation}\label{eq norm of W(W+V)inverse inverse} \begin{aligned} \left\| \left( I + V W^{-1} \right)^{-1} \right\| =  \left\|W (W + V)^{-1} \right\| =  \left\|(W + V)^{-1}  \right\|\leq \frac1K.
\end{aligned}
\end{equation} 

Finally, by \eqref{eq norm of W(W+V)inverse inverse},
$$\begin{aligned}
K \| W-V\|  &= K \| V W^{-1} -I\| =  K \| (V W^{-1} -I) (V W^{-1} -I) (V W^{-1} -I)^{-1}\| \\
&\leq   K  \| (V W^{-1} -I)^{-1}\| \| (V W^{-1} -I) (V W^{-1} -I)\|
\\
&\leq  \| (V W^{-1} -I) (V W^{-1} -I)\|  = \| V W^{-1} V W^{-1} -I\| \\
&= \| V W^{-1} V  -W\|, 
\end{aligned}$$ which concludes the proof of property $B(U,V)$, and hence the arguments in this proof. 	
\end{proof}

\begin{remark}\label{remark unitaries satisfy B and C1} Theorem~\ref{t Theorem 6 Haotri Molnar 2012} above implies the following result in the setting of unital C$^*$-algebras: Let $A$ be a unital C$^*$-algebra, and let $\mathcal{U} (A)$ be the group of unitary elements in $A$. Since, by the Gelfand--Naimark theorem $A$ can be regarded as a C$^*$-subalgebra of some $B(H)$ and the unit of $A$ coincides with the identity on $H$. Since $\mathcal{U} (A)$ can be regarded as a subgroup of Iso$(H)$, Theorem~\ref{t Theorem 6 Haotri Molnar 2012} implies that $\mathcal{U} (A)$ satisfies $C_1(a,b)$ for all $a,b\in \mathcal{U} (A)$, and $B(a,b)$ for all $a,b\in \mathcal{U} (A)$ with $\|a-b\|<\frac12$. Furthermore, if $B$ is another unital C$^*$-algebra and $\Delta: \mathcal{U} (A) \to \mathcal{U} (B)$ is a surjective isometry, we have $\Delta (a b^{-1} a) = \Delta (a) \Delta (b)^{-1}  \Delta( a)$ for all $a,b\in \mathcal{U} (A)$ with $\|a-b\|<\frac12$.
\end{remark}

\begin{remark}\label{r condition d(u,v)<1/2 in theorem 6 is necessary} As observed by Hatori and Moln{\'a}r in \cite[page 2130]{HatMol2012}, the hypothesis $d (U,V)<\frac12$ is absolutely necessary for the final conclusions in Theorem~\ref{t Theorem 6 Haotri Molnar 2012}. Consider, for example, a compact Hausdorff space $\Omega$ with at least two points and the Banach space $C(\Omega)$. Let
$\mathcal{G} = \hbox{Aut} (C(\Omega))$ denote the subgroup of Iso$C(\Omega)$ all composition operators on $C(\Omega)$ of the form $T_{\varphi} (f) = f\circ \varphi$ ($f \in C(\Omega)$), where $\phi$ runs in the set of all self-homeomorphisms on $\Omega$ (i.e., the group of all $^*$-automorphisms on $C(\Omega)$). By applying Urysohn’s lemma it can be checked that $\|T_{\varphi} - T_{\phi} \| = 0$ if $\phi= \varphi$ and  $\|T_{\varphi} - T_{\phi} \| =2$, otherwise. Therefore, every bijective
transformation of $\mathcal{G}$ onto itself is a surjective isometry. 
\end{remark}

It is worth commenting at this point that finding an extension of every surjective isometry between the unitary groups of two unital C$^*$-algebras is not always possible (that is, the final conclusion in Theorem~\ref{t Hatori Molnar surjecitve isometries between unitaries vN} is not, in general valid when von Neumann algebras are replaced with mere unital C$^*$-algebras). A counter example can be found even in the setting of commutative unital C$^*$-algebras. A complete characterization of those surjective  isometries between the unitary groups of two commutative unital C$^*$-algebras was established by Hatori and Moln{\'a}r in \cite[Corollary 8]{HatMol2014}. We shall provide concrete methods to produce counterexamples in the next subsection.    

\subsection{Unitaries in unital JB$^*$-algebras and the principal component}\label{subsec: connected components of the unitay set}\ \smallskip

We recall that an element $u$ in a unital JB$^*$-algebra $M$ is called unitary if $u$ is invertible with $u^* = u^{-1}$, which is equivalent to say to say that $u$ is a unitary tripotent in $M$ when the latter is regarded as JB$^*$-triple (cf. \cite[Proposition 4.3]{BraKaUp78} or the commentes in page \pageref{comments unitaries as JB*-algebra an triple coincide}).\smallskip 

The Jordan product of two unitaries need not be, in general, a unitary. However, given $u,v\in \mathcal{U} (M)$, the element $U_{u} (v)$ lies in $\mathcal{U} (M)$ (cf. Lemma~\ref{l first properties invertible Jordan}). As in the case of C$^*$-algebras, the set $\mathcal{U} (M)$ is not, in general, connected. The connected component of $ \mathcal{U} (M)$ containing the unit element is called the \emph{principal connected component} or simply the \emph{principal component} and will be denoted by $\mathcal{U}^{0}(M)$.\smallskip

Along this note we shall employ on several occasions that for each $u\in \mathcal{U} (M)$ the mapping $U_u$ is a surjective linear isometry (cf. \cite[Theorem 4.2.28$(vii)$]{Cabrera-Rodriguez-vol1}).\smallskip  

As in previous sections, we shall present the results in the setting of JB$^*$-algebras in an analogy with the corresponding versions for C$^*$-algebras, to allow the reader to taste similarities and divergences.\smallskip

Let $u$ be a unitary element in a unital C$^*$-algebra $A$. A well known result in C$^*$-algebra theory shows that $\|\11b{A} -u\|<2$ is a sufficient condition to guarantee that $u$ admits a ``logarithm'', that is, $u = e^{ih}$ for some $h\in A_{s}$ (see \cite[Exercise 4.6.6]{KR1}). The value $2$ in the inequality is optimal; actually the conclusion is not true without this assumption {\rm(}see the discussion preceding Proposition 4.4.10 in \cite{KR1}{\rm)}. We present next an  analogue of this result in the setting of JB$^*$-algebras. Along the next technical results, the Shirshov-Cohn theorem will play a central role.\smallskip

In the next lemma we shall only need the results affirming that a JB$^*$-subalgebra generated by a single self-adjoint element and the unit element is isometrically isomorphic to a commutative  unital C$^*$-algebra \cite[3.2.4. The spectral theorem]{HOS}.  Let $M$ be a unital JB$^*$-algebra, and let $u$ be a unitary in $M$. To avoid confusions, when we compute the exponential of an element $a$ in the $u$-isotope JB$^*$-algebra $M_{(u^*)}$  we shall denote $exp_{(u)} (a)$ or simply $\exp(a) = e^a$ if no confusion arises.  

\begin{lemma}\label{l untaries at short distance in a unital JB$^*$-algebra} Let $M$ be a unital JB$^*$-algebra $M$. Let  $u$ and $v$ be unitaries in $M$ satisfying $\|u-v\|= \eta <2.$ Then the following statements hold:\begin{enumerate}[$(a)$] \item There exists a self-adjoint element $h$ in the $u$-isotope JB$^*$-algebra $M(u)= M_{(u^*)}= (M,\circ_u,*_u)$ such that $v = \exp_{(u)}(i h),$ where the exponential is computed in the JB$^*$-algebra $M(u)$. In particular $u$ and $v$ lie in the same connected component of $\mathcal{U} (M)$.
\item There exists a unitary $w$ in $M$ satisfying $U_w (u^*) = v$.
\item The equality $v^2 = u$ is equivalent to $U_{v^*} (u) = \11bM$.
\item  If $U_{v} (u^*) = \{v,u,v\} = u$  we have $w = u$.
	\end{enumerate} Moreover, if $\|u-v\|= \eta = \left| 1- e^{it_0} \right| = \sqrt{2}\sqrt{1-\cos(t_0)}$ for some $t_0\in (-\pi,\pi),$ we can further assume that $\|w-u\|,\|w-v\| \leq \sqrt{2}\sqrt{1-\cos(\frac{t_0}{2})}.$
\end{lemma}

\begin{proof} $(a)$ We consider the unital JB$^*$-algebra $M(u)= M_{(u^*)}=(M,\circ_u,*_u)$. Let $\mathcal{C}$ denote the JB$^*$-subalgebra of $M(u)$ generated by $v$ and its unit --i.e. $u$--. Let us observe that the product and involution on $\mathcal{C},$ which will be denoted by ``$\cdot$'' and $*_{u}$, respectively, are precisely $\circ_u|_{\mathcal{C}\times \mathcal{C}}$ and $*_u|_{\mathcal{C}}$, respectively. Since $v$ also is unitary in $M(u)$ (cf. Lemma~\ref{l unitaries in the homotope coincide}), the JB$^*$-subalgebra $\mathcal{C}$ must be isometrically Jordan $^*$-isomorphic to a unital and commutative C$^*$-algebra, that is, to some $C(\Omega)$ for an appropriate compact Hausdorff space $\Omega,$ and under this identification, $u$ corresponds to the unit  (cf. \cite[3.2.4. The spectral theorem]{HOS} or \cite[Proposition 3.4.2 and Theorem 4.1.3$(v)$]{Cabrera-Rodriguez-vol1}).\smallskip
	
The hypothesis $\|u-v\| <2$ is of pure geometric nature, and so it holds in $\mathcal{C}$. Since $v$ is a unitary in $\mathcal{C}$ and $u$ is the unit element of the C$^*$-algebra $\mathcal{C} \equiv C(\Omega)$, we can find a self-adjoint element $h\in \mathcal{C}_{sa}$ such that $v = \exp_{(u)}({i h}),$ where the exponential is computed in $\mathcal{C}$ or equivalently in $M(u) = M_{(u^*)}$ (cf. \cite[Exercise 4.6.6]{KR1}).\smallskip

To see the final statement, we simply observe that the curve $\gamma(t) = \exp_{(u)}({i t h})$ ($t\in[0,1]$), where the exponential is computed in $M(u)$, takes values in $\mathcal{U} (M)$ with $\gamma(0) = u$ and $\gamma (1) = v$.\smallskip

$(b)$ We keep the notation employed in the proof of $(a)$. We observe that $w= \exp_{(u)}({i \frac{h}{2}})$ is a unitary element in $\mathcal{C}$ satisfying $w\cdot u^{*_u} \cdot w = w\cdot u \cdot w = v$ (let us observe that the involution of $\mathcal{C}$ is precisely the restriction of $*_{u}$ to $\mathcal{C}$). Let $U_a^{(u)}$ denote the $U$-operator on the unital JB$^*$-algebra $M(u)=M_{(u^*)}$ associated with an element $a$. Since $\mathcal{C}$ is a unital JB$^*$-subalgebra of $M(u)$, we compute the corresponding operations in $\mathcal{C}$ to deduce, thanks to the uniqueness of the triple product, that $$v= w\cdot u^{*_u}\cdot w = U_w^{(u)} (u^{*_u}) = \{w,u,w\} = U_{w} (u^*),$$ and clearly $w$ is a unitary in $M$ because it is a unitary in $M(u)$ (cf. Lemma~\ref{l unitaries in the homotope coincide}).\smallskip

$(c)$ If $v^2 =u$, it is easy to see that $U_{v^*} (u) = U_{v^*} (v^2) = U_{v^*} U_v (\11bM) = \11bM.$ Reciprocally, if  $U_{v^*} (u) = \11bM,$ we get $$v^2 = U_v (\11bM) = U_v U_{v^*} (u) = u.$$ 

$(d)$ We consider in this case the JB$^*$-subalgebra $\mathcal{D}$ of $M(v)=M_{(v^*)}=(M,\circ_v,*_v)$ generated by $u$ and the unit ($v$ in this case). As before, we identify $\mathcal{D}$ with an appropriate $C(\Omega)$ in such a way that $u$ is a unitary in this commutative unital C$^*$-algebra with $u^{*_{v}} = U_{v} (u^*) = u$ (i.e. $u$ self-adjoint in $\mathcal{D}$). The hypothesis $\|u-v\|<2$ implies that $v=u$.\smallskip

To see the final conclusion we employ, once again, the identification of $\mathcal{C}$ with $C(\Omega)$ in which $u$ corresponds to the unit and $v$ and $w$ with $\exp_{(u)}({i h})$ and $\exp_{(u)}({i \frac{h}{2}}),$ respectively.\smallskip
\end{proof}

The next observation is a straightforward consequence of Lemma~\ref{l untaries at short distance in a unital JB$^*$-algebra}$(a)$. For the proof we simply observe that for any unitary $u$ in a unital JB$^*$-algebra $M$, the set $\mathbb{T} u$ lies in the same connected component of $\mathcal{U} (M)$ that $u$. 

\begin{corollary}\label{c distance in different connected components is 2} Let $u$ and $v$ be unitaries in a unital JB$^*$-algebra. Then if $u$ and $v$ lie in two different connected components of $\mathcal{U} (M)$ we have $\|u\pm v\| =2$.
\end{corollary}

Thanks to the previous corollary, we can now exhibit examples of surjective isometries between sets of unitaries in unital JB$^*$-algebras which cannot be extended to surjective real-linear isometries between the corresponding algebras. 

\begin{remark}\label{r unitaries which are squares}\cite{CuEnHirMiuPe2022} Let $u$ be a unitary element in a unital JB$^*$-algebra $M$. We cannot always guarantee that $u = e^{i h}$ for some $h\in M_{sa}$ {\rm(}cf. \cite[Exercise 4.6.9]{KR1}{\rm)}. Suppose that $u$ satisfies the weaker condition that $u = v^2$ for some unitary $v\in M$ --this automatically holds if $u = e^{i h}$. Let $\mathcal{U}(M)^c$ denote the connected component of $\mathcal{U}(M)$ containing $u$. Since the mapping $U_v: M\to M(u)$ is a unital surjective linear isometry --and hence a Jordan $^*$-isomorphism between these two algebras (see Theorem~\ref{t Kaup-Banach-Stone})-- $U_{v} (\mathcal{U}^{0}(M))$ is a connected component of $\mathcal{U} (M) =\mathcal{U} (M(u))$ (see Lemma~\ref{l unitaries in the homotope coincide}) and contains $u= U_{v} (\mathbf{1})$. It then follows $ \mathcal{U}^c (M) = U_{v} (\mathcal{U}^{0}(M)).$\smallskip
	
The existence of a unitary square root does not necessarily hold for any unitary element, even in the setting of commutative C$^*$-algebras. Let $\mathbb{T}$ denote the unit sphere of the complex plane. The principal component of $\mathcal{U}(C(\mathbb{T}))$ is precisely the subgroup ${\rm\exp}(i C_{\mathbb{R}}(\mathbb{T}))$, which is the set of all functions $u: \mathbb{T}\to \mathbb{T}$ which are deformable or homotopic to the unit element {\rm(}cf. \cite[Exercise 4.6.7]{KR1}{\rm)}. The quotient group $\mathcal{U}(C(\mathbb{T}))/ \mathcal{U}^0(C(\mathbb{T})) = \mathcal{U}(C(\mathbb{T}))/ \exp(i C_{\mathbb{R}}(\mathbb{T}))$ --known as the Bruschlinsky group-- identifies with the ring $\mathbb{Z}$ of integers {\rm(}see \cite[\S II.3]{Hu59} or \cite[Exercise I.11.3]{Tak}{\rm)}. It is known that the maps $u_1, u_2:\mathbb{T}\to \mathbb{T},$ $u_1(\lambda)= \lambda$ and $u_2(\lambda)= \lambda^2$ are not in the principal component of $\mathcal{U}(C(\mathbb{T}))$, they are actually in two different connected components {\rm(}cf. \cite[\S II, Lemma 3.2]{Hu59}{\rm)}. Therefore, $u_2\notin \exp(i C_{\mathbb{R}}(\mathbb{T}))$ but $u_2 = u_1^2$ admits a unitary square root.
\end{remark}

We can now present an example of a surjective isometry between the unitary sets of two unital JB$^*$-algebras which cannot be extended to a surjective real-linear isometry. 
 
\begin{example}\label{example existence of non-extendible surjective isometries between unitaries} Let $M$ be a unital JB$^*$-algebra whose set of unitaries is not connected. Let $\mathcal{U}^{c} (M)$ be a connected component of $\mathcal{U} (M)$ such that $\mathcal{U}^{0} (M)\neq \mathcal{U}^{c} (M)$. Let us take an element $u\in \mathcal{U}^{c} (M)$ and any non-zero $h\in \left(M_{(u)}\right)_{sa}$.  Set $w= \exp_{(u)} (i h ), v= \exp_{(u)} (i \frac{h}{2} )\in \mathcal{U} (M) =\mathcal{U} (M_{(u)})$. Clearly, $w,v$ lie in the connected component of $\mathcal{U} (M)$ containing $u$, that is, $u,v,w\in \mathcal{U}^{c} (M),$ and since $U_{v} (u) = \{v,u,v\} = v\circ_{u} v =  \exp_{(u)} (i \frac{h}{2} ) \circ_{u} \exp_{(u)} (i \frac{h}{2} )  = \exp_{(u)} (i {h} ),$ the surjective linear isometry $U_v$ maps $\mathcal{U}^{c} (M)$ onto itself. Clearly, $U_{v}\neq Id_{M}$. Define a surjective mapping $\Delta: \mathcal{U} (M)\to \mathcal{U} (M)$ by $\Delta (a) = U_v (a)$ if $a\in \mathcal{U}^{c} (M)$ and $\Delta (a) = a$ otherwise. $\Delta$ is an isometry thanks to Corollary~\ref{c distance in different connected components is 2}. If there exists an extension of $\Delta$ to a surjective real-linear isometry $T: M\to M$, the conditions $\Delta(\11b{M}) = \11b{N},$ and  $\Delta (\lambda a) = \lambda a$ for all $a\in \mathcal{U} (M) \backslash \mathcal{U}^{c} (M)$ and Theorem~\ref{t Dang real-linear surjective isometries JBstar algebras} assure that $T$ must be a Jordan$^*$-homomorphism, Furthermore, since $T|_{\mathcal{U}^0(M)} = Id_{M}|_{\mathcal{U}^0(M)}$ and $\exp (i M_{sa})\subseteq \mathcal{U}^0(M),$ the mapping $T$ must be the identity, which contradicts that $\Delta|_{\mathcal{U}^{c} (M)} = U_{v} |_{\mathcal{U}^{c} (M)}.$ See \cite[Corollary 8]{HatMol2014}) and \cite[Corollary 5.1]{Hatori14} for additional examples.
\end{example}

We have described the principal component of the invertible elements in a unital Jordan--Banach algebra in subsection~\ref{subsect: connected components of invertible elements}. Now it is the turn of the principal component of the unitary elements in a unital JB$^*$-algebra. It is perhaps worth to revisit the structure in the case of unital C$^*$-algebras. The next theorem gathers some classical results and a more recent reinterpretation due to Hatori. 

\begin{theorem}\label{t principal connected component in unital C$^*$-algebras}{\rm(\cite[Exercises 4.6.6 and 4.6.7]{KR1}, \cite[Exercise 2, page 56]{Tak}, and \cite[Lemma 3.2]{Hatori14})} Let $A$ be a unital C$^*$-algebra, and let $\mathcal{U}^0(A)$ denote the principal component of the set $\mathcal{U}(A)$ of all unitary elements in $A$. Then  \begin{equation}\label{eq algebraic descirption of the principal component for unital Cstar algebras} \begin{aligned}
			\mathcal{U}^0(A) &= \{ e^{i h_1} \cdot \ldots \cdot e^{i h_n} : n\in \mathbb{N}, \ h_1, \ldots, h_n\in A_{sa} \}\\ &= \{ e^{i h_1} \cdot \ldots \cdot e^{i h_n} \11b{A} e^{i h_n} \cdot \ldots \cdot e^{i h_1} : n\in \mathbb{N}, \ h_1, \ldots, h_n\in A_{sa} \}
		\end{aligned}
	\end{equation} is open, closed, and path connected in the (relative) norm topology on $\mathcal{U}(A)$. Furthermore, by the continuity of the module mapping $x\mapsto |x|:=(x^* x)^\frac12$, we have $$A^{-1}_{\mathbf{1}} \cap \mathcal{U}(A) =\mathcal{U}^0(A),$$ where $A^{-1}_{\mathbf{1}}$ stands for the principal component of the invertible elements in $A$.  
\end{theorem}

The final conclusion in the previous theorem is usually obtained via ``polar decompositions'', this tool is not, in general, available in the setting of JB$^*$-algebras.\smallskip

Our next goal is to present a recent description of the principal component of the set of unitaries in a unital JB$^*$-algebra. To this end, let us fix a unital JB$^*$-algebra $M$. Let us consider a unitary element $u$ in $M$, and the $u$-isotope $M(u)=M_{(u^*)}=(M,\circ_u,*_u)$ as defined in previous sections. We know that $\mathcal{U}(M) =\mathcal{U}(M(u))$ (cf. Lemma~\ref{l unitaries in the homotope coincide}). Since the connected components of a topological space are disjoint, we can affirm that
\begin{equation}\label{eq princ components of isotopes coincide}
	\mathcal{U}^0(M(u))=\mathcal{U}^0(M), \text{ for every $u\in \mathcal{U}^0(M)$.}
\end{equation}
 
Our arguments below require the employment of another tool coming from JB$^*$-triple theory. The subtriple generated by a single element admits an useful concrete description. This is not natural for single generated C$^*$-sualgebras, which cannot be identified, via Gelfand theory, with a commutative C$^*$-algebra. The JB$^*$-subtriple of a JB$^*$-triple $E$ generated by a single element $a\in E$ (denoted by $E_a$) coincides with the norm closure of the linear span of the odd powers of $a$ defined as $a^{[1]}= a$, $a^{[3]} = \{ a,a,a\}$, and $a^{[2n+1]} := \{ a,a,{a^{[2n-1]}}\},$ $(n\in \mathbb{N})$. The Gelfand theory for commutative JB$^*$-triples shows that $E_a$ is isometrically JB$^*$-triple isomorphic to some $C_0 (\Omega_{a})$ for a unique compact Hausdorff space $\Omega_{a}$ contained in the set $[0,\|a\|],$ such that $0$ cannot be an isolated point in $\Omega_a$, where, $C_0 (\Omega_{a})$ stands for the Banach space of all complex-valued continuous functions on $\Omega_a$ vanishing at $0$ in case that $0\in \Omega_a$, and all continuous functions on $\Omega_{a}$ otherwise. We can further assume that under the corresponding isomorphism 
$a$ identifies with the inclusion function of $\Omega_a$ in $\mathbb{C}$ (cf. \cite[Lemma 3.2, Corollary 3.4 and Proposition 3.5]{Ka96}, \cite[Corollary 1.15]{Ka83}, see also \cite[Theorem 4.2.9]{Cabrera-Rodriguez-vol1}). The set $\Omega_a$ is called \emph{the triple spectrum} of $a$ (in $E$), and it does not change when computed with respect to any JB$^*$-subtriple $F\subseteq E$ containing $a$ \cite[Proposition 3.5]{Ka96}. The reader should be warned that in some reference, like \cite[Theorem 4.2.9]{Cabrera-Rodriguez-vol1}, the triple spectrum of an element $a$ is defined as the locally compact set $\Omega \backslash \{0\}$. Here we follow the notation in \cite{Ka96}, where, in a more natural terminology, the triple spectrum is a compact set.
\smallskip

We can naturally define the \emph{continuous triple functional calculus} at the element $a$. Let us consider, via Gelfand theory, the triple isomorphism $\Psi : E_a \to C_{0}(\Omega_a)$. For each continuous function $f\in C_{0}(\Omega_a)$, we set $f_t(a) :=\Psi^{-1} (f)$, and we call it the \emph{continuous triple functional calculus} of $f$ at $a$. If $p(\lambda) = \alpha_1 \lambda + \alpha_3 \lambda^3 + \ldots + \alpha_{2 n -1} \lambda^{2 n -1}$ is an odd polynomial with complex coefficients, it is easy to check that $p_t (a) = \alpha_1 a + \alpha_3 a^{[3]} + \ldots + \alpha_{2 n -1} a^{[2 n -1]}$. However, for more general functions in $C_0(\Omega_a)$ the continuous triple functional calculus is not so obvious. For example if $a$ is a non-self-adjoint element in a C$^*$-algebra $A,$ the element $a^{[2]}= h_t(a)$ with $h(t) = t^2$ does not coincide with $a^2$ in $A$.\smallskip

 Another example: the function $g(t) =\sqrt[3]{t}$ produces $g_t(a) = a^{[\frac{1}{3}]}$. Let us note that $g_t(a) = a^{[\frac{1}{3}]}$ is the unique \emph{cubic root} of $a$ in $E_a$, i.e. the unique element $a^{[\frac13 ]}\in E_a$ satisfying \begin{equation}\label{eq existence of cubic roots} \{ {a^{[\frac13]}},{a^{[\frac13 ]}},{a^{[\frac13 ]}}\}=a.
\end{equation} The sequence $(a^{[\frac{1}{3^n}]})_n$ can be recursively defined by $a^{[\frac{1}{3^{n+1}}]} = \left(a^{[\frac{1}{3^{n}}]}\right)^{[\frac 13]}$, $n\in \mathbb{N}$. It is easy to see that the sequence $(a^{[\frac{1}{3^n}]})_n$ need not be, in general, norm convergent in $E$ (consider, for example, the element $a(t) =t$ in $E=C([0,1])$). We shall need another topology to assure the convergence of the latter sequence.  

If we strength the assumptions on $E$, for example if $a$ is an element in a JBW$^*$-triple $W$ (in particular when $E$ is regarded inside $E^{**}$), the sequence $(a^{[\frac{1}{3^n}]})_n$ converges in the weak$^*$-topology of $W$ to a (unique) tripotent denoted by $r(a)$. The tripotent $r(a)$ is called the \emph{range tripotent} of $a$ in $W$. This range tripotent $r(a)$ can be characterized as the smallest tripotent $e\in W$ satisfying that $a$ is positive in the JBW$^*$-algebra $W_{2} (e)$ (compare \cite[Lemma 3.3]{EdRu96}). When $E$ is just JB$^*$-triple, we shall denote by $r(a)$ the range tripotent of $a$ in $E^{**}$.\smallskip

The range tripotents of von Neumann regular elements in a JB$^*$-triple $E$ always provide elements in $E$.  We recall that an element $a$ in a JB$^*$-triple $E$ is called \emph{von Neumann regular} if $a\in Q(a) (E)$; if $a\in Q(a)^2 (E)$ we say that $a$ is \emph{strongly von Neumann regular} (cf.  \cite{BurKaMoPeRa,FerGarSanSi,Ka96,JamPerSidTah2015}). For each von Neumann regular element $a$ in a JB$^*$-triple $E$ there might exist many different elements $c$ in $E$ satisfying $Q(a)(c) =a$. According to \cite[Theorem 1]{FerGarSanSi}, 
\cite[Lemma 4.1, Lemma 3.2 and comments after its proof]{Ka96}, \cite[Theorem 2.3 and Corollary 2.4]{BurKaMoPeRa}, the following five statements are equivalent reformulations:
\begin{enumerate}[$(a)$]\item $a$ is von Neumann regular;
	\item There exists $b\in E$ such that $[Q(a),Q(b)]=Q(a) Q(b)$ $ - Q(b)\, Q(a)=0$, $Q(a) (b) =a,$ and  $Q(b) (a) =b$;
	\item $0$ does not belong to the triple spectrum, $\Omega_a,$ of $a$;
	\item $Q(a)$ has norm-closed range;
	\item The range tripotent $r(a)$ of $a$ lies in $E$ and $a$ is positive and invertible in the JB$^*$-algebra $E_2 (r(a))$.
\end{enumerate} The element $b$ appearing in statement $(b)$ above is unique. It is called the \emph{generalized inverse} of $a$ in $E$  and is denoted by $a^{\dag}$. The symbol  Reg$_{_{vN}}(E)$ will stand for the set of all von Neumann regular elements in $E$. Another interesting property assures that $$ L(a,a^{\dag}) = L(a^{\dag},a) = L(r(a),r(a)), \hbox{ and } Q(a) Q(a^{\dag}) = Q(a^{\dag}) Q(a) = P_2(r(a))$$ (see, for example, \cite[page 589]{JamPerSidTah2015} or \cite[comments in page 192]{BurKaMoPeRa}).\smallskip

As a illustrative example, we mention that in a unital JB$^*$-algebra $M$, we have $M^{-1}\subseteq $Reg$_{vN} (M)$, and for each $z\in M^{-1},$ its range tripotent, $r(z),$ lies in $\mathcal{U}(M)$: \begin{equation}\label{eq range of invertible is unitary} z\in M^{-1} \Rightarrow r(z)\in \mathcal{U} (M) \ \ \hbox{ (cf. \cite[Lemma 2.2 and Remark 2.3]{JamPerSidTah2015}).}
\end{equation} A brief argument reads as follows: $$\begin{aligned}Q(z) ((z^{-1})^*) &= \{z, (z^{-1})^*, z\} = U_{z} (z^{-1}) = z; \\
	Q((z^{-1})^*) (z) &= \{(z^{-1})^*,z,(z^{-1})^*\} = U_{(z^{-1})^*} (z^*) = (z^{-1})^*;\\
	Q(z) Q((z^{-1})^*) (x) &=  U_z \left( U_{(z^{-1})^*} (x^*) \right)^* = U_z U_{z^{-1}} (x) =  U_z U_{z}^{-1} (x) = x; \\
	Q((z^{-1})^*) Q(z) (x) &= U_{(z^{-1})^*} \left( U_z (x^*) \right)^* = U_{(z^*)^{-1}} U_{z^*} (x) = U_{z^*}^{-1} U_{z^*} (x) = x,
\end{aligned}$$ for all $x\in M$, which proves that $z$ is von Neumann regular with $z^{\dag} = (z^{-1})^* = (z^*)^{-1},$ and since $P_2(r(z)) = Q(z) Q(z^{\dag}) = Id$, the tripotent $r(z)$ is a unitary in $M$.\smallskip

In principle the continuous triple functional calculus seems to be limited to odd powers of the involved elements. Fix an element $a$ in a JB$^*$-triple $E$ and a triple isomorphism $\Psi: E_a\to C_0(\Omega_a)$ satisfying $\Psi (a) (t) =t $ ($t\in \Omega_a$). The function $g(t) = t^2$ lies in $C_0(\Omega_a)$, and hence $g_t(a) = a^{[2]}$ is an element in $E_a$. The reader should be warned that $a^{[2]}$ is only a square for a ``local binary product'' defined only on $E_a$, depending on the element $a$.\smallskip

If $a$ is von Neumann regular (i.e. $\Omega_a\subset \mathbb{R}^+$),  the function $h(t) = \frac1t$ lies in $C_0(\Omega_a)$ and $a^{[-1]} = h_t(a)\in E_a$. It is easy to see that $$\{a,a, a^{[-1]}\} = a, \{a,r(a), a^{[-1]}\} = r(a),\hbox{ and } \{r(a), a^{[-1]}, r(a)\} = a^{[-1]}.$$
By definition $a^{[2]}$ belongs to the subtriple of $E$ generated by $a$. Since the range tripotent of $a$ satisfies that $a\in E_2(r(a))$, it follows that $E_a\subseteq E_2(r(a))$ because the latter is a JB$^*$-subtriple of $E$. The elements $a, a^{[-1]},$ $a^{[2]}$ and $a^{\dag}$ belong to the unital JB$^*$-algebra $E_2 (r(a))$. By Kaup's Banach--Stone theorem (Theorem~\ref{t Kaup-Banach-Stone}) the restriction of the triple product of $E$ to $E_2(r(a))$ coincides with
$$\{x,y,z\} = (x\circ_{r(a)} y^{*_{r(a)}}) \circ_{r(a)} z + (z\circ_{r(a)} y^{*_{r(a)}})\circ_{r(a)} x -
(x\circ_{r(a)} z)\circ_{r(a)} y^{*_{r(a)}}, $$ for all $x,y,z\in E_2(r(a)).$ We consequently have $$\begin{aligned} r(a) &= \{a^{\dag}, a, r(a)\} = L(a^{\dag},a) (r(a)) = a^{\dag}\circ_{r(a)} a,\\ a &= L(a^{\dag},a) (a)  = \{a^{\dag}, a, a\} \\
	&= (a^{\dag}\circ_{r(a)} a^{*_{r(a)}}) \circ_{r(a)} a + (a\circ_{r(a)} a^{*_{r(a)}})\circ_{r(a)} a^{\dag} - (a^{\dag}\circ_{r(a)} a)\circ_{r(a)} a^{*_{r(a)}}\\
	&= (a\circ_{r(a)} a)\circ_{r(a)}  a^{\dag},\end{aligned}$$ which proves that $ a^{\dag} = a^{[-1]}$.\smallskip

Among the consequences of the basic properties of the holomorphic functional calculus for complex Jordan--Banach algebras reviewed in Theorem~\ref{t holomorphic functional calculus}, we see that for each holomorphic function $f$ on an open set $\Omega\subseteq \mathbb{C}$, the mapping $a\mapsto f(a)$ is a continuous application from the set of all elements in a complex Jordan--Banach $M$ whose spectrum in inside $\Omega$ into $M$. We insert now a result assuring a similar property for the continuous triple functional calculus (cf. the arguments in the proof of \cite[Proposition I.4.10]{Tak}).

\begin{proposition}\label{p continuity of the square}\cite[Proposition 2.1]{CuEnHirMiuPe2022} Let $\Omega\subset \mathbb{R}_0^+$ be a compact set, and let $E$ be a JB$^*$-triple. Suppose  $\mathcal{E}_{\Omega}$ denotes the set of all elements $a\in E$ whose triple spectrum is inside $\Omega$, that is $$\mathcal{E}_{\Omega} = \{a\in E  : \Omega_a\subseteq \Omega \}.$$ If $f:\Omega\to \mathbb{C}$ is a continuous function, then the mapping defined by the continuous triple functional calculus: $$\begin{aligned}
\mathcal{E}_{\Omega}& \longrightarrow E, \  a\mapsto f_t(a)
	\end{aligned}$$ is {\rm(}norm{\rm)} continuous. In particular, the mapping $a\in E\mapsto a^{[2]}\in E$ is continuous.
\end{proposition}

\begin{proof} Let $\iota: \Omega\to\mathbb{C}$ denote the natural inclusion mapping. By the sharp refinement of Stone-Weierstrass theorem in \cite{Ka83} and \cite[Lemma 4.2.8]{Cabrera-Rodriguez-vol1} the set $$\{p(\iota^2) \iota : p \hbox{ a polynomial with complex coefficients}\}$$ is dense in $C_0 (\Omega)$. Therefore, given $\varepsilon>0$ there exists a polynomial $q$ with complex coefficients such that for $q(\lambda) = p(\lambda^2) \lambda$ we have $\sup_{t\in \Omega} |q(t) - f(t) |<\varepsilon$. The norm continuity of $q$ follows from the norm continuity of the triple product of $E$ (cf. \eqref{eq triple product non-expansive}). We can therefore find $\delta > 0$ such that $\| q_t(a)- q_t(b)\| < \varepsilon$ for all $a,b\in E$ with $\|a-b\|<\delta$ and $\|a\|,$ $\|b\|\leq \max \Omega$. Taking $a,b\in \mathcal{E}_{\Omega}$ with $\|a-b\|<\delta$ we have $$\|f_t(a) -f_t(b) \|\leq \|f_t(a) -q_t(a) \| + \|q_t(a) -q_t(b) \| + \|q_t(b) -f_t(b) \|< 3\varepsilon.$$ The rest is standard.
\end{proof}

Despite the lacking of polar decompositions in the setting of JB$^*$-algebras, the previous result assures the continuity of the mapping which sends invertible elements in a unital JB$^*$-algebra to their corresponding range tripotents.

\begin{corollary}\label{c continuity of the range tripotent}\cite[Corollary 2.2]{CuEnHirMiuPe2022} Let $M$ be a unital JB$^*$-algebra. Then the mapping $$M^{-1}\longrightarrow  \mathcal{U} (M), \ a \mapsto r(a) \hbox{ is continuous.}$$
\end{corollary}

\begin{proof} By considering the JB$^*$-subtriple of $M$ generated by an element $a$, and the properties of the continuous triple functional calculus, we see that $r(a) = Q(a^{\dag}) (a^{[2]})$.\smallskip
	
The set $M^{-1}$ is open, and the mapping $a\in M^{-1}\mapsto a^{-1}\in M$ is continuous (cf.  Lemma~\ref{l first properties invertible Jordan}). Therefore, the mapping $a\in M^{-1}\mapsto a^{\dag} = (a^{-1})^*\in M^{-1}$ is continuous. Consequently, by the continuity of the triple product (see \eqref{eq triple product non-expansive}), the mapping $$a\in M^{-1} \mapsto r(a) = Q(a^{\dag}) (a^{[2]})$$ is continuous. Let us finally note that $r(a)\in \mathcal{U} (M)$ for all $a\in M^{-1}$ (cf. \eqref{eq range of invertible is unitary}).
\end{proof}

The pursued algebraic characterization of the principal component of the set of unitaries in a unital JB$^*$-algebra can be now established.

\begin{theorem}\label{t principal component as product of exp}\cite[Theorem 2.3.]{CuEnHirMiuPe2022} Let $M$ be a unital JB$^*$-algebra, let $\mathcal{U}^0(M)$ denote the principal component of the set of unitaries in $M$ and let $u$ be an element in $\mathcal{U}(M)$. Then the following statements are equivalent:
	\begin{enumerate}[$(a)$]\item $u\in M^{-1}_{\mathbf{1}}\cap\mathcal{U} (M)$;
		\item There exists a continuous path $\Gamma : [0,1]\to \mathcal{U} (M)$ with $\Gamma (0) = \11b{M}$ and $\Gamma(1) = u$;
		\item $u\in \mathcal{U}^0(M)$;
		\item $u = U_{e^{i h_n}} \cdots U_{e^{i h_1}}(\11b{M}),$ for some $n\in \mathbb{N}$, $h_1,\ldots,h_n\in M_{sa}$;
		\item There exists $w\in \mathcal{U}^0(M)$ such that $\|u-w\|<2$.
	\end{enumerate}
	Consequently,
	\begin{equation}\label{eq algebraic charact principal component unitaries JBstar}\begin{aligned} \mathcal{U}^0(M) &=  M^{-1}_{\mathbf{1}}\cap\mathcal{U} (M)\\
			&=\left\lbrace  U_{e^{i h_n}}\cdots U_{e^{i h_1}}(\11b{M}) \colon n\in \mathbb{N}, \ h_j\in M_{sa} \ \forall\ 1\leq j \leq n  \right\rbrace \\
			&= \left\lbrace  u\in \mathcal{U} (M) : \hbox{ there exists } w\in \mathcal{U}^0(M) \hbox{ with } \|u-w\|<2 \right\rbrace
	\end{aligned}\end{equation} is analytically arc-wise connected.
\end{theorem}

\begin{proof} $(a)\Rightarrow (b)$ Suppose $u\in M^{-1}_{\mathbf{1}}\cap\mathcal{U} (M)$. Since by Theorem~\ref{t Loos thm characterization of the principal component}, the set $M^{-1}_{\mathbf{1}}$ is analytically arc-wise connected, there exists a continuous path $\gamma : [0,1]\to M^{-1}_{\mathbf{1}}$ satisfying $\gamma(0) = \11b{M}$ and $\gamma(1) = u$. Corollary~\ref{c continuity of the range tripotent} assures that the mapping $\Gamma : [0,1]\to \mathcal{U} (M),$ $\Gamma (t) := r(\gamma(t))$ is continuous with $\Gamma (0) = \11b{M}$ and $\Gamma(1) = u$. \smallskip

The implications $(b)\Rightarrow (c)$ and $(c)\Rightarrow (e)$ are clear.\smallskip
	
	$(c)\Rightarrow (a)$ It follows from $(a)\Rightarrow (b)\Rightarrow (c)$ that $M^{-1}_{\mathbf{1}}\cap\mathcal{U} (M)\subseteq \mathcal{U}^0(M)$. On the other hand,  $M^{-1}_{\mathbf{1}}$ is a clopen subset of $M^{-1}$ (cf. Theorem~\ref{t Loos thm characterization of the principal component}), and thus $M^{-1}_{\mathbf{1}}\cap\mathcal{U} (M)$ is a clopen subset of $\mathcal{U} (M)$ which contains the unit and is contained in $\mathcal{U}^0(M)$. Having in mind that $\mathcal{U}^0(M)$ is a connected set, we deduce that $M^{-1}_{\mathbf{1}}\cap\mathcal{U} (M) = \mathcal{U}^0(M).$\smallskip
	
The arguments up to this point assure that statements $(a)$, $(b)$ and $(c)$ are equivalent, and moreover, $\mathcal{U}^0(M)= M^{-1}_{\mathbf{1}}\cap\mathcal{U} (M)$ is analytically arc-wise connected.\smallskip

	$(b)\Rightarrow (d)$  Suppose there exists a continuous path $\Gamma : [0,1]\to \mathcal{U} (M)$ with $\Gamma (0) = \11b{M}$ and $\Gamma(1) = u$. Let us find, by continuity and compactness, $0=t_0 < t_1 < \ldots < t_n < 1= t_{n+1}$ such that $\|\Gamma (t_i) - \Gamma (t_{i+1})\|<2$ for all $0\leq i\leq n$.  Set $u_i=\Gamma(t_i)$ for all $0\leq i \leq n+1$. By applying Lemma~\ref{l untaries at short distance in a unital JB$^*$-algebra}$(a)$ to the element $u_1$ --which satisfies $\|u_1 -\11b{M}\|<2$-- we deduce the existence of $h_1\in M_{sa}$ such that $U_{e^{i {h_{1}}}} (\11b{M}) = u_1$.  Since $2> \|u_1 - u_2\| = \|\11b{M} - U_{ e^{-i {h_{1}}}} ( u_2)\|$, a new application of Lemma~\ref{l untaries at short distance in a unital JB$^*$-algebra}$(a)$ shows the existence of $h_{2}\in M_{sa}$ satisfying $  U_{ e^{-i {h_{1}}}} ( u_2) =  U_{ e^{i {h_2}}} ( \11b{M}), $ or equivalently, $ u_2  = U_{ e^{i {h_{1}}}}  U_{ e^{i {h_2}}} ( \11b{M})$. Suppose, by an induction argument on $n$, that we have found $h_1,\ldots,h_n\in M_{sa}$ with $u_{n} = U_{e^{i h_1}} \cdots U_{e^{i h_n}}(\11b{M}).$ As before, the condition $$2>\| u_{n+1} - u_{n} \| =  \| U_{e^{-i h_n}} \cdots  U_{e^{-i h_1}} ( u_{n+1} ) - \11b{M} \|,$$ implies, via Lemma~\ref{l untaries at short distance in a unital JB$^*$-algebra}$(a)$, the existence of $h_{n+1}\in M_{sa}$ such that $$u = u_{n+1} =  U_{e^{i h_{1}}}  U_{e^{i h_2}} \cdots U_{e^{i h_{n+1}}}(\11b{M}),$$ which gives $(d)$.\smallskip
	
	The implication $(d)\Rightarrow (c)$ is also easy because given $u = U_{e^{i h_n}} \cdots U_{e^{i h_1}}(\11b{M}),$ with $n\in \mathbb{N}$, $h_1,\ldots,h_n\in M_{sa},$ the mapping $\Gamma: [0,1]\to \mathcal{U} (M) \subseteq M^{-1}$, $\Gamma (t) = U_{e^{i t h_n}} \cdots U_{e^{i t h_1}}(\11b{M})$ is an analytic curve with $\Gamma (0) = \11b{M}$ and $\Gamma (1) = u$. This also proves the second equality in \eqref{eq algebraic charact principal component unitaries JBstar}. Therefore, $(d)$ also is equivalent to the first three statements.\smallskip
	
	
Suppose finally that statement $(e)$ holds, that is, there exists $w\in \mathcal{U}^{0} (M)$ satisfying $\|w-u\|<2$. By the first two equalities in \eqref{eq algebraic charact principal component unitaries JBstar}, $w = U_{e^{i h_n}} \cdots U_{e^{i h_1}}(\11b{M})$ for some $h_1, \ldots, h_{n}$ in $M_{sa}$. A similar argument to that in the proof of $(b)\Rightarrow (d)$ implies that $u = U_{e^{i h_{n+1}}} U_{e^{i h_n}} \cdots U_{e^{i h_1}}(\11b{M})\in \mathcal{U}^{0} (M)$.
\end{proof}

In a unital C$^*$-algebra $A$ all connected components of $\mathcal{U} (A)$ are determined by the principal component. Namely, for each $u\in \mathcal{U} (A)$, the left multiplication operator $L_u: A\to A$ is a surjective linear isometry mapping the unit to $u$. Consequently, $L_u (\mathcal{U}^0(A)) = u  \mathcal{U}^0(A)$ is precisely the connected component of $\mathcal{U}(A)$ containing $u$. However, the left multiplication operator admits no direct Jordan analogue. The good behavior exhibited by unital C$^*$-algebras fails in the setting of unital JB$^*$-algebras, not only due to the lacking of a left (or right) multiplication operator, but for the possibility that two different connected components of the unitary set of a unital JB$^*$-algebra are not isometric as metric spaces (see Remark~\ref{r existence of unitaries with non isomorphic nor isometric connected components in the Jordan setting} below).


\begin{remark}\label{r connected components which are not the pricipal one} Let $w$ be a unitary element in a unital JB$^*$-algebra $M$. Suppose $\mathcal{U}^c(M)$ is a connected component of $\mathcal{U} (M)$, $w\in \mathcal{U}^c(M)$ and $M^{-1}_{w}$ is the connected component of $M^{-1}$ containing $w$. Let $M(w)$ denote the $w$-isotope admitting $w$ as unit element. Since $\mathcal{U}(M) = \mathcal{U}(M(w))$ and $M^{-1} = (M(w))^{-1}$ as sets {\rm(}cf. Lemma~\ref{r invertible in the c-isotope}{\rm)}, $\mathcal{U}^c(M)$ and $M^{-1}_{w}$ are the principal components of $\mathcal{U}(M(w))$ and $(M(w))^{-1}$, respectively. Therefore the previous Theorem~\ref{t principal component as product of exp} applied to the $w$-isotope $M(w)$ gives \begin{equation}\label{eq algebraic charact secondary component unitaries JBstar} \begin{aligned} \mathcal{U}^c (M)&=  M^{-1}_{w}\cap\mathcal{U} (M)\\
			&=\left\lbrace  U^{(w)}_{\exp_{(w)} ({i h_n})}\cdots U^{(w)}_{\exp_{(w)}({i h_1})}(w) \colon n\in \mathbb{N}, h_j\in (M(w))_{sa}, \forall\ 1\leq j \leq n  \right\rbrace \\
			&= \left\lbrace  u\in \mathcal{U} (M) : \hbox{ there exists } v\in \mathcal{U}^c(M) \hbox{ with } \|u-v\|<2 \right\rbrace
		\end{aligned}
	\end{equation} is analytically arc-wise connected. We can actually take any $\widetilde{w}\in \mathcal{U}^c(M)$ and the same conclusion holds for the connected component $\mathcal{U}^c(M)$ and the $\widetilde{w}$-isotope  $M(\widetilde{w})$.
\end{remark}

The existence of a square root for a unitary $u$ in a unital JB$^*$-algebra $M$ facilitates the arguments to understand the connected component of $\mathcal{U} (M)$ containing $u$.

\begin{remark}\label{r unitaries which are squares 2} Let $u$ be a unitary element in a unital JB$^*$-algebra $M$. Suppose that $u$ satisfies that $u = v^2$ for some unitary $v\in M$. Let $\mathcal{U}(M)^{(c)}$ denote the connected component of $\mathcal{U}(M)$ containing $u$. Since the mapping $U_v: M\to M(u)=M_{(u^*)}$ is a unital surjective linear isometry --and hence a Jordan $^*$-isomorphism between these two algebras-- $U_{v} (\mathcal{U}^{0}(M))$ is a connected component of $\mathcal{U} (M) =\mathcal{U} (M(u))$ and contains $u$. It then follows from Theorem~\ref{t principal component as product of exp} that $$ \mathcal{U}^c (M) = U_{v} (\mathcal{U}^{0}(M)) = \left\lbrace U_{v} U_{e^{i h_n}}\cdots U_{e^{i h_1}}(\11b{M}) \colon n\in \mathbb{N}, \ h_j\in M_{sa} \ \forall\ 1\leq j \leq n  \right\rbrace.$$
	
\end{remark}

Previous studies have searched for sufficient conditions to assure  that a unitary element $u$ in a JB$^*$-algebra $M$ admits a square root (see, for example, \cite[Problem 5.1 and Lemma 5.2]{BraKaUp78}, \cite{HatMi2000}, \cite{ChiKasKawValov2005} and \cite{MatNagYam2008}.\smallskip

A consequences of the previous Theorem~\ref{t principal component as product of exp} allows us to describe the principal component of $\mathcal{U} (M)$ in terms of quadratic subsets.

\begin{proposition}\label{p self adjoint quadratic subset princ comp}\cite[Proposition 2.6]{CuEnHirMiuPe2022} Let $M$ be a unital JB$^*$-algebra. Then $\mathcal{U}^{0} (M)$ is a self-adjoint quadratic subset of $\mathcal{U} (M)$, that is, for all $u,w\in \mathcal{U}^0 (M)$ the elements $w^*$ and $U_w (u)$ belong to $\mathcal{U}^0 (M)$. Furthermore, $\mathcal{U}^{0} (M)$ is the smallest quadratic subset of $\mathcal{U} (M)$ containing the set $e^{i M_{sa}}$. In particular, $$U_{\mathcal{U}^{0} (M)} \left( \mathcal{U}^{0} (M) \right)\! = \mathcal{U}^{0} (M) =\left\lbrace  U_{e^{i h_n}}\cdots U_{e^{i h_1}}(e^{i h_0}) \colon n\in \mathbb{N}, \ h_j\in M_{sa} \ \forall\ 0\leq j \leq n  \right\rbrace.$$
\end{proposition}

\begin{proof} Since $u,w\in \mathcal{U}^0(M),$ Theorem~\ref{t principal component as product of exp} guarantees the existence of hermitian elements $h_m,\dots,h_1, k_n, \ldots, k_1\in M_{sa}$ such that $w=U_{e^{i h_m}}\cdots U_{e^{i h_1}}(\11b{M})$ and $u=U_{e^{i k_n}}\cdots U_{e^{i k_1}}(\11b{M})$ with $n,m\in \mathbb{N}$. By applying that the involution is a conjugate linear Jordan $^*$-isomorphism we get 
	$$ w^*= (U_{e^{i h_m}}\cdots U_{e^{i h_1}}(\11b{M}))^* = U_{e^{-i h_m}}\cdots U_{e^{-i h_1}}(\11b{M}),$$  which, by virtue of Theorem~\ref{t principal component as product of exp}, proves that $w^*\in \mathcal{U}^0(M)$.\smallskip
	
On the other hand, by the fundamental identity \eqref{eq fundamental identity n elements} we arrive to 
	$$\begin{aligned}
		U_w (u)&=U_w\left( U_{e^{i k_n}}\cdots U_{e^{i k_1}}(\11b{M}) \right) = U_{U_{e^{i h_m}}\cdots U_{e^{i h_1}}(\11b{M})}\left( U_{e^{i k_n}}\cdots U_{e^{i k_1}}(\11b{M}) \right) \\
		&= U_{e^{i h_m}}\cdots U_{e^{i h_1}} U_{\11b{M}} U_{e^{i h_1}}\cdots U_{e^{i h_m}}\left( U_{e^{i k_n}}\cdots U_{e^{i k_1}}(\11b{M}) \right)\\
		&=U_{e^{i h_m}}\cdots U_{e^{i h_1}} U_{e^{i h_1}}\cdots U_{e^{i h_m}} U_{e^{i k_n}}\cdots U_{e^{i k_1}}( \11b{M})\in \mathcal{U}^0(M),
	\end{aligned}$$ by Theorem~\ref{t principal component as product of exp}. We have therefore proved that $\mathcal{U}^0(M)$ is a self-adjoint quadratic subset of $M$. \smallskip

	
Suppose now that $\mathcal{V}$ is a quadratic subset of $\mathcal{U} (M)$ containing the set $e^{i M_{sa}}$. Clearly for each $h_1,\ldots, h_n$ in $M_{sa}$ an induction argument shows that $$U_{e^{i h_n}} \ldots U_{e^{i h_1}} (\11b{M}) \in \mathcal{V},$$ and thus, by Theorem~\ref{t principal component as product of exp}, $\mathcal{U}^0(M)\subseteq \mathcal{V}$. The rest is clear.
\end{proof}

The next remark leads to a natural open problem.

\begin{remark}\label{r connected component building from u as an pen question} Let $M$ be a unital JB$^*$-algebra and let $u$ be an element in $\mathcal{U} (M)\backslash \mathcal{U}^0(M)$. Let $\mathcal{U} (M)^{(u)}$ denote the connected component of $\mathcal{U} (M)$ containing $u$. The set $$\mathcal{E}_u := \left\{ U_{e^{i h_1}}\ldots U_{e^{i h_n}} (u) : \begin{array}{c}
		n\in \mathbb{N}, \ h_j\in M_{sa} \\
		\hbox{ for all } 1\leq j \leq n
	\end{array}
	\right\}$$ is arc-wise connected and contains $u$. Thus, $\mathcal{E}_u\subseteq \mathcal{U} (M)^{(u)}$. We do not know if the reverse inclusion holds, in general. If $M= A$ is a unital C$^*$-algebra both inclusions hold, a fact which is essentially due to the existence of left and right multiplicaiton operators. Namely, since $\mathcal{U}^0(A) = \{ e^{i h_1} \cdot \ldots \cdot e^{i h_n} : n\in \mathbb{N}, \ h_1, \ldots, h_n\in A_{sa} \}$ {\rm(}cf. \eqref{eq algebraic descirption of the principal component for unital Cstar algebras}{\rm)}, and the right multiplication operator $R_u: A\to A$, $x\mapsto x u$ is a surjective linear isometry mapping the unit to $u$, we deduce that $$\left\{ e^{i h_1} \cdot \ldots \cdot e^{i h_n} u : n\in \mathbb{N}, \ h_1, \ldots, h_n\in A_{sa} \right\}= R_u (\mathcal{U}^0(A)) = \mathcal{U}(A)^{(u)}.$$
	
	Similarly, the products of the form $w e^{i h_n} \cdot \ldots \cdot e^{i h_1}$ lie in $\mathcal{U} (A)^{(u)}$ for all $w\in \mathcal{U}(A)^{(u)}$, $h_1, \ldots, h_n\in A_{sa}$. So, given $h_1, \ldots, h_n\in A_{sa}$, the mapping $R_{e^{i h_1} \cdot \ldots \cdot e^{i h_n}}$ is a surjective linear isometry on $A$ mapping $\mathcal{U}(A)^{(u)}$ onto itself, and hence $$\begin{aligned}
		\mathcal{U}(A)^{(u)}  &= \left\{ e^{i h_1} \cdot \ldots \cdot e^{i h_n} u e^{i h_n} \cdot \ldots \cdot e^{i h_1} : n\in \mathbb{N}, \ h_1, \ldots, h_n\in A_{sa} \right\} \\
		& = \left\{ U_{e^{i h_1}}\ldots U_{e^{i h_n}} (u) : \begin{array}{c}
			n\in \mathbb{N}, \ h_j\in A_{sa} \\
			\hbox{ for all } 1\leq j \leq n
		\end{array}
		\right\}.
	\end{aligned} $$\smallskip
	
Fortunately, the inclusion $\mathcal{E}_u\subseteq \mathcal{U}^c (M)$ is enough to prove the following rule concerning the $U$-operators \begin{equation}\label{eq inner main other main} U_{\mathcal{U}^0 (M)} \left( \mathcal{U} (M)^{(u)}\right) = \mathcal{U} (M)^{(u)} 
	\end{equation} for any connected component  $\mathcal{U} (M)^{(u)}\subseteq \mathcal{U} (M)$. Indeed, by Theorem~\ref{t principal component as product of exp}, every $w\in \mathcal{U}^{0}(M)$ is of the form $w = U_{e^{i h_n}}\cdots U_{e^{i h_1}}(\11b{M})$ with $n\in \mathbb{N}$, $h_j\in M_{sa}$. Therefore, given $\tilde{u}\in \mathcal{U} (M)^{(u)}$, it follows from \eqref{eq fundamental identity n elements} that $$ U_w (\tilde{u}) = U_{U_{e^{i h_n}}\cdots U_{e^{i h_1}}(\11b{M})} (\tilde{u}) = U_{e^{i h_n}} \ldots U_{e^{i h_1}} U_{e^{i h_1}} \ldots U_{e^{i h_n}} (\tilde{u})\in \mathcal{E}_{\tilde{u}}\subseteq \mathcal{U} (M)^{(\tilde{u})} = \mathcal{U} (M)^{(u)}.$$
	
\noindent We finally observe that the equality $\displaystyle \bigcup_{w\in \mathcal{U}(M)^{(u)}} \mathcal{E}_w = \mathcal{U}(M)^{(u)}$ follows trivially from the above arguments.
\end{remark}

The subsection about properties of the unitary set in a unital JB$^*$-algebras is not complete without a short discussion on one-parameter unitary quadratic maps.\smallskip
 
In the setting of operators, a  \emph{one-parameter group} of bounded linear operators on a Banach space $Z$ is a mapping $\mathbb{R} \to B(Z),$ $t\mapsto E(t)$ satisfying $$E(0) =I, \hbox{ and } E({t+s}) = E({s}) E({t}), \hbox{ for all } s,t\in \mathbb{R}.$$ There are two notions of continuity. In the first one, we say that a one-parameter group $\{ E(t) : t\in \mathbb{R}\}$ is uniformly continuous at the origin if $\displaystyle \lim_{t\to 0} \| E(t) -I\| =0$. This is the strongest assumption, and it is further known that this is equivalent to the existence of a bounded linear operator $R\in B(Z)$ such that $E(t) =e^{t R}$ for all $t\in \mathbb{R}$, where the exponential is computed in the Banach algebra $B(Z)$ (cf., for example, \cite[Proposition 3.1.1]{BratRob1987}). The one-parameter group $\{E(t) : t\in \mathbb{R}\}$ on $Z$ is called \emph{strongly continuous} if for each $z$ in $Z$ the mapping $t\mapsto  E({t}) (z)$ is (norm) continuous (\cite[Definition 5.3, Chapter X]{Conway}). A \emph{one-parameter unitary group} on a complex Hilbert space $H$ is a one-parameter group on $H$ such that $E(t)$ is a unitary element for each $t\in \mathbb{R}$. \smallskip

The celebrated Stone's one-parameter theorem affirms that for each strongly continuous one-parameter unitary group $\{E(t) : t\in \mathbb{R} \}$ on a complex Hilbert space $H$ there exists a self-adjoint operator $h\in B(H)$ such that $E(t)=e^{i t h}$, for every $t\in \mathbb{R}$ (\cite[5.6, Chapter X]{Conway}).\smallskip

Let us recall the existence of exceptional JB$^*$-algebras which cannot be represented inside $B(H)$ makes hopeless a direct application of Stone's one-parameter theorem. Uniformly continuous one-parameter groups of surjective isometries (i.e. triple isomorphisms compare Kaup's Banach--Stone theorem~\ref{t Kaup-Banach-Stone}), Jordan $^*$-isomorphisms and orthogonality preserving operators on JB$^*$-algebras have been recently recently described in terms of triple derivations in \cite{GarPe21Cstar,GarPe21JBstar}. \smallskip

It should be recalled that a \emph{triple derivation} on a JB$^*$-triple $E$ is a linear mapping $\delta: E\to E$ satisfying a ternary version of Leibniz' rule $$\delta \{a,b,c\} = \{\delta(a),b,c\} + \{a, \delta(b), c\} + \{a,b,\delta(c)\}, \ \ (a,b,c\in E).$$ Triple derivations on JB$^*$-triples are automatically continuous (cf. \cite[Corollary 2.2]{BarFri90} or \cite{HoMarPeRu,PeRu2014}). Another interesting property affirms that if $\delta: M\to M$ is a triple derivation on a unital JB$^*$-algebra, we have $\delta(\11bM)^* = - \delta(\11bM)$ (cf. \cite[Proof of Lemma 1]{HoMarPeRu}).\smallskip

The natural limitations of the set of unitary elements in a unital JB$^*$-algebra  (e.g. the Jordan product of two unitary elements is not, in general, a unitary) produces what seems, a priori, an obstacle for a literal transcription of Stone's one-parameter theorem to the setting of unital JB$^*$-algebras. The next appropriate Jordan version of Stone's one-parameter theorem has been taken from \cite{CuPe2022}.

\begin{theorem}\label{t Jordan unitary groups version of Stone's theorem}\cite[Theorem 3.1]{CuPe2022} Let $M$ be a unital JB$^*$-algebra. Suppose that $\{u(t):t\in \mathbb{R}\}$ is a family in $\mathcal{U} (M)$ such that $u(0) =\11b{M}$ and the correspondence $t\mapsto u(t)$ is (norm) continuous. Then the following statements are equivalent:
	\begin{enumerate}[$(a)$]\item There exists $h\in M_{sa}$ such that $u(t) = e^{i t h}$ for all $t\in \mathbb{R}$;
		\item $u(s+t) = u(s) \circ u(t)$,  for all $t,s\in \mathbb{R}$;
		\item $U_{u(t)} (u(s)) = u(2 t +s),$ for all $t,s\in \mathbb{R}$.
	\end{enumerate}	
\end{theorem}

\begin{proof} Since the JB$^*$-subalgebra generated by the element $h$ is a commutative (and associative) Banach algebra, and the exponential is taken in this subalgebra, the implication  $(a) \Rightarrow (b)$ is clear since $u(s+t) = e^{i(s+t) h} = e^{i s h} e^{i t h} = u(s) \circ u(t)$. The implications $(b) \Rightarrow (c)$ is also clear since, by $(b)$ we have $$\begin{aligned}
U_{u(t)} (u(s)) &= 2 (u(t)\circ u(s))\circ u(t) - u(t)^2 \circ u(s) \\ 
&= 2 u(t+ s)\circ u(t) - u(2 t)\circ u(s)= u(2 t +s).
	\end{aligned}$$  
	
$(c)\Rightarrow (a)$. In order to apply the results on one-parameter group of surjective linear isometries in \cite{GarPe21Cstar,GarPe21JBstar}, we consider the family of surjective linear isometries on $M$ defined by $\{U_{u(t)} : t\in \mathbb{R}\}$. Clearly, the mapping $t\mapsto U_{u(t)}$ is (norm) continuous with $U_{u(0)} = Id_{M}$ by hypotheses. The fundamental identity \eqref{eq fundamental identity UaUbUa} asserts that $$ U_{u(t)} U_{u(s)} U_{u(t)} = U_{U_{u(t)}(u(s))} = U_{u(2t+s)} ,$$ for all $s,t\in \mathbb{R}$. It then follows that 
$$U_{u(t)}^2 = U_{u(t)} U_{u(0)} U_{u(t)} = U_{u(2 t)},\ U_{u(t)}^3 = U_{u(t)} U_{u(t)} U_{u(t)} = U_{u(3 t)},$$ and if we assume that $U_{u(k t)} = U_{u(t)}^k$ for all $k<n$ and $n\geq 3$ we get $$U_{u(n t)}= U_{u(2 t + (n-2)t)}= U_{u(t)} U_{u((n-2)t)} U_{u(t)}= U_{u(t)} U_{u(t)}^{n-2} U_{u(t)} = U_{u(t)}^n.$$ An induction arguments proves that $$U_{u(n t)} = U_{u(t)}^n \hbox{ for all } t\in \mathbb{R},\  n\in \mathbb{N}.$$

Furthermore, for each real $t$, $U_{u(t)} U_{u(-t)} U_{u(t)} = U_{u(t)}$ and $U_{u(t)}\in \hbox{Iso} (M),$ and hence $U_{u(-t)} = U_{u(t)}^{-1}$ for all $t\in \mathbb{R}$. It follows that, for a negative integer $n$ and each real $t$, we have $$U_{u(n t)} = U_{u((-n) (-t))}= U_{u(-t)}^{-n} = U_{u(t)}^n,$$ an identity which then holds for all $n\in \mathbb{Z}$.\smallskip
	
We shall prove that the family $\{U_{u(t)} : t\in \mathbb{R}\}$ is a one-parameter group of surjective linear isometries on $M$. To this end, let us fix two rational numbers $\frac{n}{m}, \frac{n'}{m'}$ with $m,m'\in \mathbb{N}$, $n,n'\in \mathbb{Z}$. The properties proved above imply that $$ U_{u\left(\frac{n}{m} t + \frac{n'}{m'} t\right)} = U_{u\left(\frac{n m' + n' m}{m m'} t\right)} =
U_{u\left(\frac{1}{m m'} t\right)}^{n m' + n' m} = 
U_{u\left(\frac{1}{m m'} t\right)}^{n m'} U_{u\left(\frac{1}{m m'} t\right)}^{ n' m} = 
U_{u\left(\frac{n}{m} t\right)} U_{u\left(\frac{n'}{m'} t\right)}, $$ and thus it follows by continuity that $$U_{u(t+s)} = U_{u(t)} U_{u(s)}, \hbox{ for all } t,s\in \mathbb{R},$$ that is, the family $\{U_{u(t)} : t\in \mathbb{R}\}$ is a uniformly continuous one-parameter group of surjective isometries on $M$. By \cite[Lemma 3.2]{GarPe21JBstar} there exists a triple derivation $\delta: M\to M$ satisfying \begin{equation}\label{eq power series for Phi} \hbox{$\displaystyle U_{u(t)} = e^{t \delta}= \sum_{n=0}^{\infty} \frac{t^n}{n!} \delta^n$ for all $t\in \mathbb{R}$,}
	\end{equation} where the exponential is computed in the Banach algebra $B(M)$.\smallskip

We continue with the next observation: \begin{equation}\label{eq power series for ut} e^{\frac{t}{2} \delta} (\11b{M}) =U_{u\left(\frac{t}{2}\right)} (\11b{M})  = u\left(\frac{t}{2}\right)^2 = U_{u\left(\frac{t}{2}\right)} (u(0))= u( t),
\end{equation} for all $t\in \mathbb{R}$, and
$u(t)^2= U_{u(t)} (\11b{M}) =U_{u(t)} (u(0)) =  u(2 t).$ It follows that 
\begin{equation}\label{eq power series}
\left(\sum_{k=0}^{\infty} \frac{t^n}{2^n n!} \delta^n (\11b{M}) \right)^2= \left(e^{\frac{t}{2} \delta} (\11b{M}) \right)^2 = \left(u(t) \right)^2 = u(2 t ) = e^{t \delta} (\11b{M}) = \sum_{k=0}^{\infty} \frac{t^n}{n!} \delta^n (\11b{M}).
\end{equation} By comparing both extremes in the previous identity we arrive to $$ \frac{1}{2^2} \delta (\11b{M})^2 + 2  \frac{1}{2^3} \delta^2 (\11b{M}) = \frac{1}{2} \delta^2 (\11b{M}),$$ which assures that $\delta (\11b{M})^2 = \delta^2 (\11b{M})$. If we combine this identity with that in \eqref{eq power series} it follows that 
$$\frac{1}{2^2 3!} \delta^3 (\11b{M}) +   \frac{1}{2^3}  \delta (\11b{M})^3=  2  \frac{1}{2^3 3!} \delta^3 (\11b{M}) +  2 \frac{1}{2^3} \frac{1}{2} \delta^2 (\11b{M})  \delta (\11b{M})=  \frac{1}{3!} \delta^3 (\11b{M}),$$ and hence $\delta (\11b{M})^3 = \delta^3 (\11b{M})$. An easy induction argument leads to  \begin{equation}\label{eq delta1n equals deltan1} \delta(\11b{M})^n = \delta^n(\11b{M}), \hbox{ for all natural } n.
\end{equation} 

It follows from \eqref{eq power series for Phi}, \eqref{eq power series for ut}, and \eqref{eq delta1n equals deltan1} that $$u(t) = e^{\frac{t}{2} \delta} (\11b{M}) = \sum_{n=0}^{\infty} \frac{t^n}{2^n n!} \delta^n (\11b{M}) = \sum_{n=0}^{\infty} \frac{t^n}{ 2^n n!} \delta(\11b{M})^n = e^{t \frac{\delta(\11b{M})}{2} },$$ for all $t\in \mathbb{R}$.  As we commented above $\delta(\11b{M})^* = -\delta(\11b{M})$ in $M$  (cf. \cite[Proof of Lemma 1]{HoMarPeRu}), and hence $u(t) = e^{i t h}$ with $h = \frac{\delta(\11b{M})}{2 i}\in M_{sa}$.
\end{proof}

\subsection{Surjective isometries between connected components of unitary sets}\label{subsect: isometries between connected components}\ \smallskip

In the previous subsection~\ref{subsec: connected components of the unitay set} we have studied the properties of the unitary set of a unital JB$^*$-algebra and its principal connected component. We have seen that the existence of multiple connected components can be applied to find surjective isometries between unitary sets which cannot be extended to surjective real-linear isometries between the algebras (cf. Example~\ref{example existence of non-extendible surjective isometries between unitaries}). In this subsection we survey the results studying the extensibility of a surjective isometry between two connected components of the unitary sets of two unital JB$^*$-algebras. \smallskip

We would like to apply the geometric tools in section~\ref{sec: geometric tool of Hatori Hirasawa Miura Molnar}. For this purpose we need to check appropriate Jordan versions of  properties $B(a,b)$, $C_1(a,b)$ and $C_2(a,b)$. We should observe that the existence of exceptional JB$^*$-algebras which cannot be embedded inside $B(H)$ makes impossible to apply Theorem~\ref{t Theorem 6 Haotri Molnar 2012} as we did in Remark~\ref{remark unitaries satisfy B and C1}. \smallskip

Our first result gathers ideas from \cite{CuPe2022} and \cite{CuEnHirMiuPe2022}.

\begin{theorem}\label{t Delta preserves inverted triple products in U0M}{\rm(\cite[Theorem 2.9]{CuPe2022} and \cite[Theorem 3.1]{CuEnHirMiuPe2022})} Let $\Delta: \mathcal{U}^0 (M)\to \mathcal{U}^0 (N)$ be a surjective isometry, where $M$ and $N$ are unital JB$^*$-algebras. Suppose $u,v\in \mathcal{U}^0 (M)$ with $\|u-v\|<\frac12$. Then the following statements hold:
	\begin{enumerate}[$(1)$]\item For all $x,y\in \mathcal{U}^0 (M)$ we have $$\left\|U_{v}(x^{-1})- U_{v} (y^{-1})\right\|=\left\|U_{v}(x^{*})- U_{v} (y^{*})\right\|= \|x^*-y^*\| = \|x-y\|.$$
		\item The constant $K= 2 - 2 \|u-v\|>1$ satisfies that $$ \left\| U_v ( w^*)- w \right\|=\left\| U_v ( w^{-1})- w \right\|\geq K \| w- v\|,$$ for all $w$ in the set $$L^0_{u,v} =\{w\in \mathcal{U}^0 (M)  : \|u-w\| = \| U_{v} (u^{-1})- w\|= \| U_{v} (u^*)- w\| = \|u-v\|\}.$$
		\item $\mathcal{U}^0 (N)$ satisfies a Jordan version of $C_2(\Delta(u),\Delta(U_u (v^{*})))$, that is,  there exists $w$ in $\mathcal{U}^0 (N)$ satisfying \begin{equation}U_w (\Delta(u)^*) = \Delta\left( U_v (u^*) \right), \ \|w- \Delta(v) \| <2,\nonumber
		\end{equation}
		and
		$$\| U_w (x^{*}) -U_w (y^{*})\|  = \|x-y\|, \quad \forall x,y\in \mathcal{U}^0 (N).$$
		\item The equalities $\Delta(U_{v} (u^{*})) =\Delta(U_{v} (u^{-1})) =U_{\Delta(v)} (\Delta(u)^{-1}) =U_{\Delta(v)} (\Delta(u)^*)$ hold.
	\end{enumerate}
\end{theorem}

\begin{proof} $(1)$ The desired conclusion follows from the fact that $U_v$ and the involution are isometries on $M$ (cf. \cite[Theorem 4.2.28$(vii)$]{Cabrera-Rodriguez-vol1}).\smallskip

$(2)$ Let us consider the $u$-isotope JB$^*$-algebra $M(u)= (M,\circ_u,*_u)$ of $M$. The $U$-operator on $M(u)$ will be denoted by $U^u$. We fix an element $w\in L_{u,v}$. Since $\|u-w\| = \|u-v\|< \frac12,$ we deduce from Lemma~\ref{l untaries at short distance in a unital JB$^*$-algebra} the existence of two self-adjoint elements $h_1,h_2\in M(u)$ such that $v= e^{ih_1}$ and $w= e^{ih_2}$. Let $\mathcal{B}$ denote the JB$^*$-subalgebra of $M(u)$ generated by $u,h_1,h_2$ --$u$ is the unit here--. The Shirshov-Cohn theorem (see Theorem~\ref{t Shirshov-Cohn}) assures the existence of a complex Hilbert space $H$ such that $\mathcal{B}$ embeds as a JB$^*$-subalgebra of $B(H)$ and both share the same unit $u$ (the product of $B(H)$ will be denoted by mere juxtaposition and the involution by $\sharp$). Obviously $u,v,w\in \mathcal{U} (\mathcal{B})\subseteq \mathcal{U} (B(H))\subseteq \hbox{Iso}(H)$ with $\|u-w\|_{\mathcal{B}} = \|u-w\|_{M}= \|u-v\|_{M}= \|u-v\|_{\mathcal{B}}$. Let us compute $$\begin{aligned} U_{v}^u (u) &= 2 (v\circ_u u)\circ_u v - (v\circ_u v)\circ_u u 
		\\
		&= 2 v\circ_u v - v\circ_u v = v\circ_u v = \{v,u,v\} = U_v(u^*),
	\end{aligned}$$
	and
	$$\begin{aligned} \| v u^{-1} v- w\|_{_{B(H)}}&= \| v u^{-1} v- w\|_{\mathcal{B}} = \| U_{v}^u (u^{\sharp}) - w\|_{\mathcal{B}} = \| U_{v}^u (u^{*_u}) - w\|_{\mathcal{B}} \\
		&= \| U_{v}^u (u) - w\|_{\mathcal{B}} = \| U_{v} (u^*)- w\|_{\mathcal{B}}= \| U_{v} (u^*)- w\|_{M}.
	\end{aligned} $$
	
Remark~\ref{remark unitaries satisfy B and C1} proves that for $K= 2 - 2 \|u-v\|>1$ we have $$\left\| U^u_v ( w^{*_u})- w \right\|_{\mathcal{B}} = \left\| U^u_v ( w^{-1})- w \right\|_{\mathcal{B}} = \| v w^{-1} v - w \|_{_{B(H)}} \geq K \|w-v\|_{\mathcal{B}} = K \|w-v\|.$$
	On the other hand, by the uniqueness of the triple product (see Theorem~\ref{t Kaup-Banach-Stone}) we have
	$$ U^u_v ( w^{*_u}) = \{v,w,v\}_{\mathcal{B}} =  \{v,w,v\}_{M(u)} =  \{v,w,v\}_{M} = U_v(w^*).$$ All together gives $$\begin{aligned}\left\| U_v ( w^{-1})- w  \right\|&=\left\|U_v ( w^*)- w \right\|_{M} =  \left\| U_v ( w^*)- w \right\|_{\mathcal{B}} \\
		&= \left\| U^u_v ( w^{*_u})- w \right\|_{\mathcal{B}}\geq K \| w- v\|,
	\end{aligned}$$ which completes the proof of $(2)$.\smallskip
		
$(3)$ Let us consider the $u$-isotope $M(u)$, and observe that $v\in\mathcal{U} (M(u))$ (cf. Lemma~\ref{l unitaries in the homotope coincide}).  Let $\mathcal{B}$ denote the JB$^*$-subalgebra of $M(u)$ generated by $v$ and the unit of $M(u)$--$u$ in this case--. The Shirshov-Cohn theorem (see Theorem~\ref{t Shirshov-Cohn}) assures that $\mathcal{B}$ is a JB$^*$-subalgebra of $B(H)$, for some complex Hilbert space $H$, and we can further assume that $u$ is the unit in $B(H)$ (cf.  \cite[Theorem 7.2.5]{HOS} and \cite[Corollary 2.2]{Wri77}). Moreover, $u$ can be identified with the unit of $B(H)$. We shall use juxtaposition for the product in $B(H)$, and the symbol $\sharp$ for its involution.\smallskip
	
	It is worth noting that when $M$ and $\mathcal{B}$ are regarded as JB$^*$-triples, the triple product in $\mathcal{B}$ is precisely  $\{{\cdot},{\cdot},{\cdot}\}_M|_{\mathcal{B}\times \mathcal{B}\times \mathcal{B}}$. At the same time, the triple product in $\mathcal{B}$ can be expressed in terms of the associative product of $B(H)$. By Proposition~\ref{p self adjoint quadratic subset princ comp} the element $U_v (u^*)$ lies in $\mathcal{U}^0(M)$, and hence we can apply $\Delta$ at this element and compute the next distances:
	$$\begin{aligned} \| \Delta(u) - \Delta(U_v (u^*))\|_N &=  \| u - U_v (u^*)\|_M =   \| u - \{ v, u,v\}_M\|_M = \| u - \{ v, u,v\}_{\mathcal{B}}\|_{\mathcal{B}}  \\
		&= \| u - \{ v, u,v\}_{_{B(H)}}\|_{_{B(H)}} =  \| u - v u^{\sharp} v\|_{_{B(H)}} \\
		&=  \| u - v v\|_{_{B(H)}} =  \| v^{\sharp} - v\|_{_{B(H)}} \leq   \| v^{\sharp} - u \|_{_{B(H)}} + \| u- v\|_{_{B(H)}} \\
		&= 2 \| u - v\|_{_{B(H)}}= 2 \| u - v\|_{\mathcal{B}} =  2 \| u - v\|_M < 1.
	\end{aligned}.$$ 	
	
	
\noindent We are in a position to apply Lemma~\ref{l untaries at short distance in a unital JB$^*$-algebra}$(b)$ to deduce the existence of a unitary $w\in \mathcal{U} (N)$ satisfying $$U_w ( \Delta(u)^* ) = \Delta\left( U_v (u^*) \right).$$ However, by hypotheses $\| \Delta(u) - \Delta(v)\| =\|u-v\|<\frac12$, and consequently, by the final statement in Lemma~\ref{l untaries at short distance in a unital JB$^*$-algebra} we may assume that $$\|w- \Delta(v) \|\leq \|w- \Delta(u) \|+\|\Delta(u)-\Delta(v) \|\leq \| \Delta(u) - \Delta(U_v (u^*))\| +\frac12  < 1 +\frac12 <2.$$ Since $\Delta(v)\in \mathcal{U}^0 (N)$, the characterization of the principal component given in Theorem~\ref{t principal component as product of exp} implies that $w\in \mathcal{U}^0 (N)$.\smallskip
	
Now, it follows from the fact that the mapping $U_w$ is an isometry that 
	$$\| U_w (x^{*}) -U_w (y^{*})\|=\| U_w (x^{-1}) -U_w (y^{-1})\| = \|x^{-1}- y^{-1}\| = \|x^* - y^* \| = \|x-y\|,$$ for all $x,y\in \mathcal{U}^0 (N)\subseteq \mathcal{U}(N)$. We can hence conclude that the Jordan version of $C_2 (\Delta(u),\Delta(U_u (v^{*})))$ holds for $\mathcal{U}^0 (N)$.\smallskip
	
	
$(4)$ Let $w\in \mathcal{U}^0 (N)$ be the element given by $(3)$. We consider the mappings defined by
	$$\varphi (x) := U_{v} (x^{-1}) = U_{v} (x^{*}), \ \ x\in \mathcal{U}^0(M),$$ and
	$$\psi (y) := U_{w} (y^{-1}) = U_{w} (y^{*}), \ \ y\in \mathcal{U}^0(N).$$ Proposition~\ref{p self adjoint quadratic subset princ comp} assures that $\varphi \left(\mathcal{U}^0 (M)\right)= \mathcal{U}^0 (M)$ and $\psi\left( \mathcal{U}^0 (N)\right) = \mathcal{U}^0 (N)$.\smallskip
	
We aim to apply the main geometric tool presented in Theorem~\ref{t main geometric tool} to $\varphi, \psi$ and $\Delta$ at $a= u$ and $c = v$. For this purpose let us observe that, since $v$ and $w$ are unitaries lying in the principal components of $\mathcal{U}(M)$ and $\mathcal{U} (N)$, respectively, $\varphi$ and $\psi$ are distance preserving bijections because the $U$-maps the involutions are surjective isometries.\smallskip
	
It is easy to check from the definitions that  $$\varphi\circ \varphi(x)=U_{v} ((U_{v} (x^{-1}))^{-1})=U_{v} ((U_{v} (x^*))^*)= U_{v} (U_{v^*} (x))=U_{v} (U_{v^{-1}} (x))= x,$$ for every $x\in \mathcal{U}^0(M)$. That is, $\varphi\circ \varphi$ is the identity mapping on $\mathcal{U}^0 (M)$. By applying that $U_w ( \Delta(u)^* ) = \Delta\left( U_v (u^*) \right)$ we have $$\begin{aligned}
		\psi (\Delta(u)) &= U_w (\Delta(u)^*) = \Delta\left( U_v (u^*) \right) =\Delta\left( \varphi (u) \right), \hbox{ and }\\
		\psi \left(\Delta\left( \varphi (u) \right)\right) &= U_{w} \left((\Delta \left( U_{v} (u^*)\right))^* \right)= \left( U_{w^*} \left(\Delta \left( U_{v} (u^*)\right)\right)\right) ^* = \Delta(u)^{**} = \Delta (u).
	\end{aligned}$$
	
Therefore the mappings  $\varphi, \psi,$ $\Delta,$ $a= u$ and $c=v$ satisfy all hypotheses in Theorem~\ref{t main geometric tool}, and hence the mentioned result implies that $U_{w} (\Delta(v)^*) = \Delta(v).$ We also know from $(3)$ that  $\|w -\Delta(v)\|<2,$  and hence Lemma~\ref{l untaries at short distance in a unital JB$^*$-algebra}$(d)$ gives $w=\Delta(v).$ Having in mind that $U_w ( \Delta(u)^* ) = \Delta\left( U_v (u^*) \right)$ we get $\Delta(U_{v} (u^{*}))   = U_{\Delta(v)} (\Delta(u)^*),$ as desired.
\end{proof}

Before presenting the result describing surjective isometries between the principal components of the unitaries in two unital JB$^*$-algebras, we state a couple of technical tools.\smallskip

The statement of the next result has been borrowed from \cite[Lemma 3.2]{CuEnHirMiuPe2022} while the proof is coming from \cite[Lemma 3.3]{CuPe2022}.

\begin{lemma}\label{l HM Lemma 7 Jordan}\cite[Lemma 3.3]{CuPe2022} Let $M$ and $N$ be two unital JB$^*$-algebras. Let $\{u_{k}: 0\leq k\leq 2^n\}$ be a subset of $\mathcal{U}^0(M)$ {\rm(}with $n \in \mathbb{N}${\rm)} and let $\Phi: \mathcal{U}^0(M)\to \mathcal{U}^0(N)$ be a mapping such that $$U_{u_{k+1}} (u_k^*) = u_{k+2}, \hbox{ and } \Phi (U_{u_{k+1}} (u_k^*)) = U_{\Phi(u_{k+1})} (\Phi (u_k)^*),$$ for all $0\leq k\leq 2^n -2$. Then $U_{u_{2^{n-1}}} (u_0^*) = u_{2^n}$ and $$\Phi\left(U_{u_{2^{n-1}}} (u_0^*) \right) = U_{\Phi(u_{2^{n-1}})} \left(\Phi(u_0)^*\right). $$
\end{lemma} 

\begin{proof} We shall argue by induction on $n$. The statement is clear for $n=1$. Suppose that our statement is true for every family with $2^n+1$ elements satisfying the conditions above. Let $\{w_{k}: 0\leq k\leq 2^{n+1}\}$ be a subset of $\mathcal{U}^{0}(M)$ with $2^{n+1} +1$ elements such that $U_{w_{k+1}} (w_k^*) = w_{k+2}$ and $$\Phi (U_{w_{k+1}} (w_k^*)) = U_{\Phi(w_{k+1})} (\Phi (w_k)^*),$$ for all $0\leq k\leq 2^{n+1} -2$. We consider the family $\{u_{k} = w_{2k} : 0\leq K\leq 2^{n} \} \subseteq \mathcal{U}^0 (M)$.  We shall show that the latter family satisfies the desired properties in the hypotheses of the lemma, and hence  we can apply the induction hypothesis to it.\smallskip    
	
Fix $0\leq k\leq 2^n -2$, by an application of the fundamental identity \eqref{eq fundamental identity UaUbUa} we get
	$$\begin{aligned} u_{k+2}  &= w_{2 k +4} = U_{w_{2k +3}} (w_{2k +2}^*) = U_{U_{w_{2k +2}} (w_{2k +1}^*) } (w_{2k +2}^*) \\
		&= U_{w_{2k +2}} U_{w_{2k +1}^*} U_{w_{2k +2}} (w_{2k +2}^*) = U_{w_{2k +2}} U_{w_{2k +1}^*} (w_{2k +2})\\
		&= U_{w_{2k +2}} U_{w_{2k +1}^*} U_{w_{2k +1}} (w_{2k}^*) = U_{w_{2k +2}}  (w_{2k}^*) = U_{u_{k +1}}  (u_{k}^*).
	\end{aligned} $$ 
Actually, the previous identities and the induction hypothesis combined with the fundamental identity \eqref{eq fundamental identity UaUbUa} lead to $$\begin{aligned}\Phi (U_{u_{k+1}} (u_k^*))& = \Phi (U_{w_{2k+2}} (w_{2k}^*)) = \Phi ( U_{w_{2k +3}} (w_{2k +2}^*) ) = U_{\Phi(w_{2 k+3})} (\Phi (w_{2k+2})^*) \\
		&= U_{\Phi( U_{w_{2 k+2}} (w_{2k+1}^*))} (\Phi (w_{2k+2})^*) = U_{ U_{\Phi(w_{2 k+2})} (\Phi(w_{2k+1})^*)} (\Phi (w_{2k+2})^*)\\
		& = U_{\Phi(w_{2 k+2})} U_{\Phi(w_{2k+1})^*} U_{\Phi(w_{2 k+2})} (\Phi (w_{2k+2})^*) \\
		&= U_{\Phi(w_{2 k+2})} U_{\Phi(w_{2k+1})^*} (\Phi(w_{2 k+2})) \\
		&= U_{\Phi(w_{2 k+2})} U_{\Phi(w_{2k+1})^*} (\Phi(U_{w_{2 k+1}} (w_{2k}^*))) \\
		&= U_{\Phi(w_{2 k+2})} U_{\Phi(w_{2k+1})^*} U_{\Phi(w_{2 k+1})} (\Phi(w_{2k})^*)\\
		&= U_{\Phi(w_{2 k+2})} (\Phi(w_{2k})^*) = U_{\Phi(u_{k+1})} (\Phi (u_k)^*).
	\end{aligned}$$ We have proved that the family formed by the $u_k$'s satisfies the properties in the statement of the lemma.\smallskip
	
The induction hypothesis applied to the family $\{u_{k} : 0\leq K\leq 2^{n} \} \subseteq \mathcal{U}^0 (M)$ gives  $$U_{w_{2^{n}}} (w_0^*) = U_{u_{2^{n-1}}} (u_0^*) = u_{2^n} = w_{2^{n+1}}$$ and $$\Phi\left(U_{w_{2^{n}}} (w_0^*) \right) =\Phi\left(U_{u_{2^{n-1}}} (u_0^*) \right) = U_{\Phi(u_{2^{n-1}})} \Phi(u_0)^*= U_{\Phi(w_{2^{n}})} \Phi(w_0)^*,$$ which concludes the proof.
\end{proof}
 
As we have seen in Proposition~\ref{p self adjoint quadratic subset princ comp}, the principal component of the set of unitary elements in a unital JB$^*$-algebra $M$ is the smallest quadratic subset of $\mathcal{U} (M)$ containing $e^{i M_{sa}}$. So, in order to understand the precise form of a surjective isometry between principal components the behaviour on finite compositions of $U$-maps of elements in $e^{i M_{sa}}$ seems to be crucial. We recall that a mapping between unital JB$^*$-algebras is called \emph{unital} if it maps the unit to the unit.

\begin{proposition}\label{p unital SI over U-products of exp}\cite[Proposition 3.3]{CuEnHirMiuPe2022} Let  $M$ and $N$ be unital JB$^*$-algebras. Let  Iso$_{1}(\mathcal{U}^0 (M), \mathcal{U}^0 (N))$ denote the set of all surjective unital isometries from $\mathcal{U}^0 (M)$ onto $\mathcal{U}^0 (N)$. Then the identity $$ \Delta_0(U_{e^{i h_n}}\cdots U_{e^{i h_1}}(e^{i h_0}))=U_{\Delta_0(e^{i h_n})}\cdots U_{\Delta_0(e^{i h_1})}(\Delta(e^{i h_0})),$$ holds for all $\Delta_0\in \hbox{Iso}_{1} (\mathcal{U}^0 (M), \mathcal{U}^0 (N))$, all natural $n$, and all $h_n,\dots,h_1, h_0\in M_{sa}$. \smallskip

\noindent Furthermore, there exists a mapping $$k:M_{sa}\times \hbox{Iso}_{1} (\mathcal{U}^0 (M), \mathcal{U}^0 (N))\to N_{sa}$$ satisfying the following properties:\begin{enumerate}[$(a)$]\item If we fix $\Delta_0\in \hbox{Iso}_{1} (\mathcal{U}^0 (M), \mathcal{U}^0 (N)),$ the mapping $k(\cdot, \Delta_0): M_{sa}\to N_{sa}$ is a surjective linear isometry;
		\item $k(\cdot, \Delta_0)^{-1} = k(\cdot, \Delta_0^{-1}),$ for every $\Delta_0\in \hbox{Iso}_{1} (\mathcal{U}^0 (M), \mathcal{U}^0 (N));$
		\item $\Delta_0(e^{i t h } )= e^{i t k(h, \Delta_0)},$ for all $\Delta_0\in \hbox{Iso}_{1} (\mathcal{U}^0 (M), \mathcal{U}^0 (N))$, $h\in M_{sa}$ and $t\in \mathbb{R}$.
	\end{enumerate}	
\end{proposition}

\begin{proof} We begin the proof by establishing a series of properties and claims. Fix an arbitrary $\Delta_0\in$ Iso$_{1}(\mathcal{U}^0 (M), \mathcal{U}^0 (N))$ and $h\in M_{sa},$ and consider the continuous mapping $ \mathbb{R}\to \mathcal{U}^0 (M)$, $t\mapsto v_t = v_t^{h}:=e^{i t h}$. In a first step we shall show that 
	\begin{equation}\label{eq n1 claim1}
		\Delta_0( U_{v_t}(v_s))=
		U_{\Delta_0(v_t)}(\Delta_0( v_s^* )^*)
	\end{equation} for every  $s,t\in\mathbb{R}$. In order to prove \eqref{eq n1 claim1}, choose a positive integer $m$ such that \begin{equation}\label{eq m} e^{\frac{\| i(s+t) h \|_{_{M}}}{2^m}}-1 < \frac{1}{2},\hbox{ and hence } \left\| e^{\frac{i(s+t) h}{2^m}}-\11b{M} \right\|_{_{M}}\leq e^{\frac{\| i(s+t) h \|_{_{M}}}{2^m}}-1 < \frac{1}{2}.
	\end{equation}
	
Given $0\leq l \leq 2^{m+1}$, we set $$u_l := v_s^* \circ e^{\frac{i l (s+t) h}{2^m}}= e^{-i s h}\circ e^{\frac{i l (s+t) h}{2^m}} =e^{\left(-i s +\frac{i l (s+t)}{2^m}\right)  h},$$ and consider the family $\{ u_l : 0\leq l \leq 2^{m+1}\}\subseteq \mathcal{U}^0(M)$. 
	
By applying the properties of the holomorphic functional calculus, the reader can check that the following identities hold:
\begin{equation}\label{eq Uuk+1 is uk+2} u_0 =v_s^*,\  u_{2^m}=v_t,\  u_{2^{m+1}}=U_{v_t} (v_s), \
	U_{u_{l+1}} ({u_l}^*)  =U_{u_{l+1}}(u_l^{-1})= u_{l+2}, 
\end{equation}  for any $0\leq l \leq 2^{m+1}-2$,  and 
$$ \begin{aligned} 
			\|u_{l+1} - u_l \| &= \left\| {v_s}^*\circ e^{\frac{i(l+1)(s+t) h}{2^m}}- {v_s}^*\circ e^{\frac{i l (s+t) h}{2^m}} \right\| \\
			&\leq  \|{v_s}^*\| \ \left\| e^{\frac{i l (s+t) h}{2^m}}\circ e^{\frac{i(s+t) h}{2^m}} - e^{\frac{i l (s+t) h}{2^m}} \right\|\\
			& \leq \left\|e^{\frac{i l (s+t) h}{2^m}}\right\|\  \left\| e^{\frac{i (s+t) h}{2^m}}-\11b{M} \right\| \\
			&\leq e^{\frac{\| i(s+t) h \|}{2^m}}-1 < \frac{1}{2},\hbox{ for all } 0\leq l \leq 2^{m+1}-1.
		\end{aligned}$$ Now, Theorem~\ref{t Delta preserves inverted triple products in U0M}$(3)$ applied to $\Delta_0$ and the unitaries $u_{l+1}$ and $u_l$ gives
	\begin{equation}\label{eq application of Th2.9 3}
		\Delta_0(U_{u_{l+1}} ({u_l}^*)) = U_{\Delta_0(u_{l+1})} (\Delta_0(u_{l})^*), \hbox{  for every $l\in \{0,\dots, 2^{m+1}-1\}$.}
	\end{equation}
	
	It follows from \eqref{eq Uuk+1 is uk+2} and \eqref{eq application of Th2.9 3} that we are in a position to apply Lemma~\ref{l HM Lemma 7 Jordan} with $\Phi = \Delta_0$ and $n= m+1$ to deduce that
	$$\begin{aligned}
		\Delta_0 (U_{v_t} ({v_s}) )&= \Delta_0 (U_{v_t} (({v_s}^*)^*)) =\Delta_0(U_{u_{2^m}} ({u_0}^*)) \\
		&= U_{\Delta_0(u_{2^m})} (\Delta_0(u_{0})^*) =U_{\Delta_0 (v_t)} (\Delta_0({v_s}^*)^*),
	\end{aligned}$$  which concludes the proof of \eqref{eq n1 claim1}.\smallskip
	
By taking $t=1$, $s=0$ and $h=h_1\in M_{sa}$ in \eqref{eq n1 claim1} we obtain
	\begin{equation}\label{eq Deltaa1 to the square}
		\Delta_0({e^{i 2 h_1}}) =\Delta_0( U_{e^{i h_1}} (\11b{M}))= U_{\Delta_0(e^{i h_1})}(\11b{N})= \Delta_0(e^{i h_1})^2.
	\end{equation} 
	
On the other hand, by replacing $t=0$ in \eqref{eq n1 claim1} we get
	\begin{equation}\label{eq Delta0 adjoint} \Delta_0(e^{i s h}) =\Delta_0(e^{-i s h})^*, \ \ \forall h\in M_{sa}, \ \forall s\in \mathbb{R}.\end{equation}
	
The most important consequence derived from \eqref{eq n1 claim1} (and \eqref{eq Delta0 adjoint}) assures that, for each $h\in M_{sa}$, the mapping $$\Gamma: \mathbb{R}\to \mathcal{U}^0 (N), \ t\mapsto \Gamma(t) =\Delta_0(e^{i t h})$$ is a uniformly continuous one-parameter family of unitaries satisfying  $$\begin{aligned} U_{\Gamma(t)} ( \Gamma(s) ) &= U_{ \Delta_0(e^{i t h}) } ( \Delta_0(e^{ i s h}) ) = U_{ \Delta_0(e^{i t h}) } ( \Delta_0(e^{- i s h})^* ) \\ 
	&=  \Delta_0( U_{e^{i t h}}(e^{i s h})) =  \Delta_0( e^{i(2 t + s) h})= \Gamma(2t +s),
\end{aligned}$$ and hence $\Gamma$ satisfies the hypotheses in the Jordan version of Stone's one-parameter theorem presented in Theorem~\ref{t Jordan unitary groups version of Stone's theorem}. Therefore there exists a unique $k(h, \Delta_0)\in N_{sa},$ depending on $h$ and $\Delta_0,$ satisfying \begin{equation}\label{eq application Stone one parameter thm}  \Delta_0(e^{i t h}) = e^{i t  k(h, \Delta_0)}, \hbox{ for all } t\in \mathbb{R}.
	\end{equation}
	
	We observe that the mapping $k(\cdot, \Delta_0): M_{sa} \to N_{sa}$ is well defined, and moreover, by \eqref{eq application Stone one parameter thm}, $\11b{N}= \Delta_0(e^{i t 0}) = e^{i t  k(0, \Delta_0)},$ for all $t\in \mathbb{R},$ which proves that $k(0, \Delta_0)=0$. \smallskip
	
We claim that \begin{equation}\label{eq k is a surjective isometry} 
	k(\cdot, \Delta_0): M_{sa} \to N_{sa} \hbox{  is a surjective linear isometry.}
\end{equation} Namely, the previous conclusion applied to $\Delta_0^{-1}: \mathcal{U}^0 (N) \to \mathcal{U}^{0} (M)$, show that  for each $k_0\in N_{sa}$ there exists a well defined mapping $h(\cdot, \Delta_0^{-1}): N_{sa}\to M_{sa}$ satisfying \begin{equation}\label{eq application Stone one parameter thm Delta inverse}  \Delta_0^{-1} (e^{i t k_0}) = e^{i t  h(k_0, \Delta_0^{-1})}, \hbox{ for all } t\in \mathbb{R}.
\end{equation} The identities \eqref{eq application Stone one parameter thm} and \eqref{eq application Stone one parameter thm Delta inverse} prove that $$\begin{aligned}e^{i t h_0} &= \Delta_0^{-1} \Delta_0(e^{i t h_0}) = \Delta_0^{-1} \left( e^{i t  k(h_0, \Delta_0)} \right) = e^{i t  h(k(h_0, \Delta_0), \Delta_0^{-1})},\ \hbox{ and }\\
e^{i t k_0} &= \Delta_0 \Delta_0^{-1} (e^{i t k_0}) = \Delta_0 \left( e^{i t  h(k_0, \Delta_0^{-1})} \right) = e^{i t  k(h(k_0, \Delta_0^{-1}), \Delta_0)},
\end{aligned}  $$ for all $t\in \mathbb{R}$, $h_0\in M_{sa}$, $k_0\in N_{sa}$, and hence $$h(k(h_0, \Delta_0), \Delta_0^{-1}) = h_0, \hbox{ and } k(h(k_0, \Delta_0^{-1}), \Delta_0) = k_0, $$ for all $h_0\in M_{sa}$, $k_0\in N_{sa}$, that is, $k(\cdot, \Delta_0)$ and $h(\cdot, \Delta_0^{-1})$ are bijections with $h(\cdot, \Delta_0^{-1}) = k(\cdot, \Delta_0)^{-1}$.\smallskip

Furthermore, given $h_0,h_0'\in M_{sa}$ and a real number $t$ we have
\begin{equation}\label{eq derivative} \|\cdot\|\hbox{-}\lim_{t\to 0^+} \frac{e^{i t h_0}-e^{ith_0'}}{t}= \|\cdot\|\hbox{-}\lim_{t\to 0^+} \frac{e^{i t h_0}-\11b{M}}{t}-\frac{e^{i t h_0'}-\11b{M}}{t}  = ih_0-ih_0'.
\end{equation}
We deduce from \eqref{eq application Stone one parameter thm} and the hypotheses on $\Delta$ that $$ \|e^{it k(h_0, \Delta_0)} - e^{it k(h_0', \Delta_0)} \|_{_{N}}=\|\Delta_0 (e^{ith_0})-\Delta_0 (e^{ith_0'}) \|_{_{N}}=\| e^{i t h_0}-e^{it h_0'} \|_{_{M}},$$ which by \eqref{eq derivative} leads to   $$\|k(h_0, \Delta_0) - k(h_0', \Delta_0) \|_{_{N}}=\|h_0-h_0' \|_{_{M}},$$ witnessing that $k(\cdot, \Delta_0) :M_{sa}\to N_{sa}$ is a surjective isometry, and hence a surjective linear isometry by the Mazur--Ulam theorem. This proves the claim in \eqref{eq k is a surjective isometry}.\smallskip

Our next goal consist in proving that 
	\begin{equation}\label{eq n1 with 2 exponentials}
		\Delta_0( U_{e^{i h_1}}(e^{i h_0}))= U_{\Delta_0 (e^{i h_1})} (\Delta_0(e^{i h_0})), \hbox{ for all } h_1,h_0\in M_{sa}.
	\end{equation}
	
	Namely, fix $h_1\in M_{sa}$ and consider the surjective isometry $\Delta_0^{h_1} :=$ $U_{\Delta_0 (e^{i h_1})^*} \Delta_0$ $U_{e^{i h_1}} : \mathcal{U}^{0} (M) \to \mathcal{U}^{0} (N),$ which is well defined by Theorem~\ref{t principal component as product of exp}. By \eqref{eq Deltaa1 to the square} we deduce that
	$$\begin{aligned} \Delta_0^{h_1} (\11b{M}) &= U_{\Delta_0(e^{i h_1})^*} \Delta_0 U_{e^{i h_1}} (\11b{M}) = U_{\Delta_0(e^{i h_1})^*} \Delta_0 (e^{ 2 i h_1}) =  U_{\Delta_0(e^{i h_1})^*} (\Delta_0 (e^{ i h_1})^2) \\
		&= U_{\Delta_0(e^{i h_1})^*} U_{\Delta_0 (e^{ i h_1})} (\11b{M})=  \11b{N},
	\end{aligned} $$ which proves that the mapping $\Delta_0^{h_1}  = U_{\Delta_0(e^{i h_1})^*} \Delta_0 U_{e^{i h_1}}$ is unital.\smallskip
	
	Let $k(\cdot, \Delta_0), k(\cdot, \Delta_0^{h_1}): M_{sa} \to N_{sa}$ be the surjective linear isometries given by \eqref{eq application Stone one parameter thm} and \eqref{eq k is a surjective isometry} for $\Delta_0$ and $\Delta_0^{h_1}$, respectively.\smallskip
	
	Choose, by continuity, $\varepsilon>0$ such that
	$$ \|e^{i h_1}-e^{-i t h}\|<\frac12, \ \forall t\in \mathbb{R}, \ h\in M_{sa} \hbox{ with } |t-1|<\varepsilon, \hbox{ and } \|h + h_1\|<\varepsilon.$$ Theorem~\ref{t Delta preserves inverted triple products in U0M}$(3)$, applied to $\Delta_0$ and the unitaries $e^{i h_1}, e^{-it h}\in \mathcal{U}^0(M)$ with $t$ in $\mathbb{R}$ and $h\in M_{sa}$ satisfying $\|h+h_1\|<\varepsilon$ and $|t-1|<\varepsilon,$ proves that
	$$\Delta_0(U_{e^{i h_1}} (e^{it h})) =U_{\Delta_0( e^{i h_1} )} ((\Delta_0(e^{-it h}))^*)=\eqref{eq Delta0 adjoint}=U_{\Delta_0(e^{i h_1})} (\Delta_0(e^{it h})).$$
This conclusion and the properties defining the mappings $k(\cdot, \Delta_0)$ and $k(\cdot, \Delta_0^{h_1})$ assure that $$e^{i t  k(h , \Delta_0)} = e^{i t  k(h, \Delta_0^{h_1})},\  \forall t\in \mathbb{R}, \ h\in M_{sa} \hbox{ with } |t-1|<\varepsilon, \hbox{ and } \|h + h_1\|<\varepsilon.$$ If we fix an arbitrary $h\in M_{sa}$ with $\|h + h_1\|<\varepsilon$ and we take a simple derivative at $t=1$ in the above equality we obtain $$ k(h , \Delta_0)= k(h, \Delta_0^{h_1}), \hbox{ for all } h\in M_{sa} \hbox{ with } \|h + h_1\|<\varepsilon.$$ We conclude from the linearity of $k(\cdot, \Delta_0),$ $k(\cdot, \Delta_0^{h_1})$  that $k(\cdot, \Delta_0)=k(\cdot, \Delta_0^{h_1})$.\smallskip
	
The conclusion in \eqref{eq n1 with 2 exponentials} now follows from the identity $$ \Delta_0(e^{i t h_0}) = e^{i t  k(h_0, \Delta_0)}   = e^{i t  k(h_0, \Delta_0^{h_1})} = \Delta_0^{h_1} (e^{i t h_0}) = U_{\Delta_0 (e^{i h_1})^*} \Delta_0 U_{e^{i h_1}} (e^{i t h_0}) ,$$ for all real $t$ and all $h_0\in M_{sa}$ (cf. \eqref{eq application Stone one parameter thm}).\smallskip

We deal now with the first statement in the proposition.  We observe that, since for each self-adjoint element $h\in M_{sa}$ we have $U_{e^{i \frac{h}{2}}}(\11b{M}) ) = e^{i h},$ the desired identity is equivalent to  \begin{equation}\label{eq ideneity at 1} \Delta_0(U_{e^{i h_n}}\cdots U_{e^{i h_1}}(\11b{M}))=U_{\Delta_0(e^{i h_n})}\cdots U_{\Delta_0(e^{i h_1})}(\11b{N})
\end{equation} for every $\Delta_0,$ $n\in\mathbb{N}$, and $h_n,\dots,h_1\in M_{sa}$ as in the hypotheses. \smallskip

The desired statement will be proved by induction on $n$. The case $n=1$ follows from \eqref{eq ideneity at 1}.\smallskip
	
Suppose, by the induction hypothesis, that \eqref{eq ideneity at 1} holds for all natural $k\leq n$. Take an arbitrary unital surjective isometry $\Delta_1 : \mathcal{U}^0 (M)\to \mathcal{U}^0 (N)$ and $h_{n+1},h_n,\dots,h_1\in M_{sa}$. Set $w_{0} = e^{2 i h_{n+1}}$ and $\widetilde{w}_{0} = \Delta_1(e^{2 i h_{n+1}})\in \mathcal{U}^0 (N).$ We consider the unital JB$^*$-algebras given by the isotopes $M(w_{0})= M_{(w_0^*)}$ and $N(\widetilde{w}_{0})= M_{(\widetilde{w}_0^*)}$ with the corresponding Jordan products $\circ_{w_{0}}$ and $\circ_{\widetilde{w}_{0}}$ and involutions $*_{w_{0}}$ and $*_{\widetilde{w}_{0}}$, respectively. Since $w_{0}\in \mathcal{U}^0 (M)$ and $\widetilde{w}_{0} \in \mathcal{U}^0 (N),$ it follows from \eqref{eq princ components of isotopes coincide} that $$\mathcal{U}^0(M(w_{0}))=\mathcal{U}^0(M), \hbox{ and } \mathcal{U}^0(N(\widetilde{w}_{0}))=\mathcal{U}^0(N),$$ and thus $\widetilde{\Delta}_1 : \mathcal{U}^0(M(w_{0}))=\mathcal{U}^0(M)\to \mathcal{U}^0(N(\widetilde{w}_{0}))=\mathcal{U}^0(N),$ $\widetilde{\Delta}_1 (a) = {\Delta}_1(a)$ is a unital (i.e., $\Delta_1(w_{0}) = \widetilde{w}_{0}$) surjective isometry.\smallskip
	
The mapping $U_{e^{i h_{n+1}}} : M\to M(w_{0})$ is a unital surjective linear isometry, and hence a Jordan $^*$-isomorphism between these two unital JB$^*$-algebras (cf. \cite[Theorem 6]{WriYo} or \cite[Theorem 1.9]{IsRo} or Kaup's Banach--Stone in  Theorem~\ref{t Kaup-Banach-Stone}). Therefore $\hat{h}_k:= U_{e^{i h_{n+1}}} (h_k) \in (M(w_{0}))_{sa}$ for all $1\leq k \leq n$. By applying the induction hypothesis to $\widetilde{\Delta}_1 : \mathcal{U}^0(M(w_{0}))\to \mathcal{U}^0(N(\widetilde{w}_{0})),$ and $\hat{h}_1, \ldots, \hat{h}_n$ we get
\begin{equation}\label{eq induction hypothesis applied to tildeDelta1}\begin{aligned} &{\Delta}_1 U^{(w_{0})}_{\exp_{(w_{0})} (i \hat{h}_{n})} \ldots U^{(w_{0})}_{\exp_{(w_{0})} (i \hat{h}_{1})} (w_{0}) = \widetilde{\Delta}_1 U^{(w_{0})}_{\exp_{(w_{0})} (i \hat{h}_{n})} \ldots U^{(w_{0})}_{\exp_{(w_{0})} (i \hat{h}_{1})} (w_{0}) \\
			&= U^{(\widetilde{w}_{0})}_{ \widetilde{\Delta}_1(\exp_{(w_{0})} (i \hat{h}_{n}))} \ldots U^{(\widetilde{w}_{0})}_{ \widetilde{\Delta}_1(\exp_{(w_{0})} (i \hat{h}_{1}))} (\widetilde{w}_{0}),
		\end{aligned}
	\end{equation} where, as in previous results, $\exp_{(w_0)}$ denotes the exponential in the unital JB$^*$-algebra $M(w_0)$.\smallskip
	
Let us carefully decode the identity in \eqref{eq induction hypothesis applied to tildeDelta1}. Having in mind that the $U$-operator in the isotope $M(w_0) = M_{(w_0^*)}$ ca be expressed in terms of the $U$-operator in $M$ (cf. \eqref{eq quadratic operator in the isotope}) we derive that \begin{equation}\label{eq 030521 a}
		U^{(w_0)}_{v_1} (v_2) = U_{v_1} U_{w_0^*} (v_2)  = U_{v_1} U_{e^{-2 i h_1}} (v_2).
	\end{equation}
	
Furthermore, since $\exp_{(w_0)}$ denotes the exponential in the unital JB$^*$-algebra $M(w_0)$, it follows from the fact that $U_{e^{i h_{n+1}}} : M\to M(w_0)$ is a Jordan $^*$-isomorphism that \begin{equation}\label{eq exponentials of hathk} \exp_{w_0} (i \hat{h}_k)= \exp_{w_0} (i U_{e^{i h_{n+1}}} (h_k)) = U_{e^{i h_{n+1}}} (e^{ih_k}), \hbox{ for all } 1\leq k\leq n,
	\end{equation} and by the fundamental identity \eqref{eq fundamental identity UaUbUa} we derive that
	\begin{equation}\label{eq exponentials and U} \begin{aligned}
			U_{\exp_{w_0} (i \hat{h}_k)}^{(w_0)} & = U_{\exp_{w_0} (i \hat{h}_k)} U_{{w_0}^*}  = U_{U_{e^{i h_{n+1}}} (e^{ih_k})} U_{{w_0}^*} \\
			&= U_{e^{i h_{n+1}}} U_{e^{ih_k}} U_{e^{i h_{n+1}}} U_{{e^{- 2 i h_{n+1}}}} = U_{e^{i h_{n+1}}} U_{e^{ih_k}} U_{e^{-i h_{n+1}}},
		\end{aligned}
	\end{equation} (in the last step we applied that $e^{i h_{n+1}}$ and $e^{-2i h_{n+1}}$ operator commute). \smallskip
	
Now, by combining \eqref{eq 030521 a} and \eqref{eq exponentials and U}, the left hand side of \eqref{eq induction hypothesis applied to tildeDelta1} writes in the form \begin{equation}\label{eq left hand side of 25} \begin{aligned} & {\Delta}_1 U^{(w_{0})}_{\exp_{(w_{0})} (i \hat{h}_{n})} \ldots U^{(w_{0})}_{\exp_{(w_{0})} (i \hat{h}_{1})} (w_{0}) \\
			&=  {\Delta}_1 U_{e^{i h_{n+1}}} U_{e^{ih_n}} U_{e^{-i h_{n+1}}} \ldots U_{e^{i h_{n+1}}} U_{e^{ih_1}} U_{e^{-i h_{n+1}}}  (e^{2i h_{n+1}}) \\
			&= {\Delta}_1 U_{e^{i h_{n+1}}} U_{e^{ih_n}}  U_{e^{i h_{n-1}}}  \ldots U_{e^{i h_{2}}} U_{e^{ih_1}} (\11b{M}).
		\end{aligned}
	\end{equation}

Concerning the right hand side of \eqref{eq induction hypothesis applied to tildeDelta1} we observe that, by \eqref{eq exponentials of hathk}, for $1\leq k\leq n$ we have
	$$\begin{aligned}
		U^{(\widetilde{w}_{0})}_{ \widetilde{\Delta}_1(\exp_{(w_{0})} (i \hat{h}_{k}))} &=  U^{(\widetilde{w}_{0})}_{{\Delta}_1(\exp_{(w_{0})} (i \hat{h}_{k}))} = U_{{\Delta}_1(U_{e^{i h_{n+1}}} (e^{ih_k}))} U_{\widetilde{w}_{0}^*} \\
		&= \hbox{(by \eqref{eq n1 with 2 exponentials})} = U_{U_{{\Delta}_1(e^{i h_{n+1}})} {\Delta}_1 (e^{ih_k})} U_{\widetilde{w}_{0}^*} = \hbox{(by \eqref{eq fundamental identity UaUbUa})} \\
		&= U_{{\Delta}_1(e^{i h_{n+1}})} U_{{\Delta}_1 (e^{ih_k})} U_{{\Delta}_1(e^{i h_{n+1}})} U_{\Delta_1(e^{2 i h_{n+1}})^*}\\
		&= \hbox{(by \eqref{eq Deltaa1 to the square})} = U_{{\Delta}_1(e^{i h_{n+1}})} U_{{\Delta}_1 (e^{ih_k})} U_{{\Delta}_1(e^{i h_{n+1}})} U_{(\Delta_1(e^{ i h_{n+1}})^2)^*}\\
		& = U_{{\Delta}_1(e^{i h_{n+1}})} U_{{\Delta}_1 (e^{ih_k})} U_{{\Delta}_1(e^{i h_{n+1}})^*} \\
		&= U_{{\Delta}_1(e^{i h_{n+1}})} U_{{\Delta}_1 (e^{ih_k})} U_{{\Delta}_1(e^{i h_{n+1}})}^{-1}.
	\end{aligned} $$ Therefore, the right hand side of \eqref{eq induction hypothesis applied to tildeDelta1} writes in the form
	\begin{equation}\label{eq right han side simplified}\begin{aligned}
			& U^{(\widetilde{w}_{0})}_{ \widetilde{\Delta}_1(\exp_{(w_{0})} (i \hat{h}_{n}))} \ldots U^{(\widetilde{w}_{0})}_{ \widetilde{\Delta}_1(\exp_{(w_{0})} (i \hat{h}_{1}))} (\widetilde{w}_{0}) \\
			&= U_{{\Delta}_1(e^{i h_{n+1}})} U_{{\Delta}_1 (e^{ih_n})} \ldots  U_{{\Delta}_1 (e^{ih_1})} U_{{\Delta}_1(e^{i h_{n+1}})}^{-1} (\Delta_1(e^{2 i h_{n+1}}))\\
			&= \hbox{(by \eqref{eq Deltaa1 to the square})} = U_{{\Delta}_1(e^{i h_{n+1}})} U_{{\Delta}_1 (e^{ih_n})} \ldots  U_{{\Delta}_1 (e^{ih_1})} U_{{\Delta}_1(e^{i h_{n+1}})}^{-1} (\Delta_1(e^{ i h_{n+1}})^2)\\
			& =  U_{{\Delta}_1(e^{i h_{n+1}})} U_{{\Delta}_1 (e^{ih_n})} \ldots  U_{{\Delta}_1 (e^{ih_1})} (\11b{N}). \end{aligned}
	\end{equation}
	
Finally, \eqref{eq induction hypothesis applied to tildeDelta1}, \eqref{eq left hand side of 25} and \eqref{eq right han side simplified} give $$ {\Delta}_1 U_{e^{i h_{n+1}}} U_{e^{ih_n}}   \ldots U_{e^{i h_{2}}} U_{e^{ih_1}} (\11b{M})=  U_{{\Delta}_1(e^{i h_{n+1}})} U_{{\Delta}_1 (e^{ih_n})} \ldots  U_{{\Delta}_1 (e^{ih_1})} (\11b{N}),$$ which concludes the induction argument proving \eqref{eq ideneity at 1}.
\end{proof}	

The main technical results to understand the precise form of a surjective isometry between the principal components of the unitary sets in two unital JB$^*$-algebras have been already established in Proposition~\ref{p unital SI over U-products of exp}.  We can now proceed to present the desired description.

\begin{theorem}\label{t surj isom principal components}{\rm\cite[Theorem 3.4]{CuEnHirMiuPe2022}}
	Let $\Delta: \mathcal{U}^0 (M)\to \mathcal{U}^0 (N)$ be a surjective isometry between the principal components of two unital JB$^*$-algebras. Then there exist self-adjoint elements $k_1,\ldots,k_n$ in $N$, a central projection $p$ in $N,$ and a Jordan $^*$-isomorphism $\Phi:M\to N$ such that $$\begin{aligned}\Delta(u) &= p\circ U_{e^{i k_n}} \ldots  U_{e^{i k_1}}  \Phi (u) + (\11b{N}-p)\circ \left( U_{e^{-i k_n}} \ldots  U_{e^{-i k_1}} \Phi(u)\right)^*,
	\end{aligned}$$ for all $u\in \mathcal{U}^0(M)$. Moreover, the elements $k_1,\ldots,k_n\in N_{sa}$ can be replaced by any finite collection of self-adjoint elements $\tilde{k}_1,\ldots,\tilde{k}_m$ such that $ \Delta(\11b{M}) = U_{e^{i \tilde{k}_m}} \ldots  U_{e^{i \tilde{k}_1}} (\11b{N})$. Under these hypotheses $M$ and $N$ are Jordan $^*$-isomorphic, and there exists a surjective real-linear isometry {\rm(}i.e., a real-linear triple isomorphism{\rm)} from $M$ onto $N$ whose restriction to $\mathcal{U}^0 (M)$ is $\Delta$.
\end{theorem}

\begin{proof} Since $\Delta(\11b{M})$ is a unitary in the principal component of $\mathcal{U} (N)$, by the algebraic description established in Theorem~\ref{t principal component as product of exp}, there exist $k_1,\ldots,k_n\in N_{sa}$ such that $ \Delta(\11b{M}) = U_{e^{i k_n}} \ldots  U_{e^{i k_1}} (\11b{N})$. Since the principal component is a quadratic subset (cf. Proposition~\ref{p self adjoint quadratic subset princ comp}), we deduce that the mapping $\Delta_0 = U_{e^{-i k_1}} \ldots  U_{e^{-i k_n}} \Delta : \mathcal{U}^0 (M)\to \mathcal{U}^0 (N)$ is a surjective isometry. Clearly $\Delta_0 (\11b{M}) = \11b{N}$.\smallskip
	
Let $k(\cdot,\Delta_0): M_{sa}\to N_{sa}$ be the surjective linear isometry given by Proposition~\ref{p unital SI over U-products of exp} satisfying $\Delta_0(e^{i t h } )= e^{i t k(h, \Delta_0)},$ for all $h\in M_{sa}$. \smallskip
	
The self-adjoint part of any JB$^*$-algebra is a JB-algebra. Thus, $k(\cdot,\Delta_0) : M_{sa}\to N_{sa}$ is a surjective linear isometry between JB-algebras. Under these circumstances we know the existence of a central symmetry $k(\11b{M},\Delta_0)$ in $N_{sa}$ and a Jordan $^*$-isomorphism $\Phi:M\to N$ such that \begin{equation}\label{eq f description}
		k(h,\Delta_0) = k(\11b{M},\Delta_0) \circ \Phi(h),
	\end{equation} for every $h\in M_{sa}$ 
(cf. \cite[Theorem 1.4 and Corollary 1.11]{IsRo}).
	\smallskip
	
Now, the arguments mix ideas from \cite{CuPe2022, HatMol2014} and \cite{CuEnHirMiuPe2022}. Since $k(\11b{M},\Delta_0)$ is a central symmetry, there exists a central projection $p$ in $N$ such that $k(\11b{M},\Delta_0)= 2p -\11b{N}= p - (\11b{N}-p)$, with $p\perp (\11b{N}-p)$. It is not hard to check that for any $n>0$ we have $$\left( 2p-\11b{N}\right)^n= \left( p-(\11b{N}-p)\right)^n= p + (-1)^n(\11b{N}-p).$$

	The mapping $\Delta_0$ can be expressed in terms of $p$ and $\Phi$ by some simple computations. Namely, given an arbitrary $h\in M_{sa}$, we have
	$$\begin{aligned}
		\Delta_0(e^{i h})&=e^{i k(h,\Delta_0)}=e^{i k(\11b{M},\Delta_0)\circ \Phi(h) }=e^{i(2p-\11b{N})\circ \Phi(h) }= \sum_{n=0}^{\infty}\frac{\left( i  (2p-\11b{N})\circ \Phi(h)\right)^n }{n!}.
	\end{aligned}$$
	We now make use of the properties of $\Phi$ as Jordan $^*$-isomorphism and the expression of $\left( 2p-\11b{N}\right)^n$ given above. It is worth noting that $k(\11b{M},\Delta_0)$ (and hence $p$) is central, and so it operator commutes with any element in $N$. Therefore,
	\begin{equation}\label{eq Delta0 in eix}\begin{aligned}
			\Delta_0(e^{ih})&=\sum_{n=0}^{\infty}\frac{ i^n  (2p-\11b{N})^n\circ \Phi(h)^n }{n!} =\sum_{n=0}^{\infty}\frac{ i^n  \left( p + (-1)^n(\11b{N}-p)\right) \circ \Phi(h^n) }{n!} \\
			&= \sum_{n=0}^{\infty}\frac{ i^n p \circ \Phi(h^n) }{n!} + \sum_{n=0}^{\infty}\frac{ i^n (-1)^n(\11b{N}-p) \circ \Phi(h^n) }{n!} \\
			&=p \circ \sum_{n=0}^{\infty}\frac{ i^n \Phi(h^n) }{n!} + (\11b{N}-p) \circ\sum_{n=0}^{\infty}\frac{ i^n (-1)^n \Phi(h^n) }{n!} \\
			&=p \circ \Phi\left( \sum_{n=0}^{\infty}\frac{ i^n h^n }{n!}\right)  + (\11b{N}-p) \circ\Phi\left( \sum_{n=0}^{\infty}\frac{ i^n (-1)^n h^n}{n!}\right) \\
			&= p \circ \Phi(e^{ih})  + (\11b{N}-p) \circ\Phi( e^{-ih})= p \circ \Phi(e^{ih})  + (\11b{N}-p) \circ\Phi( e^{ih})^*.
		\end{aligned}
	\end{equation}
	
We now deal with a unitary element $u$ in $\mathcal{U}^0(M)$. A new application of Theorem~\ref{t principal component as product of exp} assures the existence of a natural $m$ and $h_1,\dots, h_m\in M_{sa}$ such that $u=U_{e^{i h_m}}\cdots U_{e^{i h_1}}(\11b{M})$. By combining  Proposition~\ref{p unital SI over U-products of exp} and \eqref{eq Delta0 in eix} we arrive to
	$$\begin{aligned}
		\Delta_0(u)&=\Delta_0(U_{e^{i h_m}}\cdots U_{e^{i h_1}}(\11b{M})) = U_{\Delta_0(e^{i h_m})}\cdots U_{\Delta_0(e^{i h_1})}(\Delta_0(\11b{M})) \\
		&=U_{\left( p \circ \Phi(e^{i h_m})  + (\11b{N}-p) \circ\Phi( e^{i h_m})^* \right) }\cdots U_{\left( p \circ \Phi(e^{i h_1})  + (\11b{N}-p) \circ\Phi( e^{i h_1})^*\right) }(\11b{N})
	\end{aligned}$$
	
	Since $p$ is a central projection in $N$, it is not hard to see that $$U_{p\circ a+(\11b{N}-p)\circ b} = U_{p\circ a} + U_{(\11b{N}-p)\circ b} =p\circ U_{ a} +(\11b{N}-p)\circ U_{ b},$$ for all $a,b\in N$. Since $p\perp \11b{N}-p$ and $\Phi$ is a Jordan $^*$-isomorphism it follows that
	$$\begin{aligned}
		\Delta_0(u)&= p\circ \Phi (U_{e^{i h_m}}\cdots U_{e^{i h_1}} (\11b{M})) + (\11b{N}-p)\circ \Phi(U_{e^{i h_m}}\cdots U_{e^{i h_1}}(\11b{M}))^*\\
		&= p\circ \Phi (u) + (\11b{N}-p)\circ \Phi(u)^*.
	\end{aligned}$$
	
	Finally, by the definition of $\Delta_0$, we arrive at
	$$\begin{aligned}
		\Delta(u)&= U_{e^{i k_n}} \ldots  U_{e^{i k_1}} \Delta_0 (u)= U_{e^{i k_n}} \ldots  U_{e^{i k_1}} \left( p\circ \Phi (u) + (\11b{N}-p)\circ \Phi(u)^*\right),\\
		&=   p\circ U_{e^{i k_n}} \ldots  U_{e^{i k_1}}  \Phi (u) + (\11b{N}-p)\circ \left( U_{e^{-i k_n}} \ldots  U_{e^{-i k_1}} \Phi(u)\right)^*,
	\end{aligned}$$ for all $u\in \mathcal{U}^0(M),$ which concludes the proof.
\end{proof}

\begin{remark}\label{r thm principal component and unital} Let us observe that the conclusion in Theorem~\ref{t surj isom principal components} is much finer under the extra assumption that $\Delta$ is unital, since in that case there exist a central projection $p$ in $N,$ and a Jordan $^*$-isomorphism $\Phi:M\to N$ such that $$\begin{aligned}\Delta(u) &= p\circ  \Phi (u) + (\11b{N}-p) \Phi(u)^*= p\circ  \Phi (u) + (\11b{N}-p) \Phi(u^*),
	\end{aligned}$$ for all $u\in \mathcal{U}^0(M)$.
\end{remark}

The Anti-theorem~\ref{anti-theorem BKU} shows that the group of all surjective linear isometries on a unital JB$^*$-algebra $M$ does not always act transitively on $\mathcal{U} (M)$.   
The conclusion of Theorem~\ref{r thm principal component and unital} can be drastically simplified if we additionally assume that $\Delta(\11b{M})$ can be mapped to $\11b{N}$ under a concrete surjective linear isometry. The next sufficient condition which is an easy consequence of Lemma~\ref{l untaries at short distance in a unital JB$^*$-algebra}$(b)$. 

\begin{lemma}\label{l square root of v*} Let $u$ be a unitary element in a unital JB$^*$-algebra $M$. Suppose additionally that $\|\11b{M}- u\|<2$, then there exists $v\in \mathcal{U}(N)$ satisfying $v^2 = u$. Consequently, $U_{v^*} (u) = \11b{M}$. 
\end{lemma}

\begin{proof} Lemma~\ref{l untaries at short distance in a unital JB$^*$-algebra}$(b)$  implies the existence of a unitary $\omega$ in $M$ such that $U_{\omega} (u) = \11bN$. Then $u = U_{\omega^*} U_{\omega} (u) = U_{\omega^*} (\11b{M}) = (\omega^*)^2$ (cf. Lemma \ref{l untaries at short distance in a unital JB$^*$-algebra}$(c)$). The desired statement follows by just taking $v = \omega^*$.	
\end{proof}

The next corollary goes in the line of \cite[Theorem 3.6]{CuPe2022}.

\begin{corollary}\label{c surj isom principal components simplified} Let $\Delta: \mathcal{U}^0 (M)\to \mathcal{U}^0 (N)$ be a surjective isometry between the principal components of two unital JB$^*$-algebras. Suppose there exists $v\in \mathcal{U} (N)$ such that $U_v (\Delta(\11b{M})) = \11b{N}$. Then there exist a central projection $p$ in $N$ and a Jordan $^*$-isomorphism $\Phi:M\to N$ such that $$\begin{aligned}\Delta(w) &= p\circ U_{v^*}  \Phi (w) + (\11b{N}-p)\circ \left( U_{v} \Phi(w)\right)^*,
	\end{aligned}$$ for all $w\in \mathcal{U}^0(M)$. 
\end{corollary}

\begin{proof} The mapping $\tilde{\Delta} = U_{v} \Delta : \mathcal{U}^0 (M)\to \mathcal{U}^0 (N)$ is a surjective unital isometry. Theorem~\ref{t surj isom principal components} (see also Remark~\ref{r thm principal component and unital}) proves the existence of a Jordan $^*$-isomorphism $\Phi: M\to N$ such that $$ \tilde{\Delta} (w) = U_{v} \Delta (w) = p\circ  \Phi (w) + (\11b{N}-p) \Phi(w)^*,$$ for all $w\in \mathcal{U}^{0} (M).$ Since $p$ is central we have $$ \Delta (w) = p\circ  U_{v^*} \Phi (w) + (\11b{N}-p) U_{v^*} \Phi(w)^*= p\circ U_{v^*}  \Phi (w) + (\11b{N}-p)\circ \left( U_{v} \Phi(w)\right)^*,$$ for all $w\in \mathcal{U}^{0} (M).$   
\end{proof}

It is natural to ask what we can conclude in the case that we have a surjective isometry between two arbitrary connected components in two unital JB$^*$-algebras. Let $M$ be a unital JB$^*$-algebra, let $\mathcal{U}^c (M)$ be a connected component of $\mathcal{U} (M)$, and let $w$ be an element in $\mathcal{U}^c (M)$. We have already commented that $\mathcal{U}^c (M)$ is the principal component of the $w$-isotope $M(w) =M_{(w^*)}$ (cf. Remark~\ref{r connected components which are not the pricipal one}). When particularized to the corresponding isotopes, the previous Theorem~\ref{t surj isom principal components} can be re-stated as follows.

\begin{theorem}\label{t surj isom secondary components}{\rm\cite[Theorem 3.6]{CuEnHirMiuPe2022}} Let $M$ and $N$ be unital JB$^*$-algebras. Let $\Delta: \mathcal{U}^c (M)\to \mathcal{U}^c (N)$ be a surjective isometry between two connected components of $\mathcal{U} (M)$ and $\mathcal{U} (N),$ respectively. Fix $w_1\in \mathcal{U}^c (M)$ and $w_2\in \mathcal{U}^c (N)$. Then there exist $k_1,\ldots,k_n\in (N(w_2))_{sa}$, a central projection $p\in N(w_2)$ and a Jordan $^*$-isomorphism $\Phi:M(w_1)\to N(w_2)$ such that $$\begin{aligned}\Delta(u) &= p\circ_{w_2} U^{(w_2)}_{\exp_{(w_2)}({i k_n})} \ldots  U^{(w_2)}_{\exp_{(w_2)}({i k_1})}  \Phi (u) \\
		&+ (w_2-p)\circ_{w_2} \left( U^{(w_2)}_{\exp_{(w_2)}({-i k_n})} \ldots  U^{(w_2)}_{\exp_{(w_2)}({-i k_1})} \Phi(u)\right)^{*_{w_2}},
	\end{aligned}$$ for all $u\in \mathcal{U}^c (M)$. In particular, the unital JB$^*$-algebras $M(w_1)= M_{(w_1^*)}$ and $N(w_2) = N_{(w^*_1)}$ are Jordan $^*$-isomorphic, and there exists a surjective real-linear isometry {\rm(}i.e., a real-linear triple isomorphism{\rm)} from $M$ onto $N$ whose restriction to $\mathcal{U}^c (M)$ is $\Delta$.
\end{theorem}

Every unitary element of a unital JB$^*$-algebra $M$ need not be, in general, of the form $e^{ih}$ for some $h\in M_{sa}$, the statement is not true nor in the case of unital C$^*$-algebras {\rm(}cf. Remark~\ref{r unitaries which are squares} the discussion preceding Proposition 4.4.10 in \cite{KR1}{\rm)}.  However, in a von Neumann algebra, every unitary admits a logarithm (cf.  \cite[Theorem 5.2.5]{KR1}). The same conclusion holds for JBW$^*$-algebras. 

\begin{lemma}\label{l every untary in a JBWalgebra admits a logarithm} Let $u$ be a unitary element in a JBW$^*$-algebra $M$. Then there exists a self-adjoint element $h\in M$ such that $u = e^{ih}$. Therefore $\mathcal{U} (M) = \mathcal{U}^0 (M) =\exp \left({i M_{sa}}\right).$
\end{lemma}

\begin{proof} Let $\mathcal{W}$ denote the JBW$^*$-subalgebra of $M$ generated by $u,$ $u^*$ and the unit element. Clearly $\mathcal{W}$ is an associative JBW$^*$-algebra (cf. Lemma~\ref{l first properties invertible Jordan}). Therefore $\mathcal{W}$ is a commutative von Neumann algebra. Theorem 5.2.5 in \cite{KR1} implies the existence of an element $h\in \mathcal{W}_{sa}\subseteq M_{sa}$ such that $e^{i h} = u$ (in $\mathcal{W}$ and also in $M$).
\end{proof}

The next result is a consequence of Proposition~\ref{p unital SI over U-products of exp}, Theorem~\ref{t surj isom principal components} and Lemma~\ref{l every untary in a JBWalgebra admits a logarithm}. 

\begin{theorem}\label{t HM for unitary sets of unitary JBW*-algebras}{\rm\cite[Theorem 3.9]{CuPe2022}} Let $M$ and $N$ be two JBW$^*$-algebras, and let $\Delta: \mathcal{U} (M)\to \mathcal{U} (N)$ be a surjective isometry.  Then there exist a self-adjoint element $k_1$ in $N,$ a central projection $q\in N$, and a Jordan $^*$-isomorphism $\Phi:M\to N$ such that $$\begin{aligned}\Delta( u ) &= q \circ  U_{e^{i k_1}}\left(  \Phi(u) \right) + (\11b{N}-q)  \circ U_{e^{i k_1}}\left( \Phi( u)^*\right),\\
	\end{aligned}$$ for all $u\in \mathcal{U}(M)$. If $\Delta$ is unital we can take $k_1=0$.  Consequently, $\Delta$ admits a (unique) extension to a surjective real-linear isometry from $M$ onto $N$.\end{theorem}

\begin{proof} Lemma~\ref{l every untary in a JBWalgebra admits a logarithm} implies that  $\mathcal{U}(M) =\mathcal{U}^0 (M) = \exp\left(i M_{sa}\right)$ and $\mathcal{U}(N) =\mathcal{U}^0 (N) = \exp\left(i N_{sa}\right)$. Therefore  $\Delta(\11b{M}) = e^{2 i k_1}= U_{e^{i k_1}}\left(  \11b{N} \right)$ for some $k_1\in N_{sa}$.  Theorem~\ref{t surj isom principal components} gives the desired conclusion.
\end{proof}

We can now obtain a stronger version of the previous Theorem~\ref{t HM for unitary sets of unitary JBW*-algebras}.

\begin{theorem}\label{t HM for unitary sets of unitary JBW*-algebra to a JB*-algebra} Let $M$ and $N$ be two unital JB$^*$-algebras, and let $\Delta: \mathcal{U} (M)\to \mathcal{U} (N)$ be a surjective isometry. Suppose additionally that $M$ or $N$ is a JBW$^*$-algebra.  Then $M$ and $N$ are JBW$^*$-algebras and there exist a self-adjoint element $k_1$ in $N,$ a central projection $q\in N$, and a Jordan $^*$-isomorphism $\Phi:M\to N$ such that $$\begin{aligned}\Delta( u ) &= q \circ  U_{e^{i k_1}}\left(  \Phi(u) \right) + (\11b{N}-q)  \circ U_{e^{i k_1}}\left( \Phi( u)^*\right),\\
	\end{aligned}$$ for all $u\in \mathcal{U}(M)$. If $\Delta$ is unital we can take $k_1=0$.  Consequently, $\Delta$ admits a (unique) extension to a surjective real-linear isometry from $M$ onto $N$.\end{theorem}

\begin{proof} There is no loos of generality in assuming that $M$ is a JBW$^*$-algebra. In such a case, Lemma~\ref{l every untary in a JBWalgebra admits a logarithm} implies that  $\mathcal{U}(M) =\mathcal{U}^0 (M) = \exp\left(i M_{sa}\right)$. It follows that $\mathcal{U} (N) = \Delta (\mathcal{U}(M)) = \Delta (\mathcal{U}^0(M)) = \Delta (\exp(i M_{sa}))$ is connected, and hence  $\mathcal{U}(N) =\mathcal{U}^0 (N)$. By Theorem~\ref{t surj isom principal components} $M$ and $N$ are Jordan $^*$-isomorphic, and thus $N$ also is a JBW$^*$-algebra (cf. the comments on the Dixmier--Ng theorem in page~\pageref{comments on Dixmier-Ng theorem}). The desired conclusion follows from Theorem~\ref{t HM for unitary sets of unitary JBW*-algebras}.\end{proof}

We can now rediscover a slightly improved version of a result proved by Hatori and Moln{\'a}r in the setting of von Neumann algebras. 

\begin{theorem}\label{t HM for unitary sets of unitary C*-algebras}{\rm\cite[Corollary 3]{HatMol2014}} Let $A$ and $B$ be two unital C$^*$-algebras, and let $\Delta: \mathcal{U} (A)\to \mathcal{U} (B)$ be a surjective isometry. Suppose that $A$ or $B$ is a von Neumann algebra. Then there exist a central projection $q\in B$ and a Jordan $^*$-isomorphism $\Phi:A\to B$ such that $$\begin{aligned}\Delta( u ) &= \Delta(\11b{A}) \left( q   \Phi(u) + (\11b{B}-q)   \Phi( u)^* \right),\\
	\end{aligned}$$ for all $u\in \mathcal{U}(A)$. Consequently, $\Delta$ admits a (unique) extension to a surjective real-linear isometry from $A$ onto $B$.\end{theorem}

\begin{proof} The element $\Delta(\11b{A})$ is a unitary in $B$ and the mapping $\Delta(\11b{A})^* \Delta : \mathcal{U} (A) \to \mathcal{U} (B)$ is a unital surjective isometry. By Theorem~\ref{t HM for unitary sets of unitary JBW*-algebra to a JB*-algebra}, there exist a central projection $q\in B$, and a Jordan $^*$-isomorphism $\Phi:A\to B$ such that $$\Delta(\11b{A})^* \Delta( u ) =  q   \Phi(u) + (\11b{B}-q)   \Phi( u)^*$$ for all $u\in \mathcal{U} (A)$.  The rest is clear.
\end{proof}

In the same way that Theorem~\ref{t HM for unitary sets of unitary C*-algebras} was obtained from Theorem~\ref{t HM for unitary sets of unitary JBW*-algebras}, the next result, due to Hatori, follows from Theorem~\ref{t surj isom principal components} or from Theorem~\ref{t surj isom secondary components}. 

\begin{theorem}\label{t surj isom secondary components Cstat}{\rm\cite[Corollary 3.5]{Hatori14}} Let $A$ and $B$ be unital C$^*$-algebras. Let $\Delta: \mathcal{U}^c (A)\to \mathcal{U}^c (B)$ be a surjective isometry between two connected components of $\mathcal{U} (A)$ and $\mathcal{U} (B),$ respectively. Fix $w_1\in \mathcal{U}^c (A)$ and $w_2\in \mathcal{U}^c (B)$ such that $\Delta(w_1) = w_2$, $\mathcal{U}^c (A) = w_1 \mathcal{U}^0 (A) $ and $\mathcal{U}^c (B) = w_2 \mathcal{U}^0 (B)$. Then there exists a central projection $p \in B$, a Jordan $^*$-isomorphism $\Phi: A\to B$ and a unitary $u =  w_2 (p \Phi(w_1^*) + (\11b{B} - p) \Phi (w_1^*)^*)$ such that $$\Delta (w) = u (p \Phi (w) + (\11b{B} -p) \Phi(w)^*),$$ for all $w\in \mathcal{U}^c (A) = w_1 \mathcal{U}^0 (A)$. 
\end{theorem}

\begin{proof} We take part of the argument from \cite[Corollary 3.5]{Hatori14}. Consider the surjective isometry $\tilde{\Delta} : \mathcal{U}^0 (A) \to \mathcal{U}^0 (B)$ defined by $ \tilde{\Delta} (a) := w_2^*\Delta ( w_1 a)$. It follows from the definition that $\tilde{\Delta}$ is unital. By Theorem~\ref{t surj isom principal components} (see also Remark~\ref{r thm principal component and unital}) there exist a central projection $p$ in $B,$ and a Jordan $^*$-isomorphism $\Psi:A\to B$ such that $$\begin{aligned}\tilde{\Delta}(a) &= p\circ \Psi (a) + (\11b{N}-p)\circ \Psi(a)^*,
	\end{aligned}$$ for all $a\in \mathcal{U}^0(A)$. Therefore \begin{equation}\label{eq expression of tilde Delta 2308} \begin{aligned}{\Delta}(w) &= w_2 \left(p\circ \Psi (w_1^* w) + (\11b{N}-p)\circ \Psi(w_1^* w)^*\right), \hbox{ and }\\ 
		{\Delta}(w) &= w_2 \left(p\circ \Psi (w_1^* w) + (\11b{N}-p)\circ \Psi(w_1^* w)^*\right),
\end{aligned}
\end{equation} for all $w\in \mathcal{U}^c(A)= w_1 \mathcal{U}^0 (A) $. Define a mapping $\Phi:A\to B$ by $$\Phi (x) = p \Psi (w_1) \Psi (w_1^* x) +(\11b{B}-p) \Psi (w_1^* x) \Psi (w_1).$$ It is not hard to check that $\Phi$ is unital and $\Phi (w_1) = \Psi (w_1)$.\smallskip

A surprising geometric argument based on Kadison's Theorem~\ref{t Kadison 55} (or in Kaup's Theorem~\ref{t Kaup-Banach-Stone}) assures that $\Phi$ is a Jordan $^*$-isomorphism. Namely, since orthogonality implies $M$-orthogonality
$$\begin{aligned}
\left\| \Phi(x) \right\| &= \max\left\{\left\| p \Psi (w_1) \Psi (w_1^* x) \right\|, \left\| (\11b{B}-p) \Psi (w_1^* x) \Psi (w_1)\right\| \right\} \\
&= \max\left\{\left\| p \Psi (w_1^* x) \right\|, \left\| (\11b{B}-p) \Psi (w_1^* x)\right\| \right\} \\
&= \left\| p \Psi (w_1^* x) + (\11b{B}-p) \Psi (w_1^* x)\right\| = \left\| \Psi (w_1^* x)\right\| = \|w_1^* x\| = \|x\|.
\end{aligned}$$  This proves that $\Phi$ is a (unital) surjective isometry, and thus Kadison's theorem implies that $\Phi$ is a Jordan $^*$-isomorphism.\smallskip

Finally, having in mind that $p$ is central, \eqref{eq expression of tilde Delta 2308}, and the other properties, we get  
$$\begin{aligned}
&u (p \Phi (w) + (\11b{B} -p) \Phi(w)^*) = u (p \Psi (w_1) \Psi (w_1^* w) +(\11b{B}-p)  \Psi (w_1)^* \Psi (w_1^* w)^*) \\
&=  w_2 (p \Phi(w_1^*) + (\11b{B} - p) \Phi (w_1^*)^*) (p \Phi (w_1) \Psi (w_1^* w) +(\11b{B}-p)  \Phi (w_1)^* \Psi (w_1^* w)^*) \\
&= w_2 \left(p\circ \Psi (w_1^* w) + (\11b{N}-p)\circ \Psi(w_1^* w)^*\right) = \Delta (w).\end{aligned}$$ 
\end{proof} 

The next corollary shows that the metric space given by the set of unitary elements in a unital JB$^*$-algebra or by its principal component is a metric invariant.

\begin{corollary}\label{c two unital JB*algebras are isomorphic iff their unitaries are isometric iff their principal components are isometric}{\rm(\cite[Corollary 3.8]{CuPe2022} and \cite[Corollary 3.5]{CuEnHirMiuPe2022})}
	Let $M$ and $N$ be two unital JB$^*$-algebras. Then $M$ and $N$ are isometrically isomorphic as (complex) Banach spaces if, and only if, they are isometrically isomorphic as real Banach spaces if, and only if, $\mathcal{U}(M)$ and $\mathcal{U}(N)$ are isometrically isomorphic as metric spaces. Furthermore, the following statements are also  equivalent:
	\begin{enumerate}[$(a)$]
		\item $M$ and $N$ are Jordan $^*$-isomorphic;
		\item There exists a surjective isometry $\Delta: \mathcal{U}(M)\to \mathcal{U}(N)$ satisfying $\Delta (\11b{M})\in \mathcal{U}^0(N)$;
		\item There exists a surjective isometry $\Delta: \mathcal{U}^0(M)\to \mathcal{U}^0(N).$
	\end{enumerate}
\end{corollary}

\begin{proof} To see the first statement, we observe first that every surjective linear  isometry from $M$ onto $N$ is a surjective real-linear isometry. Every surjective real-linear isometry  $T: M\to N$  is a triple homomorphism (cf. Theorem~\ref{t Dang real-linear surjective isometries JBstar algebras}). Therefore, $T$ preserves triple products of the form $\{x,y,z\} = (x\circ y^*) \circ z + (z\circ y^*)\circ x -
	(x\circ z)\circ y^*$. In particular $T$ maps unitaries in $M$ to unitaries in $N$ (see Lemma~\ref{l unitaries in the homotope coincide}), and hence $T (\mathcal{U} (M))= \mathcal{U} (N)$. Therefore $\Delta = T|_{\mathcal{U} (M)} : \mathcal{U} (M) \to \mathcal{U} (N)$ is a surjective isometry. If there exists a surjective isometry $\Delta : \mathcal{U} (M) \to \mathcal{U} (N)$, the set $\Delta\left(\mathcal{U}^{0} (M)\right)$ must be a connected component (say $\mathcal{U}^{c} (N)$) of $\mathcal{U} (N)$ (cf. Corollary~\ref{c distance in different connected components is 2} and Remark~\ref{r connected components which are not the pricipal one}). Since $\Delta|_{\mathcal{U}^{0} (M)} : \mathcal{U}^{0} (M)\to \mathcal{U}^{c} (N)$ is a surjective isometry,  Theorem~\ref{t surj isom secondary components} proves that $M$ and $N$ are isometrically isomorphic as real Banach spaces. \smallskip

	We see next the equivalence of the three final statements. The implication $(a)\Rightarrow (b)$ is clear because every Jordan $^*$-isomorphism is a unital linear isometry and maps unitaries to unitaries.\smallskip
	
	$(b)\Rightarrow (c)$ Suppose there exists a surjective isometry $\Delta: \mathcal{U}(M)\to \mathcal{U}(N)$ satisfying $\Delta (\11b{M})\in \mathcal{U}^0(N)$. Since every surjective isometry maps connected components to connected components, $\Delta(\mathcal{U}^0(M))$ is a connected component of $\mathcal{U}(N)$ and contains an element in $\mathcal{U}^0(N)$, the equality $\Delta(\mathcal{U}^0(M))= \mathcal{U}^0(N)$ holds, and thus $\Delta|_{\mathcal{U}^0(M)} : \mathcal{U}^0(M)\to \mathcal{U}^0(N)$ is a surjective isometry.\smallskip
	
	Finally, the  implication $(c)\Rightarrow (a)$ is a consequence of Theorem~\ref{t surj isom principal components}.
\end{proof}

\begin{remark}\label{r unitaries is connected} If in Corollary~\ref{c two unital JB*algebras are isomorphic iff their unitaries are isometric iff their principal components are isometric} we additionally assume that $\mathcal{U} (M)$ or $\mathcal{U} (N)$ is connected the statements $(a)$ to $(c)$ are also equivalent to the existence of a surjective isometry $\Delta: \mathcal{U}(M)\to \mathcal{U}(N)$.
\end{remark}

It is already pointed out by the Kadison's theorem (see Theorem~\ref{t Kadison 55}) that two unital C$^*$-algebras are isometrically isomorphic if and only if they are Jordan $^*$-isomorphic. Since the set of all unitary elements in a unital C$^*$-algebra is a multiplicative subgroup, the conclusion in this setting is stronger (cf. \cite[Corollary 2]{HatMol2012}). 

\begin{corollary}\label{c two unital C*algebras are isomorphic iff their unitaries are isometric iff their principal components are isometric} Let $A$ and $B$ be two unital C$^*$-algebras. Then the following statements are equivalent:
	\begin{enumerate}[$(a)$]
		\item $A$ and $B$ are isometrically isomorphic as (complex) Banach spaces; \item $A$ and $B$ are isometrically isomorphic as real Banach spaces;
		\item $\mathcal{U}(A)$ and $\mathcal{U}(B)$ are isometrically isomorphic as metric spaces; 
		\item $A$ and $B$ are Jordan $^*$-isomorphic;
		\item There exist connected components $\mathcal{U}^1(A)\subseteq \mathcal{U}(A)$ and $\mathcal{U}^2(B)\subseteq \mathcal{U}(B)$ and a surjective isometry $\Delta: \mathcal{U}^1(A)\to \mathcal{U}^2(B);$
		\item There exists a surjective isometry $\Delta: \mathcal{U}^0(A)\to \mathcal{U}^0(B).$
	\end{enumerate}
\end{corollary}

\begin{proof} We only need to prove that $(e)\Rightarrow (d)$, the remaining implications have been established in Corollary~\ref{c two unital JB*algebras are isomorphic iff their unitaries are isometric iff their principal components are isometric}. Suppose there exists a surjective isometry $\Delta: \mathcal{U}^1(A)\to \mathcal{U}^2(B),$ where $\mathcal{U}^1(A)\subseteq \mathcal{U}(A)$ and $\mathcal{U}^2(B)\subseteq \mathcal{U}(B)$ are connected components. Pick $w_1\in \mathcal{U}^1(A)$ and $w_2\in \mathcal{U}^2(B)$. It is easy to check that $L_{w_2^*} \Delta L_{w_1} : \mathcal{U}^0(A)\to \mathcal{U}^0(B),$ $u\mapsto w_2^* \Delta (w_1 u )$ is a unital surjective isometry. Then $A$ and $B$ are Jordan $^*$-isomorphic by Corollary~\ref{c two unital JB*algebras are isomorphic iff their unitaries are isometric iff their principal components are isometric}.  
\end{proof}

We can now show that the conclusion in Corollary~\ref{c two unital JB*algebras are isomorphic iff their unitaries are isometric iff their principal components are isometric} is  optimal. The properties of unital JB$^*$-algebras differ from those in the strict subclass of unital C$^*$-algebra.

\begin{remark}\label{r existence of unitaries with non isomorphic nor isometric connected components in the Jordan setting} As we have seen in Anti-theorem~\ref{anti-theorem BKU}, there exists a unital JB$^*$-algebra $M$ with unit $\11b{M}$ containing a unitary element $w$ for which we cannot find a surjective linear isometry $T: M\to M$ mapping $\11b{M}$ to $w$, that is, the group of all surjective linear isometries on $M$ --equivalently, triple automorphisms-- on $M$ is not transitive on $\mathcal{U} (M)$. Furthermore,	the unital JB$^*$-algebras $M$ and $M(w)$ are not Jordan $^*$-isomorphic.\smallskip 
	
Let $\mathcal{U}^c (M)$  denote the connected component of $\mathcal{U} (M)$ containing $w$. We claim that $\mathcal{U}^0 (M)$ and $\mathcal{U}^c (M)$ are not isomorphic as metric spaces. Otherwise, Theorem~\ref{t surj isom principal components} (or Theorem~\ref{t surj isom secondary components}) implies that $M$ and $M(w)$ are Jordan $^*$-isomorphic, which is impossible.
	
\end{remark}

\subsection{Extendibility of surjective isometries between sets of unitary elements}\label{subsec: when is a surjective isometry extendable} \ \smallskip

We have seen in Example~\ref{example existence of non-extendible surjective isometries between unitaries} that a surjective isometry between the sets of unitary elements in two unital JB$^*$-algebras is not always extendable to a surjective real-linear isometry between the algebras. Remark \ref{r existence of unitaries with non isomorphic nor isometric connected components in the Jordan setting} shows an example of a unital JB$^*$-algebra $M$ containing a unitary element $w_0$ such that the principal component $\mathcal{U}^0 (M)$ and the connected component $\mathcal{U}^c (M)$ containing $w_0$ are not isomorphic as metric spaces.\smallskip 

Another geometric property of JB$^*$-triples which describes algebraically the extreme points of the closed unit ball of these Banach spaces has not been treated yet. Henceforth, let the symbol  $\partial_e({B}_X)$ stand for the set of all extreme points of the closed unit ball of a Banach space $X$. If $E$ is a JB$^*$-triple, the extreme points of its closed unit ball are precisely the complete tripotents in $E$, that is, \begin{equation}\label{eq extrem points complete tripotents} \partial_e({B}_E) = \hbox{Trip}_{max}(E)= \{ \hbox{complete tripotents in } E \}
\end{equation} (cf. \cite[Proposition 3.5]{KaUp77}).\smallskip

A result by Youngson proves that a JB$^*$-algebra is unital if and only if its closed unit ball have extreme points (see \cite{youngson1981} or \cite[Theorem 4.2.36]{Cabrera-Rodriguez-vol1}). It is also well known that for each unital JB$^*$-algebra $M$ we have $\mathcal{U} (M)\subset \partial_{e} (B_{M})$ (and the containing can be, in general, strict). The set $\partial_e({B}_M)$ seems to be an appropriate candidate to define a metric invariant for unital JB$^*$-algebras. However, the answer under these weak conditions is not always positive, since we can define a counterexample via the same method employed in Example~\ref{example existence of non-extendible surjective isometries between unitaries}.\smallskip 

The lacking or scarcity of unitary elements in the setting of real C$^*$-algebras makes hopless to conclude that every surjective isometry between sets of unitary elements in two unital real C$^*$-algebras is extendable to a surjective linear isometry between the algebras. We present next a counterexample taken from \cite[Remark 3.15]{CuePer20}.

\begin{example}\label{example counterexample real von Neumann algebras} Let us consider the real Banach space $A=\mathbb{R} \oplus^\infty \mathbb{R} = \ell_{\infty}^{2} (\mathbb{R}).$ Note that $A$ is actually a real von Neumann algebra with respect to the pointwise product in the usual sense (cf. \cite{Li2003}). The sets of all extreme points of  $B_A$ and the set of all unitary elements in $A$ coincide, however, it is a relative small set formed by the following four points $$\partial_e ({B}_A) = \mathcal{U} (A) =\{ e_1=(1,1),e_2=(1,-1),e_3=(-1,1),e_4=(-1,-1) \}.$$ The distances between elements in $\partial_e ({B}_A)$ satisfy $ \|e_i-e_j\|=2$ for every $i\neq j$ in $\{1,\ldots,4\}$. Therefore the mapping $\Delta:\partial_e ({B}_X)\to \partial_e ({B}_X) $ 
	$$\Delta (e_1)= e_2,\ \Delta(e_2)=e_3, \ \Delta(e_3)=e_4, \hbox{ and } \Delta(e_4)=e_1,$$ is a surjective isometry. We claim that $\Delta$ cannot be extended to a surjective linear isometry on $A$. Suppose, on the contrary, that $T: A\to A$ surjective real-linear isometry extending $\Delta$. In such a case $e_4 = \Delta (e_3) = T(e_3) = T (-e_2) ) =-T(e_2) = -\Delta (e_2) = -e_3$,  which is impossible. 
\end{example}

The natural question is whether there exist sufficient conditions to guarantee the extensibility of a surjective isometry between the sets of unitary elements.\smallskip
 
A Jordan variant of the Russo--Dye theorem for unital JB$^*$-algebras due to Wright and Youngson \cite{WriYou77} was employed in \cite{CuEnHirMiuPe2022} to find the following characterization of those surjective isometries between sets of unitary elements which admits an extension to a surjective real-linear isometry.

\begin{corollary}\label{c extendibility of surjective isometries}{\rm\cite[Corollary 3.8]{CuEnHirMiuPe2022}}  Let $\Delta :  \mathcal{U} (M)\to \mathcal{U} (N)$ be a surjective isometry between the unitary sets of two unital JB$^*$-algebras. Let $\{\mathcal{U}^0 (M)\}\cup \{ \mathcal{U}^{j} (M) : j\in \Lambda\}$ be the collection of all different connected components of $\mathcal{U} (M)$. For each $j\in \Lambda\cup \{0\}$, let $T_j : M\to N$ denote the surjective real-linear isometry extending the mapping $\Delta|_{\mathcal{U}^{j} (M)} : \mathcal{U}^{j} (M)\to \Delta(\mathcal{U}^{j} (M))$ whose existence is assured by Theorem \ref{t surj isom secondary components}. Then the following statements are equivalent:
	\begin{enumerate}[$(a)$]\item $\Delta$ admits an extension to a surjective {\rm(}real-linear{\rm)} isometry from $M$ onto $N$;
		\item $T_{j_1}=T_{j_2}$ for all $j_1,j_2\in \Lambda\cup \{0\}$.
	\end{enumerate}
\end{corollary}

\begin{proof} $(a)\Rightarrow(b)$ Suppose there exists a surjective isometry $T: M\to N$ such that $T|_{\mathcal{U} (M)} = \Delta$. Fix $j\in \Lambda\cup \{0\}$, by hypothesis, $T(u) = \Delta (u) = T_j (u)$ for all $u\in  \mathcal{U}^{j} (M)$.\smallskip
	
Fix $u_0\in \mathcal{U}^{j} (M)$. Let us consider the unital JB$^*$-algebra given by the $u_0$-isotope $M(u_0) =M_{(u_0^*)}$. By the Wright--Youngson--Russo--Dye theorem (see \cite[Corollary]{WriYou77}), the closed unit ball of the Banach space $M = M(u_0)$ coincides with the closed convex hull of the set
	$$\{ \exp_{(u_0)} (i h) : h\in (M(u_0))_{sa} \}.$$ By applying Remark \ref{r connected components which are not the pricipal one} we deduce that the latter set is contained in the principal component of $\mathcal{U} (M(u_0))$ which is precisely $\mathcal{U}^{j} (M)$. It follows from the real-linearity of $T$ and $T_j,$ combined with the fact that these two maps coincide on $\mathcal{U}^{j} (M)$, that $T=T_j$ on the closed unit ball of $M$, and hence on the whole $M$.\smallskip
	
	$(b)\Rightarrow(a)$ This implication is easy, it suffices to take $T= T_j$ for any $j\in \Lambda\cup \{0\}$.
\end{proof}

We continue with the following identity principle, which is of independent interest, taken from \cite{CuEnHirMiuPe2022}. 

\begin{proposition}\label{p T and T0}{\rm\cite[Proposition 3.9]{CuEnHirMiuPe2022}} Let $T,T_0 : M\to N$ be two surjective real-linear isometries between unital JB$^*$-algebras. Let us assume that $T_0$ is a real-linear Jordan $^*$-isomorphism, and suppose that $U_{T_0(u_0)} = U_{T(u_0)}$ for all $u_0\in \mathcal{U}^0 (M)$. Then there exists a central projection $p\in M$ such that $$ T|_{M\circ p\equiv M_2(p)} = T_0|_{M\circ p\equiv M_2(p)}, \hbox{ and } T|_{M\circ (\11b{M}-p)\equiv M_2(\11b{M}-p)} = - T_0|_{M\circ (\11b{M}-p)\equiv M_2(\11b{M}-p)}.$$
\end{proposition}

\begin{proof} Assume first that $T$ is unital and hence a real-linear Jordan $^*$-isomorphism (cf. Theorem~\ref{t Dang real-linear surjective isometries JBstar algebras}). It follows from the assumptions that given $h\in M_{sa}$ we have $$U_{T_0(e^{it h})} (\11b{N}) = U_{T(e^{it h})} (\11b{N}),$$ for all $t\in \mathbb{R}$. Taking derivatives at $t=0$ in both sides of the above identity we deduce that $$ 2 T_0( i h)= 2 U_{T_0( i h),\11b{N}} (\11b{N}) = 2 U_{T(i h),\11b{N}} (\11b{N})= 2 T(i h),$$ which proves that $T_0( i h) = T( i h)$ for all $h\in M_{sa}$ --recall that $T$ and $T_0$ are merely real-linear maps--. Having in mind that $T$ and $T_0$ are real-linear Jordan $^*$-isomorphisms we deduce that $$- T(h^2) = T(ih)^2 = T_0(ih)^2 = - T_0(h^2), \hbox{ for all } h\in M_{sa},$$ which proves that $T$ and $T_0$ coincide on positive elements in $M$. Since every element in $M_{sa}$ writes as the orthogonal sum of two positive elements in $M$, we conclude that $T$ and $T_0$ also coincide on $M_{sa}$. A general element $x\in M$ can be written in the form $x = h + i k$ with $h,k\in M_{sa},$ and hence $$T(x) = T(h +i k) = T(h) + T(i k) = T_0(h) + T_0(i k) = T_0(x).$$ We have therefore proved that $T= T_0$ in this case.\smallskip
	
We consider next the general case. By hypothesis, $Id= U_{T_0(\11b{M})} = U_{T(\11b{M})},$ and hence $T(\11b{M})^* = U_{T(\11b{M})}(T(\11b{M})^*) = T(\11b{M}),$ witnessing that $T(\11b{M})$ is a symmetry in $N$, that is, there exist two orthogonal projections $q_1,q_2\in N$ such that $T(\11b{M}) = q_1 - q_2$ and $q_1 + q_2=\11b{N}$. Since $U_{q_1} + U_{q_2} - 2 U_{q_1,q_2} =U_{q_1-q_2} = Id =U_{q_1+ q_2} = U_{q_1} + U_{q_2} + 2 U_{q_1,q_2},$ it can be easily deduced that $U_{q_1,q_2} =0$ and hence $N = N_2(q_1) \oplus^{\infty} N_2(q_2)$, which implies that $q_1$ and $q_2$ are central projections.\smallskip
	
Since $T_0$ is a real-linear Jordan $^*$-isomorphism, the elements $p_1 = T_0^{-1} (q_1)$ and $p_2 = T_0^{-1} (q_2)$ are two orthogonal central projections in $M$ with $p_1 + p_2 =\11b{M}$ and $M = M_2(p_1) \oplus^{\infty} M_2(p_2)$. We deduce from the above conclusions that every unitary $u$ in $M$ writes in the form $u = u_1 + u_2$ with $u_j \in M_2(p_j)$, and $\mathcal{U}^0 (M)$ identifies with $\mathcal{U}^0 (M_2(p_1))+ \mathcal{U}^0 (M_2(p_2))$. It is clear from the hypotheses that $T_0 $ is the orthogonal sum of the two real-linear Jordan $^*$-isomorphisms $T_0|_{M_2(p_1)} : M_2(p_1) \to N_2(q_1)$ and $T_0|_{M_2(p_2)} : M_2(p_2)\to N_2(q_2)$. Moreover, since $T$ is a surjective real-linear isometry (and hence a real-linear triple isomorphism, cf. Theorem~\ref{t Dang real-linear surjective isometries JBstar algebras}) with $T(\11b{M}) = q_1 - q_2$. Therefore, the mapping $T: M\to N(q_1-q_2)$ is a real-linear Jordan $^*$-isomorphism, and $T$ writes as the orthogonal sum of the real-linear Jordan $^*$-isomorphisms $T|_{M_2(p_1)} : M_2(p_1) \to N_2(q_1)$ and $T|_{M_2(p_2)} : M_2(p_2)\to N_2(-q_2)$.\smallskip
	
We shall next prove that the following pairs of (unital) maps $T_0|_{M_2(p_1)}, T|_{M_2(p_1)} : M_2(p_1) \to N_2(q_1)$ and $-T_0|_{M_2(p_2)},$ $T|_{M_2(p_2)} : M_2(p_2)\to N_2(-q_2)$ satisfy the hypothesis assumed for the pair $T_0,T$. To this end, let us take $u_j\in \mathcal{U}^0(M_2(p_j))$. The element $u = u_1 + u_2$ is a unitary in $\mathcal{U}^0 (M),$ and by our hypotheses and the fact that $M = M_2(p_1) \oplus^{\perp} M_2(p_2)$ and $N = N_2(q_1) \oplus^{\perp} N_2(q_2)$ we obtain that $$ U_{T(u_1)} + U_{T(u_2)}=  U_{T(u_1+u_2)} = U_{T_0(u_1+u_2)} = U_{T_0(u_1)} + U_{T_0(u_2)},$$ which implies that $U_{T_0(u_1)} =  U_{T(u_1)}$ and $U_{T_0(u_2)} = U_{T(u_2)}$, because these maps have orthogonal ranges and supports. \smallskip
	
We are therefore in a position to apply the conclusion in the first paragraph of this proof to the pairs $(T_0|_{M_2(p_1)}, T|_{M_2(p_1)})$ and $(-T_0|_{M_2(p_2)}, T|_{M_2(p_2)})$ to obtain that $T_0|_{M_2(p_1)}=T|_{M_2(p_1)}$ and $-T_0|_{M_2(p_2)}=T|_{M_2(p_2)}$. The desired conclusion follows by taking $p = p_1$.
\end{proof}

The extendibility of a surjective isometry $\Delta$ between the sets of unitary elements in two unital JB$^*$-algebras $M$ and $N$ is not always an easy task. Actually, the extendibility is, in general, hopeless. It is natural to ask whether an additional preservation hypothesis can be added to guarantee a linear extension. The instability of the set of unitaries under Jordan products invites to discard the preservation of the Jordan product as extra hypothesis. In view of the properties of surjective real-linear isometries between JB$^*$-algebras and the conclusion in Proposition~\ref{p surjective isometries preserving quadratic expressions}, the natural hypothesis is the following:  $$\Delta\{u,v,u\}= \Delta U_u (v^*) = U_{\Delta(u)} \left(\Delta(v)^* \right) =\{\Delta(u), \Delta(v), \Delta(u)\}, \hbox{ for all } u,v\in \mathcal{U} (M).$$ This preservation hypothesis is enough under assuming a mild extra assumption.

\begin{proposition}\label{p extra algebraic hypotheses}{\rm\cite[Proposition 3.10]{CuEnHirMiuPe2022}} Let $\Delta :  \mathcal{U} (M)\to \mathcal{U} (N)$ be a surjective isometry between the unitary sets of two unital JB$^*$-algebras. Let $\{\mathcal{U}^0 (M)\}\cup \{ \mathcal{U}^{j} (M) : j\in \Lambda\}$ be the collection of all different connected components of $\mathcal{U} (M)$. Assume that for each $j\in \Lambda$ there exists $u_j \in \mathcal{U} (M)$ with $u_j^2 \in \mathcal{U}^{j} (M)$. Suppose additionally that \begin{equation}\label{eq preservation of quadratic products extra hypo} \Delta\{u,v,u\}= \Delta U_u (v^*) = U_{\Delta(u)} \left(\Delta(v)^* \right) =\{\Delta(u), \Delta(v), \Delta(u)\},
	\end{equation} for all $u,v\in \mathcal{U} (M)$. Then $\Delta$ admits an extension to a surjective {\rm(}real-linear{\rm)} isometry from $M$ onto $N$.
\end{proposition}

\begin{proof} Up to replacing the original Jordan product and involution on $N$ by the one in the $\Delta(\11b{M})$-isotope, we can assume that $\Delta$ is unital, and hence $\Delta (\mathcal{U}^0 (M)) = \mathcal{U}^0 (N)$.  Let $T_j : M\to N$ denote the surjective real-linear isometry extending the mapping $\Delta|_{\mathcal{U}^{j} (M)} : \mathcal{U}^{j} (M)\to \Delta(\mathcal{U}^{j} (M))$ whose existence is assured by Theorem \ref{t surj isom secondary components}.\smallskip
	
	Fix $j,k\in \Lambda$ satisfying $\left(\mathcal{U}^j (M)\right)^* = \mathcal{U}^k(M)$. For each $u_j\in \mathcal{U}^j (M)$, it follows from \eqref{eq preservation of quadratic products extra hypo} that $$T_k (u_j^*) = \Delta (u_j^*) = \Delta U_{\mathbf{1}} (u_j^*) =  U_{\mathbf{1}} (\Delta(u_j)^*) = T_j (u_j)^*,$$ which proves that $T_k (u_j^*) = T_j (u_j)^*$ for all $u_j \in \mathcal{U}^j (M)$. We fix $w_j\in \mathcal{U}^j (M)$. Since, by the Wright--Youngson--Russo--Dye theorem \cite[Corollary]{WriYou77}, the convex hull of the set $\{ \exp_{(w_j)} (i h) : h\in (M(w_j))_{sa} \}(\subset \mathcal{U}^j (M) = \mathcal{U}^0 (M(w_j)))$ is norm dense in the closed unit ball of $M(w_j)$, it follows from the linearity and continuity of $T_k$ and $T_j$ that \begin{equation}\label{eq Tj and Tk coincide} T_k (x^*) = T_j (x)^* \hbox{ for all } x \in M.
	\end{equation}
	
	Take now an arbitrary $u_0 \in \mathcal{U}^0(A)$. By \eqref{eq inner main other main} in Remark~\ref{r connected component building from u as an pen question} we have $U_{u_0} (u_j^*)\in \mathcal{U}^k(M),$ and hence by the hypothesis \eqref{eq preservation of quadratic products extra hypo} we get $$\begin{aligned}
		U_{T_k(u_0)} (T_j(u_j)^*) &= U_{T_k(u_0)} (T_k(u_j^*)) = T_k U_{u_0} (u_j^*) = \Delta U_{u_0} (u_j^*) = U_{\Delta(u_0)} \left(\Delta(u_j^*)\right) \\
		&= U_{T_0(u_0)} \left( T_k(u_j^*)\right)= \hbox{\eqref{eq Tj and Tk coincide}}= U_{T_0(u_0)} \left( T_j(u_j)^*\right)
	\end{aligned},$$ for all $u_j\in \mathcal{U}^j (M)$. A new application of the Wright--Youngson--Russo--Dye theorem as in the previous paragraph --together with the surjectivity of $T_j$-- implies that $U_{T_k(u_0)} = U_{T_0(u_0)}$ for all $u_0 \in \mathcal{U}^0(M)$. Proposition \ref{p T and T0} assures the existence of a central projection $p_k\in M$ such that \begin{equation}\label{eq existence of pk} T_k = T_0|_{M_2(p_k)} \oplus -T_0|_{M_2(\11b{M}- p_k)}.
	\end{equation}

Fix an arbitrary $l\in \Lambda$. By assumptions, there exists $u_l\in \mathcal{U} (M)$ with $u_l^2 \in \mathcal{U}^{l} (M)$. We can assume that $u_l\in \mathcal{U}^{j} (M)$ and $u_l^*\in \mathcal{U}^{k} (M)$ for some $j,k\in \Lambda$. By hypotheses, \begin{equation}\label{eq Tl and Tj ul square} T_l (u_l^2) = \Delta (u_l^2) = \Delta U_{u_l} (\11b{M}) =  U_{\Delta(u_l)} \Delta (\11b{M}) = \Delta(u_l)^2 = T_j (u_l)^2.
	\end{equation} By considering the central projections $p_l$ and $p_j$ given by \eqref{eq existence of pk}, by orthogonality, we have $$ T_j (u_l)^2 = \left( T_0 (u_l\circ p_j) - T_0(u_l\circ (\11b{M}-p_j))\right)^2 = T_0 (u_l^2 \circ p_j) + T_0(u_l^2\circ (\11b{M}-p_j) ) = T_0 (u_l^2)$$ and $$ T_l (u_l^2) = T_0 (u_l^2\circ p_l) - T_0(u_l^2\circ (\11b{M}-p_l)),$$ which combined with \eqref{eq Tl and Tj ul square} give $T_0(u_l^2\circ (\11b{M}-p_l))=0,$ and hence $u_l^2\circ (\11b{M}-p_l)=0$ (by the injectivity of $T_0$). Having in mind that $p_l$ is a central projection, we have $$0 = u_l^2\circ (\11b{M}-p_l) = U_{u_l} (\11b{M}-p_l),$$ witnessing that $\11b{M}-p_l=0,$ and consequently, $T_l = T_0$. If we have in mind the arbitrariness of $l\in \Lambda$, the proof concludes by applying Corollary~\ref{c extendibility of surjective isometries}.
\end{proof}

The connections between the results in sections~\ref{sec: invertible elements as an invariant}, \ref{sec: geometric goodness of C* and JB-algebras} and \ref{sec: unitaries as invariant} and the different metric invariants studied in each one of them are much deeper than one might expect. For example, the concrete description of the surjective isometries between sets of invertible elements in two unital JB$^*$-algebras obtained in Theorem~\ref{t surjective isometries between invertible clopen JBstar} can be also deduced from Theorem~\ref{t surjective isometries between invertible clopen} and Theorem~\ref{t surj isom principal components} (Remark~\ref{r thm principal component and unital}). The argument relies on the following geometric characterization of unitary elements obtained in \cite{CuePer20unitaries}.

\begin{theorem}\label{t geometric char unitaries}{\rm\cite[Theorem 3.8]{CuePer20unitaries}}
Let $u$ be an extreme point of the closed unit ball of a unital JB$^*$-algebra $M$. Then $u$ is a unitary if and only if the set $$\mathcal{M}_{u} := \{e\in \partial_e({B}_{M}) : \|u\pm e\|\leq \sqrt{2} \}$$ contains an isolated point.
\end{theorem}

To see  how to get the conclusion in Theorem~\ref{t surjective isometries between invertible clopen JBstar}, let $\Delta: \mathfrak{M}\to \mathfrak{N}$ be a surjective isometry in the conditions of the just quoted theorem. Having in mind that every JB$^*$-algebra is semisimple, Theorem \ref{t surjective isometries between invertible clopen} proves the existence of a surjective real-linear isometry $T_0: M\to N$ whose restriction to $\mathfrak{M}$ is $\Delta$. Clearly $T_0$ maps $\partial_e(\mathcal{B}_M)$ onto $\partial_e(\mathcal{B}_N)$. Moreover, Theorem~\ref{t geometric char unitaries} assures that  $T_0(\mathcal{U}(M)) = \mathcal{U}(N)$, and hence $T_0|_{\mathcal{U}(M)} : \mathcal{U}(M)\to \mathcal{U}(N)$ is a surjective isometry. In particular, $u= T_0(\11b{M})$ is a unitary in $N$.\smallskip 

Let $\mathcal{U}(N_{(u^*)})=\mathcal{U}(N({u}))$ stand for the set of unitaries in the $u$-isotope $N({u})$ of $N$. As we observed in Lemma~\ref{l unitaries in the homotope coincide} $\mathcal{U}(N({u})) =\mathcal{U}(N)$. We can therefore conclude that $T_0|_{\mathcal{U}(M)} : \mathcal{U}(M)\to \mathcal{U}(N{(u)})$ is a unital surjective isometry. By applying Theorem~\ref{t surj isom principal components} we deduce the existence of a central projection $p$ in $M$ and an isometric Jordan $^*$-isomorphism $J:M\to N{(u)}$ such that \begin{equation}\label{eq T0 and the real Jordan isom coincide on the principal component} T_0( w ) =  J(p \circ w) + J( (\11b{M}-p)\circ w^*),
	\end{equation} for all $w\in \exp(i M_{sa})\subseteq M_{\mathbf{1}}^{-1}\subseteq \mathfrak{M}$.\smallskip
	
A new application of the Wright--Youngson--Russo--Dye theorem for unital JB$^*$-algebras (see \cite[Corollary 2.4]{WriYou77} or  \cite[Corollary 3.4.7]{Cabrera-Rodriguez-vol1}{\rm)} asserts that the closed unit ball of $M$ coincides with the closed convex-hull of the set $\exp(i M_{sa})$. Since the maps $T_0$ and $x\mapsto J(p \circ x + (\11b{M}-p)\circ x^*)$ are real-linear and continuous, we deduce from \eqref{eq T0 and the real Jordan isom coincide on the principal component} that $T_0( x ) =  J(p \circ x) + J( (\11b{M}-p)\circ x^*),$ for all $x\in M$, and consequently $$\Delta( a ) = T_0(a)=  J(p \circ a) + J( (\11b{M}-p)\circ a^*),$$ for all $a\in \mathfrak{M}$. The rest follows from the arguments given in the final part of the proof of Theorem \ref{t surjective isometries between invertible clopen JBstar}.\smallskip\smallskip

{\bf Acknowledgments.} Author partially supported by grant PID2021-122126NB-C31 funded by MCIN/AEI/10.13039/501100011033, Junta de Andaluc\'{\i}a grant FQM375, and by the IMAG--Mar\'{\i}a de Maeztu grant CEX2020-001105-M/AEI/10.13039/\linebreak501100011033.

\end{document}